\documentclass[11pt]{article}

\usepackage[margin=1in]{geometry}
\usepackage{amsmath,amssymb,amsthm,mathtools}
\usepackage{bm,mathrsfs}
\usepackage[shortlabels]{enumitem}
\usepackage{natbib}
\usepackage{xcolor}
\usepackage{microtype}
\usepackage[colorlinks=true,linkcolor=blue,citecolor=blue,urlcolor=blue]{hyperref}
\usepackage[nameinlink]{cleveref}
\usepackage[most]{tcolorbox}
\usepackage{graphicx}
\newtcolorbox{heuristicbox}{
  enhanced,
  breakable,
  colback=black!2,
  boxrule=0pt,
  borderline west={1.8pt}{0pt}{black!55},
  sharp corners,
  boxsep=0pt,
  left=10pt,
  right=3pt,
  top=8pt,
  bottom=8pt,
  before skip=1.2\baselineskip,
  after skip=1.2\baselineskip
}

\allowdisplaybreaks
\numberwithin{equation}{section}
\newcommand{\I}{\bm I}
\newcommand{\Id}{\operatorname{Id}}
\newcommand{\Tr}{\operatorname{Tr}}
\newcommand{\rmd}{\mathrm d}
\newcommand{\blk}{\operatorname{bl}}
\newcommand{\card}{\operatorname{card}}
\newcommand{\spec}{\operatorname{Spec}}
\newcommand{\ran}{\operatorname{Ran}}
\newcommand{\Ran}{\operatorname{Ran}}

\theoremstyle{plain}
\newtheorem{thm}{Theorem}[section]
\newtheorem{theorem}[thm]{Theorem}
\newtheorem{lemma}[thm]{Lemma}

\newtheorem{proposition}[thm]{Proposition}
\newtheorem{prop}[thm]{Proposition}
\newtheorem{corollary}[thm]{Corollary}
\newtheorem{coro}[thm]{Corollary}
\theoremstyle{definition}
\newtheorem{definition}[thm]{Definition}
\newtheorem{defn}[thm]{Definition}
\newtheorem{example}[thm]{Example}
\theoremstyle{remark}
\newtheorem{remark}[thm]{Remark}
\newtheorem{rmk}[thm]{Remark}

\title{Quantum Diffusion and the Block-Size Crossover in the Block Anderson Model}
\author{  Haoxuan Fu\\  {\small Department of Statistics and Data Science}\\  {\small Cornell University}\\  {\small\href{mailto:hf357@cornell.edu}{\texttt{hf357@cornell.edu}}}}
\date{}

\begin{document}
\maketitle

\begin{abstract}
We study the block Anderson Hamiltonian \(H_L=\lambda\Psi_L+\mathcal V\) on the discrete torus \((\mathbb Z/L\mathbb Z)^d\), \(d\geq2\), where \(\Psi_L\) is the nearest-neighbor Laplacian and \(\mathcal V\) consists of independent \(W^d\times W^d\) GOE blocks.  In the weak-disorder limit \(\lambda\to\infty\), we establish an entrywise local law, a spatial scaling limit for the squared resolvent entries, time-averaged dynamical delocalization, and lower bounds on the localization length of most bulk eigenfunctions.  Particularly, we allow the block size to be fixed or grow above a polylogarithmic threshold.  In sufficiently large polynomial volumes, the dynamical estimates hold for \(C_\delta(\lambda+W/\lambda)\leq T\leq \lambda^{1+\kappa_d^{\rm BA}-\delta}\), with a dimension-dependent exponent \(\kappa_d^{\rm BA}>0\), and the spreading scale is \(\sqrt{(\lambda^3+\lambda W)T}\). This identifies a transition at \(W\asymp\lambda^2\) between \(\lambda^{3/2}\sqrt T\) and \(\sqrt{\lambda WT}\).

The proof combines Gaussian moment estimates with deterministic stability of a matrix Dyson equation.  A partial Bloch transform handles fixed blocks, while scalar reduction controls growing blocks.  A Feshbach--Schur reduction of the block \(T\)-equation separates a single scalar diffusion mode from the stable internal cell modes.  A renewal expansion and exact cell counting give a positive representation of the diffusion coefficient and control the fourth moment uniformly across the block-size transition. Notably, our argument does not require an enumeration of Feynman diagrams.
\end{abstract}
\tableofcontents
\section{Introduction}
\label{sec:introduction}
Random Schr\"odinger operator, such as the Anderson tight-binding model, has been a central topics in the field of statistical physics. Since Anderson's seminal work \citep{1958PhRv..109.1492A}, it has provided a fundamental framework for describing quantum transport in disordered media, such as the motion of electrons in a metal with impurities. Mathematically, the scalar Anderson model is the random Hamiltonian
\begin{equation}
\label{eq:vanilla_anderson}
    \begin{aligned}
        H=-\Delta_{\mathbb Z^d}+\delta\mathcal V,
    \end{aligned}
\end{equation}
where $\Delta_{\mathbb Z^d}$ is the discrete Laplacian and $\mathcal V$ is a diagonal random potential with i.i.d.\ centered Gaussian entries. In dimensions $d\geq3$, the model is conjectured to undergo a metal--insulator transition, whose location depends on the disorder strength and the energy. Anderson localization is known at sufficiently large disorder or at extreme energies in every dimension, and at every nonzero disorder strength in $d=1$. Localization for one-dimensional random Schr\"odinger operators has been studied extensively; see, e.g., \citep{pure_point_1d_anderson, 10.1215/S0012-7094-02-11414-8}. A foundational breakthrough in higher dimensions was the multiscale analysis of \citet{absense_diff_low_energy}; the fractional moment method of \citet{local_diffu_absense_simp} subsequently provided an alternative approach. In $d=2$, complete Anderson localization at every disorder strength $\delta>0$ and every energy remains conjectural \citep{PhysRevLett.42.673}.

The weak-disorder regime $\delta\ll1$ is less well understood and gives rise to several outstanding conjectures in mathematical physics \citep{simon1984fifteen, 10.1063/1.533321}. In $d \geq 3$, the \textit{extended states} conjecture predicts absolutely continuous spectrum at suitable bulk energies, whereas the \textit{quantum diffusion} conjecture concerns diffusive spreading under the Schr\"odinger evolution. In $d=2$, where localization is expected at every disorder strength, diffusive behavior is still expected over a long but finite range of times. For the exploration of both conjecture, perturbative expansions have been demonstrated to be a powerful tool for analyzing scattering and recollision histories. Early work by \citet{deri_transp_equ_electr} derived a kinetic equation for the averaged momentum distribution on the timescale $t\sim\delta^{-2}$, for short kinetic times. For a continuum Gaussian random potential, \citet{https://doi.org/10.1002/(SICI)1097-0312(200006)53:6<667::AID-CPA1>3.0.CO;2-5} subsequently established the kinetic limit for the expected Wigner distribution on the whole kinetic timescale. Erd\H{o}s--Salmhofer--Yau then established diffusion for the expected Wigner distribution on timescales $t\sim\delta^{-2-\varepsilon}$, for sufficiently small $\varepsilon>0$, in the weak-coupling limit in $d\geq3$ \citep{ErdosSalmhoferYau2007, ErdosSalmhoferYau2007Recollision, ErdosSalmhoferYau2008}.

Recent developments in random matrix theory for models with spatially restricted randomness have opened another route toward understanding delocalization. Significant progress has been made for random band matrices \citep{YangYin2021, Yang_2022, YangYauYin2026}, and related methods apply to broader classes of non-mean-field models, including certain models in $d=2$ \citep{dubova2025delocalizationnonmeanfieldrandommatrices, dubova2025delocalizationtwodimensionalrandomband, yang2025delocalizationgeneralclassrandom, truong2025localizationlengthfinitevolumerandom}. Such models interpolate between the Anderson model, whose randomness is confined to the diagonal, and mean-field random matrices. Despite these advances, progress on delocalization for the Anderson model itself remains limited. Recently, \citet{black2025tailboundsdysonseries, black2025selfconsistentequationsquantumdiffusion} developed an operator-based approach to the perturbative regime. Their first work controls the Dyson series without enumerating individual Feynman diagrams; whereas in their subsequent work they derive quantum diffusion, in a time-averaged sense and for all $d \geq  2$, through self-consistent resolvent equations and dispersive estimates that control recollision effects. In this paper, we extend this framework to the block Anderson model in both fixed-block and growing-block regimes, without enumerating perturbative diagrams.

More precisely, we consider the torus $\mathbb Z_L^d\coloneqq(\mathbb Z/L\mathbb Z)^d$, with $L=nW$, decomposed into $n^d$ cubic blocks of side length $W$. We define
\[
\mathcal V=\operatorname{diag}(V_1,\ldots,V_{n^d}),
\]
where the $V_i$ are independent $W^d\times W^d$ matrices from the Gaussian orthogonal ensemble (GOE), with the standard $W^{-d}$ variance normalization. We use the large-hopping normalization
\[
H_L=\lambda\Psi_L+\mathcal V,
\qquad
\lambda^{-1}H_L=\Psi_L+\lambda^{-1}\mathcal V,
\]
where $\Psi_L=-\Delta_{\mathbb Z_L^d}$ is the nonnegative discrete Laplacian, so that $\lambda^{-1}$ plays the role of the disorder strength $\delta$ above.

The block Anderson model extends the scalar Anderson model, which is recovered when $W=1$, by replacing the diagonal potential on each block with a random matrix that couples all sites of the block. \Citet{yang2025delocalizationgeneralclassrandom} established delocalization and quantum unique ergodicity (QUE) for a general class of random block Schr\"odinger operators in $d\geq7$. Related block models in $d\geq3$, including a block Anderson model, were studied by \citet{dubova2025delocalizationnonmeanfieldrandommatrices}, who proved delocalization, QUE, and bulk universality. In low dimensions, \citet{truong2025localizationlengthfinitevolumerandom} obtained lower bounds on localization lengths in $d=1$ and delocalization and QUE in $d=2$. These results assume a power-law lower bound on the block size relative to the system size and a hopping strength of order at most one, and they exploit the many random degrees of freedom within each block. Our results instead concern the large-hopping (weak-disorder) limit $\lambda\to\infty$, and allow $W$ to remain fixed or to grow from a polylogarithmic threshold in $\lambda$, with no lower bound of the form $W\geq L^{\varepsilon}$.

These random matrix blocks introduce internal degrees of freedom and random couplings between microscopic sites within each block. The resulting deterministic approximation is matrix-valued, even though the spatial transport ultimately has a single scalar diffusion mode. Establishing this reduction to a single scalar mode, uniformly as the block size changes, is the main problem addressed here.

Our result establishes quantum diffusion, in a time-averaged sense and on a polynomial time window, for block Anderson Hamiltonians with fixed nontrivial blocks and growing blocks. More importantly, it includes a crossover phenomenon for quantum diffusion at $W\asymp\lambda^2$. They give a spatially tested resolvent profile at fixed bulk energy, high-probability time-averaged bounds on propagation, and lower bounds on the localization length of most eigenfunctions in regular bulk windows. The main analytic contribution is a local law governed by a matrix Dyson equation with block-local, globally non-flat covariance, followed by a Feshbach-Schur reduction of the matrix-valued T-equation. This reduction keeps the single nearly conserved density mode and controls the remaining internal cell modes. A positive renewal representation then quantifies the diffusion coefficient across block sizes. Together, these steps provide a strategy for studying transport in other random operators with internal cell degrees of freedom.

\subsection{Model and scaling}
Let \(\mathbb Z_L^d=(\mathbb Z/L\mathbb Z)^d\), where \(L=nW\), and partition the torus into \(n^d\) disjoint cubic blocks of side length \(W\).  The block labels form \(\widetilde{\mathbb Z}_n^d\), and \(\Pi_{[I]}\) denotes the orthogonal projection onto block \(I\).  On \(\ell^2(\mathbb Z_L^d)\), set
\[
    H_L=\lambda\Psi_L+\mathcal V,
    \qquad
    \mathcal V=\bigoplus_{I\in\widetilde{\mathbb Z}_n^d}V_I,
    \qquad \lambda\longrightarrow\infty,
\]
with $\Psi_L = 2d\I_L - \Delta_L$, where $\Delta_L$ denotes the discrete Laplacian on $\mathbb Z_L^d$ and $\I_L$ denotes the $L\times L$ identity matrix.  The matrices \(V_I\) are independent GOE matrices of dimension \(W^d\), normalized by
\begin{equation}
\label{eq:BA-GOE-block-normalization}
        \mathbb E (V_I)_{uv}(V_J)_{u'v'}
        =
        \frac{\mathbf 1_{\{I=J\}}}{W^d}
        \left(
        \mathbf 1_{\{u=u'\}}\mathbf 1_{\{v=v'\}}
        +
        \mathbf 1_{\{u=v'\}}\mathbf 1_{\{v=u'\}}
        \right).
\end{equation}
The disorder is independent between blocks and couples all microscopic coordinates within each block. In other words, $\mathcal{V}$ is a diagonal block matrix with i.i.d GOE blocks.  When \(W=1\), this reduces to a scalar Gaussian potential, with variance two under \cref{eq:BA-GOE-block-normalization}.

\begin{figure}[htbp]
    \centering
    \includegraphics[width=0.8\linewidth]{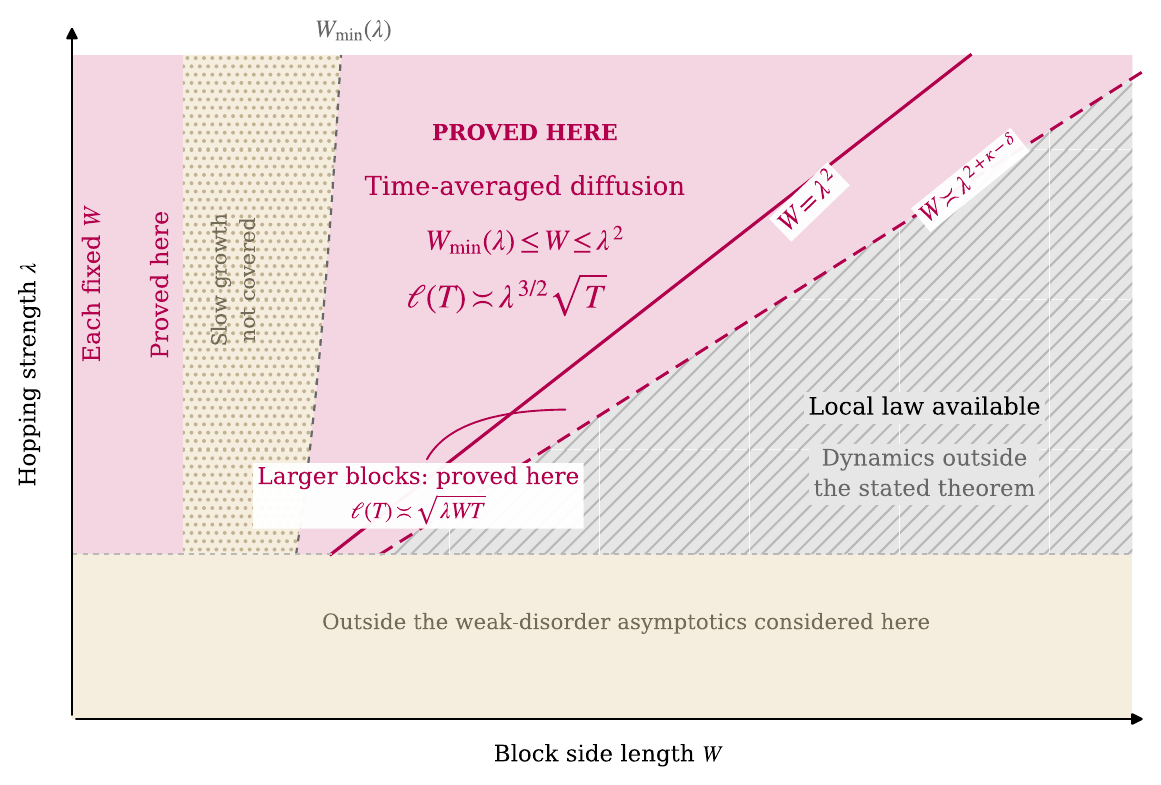}
    \renewcommand{\figurename}{Fig.}
    \caption{\textbf{Block Anderson model: coverage and crossover.}
    \textit{The random Hamiltonian studied is \(H_L=\lambda\Psi_L+\mathcal V\), with \(d\geq2\).
    Here \(\kappa=\kappa_d^{\mathrm{BA}}\) and
    \(W_{\min}(\lambda)=[C_0(1+\log\lambda)^4]^{p/d}\),
    with \(p\) as in the growing-block assumptions, and the diffusive length is denoted as $\ell(T)$.}}
    \label{fig:BA-parameter-regimes}
\end{figure}
Fig~\ref{fig:BA-parameter-regimes} illustrates the explored regime in our paper. There are two deterministic stability arguments. For fixed \(W\), we use the matrix-valued Bloch representation.  For growing blocks, we assume $W^{d/p}\geq C_0(1+\log\lambda)^4$ for a fixed \(p>p_d\), with \(p<4\) when \(d\geq7\), and use scalar reduction.  The local law has no further upper restriction on \(W\).  The dynamical window additionally requires
\[
        C_\delta(\lambda+W/\lambda)
        \leq T\leq\lambda^{1+\kappa_d^{\rm BA}-\delta}.
\]
For example, if \(W\asymp\lambda^{2+\alpha}\) and \(0<\alpha<\kappa_d^{\rm BA}\), this window is nonempty for \(0<\delta<\kappa_d^{\rm BA}-\alpha\). The term \emph{mesoscopic} will refer to \(W\leq\lambda^{2-\gamma_0}\) for a fixed \(\gamma_0>0\) when we need a vanishing ratio \(W/\lambda^2\), in particular for the explicit coefficient limit.  The transport estimates also cover the transition region and genuinely large blocks. Our choice and definition of $p_d$ and $\kappa_d^{\rm BA}$ could be found in Section~\ref{sec:model-notation} and Section~\ref{sec:diffusion-profile}.

\subsection{Main conclusions}
Write
\[
        R(z)=(H_L-z)^{-1},\qquad z=E+\mathrm i\eta\in\mathbb C^+,
\]
and let
\[
        M(z)=(\lambda\Psi_L-z-\Theta(z))^{-1},
        \qquad \Theta(z)=S[M(z)],
\]
where we call \(S\) the covariance map of the block-GOE disorder, or superoperator in short. We defer its definition to \cref{eq:glob_super_op}. Our first result compares the entries of \(R\) and its deterministic counterpart \(M\) in a regular bulk energy window.  The exponents \(p_d\), exceptional energies \(\Sigma_d^\Psi\), and stochastic domination notation \(\prec\) are defined in Section~\ref{sec:model-notation}. The constants below may depend on the fixed energy window and, in the fixed-block regime, on \(W\).
\begin{theorem}[Entrywise local law]
\label{thm:local_law_block_anderson}
Fix \(d\geq2\), \(\varepsilon>0\), and a sufficiently small \(\delta>0\), and set \(e=E/\lambda\).  Suppose one of the following holds:
\begin{enumerate}[(a)]
\item \(W^{d/p}\geq C_0(1+\log\lambda)^4\) for some fixed \(p>p_d\) (and \(p<4\) when \(d\geq7\)), and \(e\in(0,4d)\) satisfies \(\operatorname{dist}(e,\Sigma_d^\Psi)\geq\varepsilon\);
\item \(W<\infty\) is fixed and \(e\) belongs to a fixed finite-height folded bulk window in the sense of Definition~\ref{def:fh-bulk-window}.
\end{enumerate}
The following bounds hold uniformly for \(\eta\leq K\lambda\), with \(K>0\) fixed.

If \(d\in\{2,4,5,6\}\), then for \(\eta\geq\lambda^{-1-2/9+\delta}\),
\begin{align}
\label{eq:intermediate-expectation-local-law}
        \|\mathbb ER-M\|_{1\to\infty}
        &\leq
        C_\delta\lambda^{-3/2+\delta}(1 + (\lambda\eta)^{-1})^{9/4},\\
\label{eq:intermediate-fluctuation-local-law}
        \|R-\mathbb ER\|_{1\to\infty}
        &\prec
        \lambda^{-3/2}(1+(\lambda\eta)^{-1})^{7/4}.
\end{align}
If \(d=3\), then for \(\eta\geq\lambda^{-1-6/19+\delta}\),
\begin{align}
\label{eq:d3-expectation-local-law}
        \|\mathbb ER-M\|_{1\to\infty}
        &\leq
        C_\delta\lambda^{-7/4+\delta}(1+(\lambda\eta)^{-1})^{19/8},\\
\label{eq:d3_local_law_stoch}
        \|R-\mathbb ER\|_{1\to\infty}
        &\prec
        \lambda^{-7/4}(1 + (\lambda\eta)^{-1})^{15/8}.
\end{align}
If \(d\geq7\), set $\kappa_d^{\rm LL}\coloneqq\frac{2d-12}{3d-6}.$ For \(\eta\geq\lambda^{-1-\kappa_d^{\rm LL}+\delta}\),
\begin{align}
\label{eq:highd-one-four-final}
        \|R\|_{1\to4}&\prec\lambda^{-1},\\
\label{eq:highd-fluctuation-final}
        \|R-\mathbb ER\|_{1\to\infty}&\prec\lambda^{-2},\\
\label{eq:highd-expectation-final}
        \|\mathbb ER-M\|_{1\to\infty}
        &\leq C_\delta\lambda^{-2+\delta}(1 + (\lambda\eta)^{-1})^{1/2}.
\end{align}
\end{theorem}
The local law above quantifies the behavior of a single resolvent, while the transport is encoded by the two-resolvent observable \(R\mathcal A R^*\), or the resolvent sandwich, where \(\mathcal A\) is a deterministic block-diagonal test operator. This forces us to consider the correlation between two resolvents. By Gaussian integration by parts, we establish a block \(T\)-equation for its expectation.  After eliminating the stable internal cell modes, its inverse is governed at low frequency by
\[
        \frac1{s(k)}\approx
        \frac1{s(0)+d_{\rm BA}|k|^2},
        \qquad s(0)\asymp\lambda\eta.
\]
Heuristically, the diffusion coefficient $d_{\rm BA}$ satisfies
\[
        W^2d_{\rm BA}\asymp\lambda^4+\lambda^2W.
\]
Thus division by the intensity factor \(\lambda\) gives the microscopic transport scale \(\lambda^3+\lambda W\). This explains the crossover phenomenon due to the block size transition, and our next two theorems explicitly quantify such transport scale crossover. We defer the explicit exponent \(\kappa_d^{\rm BA}\) appears in the next theorems to \cref{eq:BA-dynamical-exponent}.
\begin{thm}[Time-averaged quantum diffusion]
\label{thm:BA-dynamical-delocalization}
Let \(d\geq2\) and \(0<\delta<\kappa_d^{\rm BA}\).  Suppose that one of the following two regimes holds.
\begin{enumerate}[(a)]
    \item \textit{Fixed block size regime:} The block size \(W<\infty\) is fixed.
    \item \textit{Growing block size regime:} For some fixed \(p>p_d\), with \(p<4\) when \(d\geq7\),
    \[
        W^{d/p}\geq C_0(1+\log\lambda)^4,
        \qquad W\leq\lambda^2.
    \]
\end{enumerate}
Let \(x=(0,u)\), where \(u\in[W^d]\).  There are constants \(c_\delta,C_\delta>0\) and \(\lambda_0<\infty\) such that, for every \(\lambda\geq\lambda_0\) and every deterministic \(T\) satisfying
\begin{equation}
\label{eq:BA-dynamical-time-volume-range}
       C_\delta\lambda
        \leq T
        \leq
        \lambda^{1+\kappa_d^{\rm BA}-\delta}, \quad \frac{n}{\lambda T}\geq(\log\lambda)^2
\end{equation}
the following bounds hold:
\begin{equation}
\label{eq:BA-dynamical-second-moment}
\begin{aligned}
        \frac1T\int_0^T
        \sum_{y\in\mathbb Z_L^d}|y-x|^2
        \left|(e^{-\mathrm i tH_L})_{xy}\right|^2\,\rmd t
        \geq c_\delta\lambda^3T
\end{aligned}
\end{equation}
and
\begin{equation}
\label{eq:BA-dynamical-outer-tail}
\begin{aligned}
        \frac1T\int_0^T
        \sum_{\substack{y\in\mathbb Z_L^d:\\
        |y-x|\geq C_\delta\lambda^{3/2}\sqrt T}}
        \left|(e^{-\mathrm i tH_L})_{xy}\right|^2\,\rmd t
        \leq\delta
\end{aligned}
\end{equation}
with probability at least \(1-C_\delta\lambda^{-1000}\).  Here \(|y-x|\) denotes the periodic Euclidean distance on \(\mathbb Z_L^d\), represented in \([-L/2,L/2)^d\).  In the fixed block size regime the constants may depend on \(W\) and on the finite collection of folded bulk windows.
\end{thm}

The following statement unifies the admissible growing-block regimes, including the transition \(W\asymp\lambda^2\).
\begin{thm}[Time-averaged quantum diffusion with genuine large blocks]
\label{thm:BA-genuine-large-block-dynamics}
Let \(d\geq2\), \(0<\delta<\kappa_d^{\rm BA}\), and suppose that, for some fixed \(p>p_d\), with \(p<4\) when \(d\geq7\),
\[
        W^{d/p}\geq C_0(1+\log\lambda)^4,
        \qquad W\leq\lambda^{2+\kappa_d^{\rm BA}-\delta}.
\]
Keep the model and polynomial-volume assumptions of Theorem~\ref{thm:BA-dynamical-delocalization}. There are constants \(c_\delta,C_\delta>0\) and \(\lambda_0<\infty\), independent of \(W\), such that for every \(\lambda\geq\lambda_0\), every deterministic \(T\) satisfying
\begin{equation}
\label{eq:BA-large-block-time-range}
        C_\delta\left(\lambda+\frac W\lambda\right)
        \leq T\leq\lambda^{1+\kappa_d^{\rm BA}-\delta},
        \qquad \frac n{\lambda T}\geq(\log\lambda)^2,
\end{equation}
and every fixed \(x=(0,u)\), with probability at least \(1-C_\delta\lambda^{-1000}\),
\begin{equation}
\label{eq:BA-large-block-second-moment}
        \frac1T\int_0^T\sum_{y\in\mathbb Z_L^d}|y-x|^2
        \left|(e^{-\mathrm itH_L})_{xy}\right|^2\,\rmd t
        \geq c_\delta(\lambda^3+\lambda W)T
\end{equation}
and
\begin{equation}
\label{eq:BA-large-block-outer-tail}
        \frac1T\int_0^T
        \sum_{\substack{y\in\mathbb Z_L^d:\\
        |y-x|\geq C_\delta\sqrt{(\lambda^3+\lambda W)T}}}
        \left|(e^{-\mathrm itH_L})_{xy}\right|^2\,\rmd t
        \leq\delta.
\end{equation}
Distances have the same periodic convention as in Theorem~\ref{thm:BA-dynamical-delocalization}.
\end{thm}
\begin{remark}
    Theorem~\ref{thm:BA-dynamical-delocalization}-\ref{thm:BA-genuine-large-block-dynamics} reveal how the block size changes the transport scale. A free flight has microscopic length of order $\lambda^2$. When $W\ll\lambda^2$, one flight typically traverses many blocks, whereas when $W\gg\lambda^2$, the probability that one flight crosses a block boundary is of order $\lambda^2/W$, averaged over the starting position in the block. Exact endpoint counting keeps both contributions and gives an effective microscopic diffusion coefficient of order $\lambda^3+\lambda W$ for the evolution generated by $H_L$. Equivalently, the spreading scale at time $T$ is $\sqrt{(\lambda^3+\lambda W)T}$. For $W\ll\lambda^2$ this is $\lambda^{3/2}\sqrt T$, which matches the scalar prediction $\delta^{-1}\sqrt t$ in the time variable $t=\lambda T$ in the Anderson model \cref{eq:vanilla_anderson}. A uniform fourth-moment bound allows us to pass through $W\asymp\lambda^2$ without a gap between the two regimes.
\end{remark}

The same resolvent estimates also give a spatial scaling limit and a lower bound on the localization length of most eigenfunctions. The coefficients \(d_{\rm BA}\) and \(j_u(0)\) below are defined in Subsection~\ref{subsec:BA-Fourier-profile}.

\begin{thm}[Resolvent scaling limit]
\label{thm:BA-resolvent-scaling}
Let \(d\geq3\), \(0<\kappa<\kappa_d^{\rm BA}\), and \(0<a<\frac18\min\{\kappa,\kappa_d^{\rm BA}-\kappa\}\). Suppose that one of the following holds:
\begin{enumerate}[(a)]
\item \(W<\infty\) is fixed and \(e\) belongs to a fixed finite-height folded bulk window;
\item \(W^{d/p}\geq C_0(1+\log\lambda)^4\) for a fixed \(p>p_d\), with \(p<4\) when \(d\geq7\), \(W\leq\lambda^2\), and \(\operatorname{dist}(e,\Sigma_d^\Psi)\geq\varepsilon\);
\item \(\lambda^2\leq W\leq\lambda^{2+\kappa-8a}\) and \(\operatorname{dist}(e,\Sigma_d^\Psi)\geq\varepsilon\).
\end{enumerate}
In \textnormal{(b)} and \textnormal{(c)}, assume also \(e\in(0,4d)\).  Set
\[
        z=\lambda e+\mathrm i\eta,
        \qquad \eta=\lambda^{-1-\kappa},
        \qquad
        r_\lambda=\lambda^{-a}
        \sqrt{\frac{\lambda^4+\lambda^2W}{\lambda\eta}}.
\]
Keep the polynomial-volume assumption and suppose that \(n\eta/\lambda\geq(\log\lambda)^2\). For every \(f\in C_c^\infty(\mathbb R^d)\) and every \(x=(0,u)\),
\begin{equation}
\label{eq:BA-resolvent-scaling-limit}
        \frac{W^2d_{\rm BA}}{j_u(0)r_\lambda^2}
        \sum_{y\in\mathbb Z_L^d}
        f\left(\frac{y-x}{r_\lambda}\right)|R(z)_{xy}|^2
        \longrightarrow
        \frac1{(d-2)|\mathbb S^{d-1}|}
        \int_{\mathbb R^d}f(y)|y|^{2-d}\,\rmd y
\end{equation}
in probability as \(\lambda\to\infty\).  Here \(y-x\) is represented in the centered fundamental domain, and \(d_{\rm BA}\) and \(j_u(0)\) are evaluated at \(z\). The difference in \cref{eq:BA-resolvent-scaling-limit} is \(\prec\lambda^{-a}\).

The same statement holds for \(d=2\) provided \(\int f=0\), with the right-hand side replaced by $-\frac1{2\pi}\int_{\mathbb R^2}f(y)\log|y|\,\rmd y.$ In case \textnormal{(b)}, if additionally \(W\leq\lambda^{2-\gamma_0}\) for some fixed \(\gamma_0>0\), one may replace \(W^2d_{\rm BA}/j_u(0)\) on the left-hand side by \(\lambda^4\) and multiply the right-hand side by
\begin{equation}
\label{eq:BA-mesoscopic-beta}
        \beta_e\coloneqq
        \frac{4\bigl(\pi\rho(e)\bigr)^3}
        {\displaystyle\pi(2\pi)^{-d}
        \int_{\{\Omega=e\}}
        \frac{|\partial_1\Omega(p)|^2}{|\nabla\Omega(p)|}
        \,\rmd\sigma(p)},
\end{equation}
where \(\rho\) is the free density of states and \(\rmd\mathcal{H}^{d-1}(p)\) is the $d-1$-dimensional Hausdorff measure. Both limits are almost sure along \(\lambda=m\in\mathbb N\), under any joint coupling of the admissible matrices.
\end{thm}

Here the scale \(r_\lambda\) is comparable to \(\lambda^{2+\kappa/2-a}\) for fixed blocks and growing blocks with \(W\leq\lambda^2\), and to \(\lambda^{1+\kappa/2-a}\sqrt W\) in the large-block regime.  The normalization in \cref{eq:BA-resolvent-scaling-limit} maintains the deterministic diffusion coefficient and the observation profile.  In the proof we will highlight that
\[
        \frac{W^2d_{\rm BA}}{j_u(0)}
        \asymp \lambda^4+\lambda^2W.
\]
Outside the mesoscopic range, a limit with normalization \((\lambda^4+\lambda^2W)r_\lambda^{-2}\) and a fixed prefactor additionally requires convergence of \((\lambda^4+\lambda^2W)j_u(0)/(W^2d_{\rm BA})\). Note that no such coefficient limit is needed for \cref{eq:BA-resolvent-scaling-limit}.

Finally, we state our result on the delocalization bound of most eigenfunctions. Let \((E_j,\psi_j)_{j=1}^{L^d}\) be the eigenpairs of \(H_L\), with \((\psi_j)\) orthonormal. For \(r>0\), define
\[
        \mathfrak L(r)\coloneqq
        \left\{E_j:\text{ there exists }x_0\in\mathbb Z_L^d
        \text{ such that }
        \|\psi_j\|_{\ell^2(B_r(x_0))}\geq1-r^{-4d}\right\},
\]
where balls use the periodic distance and eigenvalues are counted with multiplicity. Here $\mathfrak L(r)$ contains all eigenvalues whose eigenfunction decay exponentially at scale $\ll r$.

\begin{thm}[Delocalization bound for most eigenfunctions]
\label{thm:BA-most-eigenfunctions}
Let \(d\geq2\), \(0<\delta<\kappa_d^{\rm BA}/4\), and keeps the polynomial-volume assumptions of Section~\ref{sec:model-notation}. Suppose one of the following holds:
\begin{enumerate}[(a)]
\item \(W<\infty\) is fixed;
\item \(W\) satisfies the growing-block assumptions in part \textnormal{(b)} of Theorem~\ref{thm:BA-dynamical-delocalization};
\item \(\lambda^2\leq W\leq \lambda^{2+\kappa_d^{\rm BA}-4\delta}\).
\end{enumerate}
Fix \(\varepsilon>0\) and \(e\in(0,4d)\) such that \([e-2\varepsilon,e+2\varepsilon]\) is contained in a compact regular energy set in the applicable regime. Set
\[
        r=
        \begin{cases}
        \lambda^{2+\kappa_d^{\rm BA}/2-\delta},
                &\text{in cases (a) and (b)},\\[1mm]
        W^{1/2}\lambda^{1+\kappa_d^{\rm BA}/2-\delta},
                &\text{in case (c)}.
        \end{cases}
\]
If $n\lambda^{-2-\kappa_d^{\rm BA}+\delta}\geq(\log\lambda)^2,$ then, for all sufficiently large \(\lambda\),
\[
        \frac{\card\bigl(\mathfrak L(r)\cap
               [\lambda(e-\varepsilon),\lambda(e+\varepsilon)]\bigr)}
             {\card\bigl(\spec(H_L)\cap
               [\lambda(e-\varepsilon),\lambda(e+\varepsilon)]\bigr)}
        \leq C_{\varepsilon,\delta}\lambda^{-\delta/2}
\]
with probability at least \(1-C_{\varepsilon,\delta}\lambda^{-1000}\). The constants are uniform in \(W\) in cases \textnormal{(b)} and \textnormal{(c)}; in case \textnormal{(a)} they may depend on \(W\) and the folded bulk window.
\end{thm}

For fixed blocks and growing blocks with \(W\leq\lambda^2\), the delocalization bound is \(\lambda^{2+\kappa_d^{\rm BA}/2-\delta}\), while for genuinely large blocks it becomes \(W^{1/2}\lambda^{1+\kappa_d^{\rm BA}/2-\delta}\). The driving mechanism underneath is that we need a strict upper restriction on \(W\) to ensure this radius contains enough blocks for the diffusive behavior.
\subsection{The general proof strategy}
\begin{figure}[htbp]
    \centering
    \includegraphics[width=1.0\linewidth]{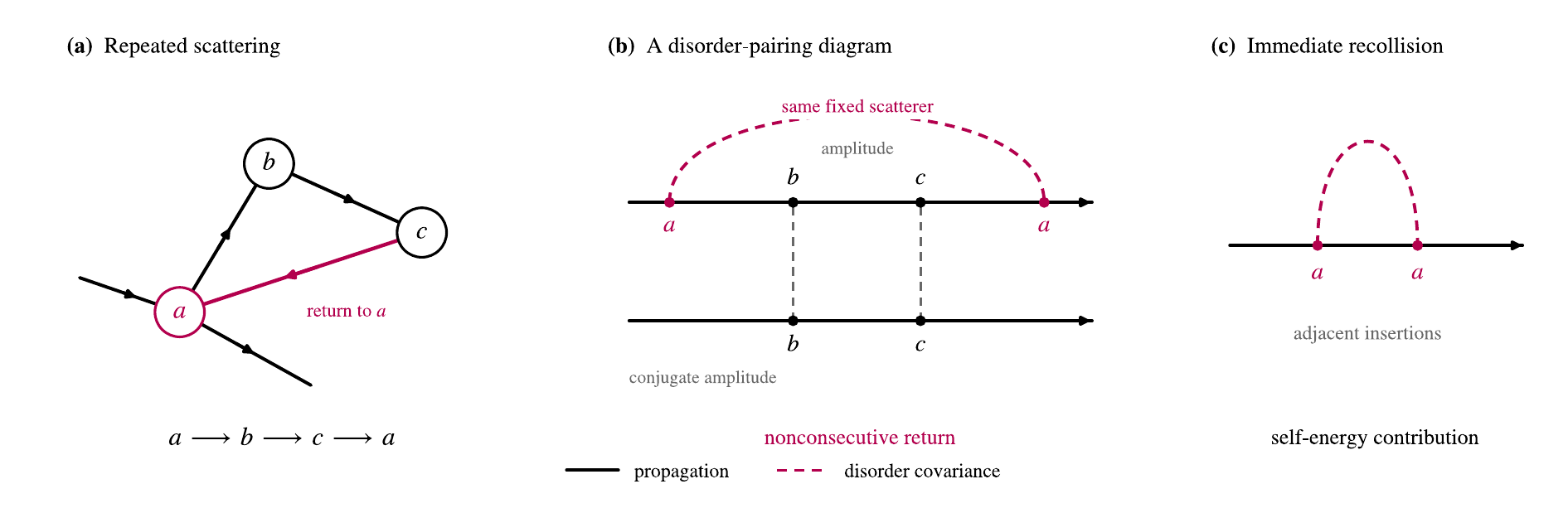}
    \renewcommand{\figurename}{Fig.}
    \caption{\textbf{Recollision and self-energy renormalization}
    \textit{Panel (a) shows a possible scattering history $a\to b\to c\to a$, with a recollision that obstructs a Markovian description; panel (b) visualizes a possible disorder pairing in the Wick expansion of the resolvent sandwich; panel (c) highlights the immediate recollisions resummed by self-energy renormalization.}}
    \label{fig:recollision_self_energy}
\end{figure}
During the random evolution, one of the most delicate issues is particle recollision. Recollisions correlate successive scattering events. See panel (a) of Fig.~\ref{fig:recollision_self_energy} for a visualization of a recollision process. Such recollisions obstruct the Markovian description of the scattering process and induce correlations even within one rail of the propagator expansion. These correlations introduce a even more complicated pairing structure when the resolvent sandwich is considered. See panel (b) for a visualization of one possible pairing. To address this, following \citet{erdos2006quantumdiffusionrandomschrodinger}, we use self-energy renormalization to account for immediate recollisions, where the renormalization is generally matrix-valued because of the block structure. See panel (c) for a visualization.

To tackle the remaining recollision contributions, we begin with a priori resolvent bounds obtained from the free Laplacian and a Dyson-series comparison between the random and free propagators. Gaussian moment inequalities then control the fluctuations of both a single resolvent and the resolvent sandwich \(R\mathcal AR^*\). Fluctuation control and the a priori dispersive resolvent bounds control the remaining recollision contributions in the relevant spectral range. Together with stability, these estimates allow the resolvent to be approximated using a self-consistent equation obtained from Gaussian integration by parts.

\emph{Step 1: self-consistent equation, stability, and the local law.} Up to controlled errors, the self-consistent equation is a matrix Dyson equation \citep{Erds2019TheMD, Ajanki_2018}. The goal is to prove that this equation remains stable against the perturbation error from the remaining recollision contributions and the error from the finite-to-infinite-volume comparison. For growing blocks, the self-energy is approximated by a scalar renormalization. For fixed blocks, the partial Bloch transform gives a family of finite-dimensional matrix equations, coupled through the self-energy.  Geometrically analyzing regular level sets of the free dispersion provides the necessary dissipative bounds and stability.  In both cases, stability holds in a fixed neighborhood of the deterministic solution, which is exactly what the continuity bootstrap in \(\eta\) requires. A continuity argument then proves the local law down to the spectral scale corresponding to the diffusive timescale.

\emph{Step 2: the block \(T\)-equation.} The local law controls a single resolvent; to obtain the diffusive profile, we need to analyze the resolvent sandwich. By Gaussian integration by parts, we get an equation for \(\mathbb E[R\mathcal A R^*]\), called the block $T$-equation.  Its kernel propagates a cell matrix from one block to another and applies the GOE covariance map.  Iterating the kernel recovers the ladder series.  Technically, we call the resolvent propagation between successive covariance contractions a single flight; its positive intensity, averaged over the starting cell, is given by ${W^{-d}}\|\widetilde M_\infty(X)\|_{\mathrm{HS}}^2$ for a block displacement \(X\), where \(\widetilde M_\infty(X)\) is the block at displacement \(X\) of the infinite-volume deterministic resolvent. The local law and fluctuation estimates control the error in replacing the random two-resolvent observable by this deterministic ladder response.

\emph{Step 3: reduction to the scalar diffusion mode.} The ladder series involves convolution in the block coordinate, so the Bloch transform naturally turns the ladder kernel into a superoperator acting on cell matrices.
To identify the unique slow mode which quantifies the diffusion coefficient, a Feshbach--Schur reduction \citep{GriesemerHasler2008} inverts the stable internal cell sector and reveals the unique nearly singular direction.
At zero frequency, this slow mode is identified, up to a controlled error, by the Ward identity. We call the scalar $s(k)$ obtained from the Schur complement the effective diffusion symbol, whose low-frequency behavior determines the diffusion.

\emph{Step 4: positivity and the block-size transition.} To determine the scale of $s(k)$, notice that the GOE covariance splits the ladder kernel into an internal matrix map and a scalar trace reset.  The internal map is strictly contracting, and summing its successive applications gives a positive renewal formula for \(s(k)\).  Differentiating this formula relates the diffusion coefficient to the second moment $W^{-d}\sum_X X_i^2 \|\widetilde M_\infty(X)\|_{\mathrm{HS}}^2.$ Cubic symmetry makes the resulting diffusion coefficient scalar. For growing blocks, exact counting of the endpoints of a free flight produces the desired quantification of the diffusion coefficient. Notably, the fourth-moment estimate \eqref{eq:BA-genuine-large-block-fourth-moment} controls the error in the quadratic approximation of \(s(k)\) uniformly on both sides of \(W\asymp\lambda^2\).

\emph{Step 5: the consequences of the resolvent profile.} We compare the finite-volume ladder response with its infinite-volume counterpart and approximate the latter by a massive elliptic Green function.  Removing the mass after spatial rescaling gives the resolvent limit.  Positive radial test functions, energy integration, and the Laplace--Plancherel identity give the time-averaged dynamical bounds.  Finally, a spectral compression estimate turns the small-ball resolvent bound into a delocalization bound for most eigenfunctions.

\subsection{Organization}

Section~\ref{sec:model-notation} fixes notation, and Section~\ref{sec:apriori-resolvent} establishes the a priori resolvent bounds.  Section~\ref{sec:self-consistent-equation} derives the self-consistent equation, Section~\ref{sec:kinetic-initialization} solves it in the kinetic range, and Section~\ref{sec:local-law} proves the local law, taking the deterministic stability inputs as given. Section~\ref{sec:diffusion-profile} contains the block \(T\)-equation, the scalar reduction, the uniform block-size comparison, and the resolvent, dynamical, and eigenfunction consequences. Section~\ref{sec:common-deterministic-tools} supplies the Bloch and Combes--Thomas tools, and Sections~\ref{sec:large-block-stability} and~\ref{sec:fixed-block-stability} prove the deterministic stability estimates for growing and fixed blocks, respectively. The appendices provide the parallel-transport and cubic-symmetry arguments.
\section{Model and notation}
\label{sec:model-notation}

We keep the notation from the introduction.  Thus \(L=nW\), \(\mathbb Z_L^d=(\mathbb Z/L\mathbb Z)^d\), and every site is written as \(x=(I,u)\), where \(I\in\widetilde{\mathbb Z}_n^d\) is a block coordinate and \(u\in[W]^d\) is the coordinate inside the block.  We identify \([W]^d\) with \([W^d]\) whenever a single cell index is more convenient.  The volume is chosen polynomially large in \(\lambda\); more precisely, \(L\asymp\lambda^{C_{\mathrm{vol}}}\) along integers compatible with \(L=nW\), where \(C_{\mathrm{vol}}\) is fixed and sufficiently large in terms of \(d\).  Only this polynomial growth and the explicit finite-volume condition in Theorem~\ref{thm:BA-dynamical-delocalization} are used below.

For \(z=E+\mathrm i\eta\in\mathbb C^+\), write $ R(z)=(H_L-z)^{-1}.$ For a matrix \(A\) on \(\ell^2(\mathbb Z_L^d)\), \(\|A\|_{p\to q}\) denotes its operator norm from \(\ell^p\) to \(\ell^q\).  If \(B\) is a cell matrix, set
\[
        \langle B\rangle_W\coloneqq\frac1{W^d}\Tr B,
        \qquad
        \mathscr S_W[B]
        \coloneqq
        \frac1{W^d}\bigl(B^\top+\Tr(B)\I_{W^d}\bigr).
\]
The global covariance map is the block-diagonal operator
\begin{equation}
\label{eq:glob_super_op}
    \begin{aligned}
         S[A]
        \coloneqq
        \sum_{I\in\widetilde{\mathbb Z}_n^d}
        \iota_I\,\mathscr S_W[\iota_I^*A\iota_I]\,\iota_I^* .
    \end{aligned}
\end{equation}

For $Y\succ0$, we also define the sandwich operator $\mathcal C_Y[B]\coloneqq Y^{-1/2}BY^{-1/2}.$ When we write order-unit norm generated by $Y$, we mean the weighted norm given by $\|\cdot\|_{Y} \coloneqq \|\mathcal{C}_Y[\cdot]\|_{2\to 2}$
\begin{definition}[Stochastic domination]
\label{def:stoch_domin}
For nonnegative random quantities \(A=A_\lambda\) and \(B=B_\lambda\), we write \(A\prec B\) if, for every \(\gamma>0\) and \(N>0\), there is \(C_{\gamma,N}<\infty\) such that
\[
        \mathbb P\bigl(A>C_{\gamma,N}\lambda^\gamma B\bigr)
        \leq C_{\gamma,N}\lambda^{-N}
\]
for all sufficiently large \(\lambda\).
\end{definition}

The Fourier symbol of \(\Psi_L\) is $\Omega(\xi)=2d-2\sum_{j=1}^d\cos\xi_j.$ We use the finite exceptional set $\Sigma_d^\Psi \coloneqq \{4j:0\leq j\leq d\}\cup\{2d\}.$ The additional value \(2d\) is convenient for the fixed-block geometric argument; we always impose separately that the rescaled energy lies in \((0,4d)\).  Define
\[
p_d =
\begin{cases}
6, & d \in \{2,4\},\\[2pt]
\dfrac{14}{3}, & d = 3,\\[4pt]
\dfrac{2d}{d-3}, & d > 4,
\end{cases}
\qquad
p_d'=\frac{p_d}{p_d-1}.
\]

For a finite-dimensional Hilbert space \(\mathcal H\), the Hilbert--Schmidt inner product and norm are
\[
        \langle A,B\rangle_{\mathrm{HS}}=\Tr(A^*B),
        \qquad
        \|A\|_{\mathrm{HS}}
        =\biggl(\sum_i\|Ae_i\|_{\mathcal H}^2\biggr)^{1/2}.
\]
For each block \(I\), let
\[
        \iota_I:\mathbb C^{W^d}\longrightarrow\mathbb C^{L^d},
        \qquad
        \iota_Ie_u=e_{(I,u)},
        \qquad
        \Pi_{[I]}=\iota_I\iota_I^*.
\]
For \(A\in\mathbb C^{L^d\times L^d}\), define
\[
        \blk_{I,J}[A]\coloneqq\iota_I^*A\iota_J,
        \qquad
        \blk_I[A]\coloneqq\blk_{I,I}[A].
\]
Thus
\[
        \Pi_{[I]}A\Pi_{[J]}
        =
        \iota_I\blk_{I,J}[A]\iota_J^*.
\]
Transpose is denoted by \({}^\top\), adjoint by \({}^*\), and \(\preceq\) denotes the Loewner order on Hermitian matrices.
\section{Free comparison and a priori resolvent bounds}
\label{sec:apriori-resolvent}

This section supplies the only probabilistic input used before the local-law bootstrap: a coarse \(L^p\)-resolvent bound.  The argument compares smoothed spectral projections of \(H_L\) with those of the free operator \(\lambda\Psi_L\). 

\subsection{Smoothed spectral projections}

Fix a real function \(\rho\in C_c^\infty(\mathbb R)\) with \(\rho=1\) on \([-1/2,1/2]\) and \(\rho=0\) outside \([-2,2]\), and set
\[
        \rho_{E,\delta}(x)
        \coloneqq
        \rho\!\left(\frac{x-E}{\delta}\right).
\]

\begin{lemma}
\label{lemma:smooth_spect_diff}
Uniformly for \(\delta\geq\lambda^{-10}\),
\begin{equation}
\label{eq:smooth-spectral-comparison}
        \|\rho_{E,\delta}(H_L)
        -\rho_{E,\delta}(\lambda\Psi_L)\|_{2\to2}
        \prec
        (\lambda\delta)^{-1/2}.
\end{equation}
\end{lemma}

\begin{proof}
Write the free propagator \(U(s)=e^{-\mathrm is\Psi_L}\) and \(\mathcal V=\sum_\alpha g_\alpha A_\alpha\), where the \(g_\alpha\) are independent standard real Gaussians and the \(A_\alpha\) are the normalized real symmetric block-GOE coordinate matrices.  Thus \(\sum_\alpha A_\alpha B A_\alpha=S[B]\).  For \(\tau\geq0\), set
\[
 T_m(\tau)=\int_{0\leq s_1\leq\cdots\leq s_m\leq\tau}
 U(s_m)^*\mathcal VU(s_m)\cdots
 U(s_1)^*\mathcal VU(s_1)\,\rmd s_1\cdots\rmd s_m.
\]
We first check that the estimate on \(T_1\) is uniform in \(W\). Translation invariance and symmetry of \(U(r)\) give the exact identity
\begin{equation}
\label{eq:BA-Dyson-covariance}
 S[U(r)]=U(r)_{00}\I+\frac1{W^d}
       \sum_I\Pi_{[I]}U(r)\Pi_{[I]}.
\end{equation}
For \(|r|\leq\lambda^{13}\) and sufficiently large \(C_{\mathrm{vol}}\), periodization of the infinite-volume propagator has error \(O(\lambda^{-A})\) for every fixed \(A\), uniformly at microscopic distances at most \(C_dW\) whenever \(|r|\geq C_dW\). Indeed, the nonzero periodization terms lie outside the propagation cone, and repeated integration by parts bounds their sum by an arbitrary power of \(|r|/L\).  The product representation of the lattice propagator and one-dimensional stationary phase therefore give
\[
 |U(r)_{00}|\leq C(1+|r|)^{-d/2}+O(\lambda^{-A}),\qquad
 |U(r)_{xy}|\leq C(1+|r|)^{-d/2}+O(\lambda^{-A})
\]
for \(x,y\) in the same block and \(|r|\geq C_dW\).  Here \(C_d\) is chosen so that every stationary point is uniformly nondegenerate. See \citep[Lemma~2.2]{black2025tailboundsdysonseries} for more details on local dispersive estimates. For smaller \(|r|\), compression of a unitary operator has norm at most one.  Consequently, the Schur test on each block gives
\begin{equation}
\label{eq:BA-Dyson-integrated-covariance}
\begin{aligned}
 \|S[U(r)]\|_{2\to2}
 &\leq C(1+|r|)^{-d/2}
       +W^{-d}\mathbf1_{\{|r|\leq C_dW\}}
       +O(\lambda^{-A}),\\
 \int_0^\tau\|S[U(r)]\|_{2\to2}\,\rmd r
 &\leq C\begin{cases}
             \log(2+\tau),&d=2,\\
             1,&d\geq3,
          \end{cases}
 \qquad 0\leq\tau\leq\lambda^{13}.
\end{aligned}
\end{equation}
The additional block contribution is bounded by \(C_dW^{1-d}\leq C_d\); this is where the GOE normalization is used.

Write \(T_1(\tau)=\sum_\alpha g_\alpha B_\alpha(\tau)\), with \(B_\alpha(\tau)=\int_0^\tau U(s)^*A_\alpha U(s)\,\rmd s\). The matrices \(B_\alpha\) are Hermitian, and
\[
 \sum_\alpha B_\alpha(\tau)^2
 =\int_0^\tau\!\int_0^\tau
       U(s)^*S[U(s-t)]U(t)\,\rmd s\,\rmd t.
\]
It follows from \cref{eq:BA-Dyson-integrated-covariance} that
\[
 \left\|\sum_\alpha B_\alpha(\tau)^2\right\|_{2\to2}
 \leq C\begin{cases}
             \tau\log(2+\tau),&d=2,\\
             \tau,&d\geq3.
          \end{cases}
\]
Matrix Gaussian concentration \citep[Theorem~A.1]{black2025tailboundsdysonseries} now applies directly in dimension \(L^d\).  Its dimensional loss is only \(\sqrt{\log(2L)}\).  In particular, \(\|\mathcal V\|_{2\to2}\prec1\), since \(S[\I]=(1+W^{-d})\I\), and
\[
 \|T_1(\tau)\|_{2\to2}\prec
 \begin{cases}
      \sqrt{\tau\log(2+\tau)},&d=2,\\
      \sqrt\tau,&d\geq3.
 \end{cases}
\]
The bound holds simultaneously for \(0\leq\tau\leq\lambda^{13}\): apply the tail estimate at dyadic times, use \(\|T_1(a+b)\|\leq\|T_1(a)\|+\|T_1(b)\|\), and use \(\|T_1(s)\|\leq s\|\mathcal V\|\) below time one.

The passage from \(T_1\) to \(T_m\) in \citep[Lemma~4.1]{black2025tailboundsdysonseries} is deterministic. It uses only the simplex decomposition, unitary conjugation, and \(\|T_m(\tau)\|\leq\tau^m\|\mathcal V\|^m/m!\), so it applies unchanged to the present matrices.  More explicitly, with \(M_\tau=\sup_{0<s\leq\tau}s^{-1/2}\|T_1(s)\|\), it gives
\[
 \|T_m(\tau)\|
 \leq\left[\frac{C}{\sqrt m}
       (\|\mathcal V\|+M_\tau)
       \sqrt{\tau\log(2+\tau)}\right]^m.
\]
Summing the Duhamel series with \(\tau=\lambda|t|\) and coupling \(\lambda^{-1}\), and using unitarity when this bound is larger than one, yields
\[
 \|e^{-\mathrm itH_L}-e^{-\mathrm it\lambda\Psi_L}\|_{2\to2}
 \prec\min\{1,(|t|/\lambda)^{1/2}\},
 \qquad |t|\leq\lambda^{12}.
\]
Finally, Fourier inversion gives
\[
 \rho_{E,\delta}(H_L)-\rho_{E,\delta}(\lambda\Psi_L)
 =\delta\int_{\mathbb R}e^{-\mathrm itE}
   \bigl(e^{\mathrm itH_L}-e^{\mathrm it\lambda\Psi_L}\bigr)
   \widehat\rho(\delta t)\,\rmd t.
\]
On \(|t|\leq\lambda^{12}\), substitution \(u=\delta t\) gives \((\lambda\delta)^{-1/2}\).  On the complement, unitarity and Schwartz decay give \(C_N(\delta\lambda^{12})^{-N}\), which is smaller than the same bound for \(\delta\geq\lambda^{-10}\) and sufficiently large \(N\).  This proves \cref{eq:smooth-spectral-comparison} uniformly in the block size.
\end{proof}

\subsection{The a priori bound}

\begin{lemma}
\label{lemma:resolv_prior_pre}
Let \(X\) be any norm on functions on \(\mathbb Z_L^d\).  For \(z=E+\mathrm i\eta\) with \(\eta\geq\lambda^{-10}\),
\begin{equation}
\label{eq:half-resolvent-comparison}
        \|(H_L-z)^{-1/2}\|_{2\to X}
        \prec
        \bigl(1+(\lambda\eta)^{-1/2}\bigr)
        \|(\lambda\Psi_L-z)^{-1/2}\|_{2\to X}.
\end{equation}
The analogous estimate holds for the adjoints.
\end{lemma}

\begin{proof}
Choose a smooth dyadic partition of unity in the spectral variable around \(E\), with shell widths \(2^j\eta\).  On overlapping shells the claim follows from functional calculus.  On nonoverlapping shells, insert and subtract the corresponding free spectral cutoff. Lemma~\ref{lemma:smooth_spect_diff}, with \(\delta=2^{\max\{j,\ell\}}\eta\), controls this difference.  Summing the \(O(\log(\lambda/\eta))\) shells gives \cref{eq:half-resolvent-comparison}. The counterpart proof in vanilla Anderson model could be seen in \citep[Lemma~3.3]{black2025selfconsistentequationsquantumdiffusion}. 
\end{proof}

\begin{proposition}
\label{prop:tomas_stein}
Let \(q>p_d\) and assume $e=E/\lambda\in(0,4d), \operatorname{dist}(e,\Sigma_d^\Psi)\geq\varepsilon$, then
\begin{equation}
\label{eq:unif_sobolev}
\begin{aligned}
        \|(\lambda\Psi_L-z)^{-1/2}\|_{2\to q}
        +
        \|(\lambda\Psi_L-z)^{-1/2}\|_{q'\to2}
        \leq
        C_{\varepsilon,q}\lambda^{-1/2}
        \bigl(1+|\log(\eta/\lambda)|\bigr).
\end{aligned}
\end{equation}
\end{proposition}

\begin{proof}
This is the finite-torus Stein--Tomas estimate for the regular level sets of the lattice dispersion; see \citep{erdos2006quantumdiffusionrandomschrodinger,cuenin2021fouriertransformsurfacecarriedmeasures}. After factoring out $\lambda^{-1/2}$, decompose the multiplier into dyadic shells around the level set $\Omega=e$.  The restriction estimate gives a uniform $2\to q$ contribution on each shell, and summing the $O(1+|\log(\eta/\lambda)|)$ shells proves the first bound.  The second follows by duality.
\end{proof}

\begin{proposition}[A priori resolvent bounds]
\label{prop:prior_resol_bd}
Let \(a\leq p_d'\), \(b\geq p_d\), and assume
\[
        \eta\geq\lambda^{-10},
        \qquad
        e=E/\lambda\in(0,4d),
        \qquad
        \operatorname{dist}(e,\Sigma_d^\Psi)\geq\varepsilon.
\]
Then
\begin{equation}
\label{eq:apriori-pq-resolvent}
        \|R(z)\|_{a\to b}
        \prec
        \lambda^{-1}\bigl(1+(\lambda\eta)^{-1}\bigr).
\end{equation}
Moreover,
\begin{equation}
\label{eq:apriori-half-resolvent}
        \|R(z)\|_{2\to b}
        +
        \|R(z)\|_{a\to2}
        \prec
        \eta^{-1/2}\lambda^{-1/2}
        \bigl(1+(\lambda\eta)^{-1/2}\bigr).
\end{equation}
\end{proposition}

\begin{proof}
Apply Lemma~\ref{lemma:resolv_prior_pre} and Proposition~\ref{prop:tomas_stein} to the two half resolvents, using duality for the left factor.  Their product gives
\[
\begin{aligned}
        \|R(z)\|_{a\to b}
        &\leq
        \|R(z)^{1/2}\|_{a\to2}
        \|R(z)^{1/2}\|_{2\to b}\\
        &\prec
        \lambda^{-1}
        \bigl(1+(\lambda\eta)^{-1/2}\bigr)^2,
\end{aligned}
\]
which is \cref{eq:apriori-pq-resolvent}.  Combining one such factor with the trivial spectral bound \(\|R(z)^{1/2}\|_{2\to2}\leq\eta^{-1/2}\) proves \cref{eq:apriori-half-resolvent}.
\end{proof}
\section{The self-consistent equation}
\label{sec:self-consistent-equation}

Gaussian integration by parts separates the averaged resolvent into a deterministic self-consistent term and a fluctuation error.  This is the point at which the covariance map \(S\) enters.

For any deterministic block-diagonal matrix \(A\), set $M_A\coloneqq(\lambda\Psi_L-z-A)^{-1}.$ The resolvent identity is
\begin{equation}
\label{equ:resol_indentity}
        R
        =
        M_A-M_A(\mathcal V+A)R.
\end{equation}
Writing \(\mathcal V=\sum_\alpha g_\alpha A_\alpha\), we have
\[
        \partial_{g_\alpha}R=-RA_\alpha R,
        \qquad
        S[B]=\sum_\alpha A_\alpha BA_\alpha.
\]
Gaussian integration by parts in \cref{equ:resol_indentity} therefore gives the exact averaged equation
\begin{equation}
\label{eq:averaged-resolvent-exact}
        \mathbb ER
        =
        M_A+
        M_A\,\mathbb E\!\left[(S[R]-A)R\right].
\end{equation}

Choose
\[
        \widetilde\Theta\coloneqq S[\mathbb ER],
        \qquad
        \widetilde M\coloneqq
        (\lambda\Psi_L-z-\widetilde\Theta)^{-1}.
\]
Then the centered part of \(S[R]\) is the only error:
\begin{equation}
\label{eq:averaged-resolvent-fluctuation-identity}
        \mathbb ER-\widetilde M
        =
        \widetilde M\,
        \mathbb E\!\left[S[R-\mathbb ER]R\right].
\end{equation}
Applying \(S\) gives the approximate self-consistent equation
\begin{equation}
\label{equ:app_self_consist_eq}
\begin{aligned}
        \widetilde\Theta
        =
        S[\widetilde M]
        +
        S\!\left[
        \widetilde M\,
        \mathbb E\!\left[S[R-\mathbb ER]R\right]
        \right].
\end{aligned}
\end{equation}

The deterministic comparison object is the solution of
\begin{equation}
\label{equ:self_consist_eq}
\begin{aligned}
        \Theta
        &=
        S[M],
        \qquad
        M=(\lambda\Psi_L-z-\Theta)^{-1},
\end{aligned}
\end{equation}
which we call the self-consistent equation. Here note that both self-energies are block-translation invariant.  We write
\[
        \widetilde\Theta=\bigoplus_I\widetilde\theta,
        \qquad
        \Theta=\bigoplus_I\theta
\]
for their common cell blocks. Thus the one-resolvent local law reduces to two tasks: control the fluctuation term in \cref{equ:app_self_consist_eq}, and prove stability of the fixed-point map
\[
        A\longmapsto
        S[(\lambda\Psi_L-z-A)^{-1}].
\]
The next section initializes this equation in the kinetic range; the sharper deterministic stability estimates needed at diffusive scales are proved in Sections~\ref{sec:large-block-stability} and \ref{sec:fixed-block-stability}.
\section{Kinetic initialization}
\label{sec:kinetic-initialization}

In the kinetic range the deterministic equation can be initialized by a direct contraction argument.  The free resolvent has columns of size \(O((\lambda\eta)^{-1/2})\); consequently, the derivative of the self-consistent map is \(O((\lambda\eta)^{-1})\).

For every matrix \(B\), the block-GOE covariance gives
\[
        S[B]\big|_{[b]}
        =W^{-d}\bigl(\Pi_{[b]}B\Pi_{[b]}\bigr)^\top
        +W^{-d}\Tr_{[b]}(B)\I_{W^d}.
\]
Schur's test and the trace bound therefore imply
\begin{equation}
\label{equ:sup_op_ineq}
        \|S[B]\|_{2\to2}
        \leq
        2\sup_{b\in\widetilde{\mathbb Z}_n^d}
        \|\Pi_{[b]}B\Pi_{[b]}\|_{1\to\infty}.
\end{equation}

Let \(R_0=(\lambda\Psi_L-z)^{-1}\).  The free local density-of-states estimate, Ward's identity, and the finite-volume comparison give
\begin{equation}
\label{eq:kinetic-comparison-free-columns}
\begin{aligned}
        \|R_0\|_{1\to\infty}\leq\frac{C_0}{\lambda},\quad
        \sup_x\bigl(
        \|R_0|x\rangle\|_2+
        \|R_0^\star|x\rangle\|_2
        \bigr)\leq C_0(\lambda\eta)^{-1/2}.
\end{aligned}
\end{equation}

\begin{lemma}
\label{lemma:contra_self_map_kinetic}
There are constants \(C_K,C<\infty\) such that, whenever \(\lambda\eta\geq C_K\), the map
\[
        \Phi(A)
        \coloneqq
        S[(\lambda\Psi_L-z-A)^{-1}]
\]
is a strict contraction on
\[
        \mathcal K
        \coloneqq
        \left\{
        A:\ A\text{ is block diagonal},\quad
        \Im A\succeq0,\quad
        \|A\|_{2\to2}\leq\frac{C}{\lambda}
        \right\}.
\]
In particular, \cref{equ:self_consist_eq} has a unique fixed point in \(\mathcal K\), and this fixed point satisfies \(\|\Theta\|_{2\to2}\leq C/\lambda\).
\end{lemma}

\begin{proof}
For \(A\in\mathcal K\), denote
\[
        M_A=(\lambda\Psi_L-z-A)^{-1}.
\]
Choose the radius of \(\mathcal K\) to be \(4C_0/\lambda\), and then choose \(C_K\) so large that \(\|R_0A\|_{2\to2}\leq C/(\lambda\eta)\leq1/2\). The Neumann expansion
\[
        M_A=(\Id-R_0A)^{-1}R_0
\]
and \cref{eq:kinetic-comparison-free-columns} yield
\begin{equation}
\label{eq:kinetic-map-column-bound}
        \sup_x\bigl(
        \|M_A|x\rangle\|_2+
        \|M_A^\star|x\rangle\|_2
        \bigr)
        \leq C(\lambda\eta)^{-1/2}.
\end{equation}
The resolvent identity also gives
\[
\begin{aligned}
        \|M_A\|_{1\to\infty}
        &\leq
        \|R_0\|_{1\to\infty}
        +\|R_0AM_A\|_{1\to\infty}\\
        &\leq
        \frac{C_0}{\lambda}
        +\frac{C}{\lambda\eta}\|A\|_{2\to2}
        \leq\frac{2C_0}{\lambda},
\end{aligned}
\]
after increasing \(C_K\).  Hence \cref{equ:sup_op_ineq} shows that \(\Phi(\mathcal K)\subseteq\mathcal K\).  The map preserves the other defining properties of \(\mathcal K\): \(S\) has block-diagonal range, and
\[
        \Im M_A=M_A(\eta\I+\Im A)M_A^\star\succeq0
\]
implies \(\Im\Phi(A)\succeq0\).

For \(A,B\in\mathcal K\), the resolvent identity and \cref{equ:sup_op_ineq,eq:kinetic-map-column-bound} give
\[
\begin{aligned}
        \|\Phi(A)-\Phi(B)\|_{2\to2}
        &\leq
        2\|M_A(A-B)M_B\|_{1\to\infty}\\
        &\leq
        \frac{C}{\lambda\eta}\|A-B\|_{2\to2}.
\end{aligned}
\]
Increasing \(C_K\) once more makes the last coefficient at most \(1/2\).  Banach's fixed-point theorem proves the claim.
\end{proof}
\begin{lemma}[The MDE solution in the kinetic range]
\label{lemma:theta_22_est_kinet}
For every \(z\in\mathbb C^+\), \cref{equ:self_consist_eq} has a unique admissible solution, characterized by \(\Im M\succeq0\).  In the kinetic range \(\lambda\eta\geq C_K\), this solution is the fixed point from Lemma~\ref{lemma:contra_self_map_kinetic}; in particular,
\[
        \|\Theta\|_{2\to2}\leq\frac{C}{\lambda},
        \qquad
        \Im\Theta\succeq0.
\]
\end{lemma}

\begin{proof}
The covariance map \(S\) preserves positivity and is self-adjoint for the Hilbert--Schmidt inner product:
\[
        \Tr(A^\star S[B])=\Tr(S[A]^\star B).
\]
The standard existence and uniqueness theorem for the operator-valued Dyson equation therefore gives a unique admissible solution; see \citet{helton2007operatorvaluedsemicircularelementssolving} and \citet{erdos2019matrixdysonequationapplications}.  When \(\lambda\eta\geq C_K\), the contraction above produces an admissible solution in \(\mathcal K\).  Uniqueness identifies the two solutions and gives the stated estimate.
\end{proof}

The same contraction controls the difference between the exact self-energy \(\widetilde\Theta=S[\mathbb ER]\) and its deterministic counterpart.  This initializes the continuity bootstrap for the local law.

\begin{lemma}
\label{lemma:kinetic-M-tildeM-comparison}
Suppose that \(\lambda\eta\geq C_K\) and, for a sufficiently small constant \(c_0>0\), $\|\mathbb ER-M\|_{1\to\infty} \leq{c_0/\lambda}.$ Then
\begin{equation}
\label{eq:kinetic-M-tildeM-comparison}
        \|\widetilde M-M\|_{1\to\infty}
        \leq
        C\|\mathbb ER-\widetilde M\|_{1\to\infty}.
\end{equation}
More precisely,
\[
        \|\widetilde M-M\|_{1\to\infty}
        \leq
        \frac{C}{\lambda\eta}
        \|\mathbb ER-\widetilde M\|_{1\to\infty}.
\]
\end{lemma}

\begin{proof}
By \cref{equ:sup_op_ineq},
\begin{equation}
\label{eq:kinetic-tilde-theta-upper}
\begin{aligned}
        \|\widetilde\Theta-\Theta\|_{2\to2}
        &=\|S[\mathbb ER-M]\|_{2\to2}\\
        &\leq2\|\mathbb ER-M\|_{1\to\infty}
        \leq\frac{2c_0}{\lambda}.
\end{aligned}
\end{equation}
Moreover, \(\Im\widetilde\Theta=S[\Im\mathbb ER]\succeq0\).  Thus \(\Theta\) and \(\widetilde\Theta\) belong to a common admissible ball of radius \(C/\lambda\).  The Neumann argument in \cref{eq:kinetic-map-column-bound} applies to both resolvents and gives
\begin{equation}
\label{eq:kinetic-comparison-columns}
        \sup_x\bigl(
        \|M|x\rangle\|_2+\|M^\star|x\rangle\|_2
        +\|\widetilde M|x\rangle\|_2
        +\|\widetilde M^\star|x\rangle\|_2
        \bigr)
        \leq C(\lambda\eta)^{-1/2},
\end{equation}

Since \(\widetilde\Theta=S[\mathbb ER]\), \(\Theta=S[M]\), and \(\Phi(A)=S[(\lambda\Psi_L-z-A)^{-1}]\), we have the exact identity
\[
        \widetilde\Theta-\Theta
        =S[\mathbb ER-\widetilde M]
        +\Phi(\widetilde\Theta)-\Phi(\Theta).
\]
The contraction estimate, with the difference factor maintained, gives
\[
\begin{aligned}
        \|\widetilde\Theta-\Theta\|_{2\to2}
        \leq{}&
        2\|\mathbb ER-\widetilde M\|_{1\to\infty}\\
        &+\frac{C}{\lambda\eta}
        \|\widetilde\Theta-\Theta\|_{2\to2}.
\end{aligned}
\]
After increasing \(C_K\), the final term is absorbed.  The resolvent identity and \cref{eq:kinetic-comparison-columns} then yield
\[
\begin{aligned}
        \|\widetilde M-M\|_{1\to\infty}
        &\leq
        \|\widetilde M\|_{1\to2}
        \|\widetilde\Theta-\Theta\|_{2\to2}
        \|M^\star\|_{1\to2}\\
        &\leq
        \frac{C}{\lambda\eta}
        \|\mathbb ER-\widetilde M\|_{1\to\infty}.
\end{aligned}
\]
This proves both assertions.
\end{proof}

\section{Local law}
\label{sec:local-law}

Set
\[
        \widetilde\Theta\coloneqq S[\mathbb ER],
        \qquad
        \widetilde M\coloneqq
        (\lambda\Psi_L-z-\widetilde\Theta)^{-1},
\]
and let
\[
        M=(\lambda\Psi_L-z-\Theta)^{-1},
        \qquad
        \Theta=S[M],
\]
be the deterministic MDE resolvent.  The local-law proof has two inputs. Gaussian concentration controls \(R-\mathbb ER\) in terms of an interpolated resolvent norm.  Deterministic stability then controls \(\mathbb ER-M\) by that fluctuation.  These two estimates form a closed continuity bootstrap in \(\eta\).

The probabilistic part does not distinguish the two block-size regimes.  The difference lies entirely in the deterministic stability input: normalized trace scalarizes the MDE for growing blocks, whereas a fibered finite-dimensional MDE is used for fixed blocks.  We record the common consequence next and defer its two proofs to Sections~\ref{sec:large-block-stability} and \ref{sec:fixed-block-stability}.

\begin{proposition}
\label{prop:deterministic-local-law-inputs}
Let \(z=E+\mathrm i\eta\) lie in one of the two bulk regimes from Theorem~\ref{thm:local_law_block_anderson}, and suppose
\[
        \|\mathbb ER-M\|_{1\to\infty}\leq\frac{c_0}{\lambda}
\]
for a sufficiently small constant \(c_0\).  Then
\begin{equation}
\label{eq:local-law-deterministic-columns}
        \sup_x\left(
        \|\widetilde M e_x\|_2+\|\widetilde M^*e_x\|_2
        \right)\leq C
\end{equation}
and
\begin{equation}
\label{eq:local-law-deterministic-comparison}
        \|\widetilde M-M\|_{1\to\infty}
        \leq
        C\|\mathbb ER-\widetilde M\|_{1\to\infty}
        +C\chi_L(\eta),
\end{equation}
where
\[
        \chi_L(\eta)
        \coloneqq
        \begin{cases}
        0,&\text{in the growing-block regime},\\[2pt]
        \eta^{-1}\exp(-cn\eta/(d\lambda)),
            &\text{in the fixed-block regime}.
        \end{cases}
\]
Moreover,
\begin{equation}
\label{eq:local-law-deterministic-max}
        \|M\|_{1\to\infty}
        \leq
        \frac{C(1+\log\lambda)^C}{\lambda}
        +C\chi_L(\eta).
\end{equation}
The constants are uniform in the growing-block regime and may depend on \(W\) and on the fixed folded bulk window in the fixed-block regime.
\end{proposition}

\begin{proof}
In the growing-block regime the explicit covariance map gives
\[
        \|\widetilde\theta_W-\theta_W\|_{2\to2}
        \leq2\|\mathbb ER-M\|_{1\to\infty}
        \leq\frac{2c_0}{\lambda}.
\]
Choose \(c_0\leq c_*/2\), where \(c_*\) is the fixed constant in Lemma~\ref{lemma:large_block_stability}.  The column bounds then follow from Corollary~\ref{coro:robust_app_bootstrap}, the comparison follows from Corollary~\ref{cor:tildeM_M_compar_large_block}, and the bound on \(M\) is Theorem~\ref{thm:main_closeness}.  In the fixed-block regime they follow from the fibered MDE trap, the one-point stability estimate, and the finite-volume comparison in Theorem~\ref{thm:finite_volume_from_infinite_trap} and Corollary~\ref{cor:qfv-unscaled-random-application}.  The Combes--Thomas comparison produces precisely \(\chi_L\).
\end{proof}

We defer the proof of Theorem~\ref{thm:local_law_block_anderson} to Subsections~\ref{subsec:local-law-bootstrap-d3}-\ref{subsec:local-law-bootstrap-highd}. We first use the following lemma to control the fluctuations of the resolvent and of the resolvent sandwich with a deterministic test operator.
\begin{lemma}
\label{lemma:concent_1_4_bdd}
For \(z\in\mathbb C^+\) and \(q\geq2\),
\begin{equation}
\label{eq:one-resolvent-concentration}
 \left(\mathbb E\|R-\mathbb ER\|_{1\to\infty}^{q}\right)^{1/q}
 \leq C\sqrt{q+\log(2L)}
       \left(\mathbb E\|R\|_{1\to4}^{2q}\right)^{1/q}.
\end{equation}
For every deterministic \(\mathcal A\),
\begin{equation}
\label{eq:two-resolvent-concentration}
\begin{aligned}
 &\left(\mathbb E
   \|R\mathcal AR^*-\mathbb E(R\mathcal AR^*)\|_{1\to\infty}^{q}
       \right)^{1/q}\\
 &\qquad\leq C\sqrt{q+\log(2L)}\,\|\mathcal A\|_{2\to2}
 \left[\mathbb E\left(
      \|R\|_{1\to4}\|R\|_{1\to2}\|R\|_{2\to4}
                  \right)^q\right]^{1/q}.
\end{aligned}
\end{equation}
The constant is independent of \(W,L,z,q\).  In the polynomial spectral ranges considered below, a deterministic stochastic upper bound for the product on the right gives the same stochastic upper bound for the fluctuation on the left, provided its reciprocal grows at most polynomially in \(\lambda\).
\end{lemma}

\begin{proof}
We first record a maximal Gaussian moment inequality, including its proof to keep track of the dimension.  If \(g\) is a vector of independent standard real Gaussians and \(F_1,\ldots,F_m\) are smooth complex-valued functions, then
\begin{equation}
\label{eq:Gaussian-maximal-moment}
 \left\|\max_{j\leq m}|F_j(g)-\mathbb EF_j(g)|\right\|_{L^q}
 \leq C\sqrt{q+\log(2m)}
       \left\|\max_{j\leq m}|\nabla F_j(g)|\right\|_{L^q},
 \qquad q\geq2.
\end{equation}
Here \(\|X\|_{L^q}=(\mathbb E|X|^q)^{1/q}\), and \(|\nabla F|^2=\sum_\alpha|\partial_\alpha F|^2\). Let \(g'\) be an independent copy of \(g\), and denote \(g_\theta=g\cos\theta+g'\sin\theta\) and \(h_\theta=-g\sin\theta+g'\cos\theta\).  Jensen's inequality, the fundamental theorem of calculus, and Minkowski's inequality give
\[
 \left\|\max_j|F_j(g)-\mathbb EF_j(g)|\right\|_{L^q}
 \leq\int_0^{\pi/2}
       \left\|\max_j
        |\langle h_\theta,\nabla F_j(g_\theta)\rangle|
       \right\|_{L^q}\,\rmd\theta.
\]
For each \(\theta\), \(h_\theta\) is a standard Gaussian vector independent of \(g_\theta\).  Conditional Gaussian tail bounds and a union bound over the \(2m\) real and imaginary parts bound the integrand by the right-hand side of \cref{eq:Gaussian-maximal-moment}.  Integration proves that inequality. In particular, its constant does not depend on the number of disorder coordinates.  Approximation gives the same statement whenever the quantities on the right are integrable.

Write \(\mathcal V=\sum_\alpha g_\alpha A_\alpha\) in normalized block-GOE coordinates.  For arbitrary complex vectors \(a,b\), the covariance identity gives
\[
 \sum_\alpha|a^*A_\alpha b|^2
 =\frac1{W^d}\sum_I
    \left(\|a_I\|_2^2\|b_I\|_2^2+|a_I^\top b_I|^2\right)
 \leq2\|a\|_4^2\|b\|_4^2.
\]
The last inequality follows from Cauchy--Schwarz over the blocks and \(\|a_I\|_2^4\leq W^d\|a_I\|_4^4\).  Since \(\partial_{g_\alpha}R=-RA_\alpha R\), it follows that
\[
 |\nabla R_{xy}|^2
 \leq2\|R^*e_x\|_4^2\|Re_y\|_4^2
 \leq2\|R\|_{1\to4}^4.
\]
We used that \(H_L\) is real symmetric, so \(R^*=\overline R\) and the corresponding operator norms agree.  Differentiating \(R\mathcal AR^*\) gives
\[
 \partial_{g_\alpha}(R\mathcal AR^*)
 =-RA_\alpha R\mathcal AR^*
  -R\mathcal AR^*A_\alpha R^*.
\]
Apply the same covariance inequality to each term and use
\[
 \|R\mathcal AR^*e_y\|_4
 \leq\|R\|_{2\to4}\|\mathcal A\|_{2\to2}\|R\|_{1\to2}.
\]
This gives
\[
 |\nabla(R\mathcal AR^*)_{xy}|^2
 \leq8\|\mathcal A\|_{2\to2}^2
      \|R\|_{1\to4}^2\|R\|_{1\to2}^2\|R\|_{2\to4}^2.
\]
Apply \cref{eq:Gaussian-maximal-moment} to the \(L^{2d}\) entries. This proves both moment bounds, including their uniformity in \(W\). For the final assertion, \(\|R\|_{2\to2}\leq\eta^{-1}\) provides a polynomial bound on the exceptional event.  Thus a deterministic stochastic upper bound gives the corresponding moment bounds with an arbitrary factor \(\lambda^\varepsilon\); Markov's inequality with sufficiently large fixed \(q\) proves the assertion.
\end{proof}

\begin{lemma}
\label{lemma:flut_exp_bdd}
Under the assumptions of Proposition~\ref{prop:deterministic-local-law-inputs},
\begin{equation}
\label{eq:fluctuation-expectation-feedback}
\begin{aligned}
        \|\mathbb ER-M\|_{1\to\infty}
        \leq
        C\eta^{-1/2}
        \left(
        \mathbb E\|R-\mathbb ER\|_{1\to\infty}^2
        \right)^{1/2}
        \left(
        \mathbb E\|R\|_{1\to\infty}
        \right)^{1/2}
        +C\chi_L(\eta).
\end{aligned}
\end{equation}
\end{lemma}

\begin{proof}
Gaussian integration by parts gives the exact identity
\[
        \mathbb ER-\widetilde M
        =
        \widetilde M\,
        \mathbb E\!\left[S[R-\mathbb ER]R\right].
\]
By \cref{eq:local-law-deterministic-columns}, \(\|S[B]\|_{2\to2}\leq2\|B\|_{1\to\infty}\), and Ward's identity,
\[
\begin{aligned}
        \|\mathbb ER-\widetilde M\|_{1\to\infty}
        &\leq
        C
        \left(
        \mathbb E\|R-\mathbb ER\|_{1\to\infty}^2
        \right)^{1/2}
        \left(
        \mathbb E\|R\|_{1\to2}^2
        \right)^{1/2}\\
        &\leq
        C\eta^{-1/2}
        \left(
        \mathbb E\|R-\mathbb ER\|_{1\to\infty}^2
        \right)^{1/2}
        \left(
        \mathbb E\|R\|_{1\to\infty}
        \right)^{1/2}.
\end{aligned}
\]
Combining this estimate with \cref{eq:local-law-deterministic-comparison} and the triangle inequality proves the lemma.
\end{proof}
\subsection{Proof of the local law for \texorpdfstring{$d=3$}{d=3}}
\label{subsec:local-law-bootstrap-d3}

Set
\[
        \eta_*\coloneqq\lambda^{-1-6/19+\delta},
        \qquad 0<\delta<6/19.
\]
The bootstrap itself is identical in the two block-size regimes.  Their only effect is through Proposition~\ref{prop:deterministic-local-law-inputs}; the finite-volume error \(\chi_L\) is present only when \(W\) is fixed.  In particular,
\begin{equation}
\label{eq:d3-deterministic-M-max}
        \|M(z)\|_{1\to\infty}
        \leq
        \frac{C(1+\log\lambda)^C}{\lambda}+C\chi_L(\eta).
\end{equation}

\begin{proof}[Proof of Theorem~\ref{thm:local_law_block_anderson} for \(d=3\)]
We follow the expectation--fluctuation bootstrap of \citet[Theorem~4.1]{black2025selfconsistentequationsquantumdiffusion}.

\emph{Step 1: concentration at the interpolated scale.} Ward's identity and Riesz--Thorin interpolation give
\begin{align*}
    \|R\|_{1\to4}
    &\leq
    \|R\|_{1\to2}^{1/8}
    \|R\|_{1\to14/3}^{7/8}\\
    &\leq
    \eta^{-1/16}
    \|R\|_{1\to\infty}^{1/16}
    \|R\|_{1\to14/3}^{7/8}.
\end{align*}
By Proposition~\ref{prop:prior_resol_bd},
\[
    \|R\|_{1\to14/3}\prec\lambda^{-1}(1+(\lambda\eta)^{-1}).
\]
Since $\eta^{-1/16}\leq\lambda^{1/16}(1+(\lambda\eta)^{-1})^{1/16}$, we obtain
\begin{equation}
\label{eq:d3-one-four-interpolated}
    \|R\|_{1\to4}
    \prec
    \lambda^{-13/16}(1+(\lambda\eta)^{-1})^{15/16}
    \|R\|_{1\to\infty}^{1/16}.
\end{equation}
Apply \cref{eq:one-resolvent-concentration} with \(q=2\). On the event furnished by Proposition~\ref{prop:prior_resol_bd}, \cref{eq:d3-one-four-interpolated} bounds the squared gradient scale by a deterministic multiple of \(\|R\|_{1\to\infty}^{1/8}\). The exceptional event contributes an arbitrarily small power of \(\lambda^{-1}\), since \(\eta^{-1}\) is polynomially bounded. Finally, Jensen's inequality gives \(\mathbb E\|R\|_{1\to\infty}^{1/4} \leq(\mathbb E\|R\|_{1\to\infty})^{1/4}\). Absorbing \(\sqrt{\log(2L)}\) into \(\lambda^{\epsilon_0}\) yields
\begin{equation}
\label{eq:d3-maximal-moment-fluctuation}
     \left(\mathbb E\|R-\mathbb ER\|_{1\to\infty}^2\right)^{1/2}
    \leq
    C_{\epsilon_0}\lambda^{\epsilon_0}
    \lambda^{-13/8}(1+(\lambda\eta)^{-1})^{15/8}(\mathbb E\|R\|_{1\to \infty})^{1/8}.
\end{equation}

\emph{Step 2: the bootstrap implication.} Choose $c_0>0$ sufficiently small for Lemma~\ref{lemma:flut_exp_bdd}, and assume at a fixed $\eta$ that
\begin{equation}
\label{eq:d3-bootstrap-assumption}
    \|\mathbb ER - M\|_{1\to \infty}\leq\frac{c_0}{\lambda}.
\end{equation}
By the triangle inequality, \cref{eq:d3-deterministic-M-max}, and \cref{eq:d3-maximal-moment-fluctuation},
\begin{equation}
\label{eq:d3-U-recursion}
    \mathbb E\|R\|_{1\to\infty}
    \leq
    C\lambda^{\epsilon_0-13/8}(1+(\lambda\eta)^{-1})^{15/8} (\mathbb E\|R\|_{1\to\infty})^{1/8}
    +\frac{C\lambda^{\epsilon_0}}{\lambda}+C\chi_L.
\end{equation}
For $\eta\geq\eta_*$,
\begin{equation}
\label{eq:d3-Q-upper}
    (1+(\lambda\eta)^{-1})\leq2\lambda^{6/19-\delta}.
\end{equation}
If $A=C\lambda^{\epsilon_0-13/8}(1+(\lambda\eta)^{-1})^{15/8}$, then
\[
    A^{8/7}
    \leq
    C\lambda^{-157/133+8\epsilon_0/7-15\delta/7}
    =o(\lambda^{-1})
\]
after choosing $\epsilon_0$ sufficiently small.  The elementary implication $u\leq Au^{1/8}+b\Rightarrow u\leq C(A^{8/7}+b)$ therefore yields
\begin{equation}
\label{eq:d3-U-closed}
    \mathbb E\|R\|_{1\to \infty}\leq C\lambda^{-1+2\epsilon_0}.
\end{equation}
Here and below logarithmic factors and the Combes--Thomas error are absorbed by $\lambda^{\epsilon_0}$. Substitution in \cref{eq:d3-maximal-moment-fluctuation} gives
\begin{equation}
\label{eq:d3-F-closed}
    \left(\mathbb E\|R- \mathbb ER\|_{1\to \infty}^2\right)^{1/2}
    \leq
    C\lambda^{-7/4+2\epsilon_0}(1+(\lambda\eta)^{-1})^{15/8}.
\end{equation}

Apply Lemma~\ref{lemma:flut_exp_bdd}.
\begin{align*}
    \|\mathbb ER - M\|_{1\to\infty}
    &\leq
    C\eta^{-1/2}
    \left(\mathbb E\|R-\mathbb ER\|_{1\to\infty}^2\right)^{1/2}
    (\mathbb E\|R\|_{1\to\infty})^{1/2}
    +C\chi_L(\eta)\\
    &\leq
    C\lambda^{-7/4+3\epsilon_0}(1+(\lambda\eta)^{-1})^{19/8}
    +C\chi_L(\eta),
\end{align*}
where we used
\begin{equation}
\label{eq:d3-eta-to-Q}
    \eta^{-1/2}\leq\lambda^{1/2}
    \bigl(1+(\lambda\eta)^{-1}\bigr)^{1/2}.
\end{equation}
Taking $3\epsilon_0<\delta$ proves the implication
\begin{equation}
\label{eq:d3-bootstrap-improvement}
    \cref{eq:d3-bootstrap-assumption}
    \quad\Longrightarrow\quad
    \|\mathbb ER - M\|_{1\to \infty}
    \leq
    C_\delta\lambda^{-7/4+\delta}(1+(\lambda\eta)^{-1})^{19/8}
    +C\chi_L(\eta).
\end{equation}
At the lower endpoint,
\begin{equation}
\label{eq:d3-bottom-gain}
    \lambda^{-7/4+\delta}(1+(\lambda\eta_*)^{-1})^{19/8}
    \leq C\lambda^{-1-11\delta/8}=o(\lambda^{-1}).
\end{equation}
This is a strict improvement of \cref{eq:d3-bootstrap-assumption}, uniformly for $\eta\geq\eta_*$.

\emph{Step 3: continuity.} Let $\eta_0=K\lambda$, where $K>4/c_0$ is fixed.  Positivity of the imaginary parts of the random resolvent and of the MDE solution gives
\[
    \|\mathbb ER(E+i\eta_0)\|_{1\to\infty}
    +\|M(E+i\eta_0)\|_{1\to\infty}
    \leq\frac{2}{K\lambda}<\frac{c_0}{\lambda}.
\]
This starts the bootstrap.  Lemma~\ref{lemma:flut_exp_bdd} is uniform in the kinetic and diffusive ranges, and the function \(\eta\mapsto\|\mathbb ER(E+\mathrm i\eta)-M(E+\mathrm i\eta)\|_{1\to\infty}\) is continuous.  If, while decreasing \(\eta\) from \(\eta_0\) to \(\eta_*\), it first reached \(c_0/\lambda\), then \cref{eq:d3-bootstrap-improvement} and \cref{eq:d3-bottom-gain} would bound it by $c_0/(2\lambda)$ there, a contradiction.  Hence \cref{eq:d3-bootstrap-assumption} holds throughout the continuation interval, and \cref{eq:d3-expectation-local-law} follows.

\emph{Step 4: stochastic domination.} The maximal moment inequality also gives the required high-probability bound.  Fix \(q\geq2\) and apply \cref{eq:one-resolvent-concentration} with this \(q\), truncating the a priori estimate as in Step~1, and using \(\|\|R\|_{1\to\infty}^{1/8}\|_{L^q}\leq (\mathbb E\|R\|_{1\to\infty}^q)^{1/(8q)}\), gives
\begin{equation}
\label{eq:local_law_R_1_inf_bdd}
 \left(\mathbb E\|R-\mathbb ER\|_{1\to\infty}^q\right)^{1/q}
 \leq C_{q,\epsilon_0}\lambda^{\epsilon_0-13/8}
       (1+(\lambda\eta)^{-1})^{15/8}(\mathbb E\|R\|_{1\to\infty}^q)^{1/(8q)}
       +O(\lambda^{-100}).
\end{equation}
The negligible remainder may be given any fixed negative power by choosing the exceptional probability sufficiently small.  By the triangle inequality and the already established expectation bound,
\[
 (\mathbb E\|R\|_{1\to\infty}^q)^{1/q}\leq C_{q,\epsilon_0}\lambda^{\epsilon_0-13/8}
       (1+(\lambda\eta)^{-1})^{15/8}(\mathbb E\|R\|_{1\to\infty}^q)^{1/(8q)}
       +C_{\epsilon_0}\lambda^{-1+\epsilon_0}.
\]
The same scalar calculation as in Step~2 therefore proves, for every fixed \(q\) and every \(\epsilon>0\),
\[
 (\mathbb E\|R\|_{1\to\infty}^q)^{1/q}\leq C_{q,\epsilon}\lambda^{-1+\epsilon},\qquad
 \left(\mathbb E\|R-\mathbb ER\|_{1\to\infty}^q\right)^{1/q}
 \leq C_{q,\epsilon}\lambda^{-7/4+\epsilon}
       (1+(\lambda\eta)^{-1})^{15/8}.
\]
Markov's inequality, first choosing \(\epsilon\) smaller than the requested stochastic loss and then \(q\) sufficiently large, proves \(\|R\|_{1\to\infty}\prec\lambda^{-1}\) and \cref{eq:d3_local_law_stoch}.
\end{proof}

\subsection{Proof of the local law for \texorpdfstring{$d = 2, 4, 5, 6$}{d = 2, 4, 5, 6}}
\label{subsec:local-law-bootstrap-intermediate}

 Denote the stopping point $\eta_*\coloneqq\lambda^{-1-2/9+\delta}.$ The a priori bound and Ward's identity give
\[
        \|R\|_{1\to6}\prec\lambda^{-1}(1+(\lambda\eta)^{-1}),
        \qquad
        \|R\|_{1\to2}
        \leq
        \eta^{-1/2}\|R\|_{1\to\infty}^{1/2}.
\]
Riesz--Thorin interpolation therefore yields
\begin{equation}
\label{eq:intermediate-one-four}
\begin{aligned}
        \|R\|_{1\to4}
        &\leq
        \|R\|_{1\to2}^{1/4}
        \|R\|_{1\to6}^{3/4}\\
        &\prec
        \lambda^{-5/8}(1+(\lambda\eta)^{-1})^{7/8}
        \|R\|_{1\to\infty}^{1/8}.
\end{aligned}
\end{equation}

\begin{proof}[Proof of Theorem~\ref{thm:local_law_block_anderson} for
\(d\in\{2,4,5,6\}\)] Assume \(\|\mathbb ER-M\|_{1\to\infty}\leq c_0/\lambda\). For fixed \(q\geq2\), Lemma~\ref{lemma:concent_1_4_bdd} and \cref{eq:intermediate-one-four} give
\[
 \left(\mathbb E\|R-\mathbb ER\|_{1\to\infty}^q\right)^{1/q}
 \leq C_{q,\epsilon_0}\lambda^{\epsilon_0-5/4}
       (1+(\lambda\eta)^{-1})^{7/4}(\mathbb E\|R\|_{1\to\infty}^q)^{1/(4q)}
       +O(\lambda^{-100}).
\]
As in the preceding proof, the exceptional part of the a priori bound is controlled by \(\|R\|_{2\to2}\leq\eta^{-1}\).  The triangle inequality and \cref{eq:local-law-deterministic-max} thus give
\[
 (\mathbb E\|R\|_{1\to\infty}^q)^{1/q}\leq A (\mathbb E\|R\|_{1\to\infty}^q)^{1/(4q)}+C_{\epsilon_0}\lambda^{-1+\epsilon_0},
 \qquad
 A=C_{q,\epsilon_0}\lambda^{\epsilon_0-5/4}
       (1+(\lambda\eta)^{-1})^{7/4}.
\]
Since \(1+(\lambda\eta)^{-1}\leq C\lambda^{2/9-\delta}\),
\[
 A^{4/3}\leq C_{q,\epsilon_0}
      \lambda^{-31/27+4\epsilon_0/3-7\delta/3}
      =o(\lambda^{-1}).
\]
The implication \(u\leq Au^{1/4}+b\Rightarrow u\leq C(A^{4/3}+b)\) proves, for every fixed \(q\) and every \(\epsilon>0\),
\[
 (\mathbb E\|R\|_{1\to\infty}^q)^{1/q}\leq C_{q,\epsilon}\lambda^{-1+\epsilon},\qquad
 \left(\mathbb E\|R-\mathbb ER\|_{1\to\infty}^q\right)^{1/q}
 \leq C_{q,\epsilon}\lambda^{-3/2+\epsilon}
       (1+(\lambda\eta)^{-1})^{7/4}.
\]
Taking \(q=2\) in these inequalities and applying Lemma~\ref{lemma:flut_exp_bdd} yields
\[
 \|\mathbb ER-M\|_{1\to\infty}
 \leq C_\delta\lambda^{-3/2+\delta}
       (1+(\lambda\eta)^{-1})^{9/4}+C\chi_L(\eta).
\]
At \(\eta=\eta_*\), the first term is bounded by \(C_\delta\lambda^{-1-5\delta/4}\); the finite-volume term is negligible by the standing volume assumption.  The continuity argument from \(K\lambda\) down to \(\eta_*\) closes the bootstrap exactly as in the \(d=3\) case.  Finally, the bounds for arbitrary fixed \(q\), followed by Markov's inequality, give
\[
 \|R\|_{1\to\infty}\prec\lambda^{-1},\qquad
 \|R-\mathbb ER\|_{1\to\infty}
 \prec\lambda^{-3/2}(1+(\lambda\eta)^{-1})^{7/4}.
\]
\end{proof}

\subsection{Proof of the local law in high dimensions}
\label{subsec:local-law-bootstrap-highd}
In high dimensions the deterministic comparison already has a bounded \(1\to4\) norm.  Iterating the concentration--feedback map then reaches the optimal fluctuation scale before the continuity argument is closed.

For $d\geq7$ set
\begin{equation}
\label{eq:highd-beta-kappa}
    \beta_d\coloneqq p_d-2=\frac{6}{d-3},
    \qquad
    \kappa_d^{\rm LL}
    \coloneqq
    \frac{1-\beta_d/2}{3/2+\beta_d/4}
    =\frac{2d-12}{3d-6},
\end{equation}
and define the stopping point of the bootstrap by $\eta_*\coloneqq \lambda^{-1-\kappa_d^{\rm LL}+\delta}.$ We also write $\ell\coloneqq 1+\log\lambda$. The following deterministic estimate is the high-dimensional input which replaces the deterministic $1\to\infty$ estimate used when $d=3$, and it is valid in both block-size regimes.

\begin{lemma}
\label{lem:highd-deterministic-one-four}
Let $d\geq7$, and suppose that the deterministic MDE bounds established above hold.  Then
\begin{equation}
\label{eq:highd-deterministic-one-four}
    \|M\|_{1\to4}
    \leq \frac{C\ell^2}{\lambda}.
\end{equation}
\end{lemma}

\begin{proof}
If $\eta\geq\lambda$, positivity of the MDE solution gives
\[
    \|M\|_{1\to4}
    \leq\|M\|_{2\to2}\leq\eta^{-1}\leq\lambda^{-1},
\]
so it remains to consider $\eta\leq\lambda$, which is the range of the uniform Sobolev input used below. Write the block-translation-invariant self-energy as $\theta=g\I+\vartheta$.  The deterministic MDE estimates in both regimes give
\[
    \Im g\geq \frac{c}{\lambda},
    \qquad
    |g|+\|\vartheta\|_{2\to2}\leq\frac{C}{\lambda},
    \qquad
    \|M\|_{1\to2}\leq C\ell.
\]
Put $R_g=(\lambda\Psi_L-z-g)^{-1}$.  Since $p_d<4$ for $d\geq7$, choose $p\in(p_d,4)$.  The uniform Sobolev estimate \cref{eq:unif_sobolev}, the inclusions $\ell^1\subset\ell^{p'}$ and $\ell^p\subset\ell^4$, and the dissipative half-resolvent estimate give
\[
    \|R_g\|_{1\to4}\leq \frac{C\ell^2}{\lambda},
    \qquad
    \|R_g\|_{2\to4}\leq C\ell.
\]
The exact resolvent identity
\[
    M=R_g+R_g\vartheta M
\]
therefore yields
\[
    \|M\|_{1\to4}
    \leq
    \|R_g\|_{1\to4}
    +\|R_g\|_{2\to4}\|\vartheta\|_{2\to2}\|M\|_{1\to2}
    \leq \frac{C\ell^2}{\lambda}.
\]
\end{proof}

\begin{proof}[Proof of Theorem~\ref{thm:local_law_block_anderson} for \(d\geq7\)]
It suffices to consider $0<\delta<\kappa_d^{\rm LL}$, since the remaining cases follow by decreasing $\delta$.  Fix $c_0>0$ small enough that Lemma~\ref{lemma:flut_exp_bdd} applies under
\begin{equation}
\label{eq:highd-bootstrap-assumption}
        \|\mathbb ER-M\|_{1\to\infty}
        \leq\frac{c_0}{\lambda}.
\end{equation}
In this proof write
\[
        Q\coloneqq1+(\lambda\eta)^{-1},
        \qquad \alpha\coloneqq1-\frac\beta2,
\]
where $\beta\in(\beta_d,2)$ is chosen so close to $\beta_d$ that
\begin{equation}
\label{eq:highd-exponent-margin}
        (\beta-\beta_d)
        \left(\frac12+\frac{\kappa_d^{\rm LL}}4\right)
        <\frac\delta{16},
        \qquad
        (2+\beta)\kappa_d^{\rm LL}<2.
\end{equation}
The second inequality is strict at $\beta=\beta_d$.

We first justify the initial scale, including its moment bounds. Set $U=\|R\|_{1\to\infty}$.  For every realization,
\[
        \|R\|_{1\to4}^2
        \leq
        \|R\|_{1\to2+\beta}^{1+\beta/2}U^\alpha.
\]
The a priori bound gives $\|R\|_{1\to2+\beta}\prec\lambda^{-1}Q$. For every fixed $r\geq2$ and $\epsilon>0$, the maximal moment inequality in Lemma~\ref{lemma:concent_1_4_bdd} consequently gives
\[
\begin{aligned}
        \bigl\|\|R-\mathbb ER\|_{1\to\infty}\bigr\|_{L^r}
        \leq
        C_{r,\epsilon}\lambda^\epsilon
        \lambda^{-1-\beta/2}Q^{1+\beta/2}
        \|U\|_{L^r}^\alpha+O(\lambda^{-100}).
\end{aligned}
\]
Here the a priori bound is used on the event $\|R\|_{1\to2+\beta}\leq\lambda^{\epsilon/4}\lambda^{-1}Q$; on its complement the deterministic bound $\|R\|_{2\to2}\leq\eta^{-1}$ and an arbitrarily large probability exponent give the displayed error.  Here we also used $\|U^\alpha\|_{L^r}\leq\|U\|_{L^r}^\alpha$ for $r > 1, \alpha < 1$.

By \cref{eq:highd-bootstrap-assumption} and \cref{eq:local-law-deterministic-max}, the triangle inequality yields
\[
        \|U\|_{L^r}
        \leq C\lambda^{-1}(1+\log\lambda)^C
        +C_{r,\epsilon}\lambda^\epsilon
        \lambda^{-1-\beta/2}Q^{1+\beta/2}
        \|U\|_{L^r}^\alpha+O(\lambda^{-100}).
\]
For $u,a,b\geq0$ and $0<\alpha<1$, $u\leq a+bu^\alpha$ implies $u\leq2a+(2b)^{1/(1-\alpha)}$. Since $Q\leq C\lambda^{\kappa_d^{\rm LL}-\delta}$, the second inequality in \cref{eq:highd-exponent-margin} shows that
\[
        \left(
        \lambda^{-1-\beta/2}Q^{1+\beta/2}
        \right)^{2/\beta}
        =o(\lambda^{-1}).
\]
Choosing the auxiliary loss sufficiently small proves, for every fixed $r$ and every $\epsilon>0$,
\begin{equation}
\label{eq:highd-initial-moments}
        \|U\|_{L^r}\leq C_{r,\epsilon}\lambda^{-1+\epsilon}.
\end{equation}
In particular, $\|R-M\|_{1\to\infty}\prec\lambda^{-1}$, which supplies the initial scale for the following iteration.

Suppose that $\|R-M\|_{1\to\infty}\prec D$ for a deterministic polynomial scale $D$.  On the event $\|R-M\|_{1\to\infty}\leq\lambda^\epsilon D$, split each row at $|M_{xy}|=\lambda^\epsilon D$.  On the first part, $|R_{xy}|\leq2|M_{xy}|$; on the second, $|R_{xy}|\leq2\lambda^\epsilon D$.  Therefore
\[
        \|R\|_{1\to4}^4
        \leq16\|M\|_{1\to4}^4
        +(2\lambda^\epsilon D)^{2-\beta}
          \|R\|_{1\to2+\beta}^{2+\beta}.
\]
Proposition~\ref{prop:prior_resol_bd} and Lemma~\ref{lem:highd-deterministic-one-four} imply
\begin{equation}
\label{eq:highd-coarse-one-four}
        \|R\|_{1\to4}^4
        \prec
        \lambda^{-4}
        +\lambda^{-(2+\beta)}Q^{2+\beta}D^{2-\beta}.
\end{equation}
The polynomial deterministic resolvent bound upgrades this estimate to every fixed moment, with an arbitrarily small power loss. Lemma~\ref{lemma:concent_1_4_bdd} thus gives the deterministic fluctuation scale
\begin{equation}
\label{eq:highd-coarse-F-map}
        F(D)\coloneqq
        \lambda^{-2}
        +\lambda^{-1-\beta/2}Q^{1+\beta/2}D^\alpha,
\end{equation}
meaning that, for every fixed $r\geq2$ and every $\epsilon>0$,
\[
        \bigl\|\|R-\mathbb ER\|_{1\to\infty}\bigr\|_{L^r}
        \leq C_{r,\epsilon}\lambda^\epsilon F(D).
\]
Together with \cref{eq:highd-initial-moments}, Lemma~\ref{lemma:flut_exp_bdd} then gives
\[
        \|\mathbb ER-M\|_{1\to\infty}
        \leq C_\epsilon\lambda^\epsilon Q^{1/2}F(D)+C\chi_L(\eta).
\]
Consequently a valid next scale is
\begin{equation}
\label{eq:highd-total-scale-map}
\begin{aligned}
        \mathcal T(D)\coloneqq
        C\lambda^{\epsilon_0}
        \left[
        \lambda^{-2}Q^{1/2}
        +\lambda^{-1-\beta/2}
          Q^{3/2+\beta/2}D^\alpha
        \right]+C\chi_L(\eta),
\end{aligned}
\end{equation}
where $0<\epsilon_0<\delta/32$ is fixed.  Each use of the preceding moment bounds takes a sufficiently small auxiliary exponent to fit within $\epsilon_0$.

Set $D_*\coloneqq\lambda^{-2+\delta/2}Q^{1/2}.$ We verify a power gain at this scale.  The linear term in $\mathcal T(D_*)/D_*$ is bounded by $C\lambda^{\epsilon_0-\delta/2}$.  For the nonlinear term its exponent is at most
\[
\begin{aligned}
        &\epsilon_0-1+\frac\beta2-\frac{\beta\delta}4
          +(\kappa_d^{\rm LL}-\delta)
             \left(\frac32+\frac\beta4\right)\\
        &\quad=
        \epsilon_0
        +(\beta-\beta_d)
           \left(\frac12+\frac{\kappa_d^{\rm LL}}4\right)
        -\delta\left(\frac32+\frac\beta2\right)
        <-\delta.
\end{aligned}
\]
We used the defining identity for $\kappa_d^{\rm LL}$ in \cref{eq:highd-beta-kappa}.  The finite-volume error is smaller than any fixed inverse power by the standing volume assumption.  Hence, with $\nu\coloneqq\delta/8$, for sufficiently large $\lambda$,
\begin{equation}
\label{eq:highd-scale-contraction}
        \frac{\mathcal T(D)}{D_*}
        \leq\lambda^{-\nu}
        \left[1+\left(\frac{D}{D_*}\right)^\alpha\right]
\end{equation}
for every scale arising in the iteration.

Start with $D_0=\lambda^{-1}$ and set $D_{j+1}=\mathcal T(D_j)$.  To make the finite number of iterations explicit, let $a_0=2$ and $a_{j+1}=\max\{0,\alpha a_j-\nu/2\}$. Since $D_0/D_*\leq\lambda^2$, \cref{eq:highd-scale-contraction} implies $D_j/D_*\leq\lambda^{a_j}$ until $a_j=0$; after this point the ratio remains at most one.  Before the first such point, $a_j=2\alpha^j- \frac\nu2\frac{1-\alpha^j}{1-\alpha}.$ Thus any integer $j$ with $\alpha^j<\frac\nu{4(1-\alpha)+\nu}$ is sufficient.  This number depends only on $d$ and $\delta$. We have proved
\[
        \|R-M\|_{1\to\infty}\prec D_*.
\]

For completeness, the nonlinear term of $F(D_*)/\lambda^{-2}$ has exponent at most
\[
        (\beta-\beta_d)
           \left(\frac12+\frac{\kappa_d^{\rm LL}}4\right)
        -\delta\left(\frac32+\frac\beta4\right)
        +\frac{\alpha\delta}2
        <-\frac\delta2.
\]
Therefore \cref{eq:highd-coarse-one-four,eq:highd-coarse-F-map} give
\[
        \|R\|_{1\to4}\prec\lambda^{-1},
        \qquad
        \|R-\mathbb ER\|_{1\to\infty}\prec\lambda^{-2}.
\]
Their fixed-moment versions, \cref{eq:highd-initial-moments}, and Lemma~\ref{lemma:flut_exp_bdd} yield the deterministic estimate
\[
        \|\mathbb ER-M\|_{1\to\infty}
        \leq C_\delta\lambda^{-2+\delta}Q^{1/2}+C\chi_L(\eta).
\]
This proves all three asserted bounds under the bootstrap assumption. At $\eta=\eta_*$, the right side is $O(\lambda^{-2+\kappa_d^{\rm LL}/2+\delta/2})=o(\lambda^{-1})$, since $\delta<\kappa_d^{\rm LL}<2/3$. The bootstrap starts at $\eta=K\lambda$ by the trivial resolvent bound.  Continuity in $\eta$, exactly as in the $d=3$ proof, stops the bootstrap at the first point where \cref{eq:highd-bootstrap-assumption} fails.
\end{proof}
\section{Diffusive profile and dynamical delocalization}
\label{sec:diffusion-profile}
We first state the dynamical consequence of the diffusive profile.  Recall from Section~\ref{sec:introduction} that the characteristic spreading length in this normalization is \(\lambda^{3/2}\sqrt t\).

Define
\begin{equation}
\label{eq:BA-dynamical-exponent}
        \kappa_d^{\rm BA}
        \coloneqq
        \begin{cases}
        \dfrac{2}{13},&d\in\{2,4,5,6\},\\[2mm]
        \dfrac{2}{9},&d=3,\\[2mm]
        \min\left\{\dfrac12,\dfrac{2d-12}{3d-6}\right\},&d\geq7.
        \end{cases}
\end{equation}
The last line is the intersection of the range in which the block \(T\)-equation can be inverted and the range of the high-dimensional local law in Theorem~\ref{thm:local_law_block_anderson}.

The estimates through Subsection~\ref{subsec:BA-resolvent-to-dynamics} are first used for fixed \(W\) and for growing blocks satisfying \(W\leq\lambda^{2-\gamma_0}\), with fixed \(\gamma_0>0\). Subsection~\ref{subsec:BA-genuine-large-blocks} then gives bounds uniform over all growing block sizes.  In particular, it proves the remaining part \(W\leq\lambda^2\) of Theorem~\ref{thm:BA-dynamical-delocalization} and the crossover statement in Theorem~\ref{thm:BA-genuine-large-block-dynamics}. The row-moment estimate below, the algebraic reductions in Lemmas~\ref{lemma:BA-trace-Feshbach} and~\ref{lemma:BA-Ward-mass}, and the hard-inverse, symmetry, and positive-renewal arguments in Lemmas~\ref{lemma:BA-low-frequency-hard}, \ref{lemma:BA-scalar-symbol-expansion}, and~\ref{lemma:BA-positive-one-flight} already hold without a mesoscopic upper bound.  These arguments use only the finite kernel moments stated below.  The same is true of the coercivity estimate in Proposition~\ref{prop:BA-cell-coercivity}.

\subsection{The block T-equation and the deterministic ladder}
\label{subsec:BA-T-equation}
We derive a deterministic approximation for the diagonal double resolvent associated with a block-diagonal observable.  Equivalently, we control quantities of the form
\[
        \sum_y a(y)|R_{xy}|^2
        =\langle e_x,R\mathcal AR^*e_x\rangle.
\]
\begin{lemma}[Block \(T\)-equation]
\label{lemma:block_T_equation}
Let \(\mathcal A\in\mathcal B(\ell^2(\mathbb Z_L^d))\) be deterministic. Then
\begin{equation}
\label{eq:block-T-equation}
\mathbb E[R\mathcal AR^*]
=
\tilde M \mathcal A\tilde M^\star
+
\tilde M S[\mathbb E(R\mathcal AR^\star)]\tilde M^\star
+
\mathfrak C_T^{\mathrm{BA}},
\end{equation}
where
\begin{align}
\label{eq:block-T-residual}
\mathfrak C_T^{\mathrm{BA}}
=
\tilde M \mathcal A(\mathbb ER-\tilde M)^\star+
\tilde M\,
\mathbb E\!\left[
S[R-\mathbb ER]\,R\mathcal AR^\star
\right]+
\tilde M\,
\mathbb E\!\left[
S[R\mathcal AR^*](R-\tilde M)^\star
\right].
\end{align}
\end{lemma}

\begin{proof}
By the definitions of \(R\) and \(\tilde M\),
\[
R^{-1}-\tilde M^{-1}
=
\mathcal V+\tilde\Theta.
\]
Hence the resolvent identity gives
\begin{equation}
\label{eq:block-renormalized-resolvent}
R
=
\tilde M-\tilde M
(\mathcal V+\tilde\Theta)R.
\end{equation}
Multiplying by \(\mathcal AR^*\) and taking expectations, we obtain
\begin{align*}
\mathbb E(R\mathcal AR^\star)
=
\tilde M \mathcal A\,\mathbb ER^\star
-
\tilde M\,\mathbb E(\mathcal V R\mathcal AR^\star) -
\tilde M\tilde\Theta\,
\mathbb E(R\mathcal AR^\star).
\end{align*}

Using
\[
\partial_{g_\alpha}R=-RA_\alpha R,
\qquad
\partial_{g_\alpha}R^\star=-R^\star A_\alpha R^\star,
\]
Gaussian integration by parts gives
\begin{align*}
\mathbb E(\mathcal V R\mathcal AR^\star)
&=
\sum_\alpha
A_\alpha
\mathbb E\!\left[
\partial_{g_\alpha}(R\mathcal AR^\star)
\right]\\
&=
-\mathbb E\!\left[
S[R]\,R\mathcal AR^\star
\right]
-\mathbb E\!\left[
S[R\mathcal AR^\star]R^\star
\right].
\end{align*}
Since $\tilde\Theta=S[\mathbb ER]$, trivially we have $S[R] = \tilde\Theta+S[R-\mathbb ER].$ The contribution containing \(\tilde\Theta\) therefore cancels the last term above, and we obtain
\begin{align*}
\mathbb E(R\mathcal AR^\star)
=\tilde M \mathcal A\,\mathbb ER^\star +
\tilde M\,
\mathbb E\!\left[
S[R-\mathbb ER]\,R\mathcal AR^\star
\right]+ \tilde M\,
\mathbb E\!\left[
S[R\mathcal AR^\star]R^\star
\right].
\end{align*}
Finally, use $R^\star=\tilde M^\star+(R-\tilde M)^\star$ in the first and third terms. Since \(S\) is linear, $\mathbb E S[R\mathcal AR^\star] = S[\mathbb E(R\mathcal AR^\star)],$ which proves \cref{eq:block-T-equation}--\cref{eq:block-T-residual}.
\end{proof}
The covariance map keeps the full matrix state inside each cell, so the closed unknown is initially matrix-valued instead of scalar.  Translation invariance nevertheless turns the equation into a convolution on the block lattice.  Write \(\mathcal A=\bigoplus_{J\in\widetilde{\mathbb Z}_n^d}\mathcal A_J\). For \(I\in\widetilde{\mathbb Z}_n^d\),
\begin{equation}
    \begin{aligned}
         \Pi_{[I]}\mathbb ER\mathcal{A}R^\star \Pi_{[I]}&= \Pi_{[I]}\tilde M \mathcal A\tilde M^\star \Pi_{[I]} + \Pi_{[I]}\tilde M S[\mathbb ER\mathcal{A}R^\star]\tilde M^\star \Pi_{[I]} + \Pi_{[I]}\mathfrak C_T^{\rm BA}\Pi_{[I]},
    \end{aligned}
\end{equation}
The deterministic resolvent \(\tilde M\) is block-translation invariant. For a block-matrix-valued function \(f=(f_I)_I\), define
\begin{equation}
\label{eq:BA-T-weighted-norm}
        \|f\|_{\infty,\Im\tilde\theta}
        \coloneqq
        \sup_{I\in\tilde{\mathbb Z}_n^d}
        \left\|
        \mathcal{C}_{\Im\tilde\theta}[f_I]
        \right\|_{2\to2}.
\end{equation} 
Set
\[
 \tilde M(I-J)\coloneqq\blk_{I,J}[\tilde M],\qquad
 \tilde K(I)[B]\coloneqq
 \mathscr S_W[\tilde M(I)B\tilde M(I)^\star],
\]
and \(g_I\coloneqq\mathscr S_W[\iota_I^\star \mathbb E(R\mathcal AR^\star)\iota_I]\).  Then
\begin{equation}
\label{eq:exact_closed_equation}
    \begin{aligned}
        g_I = (\tilde K*\mathcal{A})[I] + (\tilde K*g)[I] +\mathscr S_W[ \blk_{I}[\mathfrak C_T^{\rm BA}]],
    \end{aligned}
\end{equation}
where convolution is taken over the block coordinate.  This is the matrix-valued ladder equation.  The trace Feshbach reduction below will show that only its scalar trace component can remain slow.

Let $\tilde{\mathcal K}$ denote convolution by $\tilde K$ on the block torus.  With $\mathfrak e(I)=\mathscr S_W[\blk_I[\mathfrak C_T^{\rm BA}]]$, the ladder equation becomes $(\Id-\tilde{\mathcal K})g=\tilde{\mathcal K}\mathcal A+\mathfrak e.$ The contraction estimate below makes $\Id-\tilde{\mathcal K}$ invertible, and hence
\begin{equation}
\label{eq:exact_closed_equ_inver}
        g=(\Id-\tilde{\mathcal K})^{-1}
          (\tilde{\mathcal K}\mathcal A+\mathfrak e).
\end{equation}
\begin{lemma}
\label{lemma:BA-T-error-and-kernel}
Let $z=E+\mathrm i\eta$ be in the range of the local law from Theorem~\ref{thm:local_law_block_anderson}, and set
\begin{equation}
\label{eq:BA-T-Q-Phi}
        \Phi_d^{\rm BA}(z)
        \coloneqq
        \begin{cases}
        \lambda^{-1/2}(1+(\lambda\eta)^{-1})^{13/4},
            & d\in\{2,4,5,6\},\\[2mm]
        \lambda^{-3/4}(1+(\lambda\eta)^{-1})^{27/8},
            & d=3,\\[2mm]
        \lambda^{-1}(1+(\lambda\eta)^{-1})^2,
            & d\geq7.
        \end{cases}
\end{equation}
For every deterministic operator $\mathcal A$ satisfying $\|\mathcal A\|_{2\to2}\leq1$, for every fixed $\epsilon_0>0$, we have
\begin{equation}
\label{eq:BA-T-error-max-bound}
        \|\mathfrak C_T^{\rm BA}\|_{1\to\infty}
        \leq
        C_{\epsilon_0}\lambda^{\epsilon_0}
        \Phi_d^{\rm BA}(z).
\end{equation}
In the fixed-block regime the constant may depend on $W$ and on the finite-height folded bulk window. Moreover, if $0<\lambda\eta\leq1$ and $\lambda^{\epsilon_0}\Phi_d^{\rm BA}(z)$ is sufficiently small, then
\begin{equation}
\label{eq:BA-T-kernel-contraction}
        \|\tilde{\mathcal K}\|_{(\infty,\Im\tilde\theta)\to(\infty,\Im\tilde\theta)}
        \leq1-c(1+(\lambda\eta)^{-1})^{-1}.
\end{equation}
Consequently,
\begin{equation}
\label{eq:BA-T-kernel-inverse}
\begin{aligned}
        \left\|
        (\Id-\tilde{\mathcal K})^{-1}
        \right\|_{(\infty,\Im\tilde\theta)\to(\infty,\Im\tilde\theta)}
        &\leq C(1+(\lambda\eta)^{-1}),
        \\
        \left\|
        (\Id-\tilde{\mathcal K})^{-1}
        \right\|_{\ell^\infty(2\to2)\to\ell^\infty(2\to2)}
        &\leq C(1+(\lambda\eta)^{-1}).
\end{aligned}
\end{equation}
In particular, the smallness condition in the lemma holds, after choosing $\epsilon_0>0$ sufficiently small, on the intersection of the local-law range with
\[
        (1+(\lambda\eta)^{-1})\leq
        \begin{cases}
        \lambda^{2/13-\delta},&d\in\{2,4,5,6\},\\
        \lambda^{2/9-\delta},&d=3,\\
        \lambda^{1/2-\delta},&d\geq7,
        \end{cases}
\]
for any fixed $\delta>0$.
\end{lemma}

\begin{proof}

We record the moment estimates used below.  For a random matrix $B$, write
\[
        \|B\|_{L^p(a\to b)}
        \coloneqq
        \left(\mathbb E\|B\|_{a\to b}^p\right)^{1/p}.
\]
The moment form of the local law, Lemma~\ref{lemma:concent_1_4_bdd}, the Ward identity, and the deterministic column bounds imply, for every fixed sufficiently large $p$,
{\footnotesize
\begin{equation}
\label{eq:BA-T-input-scales}
\begin{array}{c|c|c|c}
 & d\in\{2,4,5,6\} & d=3 & d\geq7 \\ \hline
 \|R-\mathbb ER\|_{L^p(1\to\infty)}
 & \lambda^{-3/2}(1+(\lambda\eta)^{-1})^{7/4}
 & \lambda^{-7/4}(1+(\lambda\eta)^{-1})^{15/8}
 & \lambda^{-2}(1+(\lambda\eta)^{-1})^{1/2}
 \\[1mm]
 \|\mathbb ER-\tilde M\|_{1\to2}
 +\|(\mathbb ER-\tilde M)^\star\|_{1\to2}
 & \lambda^{-1/2}(1+(\lambda\eta)^{-1})^{9/4}
 & \lambda^{-3/4}(1+(\lambda\eta)^{-1})^{19/8}
 & \lambda^{-1}(1+(\lambda\eta)^{-1})
 \\[1mm]
 \bigl\|S[R\mathcal A R^\star]
       -\mathbb ES[R\mathcal A R^\star]\bigr\|_{L^p(2\to2)}
 & \lambda^{-1/2}(1+(\lambda\eta)^{-1})^{19/8}
 & \lambda^{-3/4}(1+(\lambda\eta)^{-1})^{39/16}
 & \lambda^{-1}(1+(\lambda\eta)^{-1})^{3/2}.
\end{array}
\end{equation}
}
Every entry on the right-hand side of \cref{eq:BA-T-input-scales} may carry the common factor $C_{\epsilon_0}\lambda^{\epsilon_0}$.  We suppress this factor until the end of the proof.

We make explicit the interpolation input for the third line.  The moment versions of the estimates used in the proof of the local law give
\begin{align}
\label{eq:BA-T-one-four-input}
        \|R\|_{L^{3p}(1\to4)}
        &\leq C_{\epsilon_0}\lambda^{\epsilon_0}
        \begin{cases}
        \lambda^{-3/4}(1+(\lambda\eta)^{-1})^{7/8},&d\in\{2,4,5,6\},\\
        \lambda^{-7/8}(1+(\lambda\eta)^{-1})^{15/16},&d=3,\\
        \lambda^{-1},&d\geq7,
        \end{cases}\\
\label{eq:BA-T-two-four-input}
        \|R\|_{L^{3p}(2\to4)}
        &\leq C_{\epsilon_0}\lambda^{\epsilon_0}
        \begin{cases}
        \lambda^{1/4}(1+(\lambda\eta)^{-1}),&d\in\{2,4,5,6\},\\
        \lambda^{1/8}(1+(\lambda\eta)^{-1}),&d=3,\\
        (1+(\lambda\eta)^{-1}),&d\geq7.
        \end{cases}
\end{align}
Together with $\|R\|_{L^{3p}(1\to2)}\leq C (1+(\lambda\eta)^{-1})^{1/2}$, these estimates, H\"{o}lder's inequality and Lemma~\ref{lemma:concent_1_4_bdd} imply
\begin{align*}
        \bigl\|S[R\mathcal A R^\star]
        -\mathbb ES[R\mathcal A R^\star]\bigr\|_{L^p(2\to2)}
        \leq
        C_{\epsilon_0}\lambda^{\epsilon_0}
        \|R\|_{L^{3p}(1\to4)}
        \|R\|_{L^{3p}(1\to2)}
        \|R\|_{L^{3p}(2\to4)},
\end{align*}
which is exactly the third line of \cref{eq:BA-T-input-scales}. Here \cref{eq:BA-T-one-four-input} follows by the same Riesz--Thorin interpolation used in the local-law bootstrap, while \cref{eq:BA-T-two-four-input} follows from the half-resolvent estimate in Proposition~\ref{prop:prior_resol_bd}.  Stochastic domination upgrades to these fixed $L^p$ bounds because the resolvents have a polynomial deterministic bound and $L^d$ is polynomial in $\lambda$.

For completeness, the second line of \cref{eq:BA-T-input-scales} follows from the exact averaged resolvent equation
\begin{equation}
\label{eq:BA-T-averaged-resolvent-equation}
        \mathbb ER-\tilde M
        =\tilde M\,
        \mathbb E[S[R-\mathbb ER]R].
\end{equation}
Indeed,
\[\|R\|_{1\to 2}^2 \leq \frac{1}{\eta}\sup_x\Im R_{(x,x)} \leq \frac{1}{\eta}\|R\|_{1\to \infty},\]
and \(\|R\|_{1\to\infty}\prec\lambda^{-1}\).  Therefore
\begin{equation}
\label{eq:BA-T-column-inputs}
\begin{aligned}
        \|\tilde M\|_{2\to2}&\leq C\lambda,
        \\
        \|\tilde M\|_{1\to2}
        +\|\tilde M^\star\|_{1\to2}&\leq C,
        \\
        \|R\|_{L^p(1\to2)}
        +\|R^\star\|_{L^p(1\to2)}&\leq C(1+(\lambda\eta)^{-1})^{1/2}.
\end{aligned}
\end{equation}
Together with \cref{equ:sup_op_ineq}, equation \cref{eq:BA-T-averaged-resolvent-equation} implies
\begin{equation}
\label{eq:BA-T-one-two-from-fluctuation}
        \|\mathbb ER-\tilde M\|_{1\to2}
        \leq
        C\lambda (1+(\lambda\eta)^{-1})^{1/2}
        \|R-\mathbb ER\|_{L^2(1\to\infty)},
\end{equation}
which gives precisely the second line of \cref{eq:BA-T-input-scales}.  The same equation, now estimated entrywise using the second line of \cref{eq:BA-T-column-inputs}, also gives
\begin{equation}
\label{eq:BA-T-averaged-max-bound}
        \|\mathbb ER-\tilde M\|_{1\to\infty}
        \leq
        C_{\epsilon_0}\lambda^{\epsilon_0}
        \lambda^{-1}(1+(\lambda\eta)^{-1})^{-1}\Phi_d^{\rm BA}(z).
\end{equation}
The matrices $R$, $\mathbb ER$, and $\tilde M$ are complex symmetric. Consequently,
\[
        \|(\mathbb ER-\tilde M)^\star\|_{1\to2}
        =\|\mathbb ER-\tilde M\|_{1\to2},
\]
which proves the adjoint part of the second line of \cref{eq:BA-T-input-scales}.  The use of \cref{equ:sup_op_ineq} above also shows explicitly that the block covariance map incurs no factor $W^d$.

We now estimate the three terms in \cref{eq:block-T-residual}, denoting their $1\to\infty$ norms by ${\rm I},{\rm II},{\rm III}$.  From \cref{eq:BA-T-column-inputs} and the second line of \cref{eq:BA-T-input-scales},
\begin{equation}
\label{eq:BA-T-error-I}
\begin{aligned}
        {\rm I}
        &\leq
        \|\tilde M^\star\|_{1\to2}
        \|\mathcal A\|_{2\to2}
        \|\mathbb ER-\tilde M\|_{1\to2}
        \\
        &\leq
        C_{\epsilon_0}\lambda^{\epsilon_0}
        \Phi_d^{\rm BA}(z).
\end{aligned}
\end{equation}
For the second term, H\"older's inequality, $\|R\|_{2\to2}\leq\eta^{-1}\leq\lambda (1+(\lambda\eta)^{-1})$, and \cref{eq:BA-T-column-inputs} give
\begin{equation}
\label{eq:BA-T-error-II}
\begin{aligned}
        {\rm II}
        &\leq
        C\lambda (1+(\lambda\eta)^{-1})^{3/2}
        \|R-\mathbb ER\|_{L^p(1\to\infty)}
        \\
        &\leq
        C_{\epsilon_0}\lambda^{\epsilon_0}
        \Phi_d^{\rm BA}(z).
\end{aligned}
\end{equation}
For the third term,
\begin{equation}
\label{eq:BA-T-error-III-split}
\begin{aligned}
        \mathbb E[S[R\mathcal AR^\star](R-\tilde M)^\star]
        =&
        \mathbb E\left[\left[S[R\mathcal AR^\star]-\mathbb E(S[R\mathcal AR^\star])\right](R-\mathbb ER)^\star\right]+\mathbb E[S[R\mathcal AR^\star]](\mathbb ER-\tilde M)^\star.
\end{aligned}
\end{equation}
Moreover, by \cref{eq:BA-T-column-inputs},
\begin{equation}
\label{eq:BA-T-mean-Z}
        \big\|\mathbb E[S[R\mathcal AR^\star]]\big\|_{2\to2}
        \leq C(1+(\lambda\eta)^{-1})\|\mathcal A\|_{2\to2}.
\end{equation}
Consequently, \cref{eq:BA-T-error-III-split}, H\"older's inequality, and \cref{eq:BA-T-input-scales} imply
\begin{equation}
\label{eq:BA-T-error-III}
\begin{aligned}
        {\rm III}
        \leq{}&
        C(1+(\lambda\eta)^{-1})^{1/2}
        \big\|S[R\mathcal AR^\star]-\mathbb E[S[R\mathcal AR^\star]]\big\|_{L^p(2\to2)}
        +C(1+(\lambda\eta)^{-1})\|\mathbb ER-\tilde M\|_{1\to2}
        \\
        \leq{}&
        C_{\epsilon_0}\lambda^{\epsilon_0}
        \Phi_d^{\rm BA}(z).
\end{aligned}
\end{equation}
Equations \cref{eq:BA-T-error-I}, \cref{eq:BA-T-error-II}, and \cref{eq:BA-T-error-III} prove \cref{eq:BA-T-error-max-bound}.

We turn to the deterministic kernel.  If $\|B\|_{2\to2}\leq1$, then the explicit formula for $\mathscr S_W$ yields
\begin{equation}
\label{eq:BA-T-kernel-pointwise-proof}
\begin{aligned}
        \|\tilde K(X)[B]\|_{2\to2}
        &\leq
        \frac{1}{W^d}
        \left(
        \|\tilde M(X)B\tilde M(X)^\star\|_{2\to2}
        +\left|
        \Tr\bigl(\tilde M(X)B\tilde M(X)^\star\bigr)
        \right|
        \right)
        \\
        &\leq
        \frac{1}{W^d}
        \left(
        \|\tilde M(X)\|_{2\to2}^2
        +\|\tilde M(X)\|_{\rm HS}^2
        \right).
\end{aligned}
\end{equation}

We next justify the two deterministic bounds necessary for the kernel argument. In the large-block regime, Theorem~\ref{thm:main_closeness}, Corollary~\ref{coro:robust_app_bootstrap}, and Corollary~\ref{cor:tildeM_M_compar_large_block} imply
\[
        \frac{c}{\lambda}\I_{W^d}
        \preceq\Im\theta
        \preceq\frac{C}{\lambda}\I_{W^d},
        \qquad
        \sup_x\Im M_{xx}\leq\frac C\lambda,
        \qquad
        \|\tilde M-M\|_{1\to\infty}=o(\lambda^{-1}).
\]
The same three conclusions in the fixed-block regime follow from Theorem~\ref{thm:fh-fixed-block-MDE}, Theorem~\ref{thm:finite_volume_from_infinite_trap}, Corollary~\ref{cor:qfv-unscaled-random-application}, and Lemma~\ref{lemma:fh-finite-M-max-bound}.  Moreover, the local law and \cref{equ:sup_op_ineq} give
\[
\begin{aligned}
        \left\|
        \Im\tilde\theta-\Im\theta
        \right\|_{2\to2}
        &\leq2\|\mathbb ER-M\|_{1\to\infty}
        =o(\lambda^{-1}).
\end{aligned}
\]
Therefore
\begin{equation}
\label{eq:BA-T-Y-and-M-trap}
        \frac{c}{\lambda}\I_{W^d}
        \preceq \Im\tilde\theta\preceq
        \frac{C}{\lambda}\I_{W^d},
        \qquad
        \sup_x\Im\tilde M_{xx}\leq\frac C\lambda,
        \qquad
        \Tr\operatorname{bl}_0[\Im\tilde M]
        \leq\frac{CW^d}{\lambda}.
\end{equation}
Again, set
\begin{equation}
\label{eq:BA-T-Gamma-definition}
        \Gamma\coloneqq\eta\I_{W^d}+\Im\tilde\theta.
\end{equation}
The generalized Ward identity \cref{eq:gen_ward} and block translation invariance give
\begin{equation}
\label{eq:BA-T-block-Ward}
        \sum_X
        \tilde M(X)\Gamma\tilde M(X)^\star
        =\operatorname{bl}_0[\Im\tilde M].
\end{equation}
Since $\Gamma\succeq c\lambda^{-1}\I_{W^d}$, taking the trace in \cref{eq:BA-T-block-Ward} and using \cref{eq:BA-T-Y-and-M-trap} gives
\[
        \sum_X\|\tilde M(X)\|_{\rm HS}^2
        \leq CW^d.
\]
It remains to prove the strict contraction.  Define
\begin{equation}
\label{eq:BA-T-D-and-remainder}
\begin{aligned}
        D&\coloneqq
        \mathscr S_W\!\left[
        \sum_X\tilde M(X)\tilde M(X)^\star
        \right],\quad\mathcal R_K \coloneqq
        \mathscr S_W\!\left[
        \operatorname{bl}_0[
        \Im(\tilde M-\mathbb ER)]
        \right].
\end{aligned}
\end{equation}
Applying $\mathscr S_W$ to \cref{eq:BA-T-block-Ward} and using $        \Im\tilde\theta=\mathscr S_W\!\left[\operatorname{bl}_0[\Im\mathbb ER]\right]$ gives the exact order-unit identity
\begin{equation}
\label{eq:BA-T-order-unit-identity}
        \sum_X\tilde K(X)[\Im\tilde\theta]
        =\Im\tilde\theta-\eta D+\mathcal R_K.
\end{equation}
This is the block replacement for the scalar mass identity. By \cref{equ:sup_op_ineq}, \cref{eq:BA-T-Y-and-M-trap}, and \cref{eq:BA-T-averaged-max-bound},
\begin{equation}
\label{eq:BA-T-remainder-weighted}
\begin{aligned}
        \left\|
        \mathcal{C}_{\Im\tilde\theta}\mathcal R_K
        \right\|_{2\to2}
        &\leq
        C\lambda
        \|\mathbb ER-\tilde M\|_{1\to\infty}
        \\
        &\leq
        C_{\epsilon_0}\lambda^{\epsilon_0}
        (1+(\lambda\eta)^{-1})^{-1}\Phi_d^{\rm BA}(z).
\end{aligned}
\end{equation}
Since $0<\lambda\eta\leq1$, equations \cref{eq:BA-T-Y-and-M-trap} and \cref{eq:BA-T-Gamma-definition} imply $c\lambda^{-1}\I_{W^d}\preceq\Gamma \preceq C\lambda^{-1}\I_{W^d}$.  Combining this fact with \cref{eq:BA-T-block-Ward} after applying $\mathscr S_W$ gives
\[
        \frac{c}{\lambda}D
        \preceq \Im\tilde\theta+\mathcal R_K
        \preceq\frac{C}{\lambda}D.
\]
Equation \cref{eq:BA-T-remainder-weighted} implies $-o(1)\Im\tilde\theta\preceq\mathcal R_K\preceq o(1)\Im\tilde\theta$; hence
\begin{equation}
\label{eq:BA-T-D-comparison}
        c\lambda \Im\tilde\theta\preceq D\preceq C\lambda \Im\tilde\theta.
\end{equation}
Therefore \cref{eq:BA-T-order-unit-identity} yields
\begin{equation}
\label{eq:BA-T-order-unit-contraction}
\begin{aligned}
        0\preceq\sum_X\tilde K(X)[\Im\tilde\theta]
        &\preceq
        \left(
        1-c\lambda\eta
        +C_{\epsilon_0}\lambda^{\epsilon_0}
        (1+(\lambda\eta)^{-1})^{-1}\Phi_d^{\rm BA}(z)
        \right)\Im\tilde\theta
        \\
        &\preceq
        \left(1-\frac{c}{1+(\lambda\eta)^{-1}}\right)\Im\tilde\theta,
\end{aligned}
\end{equation}
provided that $\lambda^{\epsilon_0}\Phi_d^{\rm BA}(z)$ is sufficiently small.  Here we used $\lambda\eta/(1+\lambda\eta)\asymp\lambda\eta$.

To pass from \cref{eq:BA-T-order-unit-contraction} to the norm of the full superoperator, conjugate by $\Im\tilde\theta$ and define
\begin{equation}
\label{eq:BA-T-conjugated-positive-map}
\begin{aligned}
        (\tilde{\mathcal K}_{\Im\tilde\theta} f)_I
        \coloneqq
        (\Im\tilde\theta)^{-1/2}
        \tilde{\mathcal K}
        \bigl[((\Im\tilde\theta)^{1/2}f_J(\Im\tilde\theta)^{1/2})_J\bigr]_I
        (\Im\tilde\theta)^{-1/2}.
\end{aligned}
\end{equation}
The map $\tilde{\mathcal K}_{\Im\tilde\theta}$ is positive on the finite-dimensional $C^*$-algebra $\ell^\infty(\tilde{\mathbb Z}_n^d; \mathbb C^{W^d\times W^d})$.  The norm of a positive map on a unital $C^*$-algebra equals the norm of its value at the identity.  Hence \cref{eq:BA-T-order-unit-contraction} implies
\[
        \|\tilde{\mathcal K}_{\Im\tilde\theta}\|_{\ell^\infty(2\to2)\to\ell^\infty(2\to2)}
        =\|\tilde{\mathcal K}_{\Im\tilde\theta}[\I]\|_{\ell^\infty(2\to2)}
        \leq1-\frac{c}{1+(\lambda\eta)^{-1}},
\]
where $\I$ in this display denotes the constant identity block function. This is exactly \cref{eq:BA-T-kernel-contraction}.  The Neumann series proves the first bound in \cref{eq:BA-T-kernel-inverse}.  Finally, \cref{eq:BA-T-Y-and-M-trap} gives the norm equivalence
\[
        c\lambda\|f\|_{\ell^\infty(2\to2)}
        \leq\|f\|_{\infty,\Im\tilde\theta}
        \leq
        C\lambda\|f\|_{\ell^\infty(2\to2)},
\]
and therefore also the second bound in \cref{eq:BA-T-kernel-inverse}.
\end{proof}

\begin{remark}
\label{rmk:BA-T-error-after-inversion}
\cref{equ:sup_op_ineq,eq:BA-T-error-max-bound} and \cref{eq:BA-T-kernel-inverse} imply
\begin{equation}
\label{eq:BA-T-error-after-inversion}
        \left\|
        (\operatorname{Id}-\tilde{\mathcal K})^{-1}
        \mathfrak e
        \right\|_{\ell^\infty(2\to2)}
        \leq
        C_{\epsilon_0}\lambda^{\epsilon_0}
        (1+(\lambda\eta)^{-1})\Phi_d^{\rm BA}(z).
\end{equation}
In contrast to the scalar Anderson model, one should not replace \cref{eq:BA-T-kernel-contraction} by an unweighted mass bound for the scalar envelope of $\tilde K$.  The weighted order-unit contraction is the natural statement for the superoperator-valued block kernel.
\end{remark}

\begin{thm}
\label{thm:rand_deter_ladder_apprx}
Assume the hypotheses and the smallness condition of Lemma~\ref{lemma:BA-T-error-and-kernel}.  Let $\mathcal A$ satisfy
\[
\mathcal A=
\bigoplus_{J\in\tilde{\mathbb Z}_n^d}\mathcal A_J,
\qquad
\mathcal{A}_J = \mathcal{A}_J^\star,
\qquad
\sup_J\|\mathcal{A}_J\|_{2\to2}\leq1,
\]
Then
\begin{equation}
\label{eq:BA-transformed-diffusive-profile}
\begin{aligned}
\sup_{I\in\tilde{\mathbb Z}_n^d}
\left\|
\mathscr S_W\!\left[
\operatorname{bl}_I
(R\mathcal A R^\star)
\right]
-(\Id - \tilde{\mathcal{K}})^{-1}\tilde{\mathcal{K}} \mathcal{A}(I)
\right\|_{2\to2}
\prec
\left(1+(\lambda\eta)^{-1}\right)
\Phi_d^{\rm BA}(z).
\end{aligned}
\end{equation}
Equivalently, the right-hand side of \cref{eq:BA-transformed-diffusive-profile} is
\begin{equation}
\label{eq:BA-diffusive-profile-scale}
\left(1+(\lambda\eta)^{-1}\right)
\Phi_d^{\rm BA}(z)
=
\begin{cases}
\lambda^{-1/2}
\left(1+(\lambda\eta)^{-1}\right)^{17/4},
& d\in\{2,4,5,6\},\\[2mm]
\lambda^{-3/4}
\left(1+(\lambda\eta)^{-1}\right)^{35/8},
& d=3,\\[2mm]
\lambda^{-1}
\left(1+(\lambda\eta)^{-1}\right)^3,
& d\geq7.
\end{cases}
\end{equation}

Moreover, the corresponding statement for the untransformed diagonal is
\begin{equation}
\label{eq:BA-raw-diffusive-profile}
\begin{aligned}
\sup_{I\in\tilde{\mathbb Z}_n^d}
\left\|
\blk_I[R\mathcal A R^\star]
-
\left(\mathscr S_W^{-1}\Bigl[\bigl(\Id-\tilde{\mathcal K}\bigr)^{-1}
\tilde{\mathcal K}\mathcal A(I)\Bigr]\right)
\right\|_{1\to\infty}
\prec
\left(1+(\lambda\eta)^{-1}\right)
\Phi_d^{\rm BA}(z),
\end{aligned}
\end{equation}
where the inverse fiberwise superoperator is given by
\begin{equation}
\label{eq:BA-covariance-map-inverse}
\mathscr S_W^{-1}[B]
=
W^d B^\top
-
\frac{W^d}{W^d+1}\Tr(B)\I_{W^d}.
\end{equation}
\end{thm}

\begin{proof}
By \cref{eq:exact_closed_equation,eq:exact_closed_equ_inver},
\[g_I = \tilde{\mathcal{K}}\mathcal{A}[I] + \tilde{\mathcal{K}}g[I] + \mathfrak e_I,\]
and
\[g(I) = (\Id - \tilde{\mathcal{K}})^{-1}(\tilde{\mathcal{K}}\mathcal A[I] + \mathfrak e_I).\]
By the explicit formula for $\mathscr S_W$,
\begin{equation}
\label{eq:BA-e-from-C}
\|\mathfrak e\|_{\ell^\infty(2\to2)}
\leq
2\|\mathfrak C_T^{\rm BA}\|_{1\to\infty}
\leq
C_{\epsilon_0}\lambda^{\epsilon_0}
\Phi_d^{\rm BA}(z).
\end{equation}
The weighted kernel estimate from Lemma~\ref{lemma:BA-T-error-and-kernel} gives
\begin{equation}
\left\|(\operatorname{Id}-\tilde{\mathcal K})^{-1}
        \right\|_{\ell^\infty(2\to2)\to\ell^\infty(2\to2)}
        \leq C(1+(\lambda\eta)^{-1})
\end{equation}
Thus
\begin{equation}
\label{eq:BA-expected-profile-final}
\|g-(\Id -\tilde{\mathcal{K}})^{-1}\tilde{\mathcal{K}}\mathcal{A}\|_{\ell^\infty(2\to2)}
\leq
C_{\epsilon_0}\lambda^{\epsilon_0}
\left(1+(\lambda\eta)^{-1}\right)
\Phi_d^{\rm BA}(z).
\end{equation}
We next estimate the random fluctuation.  For every fixed \(q\geq2\), Lemma~\ref{lemma:concent_1_4_bdd} and H\"older's inequality give
\begin{equation}
\label{eq:BA-RAR-concentration}
\begin{aligned}
 &\left(\mathbb E\left\|
 R\mathcal A R^\star-\mathbb E(R\mathcal A R^\star)
                   \right\|_{1\to\infty}^q\right)^{1/q}\\
 &\qquad\leq C\sqrt{q+\log(2L)}
 \|R\|_{L^{3q}(1\to4)}
 \|R\|_{L^{3q}(1\to2)}
 \|R\|_{L^{3q}(2\to4)},
\end{aligned}
\end{equation}
where \(\|\mathcal A\|_{2\to2}\leq1\) as above. The moment bounds used in the proof of Lemma~\ref{lemma:BA-T-error-and-kernel} give
\begin{equation}
\label{eq:BA-RAR-fluctuation-scales}
\left\|
R\mathcal A R^\star
-
\mathbb E(R\mathcal A R^\star)
\right\|_{1\to\infty}
\prec
\begin{cases}
\lambda^{-1/2}
\left(1+(\lambda\eta)^{-1}\right)^{19/8},
& d\in\{2,4,5,6\},\\[2mm]
\lambda^{-3/4}
\left(1+(\lambda\eta)^{-1}\right)^{39/16},
& d=3,\\[2mm]
\lambda^{-1}
\left(1+(\lambda\eta)^{-1}\right)^{3/2},
& d\geq7.
\end{cases}
\end{equation}
The supremum over the polynomially many matrix entries is absorbed by stochastic domination. Furthermore,
\begin{equation}
\label{eq:BA-transformed-fluctuation}
\begin{aligned}
\sup_I
\|\mathscr S_W[\blk_I[R\mathcal{A}R^\star]] - g_I\|_{2\to2}
&\leq
2
\left\|
R\mathcal A R^\star
-
\mathbb E(R\mathcal A R^\star)
\right\|_{1\to\infty}.
\end{aligned}
\end{equation}
Since $1+(\lambda\eta)^{-1}\geq1$, each quantity on the right-hand side of \cref{eq:BA-RAR-fluctuation-scales} is bounded by the corresponding quantity in \cref{eq:BA-diffusive-profile-scale}. Therefore \cref{eq:BA-expected-profile-final}, \cref{eq:BA-transformed-fluctuation}, and the triangle inequality prove \cref{eq:BA-transformed-diffusive-profile}.

It remains to recover the untransformed microscopic diagonal without losing a factor depending on $W$. Introduce the block-diagonal operators
\[
\mathbf h\coloneqq\bigoplus_I(\Id - \tilde{\mathcal{K}})^{-1}\tilde{\mathcal{K}}\mathcal{A}(I),
\qquad
S[\mathbb E(R\mathcal{A}R^\star)]
=
\bigoplus_I g_I,
\]
and denote
\begin{equation}
\label{eq:BA-deterministic-lift}
P
\coloneqq
\tilde M(\mathcal A+\mathbf h)\tilde M^\star.
\end{equation}
Consequently,
\begin{align}
\mathscr S_W[\operatorname{bl}_IP]
&=
\sum_J
\mathscr S_W\!\left[
\tilde M(I-J)\left(\mathcal A_J+(\Id - \tilde{\mathcal{K}})^{-1}\tilde{\mathcal{K}}\mathcal{A}(J)\right)
\tilde M(I-J)^\star
\right]
\nonumber\\
&=
\bigl(\tilde{\mathcal K}(\mathcal{A}+\mathbf h)\bigr)_I
= (\Id - \tilde{\mathcal{K}})^{-1}\tilde{\mathcal{K}}\mathcal{A}(I).
\label{eq:BA-lift-has-profile}
\end{align}
Since $\mathscr S_W$ is invertible,
\begin{equation}
\label{eq:BA-lift-block}
\operatorname{bl}_IP
=
\mathscr S_W^{-1}[(\Id - \tilde{\mathcal{K}})^{-1}\tilde{\mathcal{K}}\mathcal{A}(I)].
\end{equation}
The identity \cref{eq:BA-covariance-map-inverse} follows directly from
\[
\mathscr S_W[B]
=
\frac{1}{W^d}
\left(B^\top+\Tr(B)\I_{W^d}\right).
\]

Subtracting \cref{eq:BA-deterministic-lift} from the block T-equation gives
\begin{equation}
\label{eq:BA-raw-expectation-error}
\mathbb E[R\mathcal A R^\star]-P
=
\tilde M
(S[\mathbb E(R\mathcal{A}R^\star)]-\mathbf h)
\tilde M^\star
+
\mathfrak C_T^{\rm BA}.
\end{equation}
Using
\[
\sup_x\|\tilde M^\star|x\rangle\|_2
=
\|\tilde M^\star\|_{1\to2}
\leq C,
\]
we obtain from \cref{eq:BA-raw-expectation-error} that
\begin{align}
&
\sup_{I\in\tilde{\mathbb Z}_n^d}
\left\|
\blk_I\left[\mathbb E[R\mathcal A R^\star]\right]
- \left(\mathscr S_W^{-1}\Bigl[\bigl(\Id-\tilde{\mathcal K}\bigr)^{-1}
\tilde{\mathcal K}\mathcal A(I)\Bigr]\right)
\right\|_{1\to\infty}
\nonumber\\
&\qquad\leq
C\sup_{I}\left\|g_I - (\Id - \tilde{\mathcal{K}})^{-1}\tilde{\mathcal{K}}\mathcal{A}(I)\right\|_{2\to2}
+
\|\mathfrak C_T^{\rm BA}\|_{1\to\infty}
\nonumber\\
&\qquad\leq
C_{\epsilon_0}\lambda^{\epsilon_0}
\left(1+1/(\lambda\eta)\right)
\Phi_d^{\rm BA}(z).
\label{eq:BA-raw-expectation-final}
\end{align}
Finally, \cref{eq:BA-RAR-fluctuation-scales} controls the difference between $R\mathcal A R^\star$ and its expectation. Combining it with \cref{eq:BA-raw-expectation-final} proves \cref{eq:BA-raw-diffusive-profile}.
\end{proof}
\subsection{Reduction to the scalar diffusion mode}
\label{subsec:BA-Fourier-profile}
The target of the subsection is to illustrate the ellipticity of the expectation of resolvent sandwich by adopting a novel Schur-Feshbach reduction. In this subsection we will first replace the finite block torus by the infinite block lattice, and introduce some useful notations and simplifications for better readability. The error caused by replacing the domain is later controlled by a Combes--Thomas comparison. The overall inspiration and heuristics of Feshbach reduction will be introduced in Subsection~\ref{subsubsec:feshbach_heuristics}, whereas Subsection~\ref{subsubsect:scalar_reduct} serves as the major derivation and the proof of ellipticity.

From now on, for a matrix-valued function \(f=(f_X)_{X\in\mathbb Z^d}\), we denote
\begin{equation}
\label{eq:BA-block-Fourier-transform}
        \widehat f(k)
        \coloneqq
        \sum_{X\in\mathbb Z^d}e^{-\mathrm ik\cdot X}f_X,
        \qquad k\in\mathbb T^d.
\end{equation}
The Fourier symbol of the block convolution kernel is the superoperator
\begin{equation}
\label{eq:BA-kernel-symbol}
        \widehat{\tilde K}_\infty(k)[B]
        \coloneqq
        \sum_{X\in\mathbb Z^d}
        e^{-\mathrm ik\cdot X}\tilde K_\infty(X)[B].
\end{equation}
Thus
\begin{equation}
\label{eq:BA-kernel-multiplier}
        \widehat{\tilde{\mathcal K}_\infty f}(k)
        =
        \widehat{\tilde K}_\infty(k)[\widehat f(k)].
\end{equation}
For the fibered operator
\[
        G_{\tilde\theta}(\xi)
        \coloneqq
        (\hat\Psi(\xi)-z/\lambda-\tilde\theta/\lambda)^{-1},
\]
Parseval's identity gives the two-momentum representation
\begin{equation}
\label{eq:BA-kernel-two-momentum}
\begin{aligned}
        \widehat{\tilde K}_\infty(k)[B]
        &=
        \frac{1}{\lambda^2}\mathscr S_W\!\left[
        \int_{\mathbb T^d}
        G_{\tilde\theta}(\xi)B
        G_{\tilde\theta}(\xi-k)^\star
        \,\rmd\xi
        \right].
\end{aligned}
\end{equation}
Indeed,
\begin{equation}
    \begin{aligned}
        \widehat{\tilde{K}}_{\infty}(k)[B] &= \sum_{X \in \mathbb Z^d} e^{-i X\cdot k}\mathscr{S}_W[\tilde{M}_\infty(X) B\tilde{M}_\infty(X)^\star]\\
        &= \sum_{X \in \mathbb Z^d}\mathscr{S}_W\left[e^{-iX\cdot k}\blk_{(X,0)}[\tilde M_\infty] B\blk_{(X,0)}[\tilde M_\infty]^\star\right]\\
        &= \sum_{X, Y \in \mathbb Z^d}\mathscr{S}_W\left[e^{-iY\cdot k}\blk_{(X,0)}[\tilde M_\infty] B\blk_{(Y,0)}[\tilde M_\infty]^\star \int_{\xi \in \mathbb T^d}e^{-i(X - Y)\cdot \xi}\rmd \xi\right]\\
        &= \frac{1}{\lambda^2}\int_{\xi \in \mathbb T^d}\mathscr S_W\left[G_{\tilde\theta}(\xi) B G_{\tilde\theta}(\xi-k)^\star\right]\rmd \xi,
    \end{aligned}
\end{equation}
since $\sum_{X}e^{-iX\cdot \xi}\blk_{(X,0)}[\tilde M_\infty] = \lambda^{-1}G_{\tilde\theta}(\xi)$.  Furthermore, for \(u\in[W^d]\), define the row functional
\begin{equation}
\label{eq:BA-row-functional}
\begin{aligned}
        J_u(k)[B]
        &\coloneqq
        \sum_{X\in\mathbb Z^d}e^{-\mathrm ik\cdot X}
        \left\langle e_u,\tilde{M}_\infty(X)B\tilde{M}_\infty(X)^\star e_u\right\rangle
        \\
        &=
        \frac1{\lambda^2}
        \int_{\mathbb T^d}
        \left\langle e_u,
        G_{\tilde\theta}(\xi)B
        G_{\tilde\theta}(\xi-k)^\star e_u
        \right\rangle
        \,\rmd\xi.
\end{aligned}
\end{equation}
Denote the space of all $W^d\times W^d$ matrices over the complex field as \(\mathcal M_W=\mathbb C^{W^d\times W^d}\). We can define the scalar and trace-free projections on \(\mathcal M_W\) by
\[
\mathsf P[B]=\langle B\rangle_W\I,
\qquad
\mathsf Q[B]=B-\langle B\rangle_W\I.
\]
Then clearly
\[
\mathcal M_W
=\operatorname{Ran}\mathsf P\oplus\operatorname{Ran}\mathsf Q
=\mathbb C\I\oplus
 \{B\in\mathcal M_W:\langle B\rangle_W=0\}.
\]

\subsubsection{Heuristics on Feshbach reduction}
\label{subsubsec:feshbach_heuristics}
The block T-equation indeed provides more information than the standard T-equation, for allowing the block structure, and we can obtain the uniform $1\to \infty$ norm bound for all those "extra" blocks. However, the standard Fourier transform fails to reduce the block T-equation, and for that purpose the natural candidate would be partial Bloch transform. The problem of this method is that, Bloch symbol exhibits multiple different modes, and we must identify the unique mode which account for the elliptic behavior of the untransformed diagonal. To do this, we introduce the Feshbach reduction. Such reduction parallels hydrodynamic projection onto a conserved density: the rapidly relaxing internal degrees of freedom are eliminated, leaving a single density channel capable of producing a diffusion pole \citep{M_Fuchs_2002, 10.1143/PTP.33.423}.

To construct this channel, we adopt normalized trace to manifest the scalar size of the slow mode.  For the purpose of analyzing the quantum diffusion, we might simply choose the scalar source \(\I_{W^d}\), obtained after the Bloch reduction. Notably, identity matrix still generally excites trace-free components when acted on by \(\widehat{\tilde K}(k)\). We therefore seek a correction \(C(k)\in\operatorname{Ran}\mathsf Q\) such that the slow mode profile \(V(k):=\I_{W^d}+C(k)\) has no residual in the trace-free sector:
\[
\mathsf Q(\Id-\widehat {\tilde K}(k))[V(k)]=0,
\qquad \langle V(k)\rangle_W=1.
\]
Solving the equation and the uniform inverse on \(\Ran\mathsf Q\) determines this correction uniquely:
\[
C(k)=
\bigl(\Id_{\operatorname{Ran}\mathsf Q}
      -\mathsf Q\widehat{\tilde K}(k)\mathsf Q\bigr)^{-1}
      \mathsf Q\widehat{\tilde K}(k)[\I_{W^d}].
\]
Here the correction \(C(k)\) records the internal matrix structure required by the ladder equation. The remaining residual is now scalar:
\[
(\Id-\widehat{\tilde K}(k))[V(k)]=s(k)\I_{W^d}.
\]
This Feshbach reduction framework allows the diffusion pole to be carried by the small denominator \(s(k)\), with \(V(k)\) giving its slow mode profile, but it does not yet inform us about the slow mode candidate. Luckily, the natural choice at zero momentum of such $V(0)$ is revealed by the Ward identity. Writing its controlled remainder as \(\mathcal R_0\), Ward identity tells us
\[
(\Id-\widehat{\tilde K}(0))[\Gamma]=\eta\I_{W^d}-\mathcal R_0.
\]
After normalization, the trace-free residual of \(\Gamma/\langle\Gamma\rangle_W\) is determined solely by \(\mathcal R_0\), since \(\mathsf Q[\eta\I_{W^d}]=0\). Consequently,
\[
V(0)-\frac{\Gamma}{\langle\Gamma\rangle_W}
=
\bigl(\Id_{\operatorname{Ran}\mathsf Q}
      -\mathsf Q\widehat{\tilde K}(0)\mathsf Q\bigr)^{-1}
\frac{\mathsf Q\mathcal R_0}{\langle\Gamma\rangle_W}.
\]
A hard-sector bound and a Ward-remainder estimate therefore show that the corrected profile \(V(0)\) aligns with the normalized Ward profile. In short, Ward identifies a nearly conserved direction, the Feshbach reduction provides the framework to capture the unique slow mode, while hard-sector invertibility excludes all other slow internal directions. The remaining scalar analysis gives \(s(0)\asymp\lambda\eta\) and the quadratic spatial behavior \(s(k)\approx s(0)+d_{\rm BA}|k|^2\).
\subsubsection{The scalar reduction}
\label{subsubsect:scalar_reduct}
\begin{lemma}[Trace Feshbach reduction of the block symbol]
\label{lemma:BA-trace-Feshbach}
Fix a real \(k\) such that \(\Id-\widehat{\tilde K}_\infty(k)\) is invertible and \(\langle(\Id-\widehat{\tilde K}_\infty(k))^{-1} [\I_{W^d}]\rangle_W\neq0\), and define
\begin{equation}
\label{eq:BA-scalar-response-definitions}
\begin{aligned}
        s(k)
        &\coloneqq
        \frac1{
        \langle(\Id-\widehat{\tilde K}_\infty(k))^{-1}
        [\I_{W^d}]\rangle_W},
        \\
        C(k)
        &\coloneqq
        s(k)(\Id-\widehat{\tilde K}_\infty(k))^{-1}
        [\I_{W^d}]-\I_{W^d}.
\end{aligned}
\end{equation}
Then
\begin{equation}
\label{eq:BA-normalized-cell-equation}
        \langle C(k)\rangle_W=0,
        \qquad
        (\Id-\widehat{\tilde K}_\infty(k))
        [\I_{W^d}+C(k)]=s(k)\I_{W^d}.
\end{equation}
The second identity can equivalently be written as
\begin{equation}
\label{eq:BA-cell-equation-K-form}
        \widehat{\tilde K}_\infty(k)[\I_{W^d}+C(k)]
        =
        (1-s(k))\I_{W^d}+C(k).
\end{equation}
Moreover,
\begin{equation}
\label{eq:BA-trace-hard-resolvent}
        H(k)
        \coloneqq
        \left(
        \Id_{\operatorname{Ran}\mathsf Q}
        -
        \mathsf Q\widehat{\tilde K}_\infty(k)\mathsf Q
        \right)^{-1}
\end{equation}
is well defined.  The normalized cell is therefore characterized by
\begin{equation}
\label{eq:BA-cell-corrector-and-scalar-symbol}
\begin{aligned}
        C(k)
        &=
        H(k)\mathsf Q\widehat{\tilde K}_\infty(k)[\I_{W^d}],
        \\
        s(k)
        &=
        1-\left\langle
        \widehat{\tilde K}_\infty(k)[\I_{W^d}+C(k)]
        \right\rangle_W.
\end{aligned}
\end{equation}
For an arbitrary scalar \(\widehat a(k)\),
\begin{equation}
\label{eq:BA-block-scalar-response}
        (\Id-\widehat{\tilde K}_\infty(k))^{-1}
        [\widehat a(k)\I_{W^d}]
        =
        \frac{\widehat a(k)}{s(k)}
        \bigl(\I_{W^d}+C(k)\bigr).
\end{equation}
\end{lemma}

\begin{proof}
The identities in \cref{eq:BA-normalized-cell-equation} follow immediately from \cref{eq:BA-scalar-response-definitions}. Rearranging the second identity gives \cref{eq:BA-cell-equation-K-form}.

We next show that the operator defining \(H(k)\) is injective. Suppose \(X\in\operatorname{Ran}\mathsf Q\) satisfies
\[
        \bigl(
        \Id_{\operatorname{Ran}\mathsf Q}
        -
        \mathsf Q\widehat{\tilde K}_\infty(k)\mathsf Q
        \bigr)[X]=0.
\]
Since \(\mathsf Q[X]=X\), this implies
\[
        \mathsf Q(\Id-\widehat{\tilde K}_\infty(k))[X]=0.
\]
The kernel of \(\mathsf Q\) consists of scalar multiples of the identity.  Hence, for some \(c\in\mathbb C\),
\[
        (\Id-\widehat{\tilde K}_\infty(k))[X]=c\I_{W^d},
        \qquad
        X=c(\Id-\widehat{\tilde K}_\infty(k))^{-1}[\I_{W^d}].
\]
Taking normalized traces gives
\[
        0=\langle X\rangle_W
        =
        c\left\langle
        (\Id-\widehat{\tilde K}_\infty(k))^{-1}[\I_{W^d}]
        \right\rangle_W.
\]
The nonzero-trace assumption therefore implies \(c=0\), and thus \(X=0\).  Since \(\operatorname{Ran}\mathsf Q\) is finite dimensional, injectivity implies invertibility, proving that \(H(k)\) is well defined.

Applying \(\mathsf Q\) to the second identity in \cref{eq:BA-normalized-cell-equation}, and using \(\mathsf Q[C(k)]=C(k)\), yields
\[
        \bigl(
        \Id_{\operatorname{Ran}\mathsf Q}
        -
        \mathsf Q\widehat{\tilde K}_\infty(k)\mathsf Q
        \bigr)[C(k)]
        =
        \mathsf Q\widehat{\tilde K}_\infty(k)[\I_{W^d}].
\]
Applying \(H(k)\) proves the first identity in \cref{eq:BA-cell-corrector-and-scalar-symbol}. Taking the normalized trace in \cref{eq:BA-cell-equation-K-form} proves the second identity.

Finally, linearity and
\[
        (\Id-\widehat{\tilde K}_\infty(k))^{-1}[\I_{W^d}]
        =
        s(k)^{-1}\bigl(\I_{W^d}+C(k)\bigr)
\]
give \cref{eq:BA-block-scalar-response}.
\end{proof}
For simplicity, we shall use the abbreviation
\begin{equation}
\label{eq:BA-fixed-row-symbol}
        j_u(k)\coloneqq J_u(k)[\I_{W^d}+C(k)].
\end{equation}
The explicit GOE covariance admits the following decomposition into a strictly contracting matrix-valued part and a scalar multiple of the identity:
\begin{equation}
\label{eq:BA-renewal-split}
\begin{aligned}
        \tilde K_\infty(X)[B]&=
        \mathcal E_X[B]+\varphi_X(B)\I_{W^d},\\
        \mathcal E_X[B]
        &\coloneqq
        \frac1{W^d}
        \left(
        \tilde M_\infty(X)B
        \tilde M_\infty(X)^*
        \right)^\top,\\
        \varphi_X(B)
        &\coloneqq
        \frac1{W^d}
        \Tr\!\left(
        \tilde M_\infty(X)B
        \tilde M_\infty(X)^*
        \right).
\end{aligned}
\end{equation}
Write
\[
        \mathcal E(k)=\sum_Xe^{-\mathrm ik\cdot X}\mathcal E_X,
        \qquad
        \varphi_k=\sum_Xe^{-\mathrm ik\cdot X}\varphi_X.
\]
For \(B\succeq0\),
\begin{equation}
\label{eq:BA-renewal-order}
        0\preceq\mathcal E_X[B]
        \preceq\varphi_X(B)\I_{W^d}.
\end{equation}
The map \(\mathcal E_X\) maintains the complete matrix-valued cell profile, whereas the second term resets that profile to the scalar trace direction.  Ward's identity controls all matrix-valued profiles simultaneously through the following contraction. From now on, we write $\mathcal{R}_{\mathcal{E}}^{(k)}\coloneqq (\Id - \mathcal{E}(k))^{-1}$ and abbreviate $\mathcal{R}_{\mathcal{E}}^{(0)}$ to $\mathcal{R}_{\mathcal{E}}$.

\begin{lemma}[Ward mass and contraction of the cell sector]
\label{lemma:BA-Ward-mass}
Denote $\Gamma\coloneqq \Im\tilde\theta+\eta\I_{W^d}, \widehat\Gamma \coloneqq \Gamma/{\langle\Gamma\rangle_W}.$
then there exists a fixed \(\varrho<1\) such that
\begin{equation}
\label{eq:BA-internal-contraction}
        \mathcal E(0)[\widehat\Gamma]
        \preceq
        \varrho\widehat\Gamma.
\end{equation}
Consequently, $\mathcal R_{\mathcal E}= (\Id-\mathcal E(0))^{-1}$ is uniformly bounded in the order-unit norm generated by \(\widehat\Gamma\), and hence also in the matrix operator norm.

Let $U_0\coloneqq \mathcal R_{\mathcal E}[\I_{W^d}]$. Then
\begin{equation}
\label{eq:BA-renewal-U-comparison}
        c\I_{W^d}\preceq U_0\preceq C\I_{W^d},
        \qquad
        c\leq\langle U_0\rangle_W\leq C.
\end{equation}
For every \(B\) with \(\langle B\rangle_W=0\), the hard inverse in Lemma~\ref{lemma:BA-trace-Feshbach} has the exact formula
\begin{equation}
\label{eq:BA-renewal-hard-inverse}
        H(0)[B]
        = \mathcal R_{\mathcal E}[B]
        -U_0\frac{\langle\mathcal R_{\mathcal E}[B]\rangle_W
        }{\langle U_0\rangle_W}.
\end{equation}
In particular,
\begin{equation}
\label{eq:BA-derived-hard-resolvent}
        \|H(0)\|_{(2\to2)\to(2\to2)}
        \leq C.
\end{equation}

Define the infinite-volume Ward remainder by
\begin{equation}
\label{eq:BA-infinite-Ward-remainder}
        \mathcal R_{K,\infty}
        \coloneqq
        \mathscr S_W\!\left[
        \operatorname{bl}_{0,0}[\Im\tilde M_\infty]
        \right]-\Im\tilde\theta
\end{equation}
and denote
\[
        E_0
        \coloneqq
        (\Id-\widehat{\tilde K}_\infty(0))^{-1}
        [\mathcal R_{K,\infty}].
\]
Then
\begin{equation}
\label{eq:BA-infinite-Ward-cell-identity}
        \widehat{\tilde K}_\infty(0)[\Gamma]
        =
        \Im\tilde\theta+\mathcal R_{K,\infty}
\end{equation}
and
\begin{equation}
\label{eq:BA-exact-Ward-mass-and-cell}
\begin{aligned}
        s(0) = \frac{\eta}{
        \langle\Gamma\rangle_W+\langle E_0\rangle_W},\quad
        \I_{W^d}+C(0)=
        \frac{\Gamma+E_0}{
        \langle\Gamma\rangle_W+\langle E_0\rangle_W}.
\end{aligned}
\end{equation}
Equivalently,
\begin{equation}
\label{eq:BA-zero-renewal-formula}
        \I_{W^d}+C(0)
        =
        \frac{U_0}{\langle U_0\rangle_W},
        \qquad
        s(0)
        =
        \frac{1-\varphi_0(U_0)}
        {\langle U_0\rangle_W}.
\end{equation}
In particular,
\begin{equation}
\label{eq:BA-zero-cell-comparison}
        c\I_{W^d}
        \preceq
        \I_{W^d}+C(0)
        \preceq
        C\I_{W^d},
\end{equation}
and, in the diffusive window \(\eta\lesssim\lambda^{-1}\),
\begin{equation}
\label{eq:BA-Ward-soft-mass}
        s(0)\asymp\lambda\eta.
\end{equation}
\end{lemma}

\begin{proof}
The generalized Ward identity gives
\[
        \sum_X
        \tilde M_\infty(X)\Gamma
        \tilde M_\infty(X)^\star
        =
        \operatorname{bl}_{0,0}[\Im\tilde M_\infty].
\]
Applying \(\mathscr S_W\) proves \cref{eq:BA-infinite-Ward-cell-identity}.  We use this identity first to obtain the cell-sector contraction and then again to identify the soft mass.

For \(B\succeq0\), the matrix \(\tilde M_\infty(X)B\tilde M_\infty(X)^\star\) is positive semidefinite and is bounded above by its trace times the identity.  This proves \cref{eq:BA-renewal-order}.  In particular,
\[
        2\mathcal E(0)[\widehat\Gamma]
        \preceq
        \widehat{\tilde K}_\infty(0)[\widehat\Gamma].
\]
On the other hand, \cref{eq:BA-infinite-Ward-cell-identity} and the relative Ward-remainder estimate give
\[
\begin{aligned}
        \widehat{\tilde K}_\infty(0)[\widehat\Gamma]
        &=
        \frac{
        \Im\tilde\theta+\mathcal R_{K,\infty}
        }{
        \langle\Gamma\rangle_W
        }\\
        &\preceq
        \frac{(1+o(1))\Im\tilde\theta}{
        \langle\Gamma\rangle_W
        }
        \preceq
        (1+o(1))\widehat\Gamma.
\end{aligned}
\]
Thus $\mathcal E(0)[\widehat\Gamma]\preceq\left(1/2+o(1)\right)\widehat\Gamma,$ which proves \cref{eq:BA-internal-contraction}.

Conjugating \(\mathcal E(0)\) by \(\widehat\Gamma^{1/2}\) produces a positive map whose value at the identity is bounded by \(\varrho\I_{W^d}\).  The norm of a positive map on a matrix algebra equals the norm of its value at the identity.  Therefore the Neumann series for \(\mathcal R_{\mathcal E}\) converges in the corresponding order-unit norm.  Since \(\widehat\Gamma\asymp\I_{W^d}\), it is also uniformly bounded in the matrix operator norm.  Moreover, $U_0=\sum_{m\geq0}\mathcal E(0)^m[\I_{W^d}]\succeq\I_{W^d},$ while $\mathcal E(0)^m[\I_{W^d}]\preceq C\varrho^m\widehat\Gamma.$ This proves \cref{eq:BA-renewal-U-comparison}.

Fix \(\langle B\rangle_W=0\), and denote the right-hand side of \cref{eq:BA-renewal-hard-inverse} by \(X\).  Denote $a\coloneqq\langle\mathcal R_{\mathcal E}[B]\rangle_W/\langle U_0\rangle_W.$ Then
\[
        \langle X\rangle_W=0,
        \qquad
        (\Id-\mathcal E(0))[X]
        =
        B-a\I_{W^d}.
\]
Taking normalized traces in the second identity gives \(\langle\mathcal E(0)[X]\rangle_W=a\).  Since the scalar trace-reset term is annihilated by \(\mathsf Q\),
\[
\begin{aligned}
        \bigl(
        \Id-\mathsf Q\widehat{\tilde K}_\infty(0)\mathsf Q
        \bigr)[X]
        &=
        X-\mathsf Q\mathcal E(0)[X]\\
        &=
        X-\mathcal E(0)[X]
        +\langle\mathcal E(0)[X]\rangle_W\I_{W^d}\\
        &=B.
\end{aligned}
\]
This proves \cref{eq:BA-renewal-hard-inverse}; its norm bound follows from \cref{eq:BA-renewal-U-comparison}.

The Ward identity is equivalent to $(\Id-\widehat{\tilde K}_\infty(0))[\Gamma]=\eta\I_{W^d}-\mathcal R_{K,\infty}.$ The transferred order-unit contraction from Lemma~\ref{lemma:BA-T-error-and-kernel} gives the required inverse of \(\Id-\widehat{\tilde K}_\infty(0)\).  Its positive Neumann-series representation also gives $(\Id-\widehat{\tilde K}_\infty(0))^{-1}[\I_{W^d}]\succ0,$ so the normalization in Lemma~\ref{lemma:BA-trace-Feshbach} is valid. Applying the inverse to the last identity yields $\eta(\Id-\widehat{\tilde K}_\infty(0))^{-1}[\I_{W^d}]=\Gamma+E_0.$ Taking normalized traces and using the definitions of \(s(0)\) and \(C(0)\) proves \cref{eq:BA-exact-Ward-mass-and-cell}.

More explicitly, the deterministic remainder estimates give
\[
        \left\|
        \Gamma^{-1/2}\mathcal R_{K,\infty}\Gamma^{-1/2}
        \right\|_{2\to2}
        \leq
        o(1)\bigl(1+(\lambda\eta)^{-1}\bigr)^{-1},
\]
whereas the inverse costs at most \(C(1+(\lambda\eta)^{-1})\).  Hence $\left\|\Gamma^{-1/2}E_0\Gamma^{-1/2}\right\|_{2\to2}=o(1).$ Since \(\Gamma\asymp\lambda^{-1}\I_{W^d}\), this proves \cref{eq:BA-zero-cell-comparison,eq:BA-Ward-soft-mass}. Finally,
\[
        (\Id-\widehat{\tilde K}_\infty(0))[U_0]
        =
        \bigl(1-\varphi_0(U_0)\bigr)\I_{W^d}.
\]
After normalization, uniqueness in Lemma~\ref{lemma:BA-trace-Feshbach} gives \cref{eq:BA-zero-renewal-formula}.
\end{proof}

\begin{lemma}
\label{lemma:BA-fixed-row-moments}
Assume the deterministic bulk bounds and the hypotheses of Lemma~\ref{lemma:BA-Ward-mass} in either block-size regime.  Then, uniformly in \(u\in[W]^d\) and \(0\leq r\leq4\),
\begin{equation}
\label{eq:BA-fixed-row-moments}
        \sum_{X\in\mathbb Z^d}|X|^r
        \left\|
        \tilde M_\infty(X)^\star e_u
        \right\|_2^2
        \leq
        \begin{cases}
        C_{W,I}\lambda^{2r},
        & \text{in the fixed block size regime},\\[2mm]
        C\left(1+\lambda^2/W\right)^r,
        & \text{in the growing block size regime}.
        \end{cases}
\end{equation}
The constants in the large block size regime are uniform in \(W\).  Moreover,
\begin{equation}
\label{eq:BA-fixed-row-zero-mass}
        c\leq j_u(0)\leq C,
\end{equation}
where in the fixed block size regime the constants may depend on \(W\) and the folded bulk window \(I\).
\end{lemma}

\begin{proof}
The deterministic bounds already proved in the two regimes give
\begin{equation}
\label{eq:BA-row-moment-inputs}
        \|\Theta(\tilde\theta)\|_{2\to2}\leq\frac C\lambda,
        \qquad
        \eta\I+\Im\Theta(\tilde\theta)
        \succeq\frac c\lambda\I,
        \qquad
        \sup_x\left(
        \|\tilde M_\infty e_x\|_2
        +\|\tilde M_\infty^\star e_x\|_2
        \right)\leq C.
\end{equation}
For the growing block size regime, these are consequences of Corollary~\ref{coro:robust_app_bootstrap} and the same-self-energy finite-to-infinite comparison; for fixed \(W\), they follow from the fixed-height MDE bounds and the two generalized Ward identities.

Fix a microscopic site \(x\) and, for \(R>0\), let
\[
        \rho_{x,R}(y)\coloneqq
        \min\{|x-y|_1,R\},
        \qquad
        D_{x,R}e_y\coloneqq
        e^{\alpha\rho_{x,R}(y)}e_y,
        \qquad
        \alpha\coloneqq\frac{c_*}{\lambda^2+W}.
\]
The constant \(c_*>0\) will be chosen sufficiently small.  Since the hopping has finite range,
\[
        \left\|
        D_{x,R}\lambda\Psi D_{x,R}^{-1}-\lambda\Psi
        \right\|_{2\to2}
        \leq C\lambda\alpha.
\]
Indeed,
\[
        (D_{x,R}\lambda\Psi D_{x,R}^{-1}-\lambda\Psi)_{uv}
        =\lambda\Psi_{uv}
          \bigl(e^{\alpha(\rho_{x,R}(u)-\rho_{x,R}(v))}-1\bigr).
\]
This vanishes unless $|u-v|_1=1$, and in that case its absolute value is at most $C\lambda\alpha$.  The row and column sums are therefore bounded by $C\lambda\alpha$, and Schur's test proves the claim.

Since \(\Theta(\tilde\theta)\) is block diagonal and \(\|\Theta(\tilde\theta)\|_{2\to2}\leq C/\lambda\), it follows that
\[
        \left\|
        D_{x,R}\Theta(\tilde\theta)D_{x,R}^{-1}
        -\Theta(\tilde\theta)
        \right\|_{2\to2}
        \leq
        \frac{C\alpha W}{\lambda}e^{C\alpha W}.
\]
For fixed $x$ and a block $b$, define
\[
        \omega_b
        =\sup_{y\in[b]}\rho_{x,R}(y)-\inf_{y\in[b]}\rho_{x,R}(y),
        \qquad
        m_b
        =\frac{\sup_{y\in[b]}\rho_{x,R}(y)+\inf_{y\in[b]}\rho_{x,R}(y)}2,
\]
and
\[
        H_b=\operatorname{diag}_{y\in[b]}(\rho_{x,R}(y)-m_b).
\]
Then $\omega_b\leq d(W-1)$ and the scalar factor $e^{\alpha m_b}$ cancels under conjugation, so
\[
        D_{x,R}\big\rvert_{[b]}\tilde\theta
        D_{x,R}^{-1}\big\rvert_{[b]}
        =e^{\alpha H_b}\tilde\theta e^{-\alpha H_b}.
\]
By the fundamental theorem of calculus,
\[
        e^{\alpha H_b}\tilde\theta e^{-\alpha H_b}-\tilde\theta
        =\int_0^\alpha e^{tH_b}[H_b,\tilde\theta]e^{-tH_b}\,\rmd t.
\]
Since $\|H_b\|\leq\omega_b/2$, we obtain
\[
        \|e^{\alpha H_b}\tilde\theta e^{-\alpha H_b}-\tilde\theta\|
        \leq\|\tilde\theta\|(e^{\alpha\omega_b}-1).
\]
Taking the supremum over the blocks gives the corresponding direct-sum bound.  Consequently,
\[
        \|D_{x,R}(\lambda\Psi - z - \Theta(\tilde\theta)) D_{x,R}^{-1}-(\lambda\Psi - z - \Theta(\tilde\theta))\|_{2\to2}\leq\frac{Cc_*}{\lambda}.
\]
The dissipative bound in \cref{eq:BA-row-moment-inputs} gives \(\|\tilde M_\infty\|_{2\to2}\leq C\lambda\).  Choosing \(c_*\) sufficiently small, uniformly in \(R\), this yields
\[
    \left\|\tilde M_\infty \left(D_{x,R}(\lambda\Psi - z - \Theta(\tilde\theta)) D_{x,R}^{-1}-(\lambda\Psi - z - \Theta(\tilde\theta))\right)\right\|_{2\to2}\leq\frac{1}{2}.
\]
Using
\begin{equation}
    \begin{aligned}
        D_{x,R}\tilde M_\infty D_{x,R}^{-1}
        & = \left(D_{x,R}\tilde M_\infty^{-1}D_{x,R}^{-1}\right)^{-1}\\
        & = \left(\tilde M_\infty^{-1} + D_{x,R}(\lambda\Psi-z-\Theta(\tilde\theta))D_{x,R}^{-1}
        -(\lambda\Psi-z-\Theta(\tilde\theta))\right)^{-1}\\
        & = \left(\I+\tilde M_\infty\left[D_{x,R}(\lambda\Psi-z-\Theta(\tilde\theta))D_{x,R}^{-1}
        - (\lambda\Psi-z-\Theta(\tilde\theta))
        \right]
        \right)^{-1}\tilde M_\infty.
    \end{aligned}
\end{equation}
and \(D_{x,R}e_x=e_x\), we obtain from \cref{eq:BA-row-moment-inputs}
\[
        \|D_{x,R}\tilde M_\infty e_x\|_2
        +
        \|D_{x,R}\tilde M_\infty^\star e_x\|_2
        \leq C.
\]
The second term follows by applying the same argument to \((\lambda\Psi - z - \Theta(\tilde\theta))^\star\). Letting \(R\to\infty\),
\begin{equation}
\label{eq:BA-fixed-row-exponential-moment}
        \sum_y
        e^{2\alpha|x-y|_1}
        \left(
        |(\tilde M_\infty)_{yx}|^2
        +|(\tilde M_\infty)_{xy}|^2
        \right)
        \leq C.
\end{equation}

Take \(x=(0,u)\) and write \(y=(X,v)\).  The elementary cell geometry gives
\[
        \|X\|_2
        \leq
        \frac{\|WX+v-u\|_2}{W}+\sqrt d.
\]
Since \(t^r\leq C_r\alpha^{-r}e^{\alpha t}\), \cref{eq:BA-fixed-row-exponential-moment} implies
\begin{equation}
\label{eq:BA-fixed-row-moment-before-regimes}
\begin{aligned}
        \sum_X|X|^r
        \|\tilde M_\infty(X)^\star e_u\|_2^2
        &\leq
        C_r\left(1+(\alpha W)^{-r}\right)
        \\
        &\leq
        C_r\left(1+\frac{\lambda^2}{W}\right)^r.
\end{aligned}
\end{equation}
For fixed \(W\), the last expression is bounded by \(C_{W,I}\lambda^{2r}\).  This proves \cref{eq:BA-fixed-row-moments} for every growing block size. When \(W\leq\lambda^2\), the factor \(1+\lambda^2/W\) is comparable with \(\lambda^2/W\).

It remains to prove \cref{eq:BA-fixed-row-zero-mass}.  The normalized zero-momentum cell is already comparable with the identity by \cref{eq:BA-zero-cell-comparison}.  It therefore remains only to bound \(J_u(0)[\I_{W^d}]\) from below.  In the growing block size regime, Corollary~\ref{coro:robust_app_bootstrap} and the same-self-energy finite-to-infinite comparison give
\begin{equation}
\label{eq:BA-rowwise-imaginary-diagonal}
        \frac c\lambda
        \leq
        \bigl(
        \operatorname{bl}_{0,0}[\Im\tilde M_\infty]
        \bigr)_{uu}
        \leq
        \frac C\lambda.
\end{equation}
For fixed \(W\), the same estimate follows from the isolated chart in condition~\textnormal{\rm(C3)}.  Indeed, the free isolated projection has every diagonal entry equal to \(W^{-d}\).  Lemma~\ref{lemma:fh-isolated-persistence} shows that the perturbed Riesz projection \(P(s,\omega)\) differs from it by \(O_{W,I}(\lambda^{-2})\).  Thus
\[
        \langle e_u,P(s,\omega)e_u\rangle
        \geq\frac1{2W^d}
\]
for all sufficiently large \(\lambda\).  If \(\mu(s,\omega)\) denotes the corresponding persisted eigenvalue, the lower Loewner bound in \cref{eq:fh-Poisson-Loewner} gives
\[
        \left\langle e_u,
        \Im G_{\tilde\theta}(\Xi_e(s,\omega))e_u
        \right\rangle
        \geq
        \frac{c_{W,I}\lambda^{-2}}
        {\mu(s,\omega)^2+C_{W,I}\lambda^{-4}}.
\]
Moreover, Lemma~\ref{lemma:fh-isolated-persistence} gives \(|\partial_s\mu-1|\leq1/2\), and the range of \(\mu\) contains a fixed neighborhood of zero.  Changing variables from \(s\) to \(\mu(s,\omega)\), we obtain
\[
        \int_{-\delta_{\rm iso}}^{\delta_{\rm iso}}
        \left\langle e_u,
        \Im G_{\tilde\theta}(\Xi_e(s,\omega))e_u
        \right\rangle\,\rmd s
        \geq c_{W,I}.
\]
The lower Jacobian bound and \(|V_e|\geq v_{\rm iso}\), together with
\[
        \bigl(
        \operatorname{bl}_{0,0}[\Im\tilde M_\infty]
        \bigr)_{uu}
        =
        \frac1\lambda
        \int_{\mathbb T^d}
        \left\langle e_u,
        \Im G_{\tilde\theta}(\xi)e_u
        \right\rangle\,\rmd\xi,
\]
prove the lower bound in \cref{eq:BA-rowwise-imaginary-diagonal}; the upper bound is part of the deterministic trap.  Finally, the generalized Ward identity and \(\eta\I+\Im\tilde\theta\asymp\lambda^{-1}\I_{W^d}\) give
\[
        J_u(0)[\I_{W^d}]
        =
        \sum_X
        \|\tilde M_\infty(X)^\star e_u\|_2^2
        \asymp
        \lambda
        \bigl(
        \operatorname{bl}_{0,0}[\Im\tilde M_\infty]
        \bigr)_{uu}
        \asymp1.
\]
Applying the positive functional \(J_u(0)\) to \cref{eq:BA-zero-cell-comparison} proves \cref{eq:BA-fixed-row-zero-mass}.
\end{proof}
Define the $r$-th moment of the ladder kernel by
\[
        \mathfrak M_r
        \coloneqq
        \sum_{X\in\mathbb Z^d}|X|^r
        \|\tilde K_\infty(X)\|_{(2\to2)\to(2\to2)},
\]
This quantity enters the elliptic estimates below.
\begin{corollary}[Moments of the ladder kernel]
\label{coro:r_moment_bdd_2_2_K_tilde}
For \(0\leq r\leq4\),
\[
        \mathfrak M_r\leq C_{W,I}\lambda^{2r}
        \quad\text{if \(W\) is fixed},
        \qquad
        \mathfrak M_r\leq C(1+\lambda^2/W)^r
        \quad\text{if \(W\) grows}.
\]
\end{corollary}
\begin{proof}
    For any \(\|B\|_{2\to2}\leq1\),
    \begin{equation}
        \begin{aligned}
            \left\|\tilde K_\infty(X)[B]\right\|_{2\to 2} &= \left\|\mathscr S_W[\tilde M_\infty(X) B\tilde M_\infty(X)^\star]\right\|_{2\to 2}\\
            &\leq \frac{2}{W^d}\|\tilde M_\infty(X)\|_{\rm HS}^2.
        \end{aligned}
    \end{equation}
    Moreover,
    \[\|\tilde M_\infty(X)\|_{\rm HS}^2 = \sum_{1\leq u \leq W^d}\|\tilde M_\infty(X)^\star e_u\|_2^2,\]
    and therefore
    \begin{equation}
        \begin{aligned}
            \mathfrak M_r \leq \frac{2}{W^d}\sum_{1\leq u \leq W^d}\sum_{X \in \mathbb Z^d}|X|^r\|\tilde M_\infty(X)^\star e_u\|_2^2,
        \end{aligned}
    \end{equation}
    Lemma~\ref{lemma:BA-fixed-row-moments} proves the result.
\end{proof}
\begin{lemma}
\label{lemma:BA-low-frequency-hard}
There is a \(c_0>0\), depending only on the deterministic bulk constants, such that
\begin{equation}
\label{eq:BA-low-momentum-hard-resolvent}
        \sup_{
        |k|\leq
        c_0(1+\mathfrak M_4)^{-1/4}
        }
        \|H(k)\|_{(2\to2)\to(2\to2)}
        \leq C .
\end{equation}
The functions \(s(k)\) and \(C(k)\) from Lemma~\ref{lemma:BA-trace-Feshbach} are well defined throughout the same window.
\end{lemma}

\begin{proof}
For \(0\leq r\leq4\), H\"older's inequality and Corollary~\ref{coro:r_moment_bdd_2_2_K_tilde} give
\begin{equation}
\label{eq:BA-interpolated-kernel-moments}
\begin{aligned}
        \sum_X|X|^r\|\tilde K_\infty(X)\|
        &\leq
        \left(
        \sum_X\|\tilde K_\infty(X)\|
        \right)^{1-r/4}
        \mathfrak M_4^{r/4}
        \\
        &\leq
        C(1+\mathfrak M_4)^{r/4}.
\end{aligned}
\end{equation}
In particular,
\[
\begin{aligned}
        \|\widehat{\tilde{K}}_\infty(k)-\widehat{\tilde{K}}_\infty(0)\|
        &\leq
        \sum_X
        |e^{-\mathrm ik\cdot X}-1|
        \|\tilde K_\infty(X)\|
        \\
        &\leq
        C(1+\mathfrak M_4)^{1/4}|k|.
\end{aligned}
\]
On \(\operatorname{Ran}\mathsf Q\), the exact factorization
\[
\begin{aligned}
        \Id-\mathsf Q\widehat{\tilde{K}}_\infty(k)\mathsf Q
        =
        \bigl(\Id-\mathsf Q\widehat{\tilde{K}}_\infty(0)\mathsf Q\bigr)
        \Bigl[
        \Id
        -
        H(0)\mathsf Q(\widehat{\tilde{K}}_\infty(k)-\widehat{\tilde{K}}_\infty(0))\mathsf Q
        \Bigr]
\end{aligned}
\]
holds.  Since \(\|\mathsf Q\|\leq2\) in the matrix operator norm, \cref{eq:BA-derived-hard-resolvent} allows us to choose \(c_0\) so that the second factor differs from the identity by at most \(1/2\). It is then inverted by a Neumann series, and
\[
        H(k)
        =
        \Bigl[
        \Id
        -
        H(0)\mathsf Q(\widehat{\tilde{K}}_\infty(k)-\widehat{\tilde{K}}_\infty(0))\mathsf Q
        \Bigr]^{-1}H(0).
\]
This proves \cref{eq:BA-low-momentum-hard-resolvent}.

It remains only to justify the domain of \(s\) and \(C\).  The full operator \(\Id - \widehat{\tilde{K}}_\infty(k)\) is invertible by the order-unit input in Lemma~\ref{lemma:BA-T-error-and-kernel} and the finite-to-infinite-volume comparison, while the trace-free compression has just been proved invertible.  Relative to $\mathbb C\I_{W^d}\oplus\operatorname{Ran}\mathsf Q,$ the remaining Schur complement is a scalar.  It is nonzero because the full operator is invertible, and its inverse is exactly \(\langle(\Id - \widehat{\tilde{K}}_\infty(k))^{-1}[\I_{W^d}]\rangle_W\).  Hence this trace is nonzero at least for those \(k\) where \(\Id - \mathsf Q\widehat{\tilde K}_\infty(k)\mathsf Q\) is invertible, and \cref{eq:BA-scalar-response-definitions} applies throughout the claimed window.
\end{proof}

\begin{remark}[Size of the low-frequency window]
For \(k\neq0\), the phases \(e^{-\mathrm ik\cdot X}\) destroy positivity, so the cell-sector contraction is proved at \(k=0\) and transferred to low frequency by the norm perturbation above.  In particular, \(\mathfrak M_4\lesssim\lambda^8\) gives the fixed-\(W\) window \(|k|\lesssim\lambda^{-2}\), while \(\mathfrak M_4\lesssim(1+\lambda^2/W)^4\) gives $|k|\lesssim\min\!\left\{1,W/\lambda^2\right\}.$
\end{remark}

\begin{lemma}
\label{lemma:BA-superoperator-symbol}
Let \(a=(a_I)_{I\in\mathbb Z^d}\) be scalar and define $w[a]\coloneqq(\Id-\tilde{\mathcal K}_\infty)^{-1}(a_I\I_{W^d})_I.$ For \(x=(I,u)\), define
\begin{equation}
\label{eq:BA-direct-diagonal-lift}
        p_u[a](I)
        \coloneqq
        \sum_{X\in\mathbb Z^d}
        \left\langle e_u,
        \tilde{M}_\infty(X)w_{I-X}\tilde{M}_\infty(X)^\star e_u
        \right\rangle .
\end{equation}
Then, throughout the low-frequency window in Lemma~\ref{lemma:BA-low-frequency-hard},
\begin{equation}
\label{eq:BA-block-scalar-diagonal-transfer}
    \widehat p_u(k)=\frac{j_u(k)}{s(k)}\,\widehat a(k).
\end{equation}
\end{lemma}

\begin{proof}
Taking the block Fourier transform of the definition of \(w\) gives
\[
        \widehat w(k)
        =
        (\Id - \widehat{\tilde{K}}_\infty(k))^{-1}[\widehat a(k)\I_{W^d}].
\]
By \cref{eq:BA-block-scalar-response},
\[
        \widehat w(k)
        =
        \frac{\widehat a(k)}{s(k)}
        \bigl(\I_{W^d}+C(k)\bigr).
\]
On the other hand, Fourier transforming \cref{eq:BA-direct-diagonal-lift} and changing variables \(J=I-X\) gives
\[
\begin{aligned}
        \widehat p_u(k)
        &=
        \sum_Ie^{-\mathrm ik\cdot I}
        \sum_X
        \left\langle e_u,
        \tilde{M}_\infty(X)w_{I-X}\tilde{M}_\infty(X)^\star e_u
        \right\rangle
        \\
        &=
        \sum_Xe^{-\mathrm ik\cdot X}
        \left\langle e_u,
        \tilde{M}_\infty(X)\widehat w(k)\tilde{M}_\infty(X)^\star e_u
        \right\rangle
        \\
        &=
        J_u(k)[\widehat w(k)].
\end{aligned}
\]
Substituting the preceding formula for \(\widehat w(k)\) and using \cref{eq:BA-fixed-row-symbol} proves \cref{eq:BA-block-scalar-diagonal-transfer}.  The computation takes the \(u\)-th diagonal entry before the fibered superoperator is applied, which is precisely why no inversion of \(\mathscr S_W\) appears.
\end{proof}

\begin{lemma}
\label{lemma:BA-scalar-symbol-expansion}
There is a scalar $d_{\mathrm{BA}}\coloneqq\frac12\partial_1^2s(0)$ such that, for \(|k|\leq c_0(1+\mathfrak M_4)^{-1/4}\),
\begin{equation}
\label{eq:BA-scalar-symbol-expansion}
        s(k)
        =
        s(0)+d_{\mathrm{BA}}|k|^2
        +
        O\!\left((1+\mathfrak M_4)|k|^4\right),
        \qquad
        |d_{\mathrm{BA}}|
        \leq C(1+\mathfrak M_4)^{1/2}.
\end{equation}
\end{lemma}

\begin{proof}
Differentiate the normalized cell equation
\[
        (\Id-\widehat{\tilde K}_\infty(k))
        [\I_{W^d}+C(k)]
        =s(k)\I_{W^d},
        \qquad
        \langle C(k)\rangle_W=0.
\]
Taking the trace-free part at \(k=0\) gives the useful identity
\begin{equation}
        \partial_iC(0)
        =
        H(0)\mathsf Q
        \partial_i\widehat{\tilde K}_\infty(0)
        [\I_{W^d}+C(0)].
\end{equation}

The cubic covariance proved in \cref{eq:BA-scalar-symbol-cubic-invariance} gives \(s(gk)=s(k)\) for every signed coordinate permutation \(g\in O(d,\mathbb Z)\).  Coordinate reflections imply that every derivative containing an odd power of one coordinate vanishes at the origin.  Coordinate permutations imply that all diagonal entries of the Hessian agree.  Hence
\[
        Ds(0)=0,
        \qquad
        D^3s(0)=0,
        \qquad
        D^2s(0)=2d_{\mathrm{BA}}\I_d.
\]
The adjoint symmetry
\[
        \widehat{\tilde K}_\infty(-k)[B^*]
        =
        \widehat{\tilde K}_\infty(k)[B]^*
\]
shows that \(s(-k)=\overline{s(k)}\).  Since cubic symmetry also gives \(s(-k)=s(k)\), the scalar symbol and \(d_{\mathrm{BA}}\) are real for real \(k\) near the origin.

It remains to bound the derivatives.  From \cref{eq:BA-interpolated-kernel-moments}, for \(|\gamma|\leq4\),
\[
        \|\partial^\gamma
        \widehat{\tilde K}_\infty(k)\|
        \leq
        C_\gamma(1+\mathfrak M_4)^{|\gamma|/4}.
\]
Repeated differentiation of $H(k)=\bigl(\Id_{\ran\mathsf Q}-\mathsf Q\widehat{\tilde K}_\infty(k)\mathsf Q\bigr)^{-1}$ expresses every derivative of \(H\) as a finite sum of products of copies of \(H\) and derivatives of \(\widehat{\tilde K}_\infty\).  The low-frequency hard bound therefore gives $\|\partial^\gamma H(k)\|\leq C_\gamma(1+\mathfrak M_4)^{|\gamma|/4}.$ Using
\[
        C(k)=H(k)\mathsf Q
        \widehat{\tilde K}_\infty(k)[\I_{W^d}],
        \qquad
        s(k)=1-
        \left\langle
        \widehat{\tilde K}_\infty(k)[\I_{W^d}+C(k)]
        \right\rangle_W
\]
gives the same scale for the derivatives of \(C\) and \(s\). In particular,
\[
        |D^2s(0)|\leq C(1+\mathfrak M_4)^{1/2},
        \qquad
        \sup_k|D^4s(k)|
        \leq C(1+\mathfrak M_4)
\]
in the low-frequency window.  Fourth-order Taylor expansion completes the proof.
\end{proof}
The abstract trace Feshbach formula identifies the scalar denominator, but its differentiated form does not make the sign of the quadratic coefficient transparent.  The GOE renewal split \cref{eq:BA-renewal-split} does: every ladder step is a matrix-valued propagation followed, possibly, by a scalar trace reset.  The matrix-valued part is strictly contracting by Lemma~\ref{lemma:BA-Ward-mass}.  The next lemma uses this representation to compare the Hessian with a positive one-flight second moment.

\begin{lemma}
\label{lemma:BA-positive-one-flight}
For \(i\in\{1,\ldots,d\}\), set
\begin{equation}
\label{eq:BA-one-flight-variance}
        \mathcal M_i
        \coloneqq
        \frac1{W^d}
        \sum_{X\in\mathbb Z^d}
        X_i^2\|\tilde M_\infty(X)\|_{\mathrm{HS}}^2
        =
        \frac1{W^d\lambda^2}
        \int_{\mathbb T^d}
        \|\partial_iG_{\tilde\theta}(\xi)\|_{\mathrm{HS}}^2
        \,\rmd\xi .
\end{equation}
Under the hypotheses of Lemmas~\ref{lemma:BA-Ward-mass} and \ref{lemma:BA-scalar-symbol-expansion},
\begin{equation}
\label{eq:BA-diffusion-one-flight-comparison}
        c\mathcal M_i
        \leq
        \frac12\partial_i^2s(0)
        \leq
        C\mathcal M_i .
\end{equation}
\end{lemma}

\begin{proof}
The second identity in \cref{eq:BA-one-flight-variance} is Parseval's identity.  Write
\[
        U(k)\coloneqq
        (\Id-\mathcal E(k))^{-1}[\I_{W^d}]
        =\sum_Ze^{-\mathrm ik\cdot Z}U_Z.
\]
The contraction in Lemma~\ref{lemma:BA-Ward-mass} gives the convergent positive expansion
\begin{equation}
\label{eq:expan_phys_U_Z}
        U_Z=\mathbf1_{\{Z=0\}}\I_{W^d}
        +\sum_{m\geq1}
        \sum_{X_1+\cdots+X_m=Z}
        \mathcal E_{X_1}\circ\cdots\circ
        \mathcal E_{X_m}[\I_{W^d}].
\end{equation}
In particular $U_Z\succeq0$.  Define
\[
        b_Z\coloneqq\langle U_Z\rangle_W,
        \qquad
        p_Z\coloneqq\sum_{X+Y=Z}\varphi_X(U_Y).
\]
Both coefficient families are nonnegative and summable, and their Fourier transforms are $\langle U(k)\rangle_W$ and $\varphi_k(U(k))$, respectively.

We first establish their second moments, before differentiating. Let
\[
        \mathcal E_{2,i}\coloneqq\sum_XX_i^2\mathcal E_X,
        \qquad
        \varphi_{2,i}\coloneqq\sum_XX_i^2\varphi_X.
\]
By positivity, \cref{eq:BA-renewal-order}, and $\widehat\Gamma\asymp\I_{W^d}$,
\begin{equation}
\label{eq:BA-renewal-weighted-step}
\begin{aligned}
        \mathcal E_{2,i}[\widehat\Gamma]
        &\preceq C\mathcal M_i\widehat\Gamma,
        &\varphi_{2,i}(\widehat\Gamma)&\leq C\mathcal M_i,\\
        \mathcal E(0)^m[\widehat\Gamma]
        &\preceq\varrho^m\widehat\Gamma,
        &\varphi_0(\widehat\Gamma)&\leq C.
\end{aligned}
\end{equation}
These bounds can first be applied to truncated nonnegative spatial sums and then passed to the limit.

Let $p_Z^{(m)}$ be the contribution to $p_Z$ with $m$ internal $\mathcal E$-steps, and define $b_Z^{(m)}$ in the same way. For the $m+1$ spatial steps of a term in $p_Z^{(m)}$,
\[
        (X_{0,i}+\cdots+X_{m,i})^2
        \leq(m+1)\sum_{\ell=0}^mX_{\ell,i}^2.
\]
If the weighted step is the trace reset, summing the other steps leaves $\varphi_{2,i}\mathcal E(0)^m[\I_{W^d}]$, which is bounded by $C\mathcal M_i\varrho^m$.  If it is the internal step in position $\ell$, the corresponding sum is
\[
        \varphi_0\mathcal E(0)^\ell
        \mathcal E_{2,i}\mathcal E(0)^{m-1-\ell}[\I_{W^d}],
        \qquad 0\leq\ell<m,
\]
and \cref{eq:BA-renewal-weighted-step} bounds it by $C\mathcal M_i\varrho^{m-1}$.  Consequently,
\[
\begin{aligned}
        \sum_ZZ_i^2p_Z^{(m)}
        &\leq C\mathcal M_i(m+1)
             \bigl(\varrho^m+m\varrho^{m-1}\bigr),
             &&m\geq1,\\
        \sum_ZZ_i^2b_Z^{(m)}
        &\leq C\mathcal M_i m^2\varrho^{m-1},
             &&m\geq1.
\end{aligned}
\]
For $m=0$ the two moments are $\mathcal M_i$ and zero, respectively. Summation in $m$ therefore proves
\begin{equation}
\label{eq:BA-renewal-second-moments}
        \sum_ZZ_i^2p_Z+\sum_ZZ_i^2b_Z\leq C\mathcal M_i.
\end{equation}
This also justifies twice differentiating both scalar Fourier series by dominated convergence.

For $k$ close to zero, $\Id-\mathcal E(k)$ is invertible and $\langle U(k)\rangle_W\ne0$.  Moreover,
\[
        (\Id-\widehat{\tilde K}_\infty(k))[U(k)]
        =(1-\varphi_k(U(k)))\I_{W^d}.
\]
Uniqueness in Lemma~\ref{lemma:BA-trace-Feshbach} gives the exact identity
\begin{equation}
\label{eq:trace_renewal_identity}
        s(k)\langle U(k)\rangle_W=1-\varphi_k(U(k)).
\end{equation}
Since $Ds(0)=0$ by Lemma~\ref{lemma:BA-scalar-symbol-expansion}, twice differentiating yields
\begin{equation}
\label{eq:BA-short-positive-Hessian}
        \frac12\partial_i^2s(0)
        =\frac{\sum_ZZ_i^2p_Z+s(0)\sum_ZZ_i^2b_Z}
               {2\langle U(0)\rangle_W}.
\end{equation}
Here $\langle U(0)\rangle_W\asymp1$ and $0\leq s(0)\leq C$ by Lemma~\ref{lemma:BA-Ward-mass} and \cref{eq:BA-zero-renewal-formula}.  The contribution to $p_Z$ with no internal step is
\[
        \varphi_Z(\I_{W^d})
        =\frac1{W^d}\|\tilde M_\infty(Z)\|_{\mathrm{HS}}^2.
\]
All other contributions in \cref{eq:BA-short-positive-Hessian} are nonnegative.  This proves the lower bound; the upper bound follows from \cref{eq:BA-renewal-second-moments}.
\end{proof}
\begin{remark}
In this remark we highlight the intuition underneath the renewal indentity and why it captures the diffusion coefficient. Notice that the following equation holds
$$(\Id - \widehat{\tilde{K}}(k))^{-1}[\I_{W^d}] = (\Id - \mathcal{E}(k))^{-1}(\Id - \varphi_k(\Id - \mathcal{E}(k))^{-1})^{-1},$$
which therefore admits an expansion
    $$(\Id - \widehat{\tilde{K}}(k))^{-1}[\I_{W^d}] = \mathcal{R}^{(k)}_{\mathcal{E}}[\I_{W^d}] + \mathcal{R}^{(k)}_{\mathcal{E}}\varphi_k\mathcal{R}^{(k)}_{\mathcal{E}}[\I_{W^d}] + \mathcal{R}^{(k)}_{\mathcal{E}}\varphi_k\mathcal{R}^{(k)}_{\mathcal{E}}\varphi_k\mathcal{R}^{(k)}_{\mathcal{E}}[\I_{W^d}] + \cdots$$
This representation means that after going through every memory layer, the direction would be immediately reset to the direction of identity. Every $\varphi_k\mathcal{R}_{\mathcal{E}}[\I_{W^d}]$ generate a factor of $\varphi_k(U_k)$, and reset the direction back to identity. Repeated resets enter through
\[
        (\Id-\widehat{\tilde K}_\infty(k))^{-1}[\I_{W^d}]
        =\frac{U(k)}{1-\varphi_k(U(k))}
        =U(k)\sum_{m\geq0}\bigl(\varphi_k(U(k))\bigr)^m.
\]
For $\eta>0$ in the stated window, this series converges near zero, since $\varphi_0(U(0))<1$ by \cref{eq:BA-zero-renewal-formula}.  Each scalar factor records a complete excursion ending in a reset, after which the internal matrix profile starts again from $\I_{W^d}$.  The final factor $U(k)$ allows propagation after the last reset.  Taking the normalized trace gives \cref{eq:trace_renewal_identity}, the renewal identity. Therefore, intuitively, the coefficients $U_Z$ describe propagation through internal $\mathcal E$-steps without a trace reset, and $p_Z$ describes such an excursion followed by one reset.
\end{remark}
\begin{prop}
\label{prop:BA-cell-coercivity}
Assume the hypotheses of Lemmas~\ref{lemma:BA-Ward-mass} and \ref{lemma:BA-T-error-and-kernel}.  Let \(z=E+\mathrm i\eta\), where \(0<\eta\leq\lambda^{-1}\), and set \(e=E/\lambda\).  Suppose either
\begin{enumerate}[(a)]
\item \(W<\infty\) is fixed and \(e\) belongs to a finite-height folded bulk window \(I\); or
\item for fixed \(\varepsilon_{\mathrm{bulk}}>0\),
\begin{equation}
\label{eq:BA-cell-coercivity-large-W-range}
        e\in(0,4d),
        \quad
        \operatorname{dist}(e,\Sigma_d^\Psi)
        \geq\varepsilon_{\mathrm{bulk}},
        \quad
        W^{d/p}\geq C_0(1+\log\lambda)^4.
\end{equation}
\end{enumerate}
In case \textnormal{(a)}, and in case \textnormal{(b)} when \(W\leq\lambda^2\),
\begin{equation}
\label{eq:BA-cell-coercivity}
        c\frac{\lambda^4}{W^2}
        \leq
        d_{\mathrm{BA}}
        \leq
        C\frac{\lambda^4}{W^2}.
\end{equation}
More generally, throughout case \textnormal{(b)},
\begin{equation}
\label{eq:BA-cell-coercivity-crossover}
        c\left(\frac{\lambda^4}{W^2}+\frac{\lambda^2}{W}\right)
        \leq d_{\mathrm{BA}}\leq
        C\left(\frac{\lambda^4}{W^2}+\frac{\lambda^2}{W}\right).
\end{equation}
In case \textnormal{(a)}, the constants may depend on \(W\) and \(I\); in case \textnormal{(b)}, they are uniform in \(W\).
\end{prop}

\begin{proof}
By Lemma~\ref{lemma:BA-positive-one-flight}, it remains to estimate \(\mathcal M_i\) in the two regimes.

\emph{Fixed block size.} The fixed-block MDE trap gives
\[
        c_{W,I}\lambda^{-2}\I_{W^d}
        \preceq
        \frac{\eta\I_{W^d}+\Im\tilde\theta}{\lambda}
        \preceq
        C_{W,I}\lambda^{-2}\I_{W^d}.
\]
The dissipative resolvent bound, Ward's identity, and $ \partial_iG_{\tilde\theta}=-G_{\tilde\theta}(\partial_i\widehat\Psi_W)G_{\tilde\theta}$ give
\[
        \frac1{W^d}
        \int_{\mathbb T^d}
        \|\partial_iG_{\tilde\theta}(\xi)\|_{\mathrm{HS}}^2
        \,\rmd\xi
        \leq C_{W,I}\lambda^6.
\]

For the matching lower bound, use the isolated chart in Definition~\ref{def:fh-bulk-window}.  The scalar-pole decomposition \cref{eq:fh-patch-clean-decomposition,eq:fh-patch-zeta-regular-clean} persists for \(\tilde\theta\).  If \(s_*(\omega)\) is the zero of the real part of the scalar denominator, then on
\[
        |s-s_*(\omega)|
        \leq(c_0\lambda^2)^{-1}
\]
we have
\[
        |\zeta(s,\omega)|\leq C\lambda^{-2},
        \qquad
        |\partial_s\zeta(s,\omega)|\geq\frac12.
\]
Differentiating the rank-one pole \(P G_{\tilde\theta}P=\zeta^{-1}P\) shows that
\[
        \|\partial_sG_{\tilde\theta}
        (\Xi_e(s,\omega))\|_{\mathrm{HS}}
        \geq c_{W,I}\lambda^4.
\]
The interval has length comparable with \(\lambda^{-2}\).  The chart bounds therefore imply
\[
        \sum_{j=1}^d
        \int_{\mathbb T^d}
        \|\partial_jG_{\tilde\theta}(\xi)\|_{\mathrm{HS}}^2
        \,\rmd\xi
        \geq c_{W,I}\lambda^6.
\]
Cubic symmetry makes the \(d\) summands equal, so
\[
        \frac1{W^d}
        \int_{\mathbb T^d}
        \|\partial_iG_{\tilde\theta}(\xi)\|_{\mathrm{HS}}^2
        \,\rmd\xi
        \asymp_{W,I}\lambda^6.
\]
Consequently,
\begin{equation}
\label{eq:BA-fixed-W-one-flight-short}
        \mathcal M_i
        \asymp_{W,I}\lambda^4
        \asymp_{W,I}\frac{\lambda^4}{W^2}.
\end{equation}

\emph{Growing block size.} Write
\[
        g=\langle\tilde\theta\rangle_W,\qquad
        V=\frac{\tilde\theta-g\I_{W^d}}{\lambda},\qquad
        \zeta=\frac{z+g}{\lambda},\qquad
        G_g(\xi)=(\widehat\Psi_W(\xi)-\zeta)^{-1}.
\]
We call $G_g(\xi)$ the normalized scalar fiber resolvent. We also write $M_g(X) = \blk_{X,0}[(\lambda\Psi - z - g)^{-1}]$. We measure the block displacement from a starting point $v$ using the function $X(t, v) = \lfloor(v+t)/W\rfloor$. The average squared block displacement is given by $Q_W^{(2)}(t) = W^{-1}\sum_{v = 0}^{W-1}X(t, v)^2$. The scalarization estimates imply
\[
        \Im\zeta\asymp\lambda^{-2},\qquad
        |\Re\zeta-e|\leq C\lambda^{-2},\qquad
        \|V\|\leq\varepsilon_*\lambda^{-2},
\]
where the fixed constant \(\varepsilon_*\) can be chosen sufficiently small by increasing the lower-bound constant for \(W\). Since \(\|G_g(\xi)\|\leq C\lambda^2\), the matrices
\[
        A(\xi)=(\I-G_g(\xi)V)^{-1},\qquad
        B(\xi)=(\I-VG_g(\xi))^{-1}
\]
and their inverses are uniformly bounded.  Differentiating the resolvent identity, using the fact that \(V\) is independent of \(\xi\), gives the exact factorization
\begin{equation}
\label{eq:BA-cell-derivative-relative-factorization}
        \partial_iG_{\tilde\theta}(\xi)
        =A(\xi)\,\partial_iG_g(\xi)\,B(\xi).
\end{equation}
Consequently the Hilbert--Schmidt norms of the two first derivatives are comparable, uniformly in \(W\) and \(\xi\).  It remains to compute the scalar one-flight moment in the actual block coordinate.

Use normalized Haar measure \(\rmd\mu(p)=(2\pi)^{-d}\rmd p\) and set
\[
        m(r)=\frac1\lambda\int_{\mathbb T^d}
        \frac{e^{\mathrm i p\cdot r}}{\Omega(p)-\zeta}
        \,\rmd\mu(p),\qquad r\in\mathbb Z^d.
\]
The scalar resolvent has block entries \(M_g(X)_{uv}=m(WX+u-v)\).  For fixed \(r\), summing over the \(W^d\) possible positions of one endpoint inside its block gives the exact counting identity
\begin{equation}
\label{eq:BA-exact-cell-crossing-second-moment}
        \frac1{W^d}\sum_X X_i^2\|M_g(X)\|_{\rm HS}^2
        =\sum_{r\in\mathbb Z^d}|m(r)|^2Q_W(r_i),
\end{equation}
where, for \(|t|=qW+s\), \(q\in\mathbb N_0\) and \(0\leq s<W\),
\begin{equation}
\label{eq:BA-exact-cell-crossing-weight}
\begin{aligned}
        Q_W^{(2)}(t)
        &=\left(1-\frac{s}{W}\right)q^2
          +\frac{s}{W}(q+1)^2\\
        &=\frac{t^2}{W^2}
          +\frac{s}{W}\left(1-\frac{s}{W}\right).
\end{aligned}
\end{equation}
Indeed, the increment in the \(i\)-th block coordinate has absolute value \(q\) for \(W-s\) endpoint positions and \(q+1\) for the remaining \(s\) positions.  In particular,
\begin{equation}
\label{eq:BA-cell-crossing-weight-comparison}
        \frac12\left(\frac{t^2}{W^2}+\frac{|t|}{W}\right)
        \leq Q_W^{(2)}(t)
        \leq\frac{t^2}{W^2}+\frac{|t|}{W}.
\end{equation}
For the lower bound, \(Q_W^{(2)}(t)\geq t^2/W^2\), and also \(Q_W^{(2)}(t)\geq |t|/W\): the latter is an equality when \(q=0\), while for \(q\geq1\) it follows from \(q^2+(2q+1)s/W\geq q+s/W\).

For \(j\geq0\), denote \(\mu_{j,i}=\sum_r|r_i|^j|m(r)|^2\), with \(\mu_0\) for \(j=0\). Parseval and the coarea formula on a fixed regular bulk neighborhood give
\begin{equation}
\label{eq:BA-scalar-microscopic-flight-moments}
        \mu_0\asymp1,\qquad
        \mu_{2,i}\asymp\lambda^4,\qquad
        \mu_{4,i}\leq C\lambda^8.
\end{equation}
Here
\[
        \mu_{2,i}
        =\frac1{\lambda^2}\int_{\mathbb T^d}
        \frac{|\partial_i\Omega(p)|^2}{|\Omega(p)-\zeta|^4}
        \,\rmd\mu(p).
\]
The level-set integral of \(|\partial_i\Omega|^2/|\nabla\Omega|\) is uniformly positive on this bulk neighborhood.  This follows from cubic symmetry and \(\sum_i|\partial_i\Omega|^2=|\nabla\Omega|^2\), since the levels are regular and nonempty.  Thus the last integral is comparable to \(\lambda^{-2}(\Im\zeta)^{-3}\asymp\lambda^4\). Similarly, differentiating \((\Omega-\zeta)^{-1}\) twice bounds \(\mu_{4,i}\) by
\[
        \frac C{\lambda^2}\int_{\mathbb T^d}
        \left(|\Omega(p)-\zeta|^{-6}
              +|\Omega(p)-\zeta|^{-4}\right)\,\rmd\mu(p)
        \leq C\lambda^8.
\]
The parts away from the regular neighborhood have bounded denominators and contribute \(O(\lambda^{-2})\). Finally, Cauchy--Schwarz and H\"older give, respectively,
\[
        \mu_{1,i}\leq\mu_0^{1/2}\mu_{2,i}^{1/2},\qquad
        \mu_{2,i}\leq\mu_{1,i}^{2/3}\mu_{4,i}^{1/3},
\]
and hence \(\mu_{1,i}\asymp\lambda^2\). Combining these bounds with \cref{eq:BA-exact-cell-crossing-second-moment,eq:BA-cell-crossing-weight-comparison} proves
\[
        \frac1{W^d}\sum_X X_i^2\|M_g(X)\|_{\rm HS}^2
        \asymp\frac{\lambda^4}{W^2}+\frac{\lambda^2}{W}.
\]
The derivative factorization \cref{eq:BA-cell-derivative-relative-factorization} and Parseval transfer this comparison to \(\mathcal M_i\).  Applying Lemma~\ref{lemma:BA-positive-one-flight} completes the proof.
\end{proof}

\begin{remark}
The block-coordinate correction is a positive boundary-crossing contribution.  When \(|t|<W\), the exact weight in \cref{eq:BA-exact-cell-crossing-weight} is \(|t|/W\). Consequently the cell curvature is of order \(\lambda^4/W^2\) for \(W\ll\lambda^2\), but of order \(\lambda^2/W\) for \(W\gg\lambda^2\).  In microscopic coordinates, division by the Ward time factor \(\lambda\) gives a diffusion scale of order \(\lambda^3+\lambda W\).  The genuine large-block dynamical consequence requires the fourth-moment estimate proved in Subsection~\ref{subsec:BA-genuine-large-blocks}.
\end{remark}

\begin{lemma}
\label{lemma:BA-row-elliptic-multiplier}
Define
\[
        j_u^{\mathrm{ev}}(k)
        \coloneqq
        \frac{j_u(k)+j_u(-k)}2.
\]
In the fixed-block regime, for \(|k|\leq c_{W,I}\lambda^{-2}\),
\begin{equation}
\label{eq:BA-fixed-W-row-multiplier-comparison}
        \left|
        \frac{j_u^{\mathrm{ev}}(k)}{s(k)}
        -
        \frac{j_u(0)}
        {s(0)+d_{\mathrm{BA}}|k|^2}
        \right|
        \leq C_{W,I}.
\end{equation}
In the growing-block regime, for \(|k|\leq cW/\lambda^2\),
\begin{equation}
\label{eq:BA-growing-W-row-multiplier-comparison}
        \left|
        \frac{j_u^{\mathrm{ev}}(k)}{s(k)}
        -
        \frac{j_u(0)}
        {s(0)+d_{\mathrm{BA}}|k|^2}
        \right|
        \leq C
\end{equation}
uniformly in \(W\) and \(u\).
\end{lemma}

\begin{proof}
Combining the derivative bounds for \(J_u(k)\) from Lemma~\ref{lemma:BA-fixed-row-moments} with those for \(C(k)\) from Lemma~\ref{lemma:BA-scalar-symbol-expansion}, and applying Leibniz's rule to \(j_u(k)=J_u(k)[\I_{W^d}+C(k)]\), gives
\[
        \sup_k|D^2j_u(k)|
        \leq
        \begin{cases}
        C_{W,I}\lambda^4,&W\text{ fixed},\\[2pt]
        C\lambda^4/W^2,&W\text{ growing}.
        \end{cases}
\]
The adjoint symmetry \(j_u(-k)=\overline{j_u(k)}\) implies that the odd linear term cancels, so
\[
        j_u^{\mathrm{ev}}(k)
        =
        j_u(0)+
        \begin{cases}
        O_{W,I}(\lambda^4|k|^2),&W\text{ fixed},\\[2pt]
        O(\lambda^4W^{-2}|k|^2),&W\text{ growing}.
        \end{cases}
\]
On the same windows, Lemma~\ref{lemma:BA-scalar-symbol-expansion} and Proposition~\ref{prop:BA-cell-coercivity} give
\[
        |s(k)|
        \geq
        c\bigl(s(0)+d_{\mathrm{BA}}|k|^2\bigr)
\]
and
\[
        s(k)
        =
        s(0)+d_{\mathrm{BA}}|k|^2+
        \begin{cases}
        O_{W,I}(\lambda^8|k|^4),&W\text{ fixed},\\[2pt]
        O(\lambda^8W^{-4}|k|^4),&W\text{ growing}.
        \end{cases}
\]
Since \(j_u(0)\asymp1\), subtracting the two fractions proves the claims.
\end{proof}

\begin{prop}
\label{prop:BA-diagonal-elliptic-approximation}
Define
\begin{equation}
\label{eq:BA-elliptic-operator}
        \mathcal L_{\mathrm{BA}}
        \coloneqq
        s(0)-d_{\mathrm{BA}}\Delta_{\mathbb R^d}.
\end{equation}
Let \(f\in C_c^\infty(\mathbb R^d)\) be real and even, let \(a_I=f(I)\), and let \(p_u[a]\) be defined by \cref{eq:BA-direct-diagonal-lift}.  In the fixed-block regime,
\begin{equation}
\label{eq:BA-fixed-W-elliptic-approximation}
\begin{aligned}
        \left|
        p_u[a](0)
        -
        j_u(0)\mathcal L_{\mathrm{BA}}^{-1}f(0)
        \right|
        \leq
        C_{W,I}\left(
        \frac{\lambda^2}{s(0)}
        \|\nabla f\|_{L^\infty}
        +
        \|\widehat f\|_{L^1}
        \right).
\end{aligned}
\end{equation}
In the growing-block regime,
\begin{equation}
\label{eq:BA-growing-W-elliptic-approximation}
\begin{aligned}
        \left|
        p_u[a](0)
        -
        j_u(0)\mathcal L_{\mathrm{BA}}^{-1}f(0)
        \right|
        \leq
        C\left(
        \frac{\lambda^2}{Ws(0)}
        \|\nabla f\|_{L^\infty}
        +
        \|\widehat f\|_{L^1}
        \right).
\end{aligned}
\end{equation}
\end{prop}

\begin{proof}
Choose a real even Schwartz function \(\chi\) such that
\[
        \int\chi=1,
        \qquad
        \operatorname{supp}\widehat\chi\subset\{|k|\leq1\},
\]
and denote \(\chi_\ell(x)=\ell^{-d}\chi(x/\ell)\), with
\[
        \ell=C_{W,I}\lambda^2
        \quad\text{for fixed \(W\)},\qquad
        \ell=C\lambda^2/W
        \quad\text{for growing \(W\)}.
\]
The constants are chosen so that the Fourier support of \(\chi_\ell*f\) lies in the corresponding window of Lemma~\ref{lemma:BA-row-elliptic-multiplier}.  We have
\begin{equation}
\label{eq:BA-smoothing-error}
        \|f-\chi_\ell*f\|_{L^\infty}
        \leq
        C\ell\|\nabla f\|_{L^\infty}.
\end{equation}

Since the map \(a\mapsto p_u[a]\) is positive, therefore
\[
        \|a\mapsto p_u[a]\|_{\ell^\infty\to\ell^\infty}
        =
        \frac{j_u(0)}{s(0)}
        \leq\frac C{s(0)}.
\]
The heat-kernel representation of \(\mathcal L_{\mathrm{BA}}^{-1}\) gives
\[
        \|\mathcal L_{\mathrm{BA}}^{-1}\|_{L^\infty\to L^\infty}
        \leq\frac1{s(0)}.
\]
Thus replacing \(f\) by \(\chi_\ell*f\) costs at most \(C\ell s(0)^{-1}\|\nabla f\|_\infty\) on both sides.

For the mollified function, Poisson summation has no aliasing because its Fourier transform is supported inside the fundamental cell. Using \cref{eq:BA-block-scalar-diagonal-transfer} at the observation block \(I=0\) gives
\[
        p_u[\operatorname{Sam}(\chi_\ell*f)](0)
        =
        \frac1{(2\pi)^d}
        \int
        \frac{j_u(k)}{s(k)}
        \widehat{\chi_\ell*f}(k)\,\rmd k.
\]
Both \(s\) and \(\widehat{\chi_\ell*f}\) are even.  Hence the odd part of \(j_u\) integrates to zero, and the numerator may be replaced by \(j_u^{\mathrm{ev}}\).  Fourier inversion for the continuum elliptic operator gives
\[
        j_u(0)\mathcal L_{\mathrm{BA}}^{-1}
        (\chi_\ell*f)(0)
        =
        \frac1{(2\pi)^d}
        \int
        \frac{j_u(0)}
        {s(0)+d_{\mathrm{BA}}|k|^2}
        \widehat{\chi_\ell*f}(k)\,\rmd k.
\]
The multiplier bounds in \cref{eq:BA-fixed-W-row-multiplier-comparison,eq:BA-growing-W-row-multiplier-comparison} and \(\|\widehat{\chi_\ell*f}\|_1\leq C\|\widehat f\|_1\) bound the difference by \(C\|\widehat f\|_1\).  Combining this with the smoothing error and the two choices of \(\ell\) proves the result.
\end{proof}

\begin{remark}
\label{rmk:BA-block-versus-microscopic-elliptic-operator}
The operator in \cref{eq:BA-elliptic-operator} is written in the block coordinate.  If \(y\) denotes the continuous block coordinate, \(x=Wy\) is the microscopic coordinate, and \(a_I=F(WI)\), then
\[
        \bigl(s(0)-d_{\mathrm{BA}}\Delta_y\bigr)F(Wy)
        =
        \bigl(s(0)-W^2d_{\mathrm{BA}}\Delta_x\bigr)F(x)
        \big|_{x=Wy}.
\]
Since
\[
        s(0)\asymp\lambda\eta,
        \qquad
        W^2d_{\mathrm{BA}}\asymp\lambda^4,
\]
the physical elliptic operator has the same scaling in the fixed and growing block size regimes.  A fixed internal coordinate can produce an odd numerator, but it cancels for the even radial observables used in the dynamical proof.
\end{remark}
\subsection{From the resolvent profile to dynamics}
\label{subsec:BA-resolvent-to-dynamics}
Up to this point the argument has been deterministic, at fixed energy, and on the infinite block lattice.  It remains to make three transfers. We compare finite- and infinite-volume ladder responses, specialize the elliptic profile to small and large radial cutoffs, and finally integrate over energy.  The small cutoff shows that a fixed fraction of the wave packet has left a ball of radius \(c\lambda^{3/2}\sqrt T\); the large cutoff shows that almost no mass lies beyond \(C\lambda^{3/2}\sqrt T\).

We begin with the finite- to infinite-volume ladder comparison needed for the diagonal response.
\begin{coro}
\label{coro:BA-fixed-row-random-profile}
Assume the hypotheses of Theorem~\ref{thm:rand_deter_ladder_apprx} and Proposition~\ref{prop:BA-cell-coercivity}, and suppose that \(0<\lambda\eta\leq1\).  Let
\[
        \mathcal A
        =
        \bigoplus_{I\in\tilde{\mathbb Z}_n^d}
        a_I\I_{W^d},
        \qquad
        0\leq a_I\leq1,
\]
Assume moreover that either $a_I=0$ for every $|I|>n/10$, or $a_I=1$ for every $|I|>n/10$.  Let \(p_u\) be the infinite-volume response in \cref{eq:BA-direct-diagonal-lift}.  If
\begin{equation}
\label{eq:BA-profile-finite-volume-condition}
        \frac{n\eta}{\lambda}\geq(\log\lambda)^2,
\end{equation}
then, for every fixed \(D>0\),
\begin{equation}
\label{eq:BA-fixed-row-random-profile}
\begin{aligned}
        \sup_{|I|\leq n/10,\,u}
        \left|
        (R\mathcal A R^\star)_{(I,u),(I,u)}-p_u(I)
        \right|
        \prec
        \left(1+(\lambda\eta)^{-1}\right)
        \Phi_d^{\rm BA}(z)+\lambda^{-D}.
\end{aligned}
\end{equation}
\end{coro}

\begin{proof}
We write \(\tilde{\mathcal K}_L\) for the finite-volume kernel and set
\[
        \mathbf h_L
        =
        (\Id-\tilde{\mathcal K}_L)^{-1}
        \tilde{\mathcal K}_L\mathcal A.
\]
The elementary identity
\begin{equation}
\label{eq:BA-source-plus-ladder}
\begin{aligned}
        \mathcal A+\mathbf h_L
        &=
        \left(
        \Id+(
        \Id-\tilde{\mathcal K}_L)^{-1}
        \tilde{\mathcal K}_L
        \right)\mathcal A \\
        &=
        (\Id-\tilde{\mathcal K}_L)^{-1}
        \left(\Id-\tilde{\mathcal K}_L+
        \tilde{\mathcal K}_L\right)\mathcal A
        =
        (\Id-\tilde{\mathcal K}_L)^{-1}\mathcal A
\end{aligned}
\end{equation}
shows that the deterministic lift in \cref{eq:BA-deterministic-lift} is
\[
        P_L
        =
        \tilde M_L
        (\Id-\tilde{\mathcal K}_L)^{-1}\mathcal A
        \tilde M_L^\star.
\]
Its \((I,u)\)-th diagonal entry is the finite-volume version \(p_{u,L}(I)\) of \cref{eq:BA-direct-diagonal-lift}.  Indeed, expanding the two block convolutions gives
\[
        (P_L)_{(I,u),(I,u)}
        =
        \sum_{X\in\tilde{\mathbb Z}_n^d}
        \left\langle e_u,
        \tilde M_L(X)w_{L,I-X}
        \tilde M_L(X)^\star e_u
        \right\rangle,
        \qquad
        w_L=(\Id-\tilde{\mathcal K}_L)^{-1}\mathcal A.
\]
Thus \cref{eq:BA-raw-diffusive-profile} gives
\begin{equation}
\label{eq:BA-random-to-finite-p}
        \sup_{|I|\leq n/10,\,u}
        \left|
        (R\mathcal A R^\star)_{(I,u),(I,u)}-p_{u,L}(I)
        \right|
        \prec
        \left(1+(\lambda\eta)^{-1}\right)
        \Phi_d^{\rm BA}(z).
\end{equation}

It remains to compare the two deterministic responses.  We write \(\tilde M_\infty\) for the infinite-volume resolvent with the same self-energy as \(\tilde M_L\).  Poisson summation gives the exact periodization identity
\begin{equation}
\label{eq:BA-M-exact-periodization}
        \tilde M_L(X)
        =
        \sum_{m\in\mathbb Z^d}
        \tilde M_\infty(X+nm),
        \qquad X\in\tilde{\mathbb Z}_n^d.
\end{equation}
Indeed, inserting the Fourier series of \(\tilde M_\infty\) into the finite Bloch sum leaves precisely those Fourier coefficients whose indices are congruent to \(X\) modulo \(n\).  The Combes--Thomas bound in Lemma~\ref{lemma:CT_est_block_renorm_resolv} therefore yields, for the centered representative of \(X\),
\begin{equation}
\label{eq:BA-M-finite-infinite-difference}
\begin{aligned}
        \left\|
        \tilde M_L(X)-\tilde M_\infty(X)
        \right\|_{2\to2}
        &\leq
        \sum_{m\neq0}
        \left\|\tilde M_\infty(X+nm)\right\|_{2\to2}\\
        &\leq
        \frac C\eta
        \exp\left(-\frac{cn\eta}{d\lambda}\right).
\end{aligned}
\end{equation}

To place both ladder operators on the same space, define the periodized infinite-volume kernel on the block torus by
\begin{equation}
\label{eq:BA-periodized-infinite-kernel}
        \tilde K_\infty^{(n)}(X)[B]
        \coloneqq
        \sum_{m\in\mathbb Z^d}
        \tilde K_\infty(X+nm)[B].
\end{equation}
Expanding \cref{eq:BA-M-exact-periodization} in the definition of \(\tilde K_L\) gives the formula
\begin{equation}
\label{eq:BA-finite-kernel-cross-images}
\begin{aligned}
        (\tilde K_L-\tilde K_\infty^{(n)})(X)[B]=
        \sum_{m\neq m'}
        \mathscr S_W\!\left[
        \tilde M_\infty(X+nm)B
        \tilde M_\infty(X+nm')^\star
        \right].
\end{aligned}
\end{equation}
If \(m\neq m'\), then at least one of \(|X+nm|\) and \(|X+nm'|\) is at least \(n/2\).  Applying the Combes--Thomas estimate to both factors gives
\[
\begin{aligned}
        &\sum_{X\in\tilde{\mathbb Z}_n^d}
        \|\tilde K_L(X)-\tilde K_\infty^{(n)}(X)\|\leq
        \frac C{\eta^2}
        \sum_X\sum_{m\neq m'}
        \exp\left[-\frac{c\eta}{d\lambda}
        \left(|X+nm|+|X+nm'|\right)\right].
\end{aligned}
\]
For the factor whose argument is at least \(n/2\), we use
\[
        e^{-c\eta|X+nm|/(d\lambda)}
        \leq
        e^{-cn\eta/(4d\lambda)}
        e^{-c\eta|X+nm|/(2d\lambda)}.
\]
The map \((X,m)\mapsto X+nm\), with \(X\) in the centered fundamental domain, is a bijection onto \(\mathbb Z^d\).  Therefore after splitting according to which argument is at least \(n/2\), we may drop the restriction \(m\neq m'\) and estimate each part by
\[
\begin{aligned}
        \left(
        \sum_{X,m}e^{-c\eta|X+nm|/(2d\lambda)}
        \right)
        \left(
        \sup_X\sum_{m'}e^{-c\eta|X+nm'|/(2d\lambda)}
        \right)\leq
        \left(
        \sum_{Y\in\mathbb Z^d}
        e^{-c\eta|Y|/(2d\lambda)}
        \right)^2
        \leq C(1+\lambda/\eta)^{2d}.
\end{aligned}
\]
Together with the extracted exponential factor, this proves
\begin{equation}
\label{eq:BA-periodized-kernel-comparison}
\begin{aligned}
        \sum_{X\in\tilde{\mathbb Z}_n^d}
        \left\|
        \tilde K_L(X)-\tilde K_\infty^{(n)}(X)
        \right\|_{(2\to2)\to(2\to2)}
        \leq
        Cn^d(\lambda+\eta^{-1})^C
        \exp\left(-\frac{cn\eta}{d\lambda}\right).
\end{aligned}
\end{equation}
Thus the two inverses may now be compared on the same block torus by the exact identity
\begin{equation}
\label{eq:BA-kernel-inverse-finite-infinite}
\begin{aligned}
        &(\Id-\tilde{\mathcal K}_L)^{-1}
        -(\Id-\tilde{\mathcal K}_\infty^{(n)})^{-1}\\
        &\qquad=
        (\Id-\tilde{\mathcal K}_L)^{-1}
        (\tilde{\mathcal K}_L-
        \tilde{\mathcal K}_\infty^{(n)})
        (\Id-\tilde{\mathcal K}_\infty^{(n)})^{-1}.
\end{aligned}
\end{equation}
Both inverse norms are \(O(1+(\lambda\eta)^{-1})\).  For the periodized kernel this follows from the same order-unit calculation as in \cref{eq:BA-T-kernel-contraction}, since
\[
        \sum_{X\in\tilde{\mathbb Z}_n^d}
        \tilde K_\infty^{(n)}(X)[\Im\tilde\theta]
        =
        \sum_{Y\in\mathbb Z^d}
        \tilde K_\infty(Y)[\Im\tilde\theta].
\]
Therefore \cref{eq:BA-periodized-kernel-comparison,eq:BA-kernel-inverse-finite-infinite} give, with \(w_L=(\Id-\tilde{\mathcal K}_L)^{-1}\mathcal A\) and \(w_n=(\Id-\tilde{\mathcal K}_\infty^{(n)})^{-1}\mathcal A\),
\begin{equation}
\label{eq:BA-finite-periodized-response-comparison}
\begin{aligned}
        \|w_L-w_n\|_{\ell^\infty(2\to2)}
        \leq
        C\left(1+(\lambda\eta)^{-1}\right)^2
        n^d(\lambda+\eta^{-1})^C
        \exp\left(-\frac{cn\eta}{d\lambda}\right).
\end{aligned}
\end{equation}

We next compare the latter response with the zero-extension response in \cref{eq:BA-direct-diagonal-lift}.  This requires a tail estimate, which we now derive in the same order-unit norm as the contraction in \cref{eq:BA-T-kernel-contraction}.  If \(S\subset\mathbb Z^d\), then, for \(0\leq\mu\leq c\eta/\lambda\),
\begin{equation}
\label{eq:BA-weighted-ladder-inverse}
\begin{aligned}
        \sup_Ze^{\mu\operatorname{dist}(Z,S)}
        \left\|
        \mathcal C_{\Im\tilde\theta}
        \left[
        \bigl((\Id-\tilde{\mathcal K}_\infty)^{-1}f\bigr)_Z
        \right]
        \right\|_{2\to2}\leq
        C\left(1+(\lambda\eta)^{-1}\right)
        \sup_Ze^{\mu\operatorname{dist}(Z,S)}
        \left\|
        \mathcal C_{\Im\tilde\theta}[f_Z]
        \right\|_{2\to2}.
\end{aligned}
\end{equation}
We prove this estimate.  The explicit formula for the fibered superoperator first gives
\begin{equation}
\label{eq:BA-kernel-from-HS-row-mass}
        \|\tilde K_\infty(X)\|_{(2\to2)\to(2\to2)}
        \leq
        \frac{2}{W^d}\|\tilde M_\infty(X)\|_{\rm HS}^2.
\end{equation}
Set \(\alpha=c/(\lambda^2+W)\).  For \(x=(0,u)\) and \(y=(X,v)\), the cell geometry gives
\[
        |WX+v-u|_1\geq W|X|_1-dW.
\]
Expanding the norm over \(v\), applying \cref{eq:BA-fixed-row-exponential-moment} to the translated microscopic row, and then using the preceding inequality give
\[
        \sum_Xe^{\alpha W|X|}
        \|\tilde M_\infty(X)^\star e_u\|_2^2
        \leq Ce^{\alpha dW}.
\]
Therefore, after summing over \(u\in[W]^d\), \cref{eq:BA-kernel-from-HS-row-mass} implies, with \(\mu_0=\alpha W\),
\[
\begin{aligned}
        \sum_Xe^{\mu_0|X|}
        \|\tilde K_\infty(X)\|
        &\leq
        \frac{C}{W^d}\sum_{u=1}^{W^d}
        \sum_Xe^{\alpha W|X|}
        \|\tilde M_\infty(X)^\star e_u\|_2^2\\
        &\leq Ce^{\alpha dW}\leq C.
\end{aligned}
\]
Consequently, for \(0\leq\mu\leq\mu_0/2\),
\[
\begin{aligned}
        \sum_X(e^{\mu|X|}-1)
        \|\tilde K_\infty(X)\|
        &\leq
        C\frac{\mu}{\mu_0}.
\end{aligned}
\]
Here we used \(e^{\mu t}-1\leq\mu t e^{\mu t}\) and \(t e^{-\mu_0t/2}\leq C/\mu_0\).

For \(R>0\), denote
\[
        \rho_{S,R}(Z)
        \coloneqq
        \min\{\operatorname{dist}(Z,S),R\},
        \qquad
        (D_{\mu,S,R}f)_Z
        \coloneqq e^{\mu\rho_{S,R}(Z)}f_Z,\qquad \tilde{\mathcal{K}}_{\mu, R} \coloneqq D_{\mu, S, R}\tilde{\mathcal{K}}_\infty D_{\mu, S, R}^{-1}.
\]
Directly,
\[
        [({\tilde{\mathcal K}}_{\mu,R}-\tilde{\mathcal K}_\infty)f]_Z
        =\sum_X\tilde K_\infty(X)[f_{Z-X}]
          \left(e^{\mu(\rho_{S,R}(Z)-\rho_{S,R}(Z-X))}-1\right).
\]
The truncated multiplier is bounded on \(\ell^\infty\), the triangle inequality
\[
        \rho_{S,R}(Z)-\rho_{S,R}(Z-X)\leq|X|
\]
and the comparability \(c\lambda^{-1}\I_{W^d}\preceq\Im\tilde\theta \preceq C\lambda^{-1}\I_{W^d}\) show that

\begin{align}
\label{eq:BA-weight-conjugated-kernel-difference}
        \left\|
        \tilde{\mathcal K}_{\mu, R}
        -\tilde{\mathcal K}_\infty
        \right\|_{(\infty,\Im\tilde\theta)
        \to(\infty,\Im\tilde\theta)}\leq
        C\sum_X(e^{\mu|X|}-1)
        \|\tilde K_\infty(X)\|
        \leq C\frac\mu{\mu_0}.
\end{align}
Now
\[D_{\mu, S, R}(\Id - \tilde{\mathcal{K}}_\infty)^{-1}f = (\Id - \tilde{\mathcal{K}}_{\mu, R})^{-1}D_{\mu, S, R}f.\]
The unweighted contraction in \cref{eq:BA-T-kernel-contraction} has gap at least \(c\lambda\eta\) when \(\lambda\eta\leq1\).  For every \(W\geq1\),
\[
        \frac{\eta}{\lambda}\frac{\lambda^2+W}{W}
        =\frac{\lambda\eta}{W}+\frac{\eta}{\lambda}
        \leq2\lambda\eta.
\]
Thus choosing \(\mu=c\eta/\lambda\) sufficiently small makes \(C\mu/\mu_0\leq c\lambda\eta/2\), uniformly in \(W\).  The Neumann series in this spatially weighted order-unit norm gives a bound uniform in \(R\).
\begin{equation}
    \begin{aligned}
        \|D_{\mu, S, R}(\Id - \tilde{\mathcal{K}}_\infty)^{-1}f\|_{(\infty, \Im\tilde{\theta})}&\leq \|(\Id - \tilde{\mathcal{K}}_{\mu, R})^{-1}\|_{(\infty, \Im\tilde\theta)}\|D_{\mu, S, R}f\|_{(\infty, \Im\tilde\theta)}\\
        &\leq C(1+(\lambda\eta)^{-1})\|D_{\mu, S, R}f\|_{(\infty, \Im\tilde\theta)}.
    \end{aligned}
\end{equation}
From the definition of the weighted norm, we have
\[ \|D_{\mu, S, R}(\Id - \tilde{\mathcal{K}}_\infty)^{-1}f\|_{(\infty, \Im\tilde{\theta})} = \sup_Z e^{\mu \rho_{S,R}(Z)}\|\mathcal{C}_{\Im\tilde\theta}[((\Id - \tilde{\mathcal{K}}_\infty)^{-1}f)_Z]\|_{2\to 2},\]
and
\[\|D_{\mu, S, R}f\|_{(\infty, \Im\tilde\theta)} = \sup_Z e^{\mu \rho_{S,R}(Z)}\|\mathcal{C}_{\Im\tilde\theta}[f_Z]\|_{2\to 2}.\]
Letting \(R\to\infty\) and using monotone convergence proves \cref{eq:BA-weighted-ladder-inverse}.

Let \(a^{(0)}\) be the zero extension of the supported source and let \(a^{\rm per}\) be its \(n\)-periodic extension.  Restricting \((\Id-\tilde{\mathcal K}_\infty)^{-1} [a^{\rm per}\I_{W^d}]\) to the centered fundamental domain \(\mathcal{F}_n \coloneqq [-n/2, n/2)^d\) gives exactly \((\Id-\tilde{\mathcal K}_\infty^{(n)})^{-1} [a\I_{W^d}]\).  Let \(S\) be the support of \(a^{(0)}\) and denote
\[
        S_*\coloneqq\bigcup_{m\neq0}(S+nm).
\]
The sequence \((a^{\rm per}-a^{(0)})\I_{W^d}\) is supported in \(S_*\).  Since
\[
        \|\mathcal C_{\Im\tilde\theta}[\I_{W^d}]\|_{2\to2}
        \leq C\lambda,
        \qquad
        \|B\|_{2\to2}
        \leq\frac C\lambda
        \|\mathcal C_{\Im\tilde\theta}[B]\|_{2\to2},
\]
the two factors of \(\lambda\) cancel when \cref{eq:BA-weighted-ladder-inverse} is applied by choosing \(S=S_*\). Consequently, for every \(Z\in\mathbb Z^d\),
\begin{equation}
\label{eq:BA-periodic-zero-source-spatial-decay}
\begin{aligned}
        &\left\|
        (\Id-\tilde{\mathcal K}_\infty)^{-1}
        [(a^{\rm per}-a^{(0)})\I_{W^d}]_Z
        \right\|_{2\to2}\leq
        C\left(1+(\lambda\eta)^{-1}\right)
        e^{-\mu\operatorname{dist}(Z,S_*)}.
\end{aligned}
\end{equation}
If \(|I|\leq n/10\), then \(\operatorname{dist}(I,S_*)\geq cn\).  Taking \(\mu=c\eta/\lambda\) in the last display gives
\begin{equation}
\label{eq:BA-periodic-zero-source-comparison}
\begin{aligned}
        &\left\|
        (\Id-\tilde{\mathcal K}_\infty)^{-1}
        [(a^{\rm per}-a^{(0)})\I_{W^d}]_I
        \right\|_{2\to2}\leq
        C\left(1+(\lambda\eta)^{-1}\right)
        \exp\left(-\frac{cn\eta}{\lambda}\right).
\end{aligned}
\end{equation}
For \(|I|\leq n/10\), the triangle inequality gives
\[
        |X|+\operatorname{dist}(I-X,S_*)\geq cn.
\]
Thus \cref{eq:BA-periodic-zero-source-spatial-decay} and \cref{eq:BA-fixed-row-exponential-moment} give the explicit lift bound with $\mathcal L_\infty\coloneqq (\Id-\tilde{\mathcal K}_\infty)^{-1}$:
\begin{align}
\label{eq:BA-periodic-zero-lift-comparison}
        \sum_X
        \|\tilde M_\infty(X)^\star e_u\|_2^2
        \left\|
        \mathcal L_\infty
        [(a^{\rm per}-a^{(0)})\I_{W^d}]_{I-X}
        \right\|_{2\to2}
        &\leq
        C\left(1+(\lambda\eta)^{-1}\right)
        e^{-\mu cn}
        \sum_Xe^{\mu|X|}
        \|\tilde M_\infty(X)^\star e_u\|_2^2
        \nonumber\\
        &\leq
        C\left(1+(\lambda\eta)^{-1}\right)
        \exp\left(-\frac{cn\eta}{\lambda}\right).
\end{align}
We also compare the finite lift with the periodized lift.  First, the Ward bound and \cref{eq:BA-finite-periodized-response-comparison} give
\[
\begin{aligned}
        \left|
        \sum_X\left\langle e_u,
        \tilde M_L(X)(w_{L,I-X}-w_{n,I-X})
        \tilde M_L(X)^\star e_u
        \right\rangle
        \right|&\leq
        \|w_L-w_n\|_{\ell^\infty(2\to2)}
        \sum_X\|\tilde M_L(X)^\star e_u\|_2^2\\
        &\leq
        C\|w_L-w_n\|_{\ell^\infty(2\to2)}.
\end{aligned}
\]
Viewing \(w_n\) as an \(n\)-periodic sequence, expanding the lift constructed from \cref{eq:BA-M-exact-periodization} gives
\[
\begin{aligned}
        &\sum_X\left\langle e_u,
        \tilde M_L(X)w_{n,I-X}\tilde M_L(X)^\star e_u
        \right\rangle-
        \sum_{Y\in\mathbb Z^d}\left\langle e_u,
        \tilde M_\infty(Y)w_{n,I-Y}
        \tilde M_\infty(Y)^\star e_u
        \right\rangle\\
        &=
        \sum_X\sum_{m\neq m'}
        \left\langle e_u,
        \tilde M_\infty(X+nm)w_{n,I-X}
        \tilde M_\infty(X+nm')^\star e_u
        \right\rangle.
\end{aligned}
\]
Moreover,
\[
        \|w_n\|_{\ell^\infty(2\to2)}
        \leq
        C\left(1+(\lambda\eta)^{-1}\right).
\]
Applying Combes--Thomas to both resolvent factors, exactly as in the derivation of \cref{eq:BA-periodized-kernel-comparison}, therefore gives
\[
\begin{aligned}
        \left|
        \sum_X\sum_{m\neq m'}
        \left\langle e_u,
        \tilde M_\infty(X+nm)w_{n,I-X}
        \tilde M_\infty(X+nm')^\star e_u
        \right\rangle
        \right|
        &\leq
        \frac{C}{\eta^2}
        \left(1+(\lambda\eta)^{-1}\right)
        \\
        &\qquad{}\times
        \sum_X\sum_{m\neq m'}
        e^{-c\eta(|X+nm|+|X+nm'|)/(d\lambda)}\\
        &\leq
        C(\lambda+\eta^{-1}+n)^C
        \exp\left(-\frac{cn\eta}{d\lambda}\right).
\end{aligned}
\]
Combining this estimate with \cref{eq:BA-finite-periodized-response-comparison,eq:BA-periodic-zero-lift-comparison} therefore yields
\[
        \sup_{|I|\leq n/10,\,u}|p_{u,L}(I)-p_u(I)|
        \leq
        C(\lambda+\eta^{-1}+W+n)^C
        \exp\left(-\frac{cn\eta}{d\lambda}\right).
\]

If instead \(a_I=1\) for every \(|I|>n/10\), apply the preceding argument to \(1-a\). For the remaining constant source, first note that
\[
        \widehat{\tilde K}_\infty^{(n)}(0)
        =
        \sum_{X\in\tilde{\mathbb Z}_n^d}
        \sum_{m\in\mathbb Z^d}\tilde K_\infty(X+nm)
        =\widehat{\tilde K}_\infty(0).
\]
Let
\[
        B_L=(\Id-\widehat{\tilde K}_L(0))^{-1}[\I_{W^d}],
        \qquad
        B_\infty=(\Id-\widehat{\tilde K}_\infty(0))^{-1}[\I_{W^d}].
\]
The zero-mode version of \cref{eq:BA-kernel-inverse-finite-infinite}, followed by \cref{eq:BA-periodized-kernel-comparison}, gives
\begin{equation}
\label{eq:BA-constant-cell-finite-infinite}
        \|B_L-B_\infty\|_{2\to2}
        \leq
        C(\lambda+\eta^{-1}+n)^C
        \exp\left(-\frac{cn\eta}{d\lambda}\right).
\end{equation}
The infinite cell equation gives
\[
        B_\infty
        =\frac{\I_{W^d}+C(0)}{s(0)}.
\]
Hence the infinite response is \(p_u[1]=j_u(0)/s(0)\).  On the other hand, the finite response is
\[
        p_{u,L}[1]
        =
        \sum_X\left\langle e_u,
        \tilde M_L(X)B_L\tilde M_L(X)^\star e_u
        \right\rangle.
\]
The Ward bound and \cref{eq:BA-constant-cell-finite-infinite} control the change from \(B_L\) to \(B_\infty\).  Expanding \cref{eq:BA-M-exact-periodization} then gives
\[
\begin{aligned}
        &\sum_X\left\langle e_u,
        \tilde M_L(X)B_\infty\tilde M_L(X)^\star e_u
        \right\rangle-p_u[1]\\
        &\quad=
        \sum_X\sum_{m\neq m'}
        \left\langle e_u,
        \tilde M_\infty(X+nm)B_\infty
        \tilde M_\infty(X+nm')^\star e_u
        \right\rangle.
\end{aligned}
\]
Since \(\|B_\infty\|_{2\to2}\leq C(1+(\lambda\eta)^{-1})\), the double geometric-sum estimate used above proves
\[
        \left|p_{u,L}[1]-p_u[1]\right|
        \leq
        C(\lambda+\eta^{-1}+n)^C
        \exp\left(-\frac{cn\eta}{d\lambda}\right).
\]
Under \cref{eq:BA-profile-finite-volume-condition}, all the preceding errors are at most \(\lambda^{-D}\) for every fixed \(D\), since \(n,W,\eta^{-1}\) are polynomial in \(\lambda\) and \(e^{-c(\log\lambda)^2}\) dominates every inverse power.  Combining this estimate with \cref{eq:BA-random-to-finite-p} proves the corollary.
\end{proof}
We now pass from the fixed-energy resolvent profile to Theorem~\ref{thm:BA-dynamical-delocalization}.  Fix a radial function \(\chi\in C_c^\infty(\mathbb R^d)\) such that
\begin{equation}
\label{eq:BA-radial-cutoff}
        0\leq\chi\leq1,
        \qquad
        \chi=1\quad\text{on }B_1,
        \qquad
        \operatorname{supp}\chi\subset B_2.
\end{equation}
For \(r>0\), denote
\begin{equation}
\label{eq:BA-radial-block-observable}
        \mathcal A_{r,\eta}
        \coloneqq
        \bigoplus_{I\in\tilde{\mathbb Z}_n^d}
        \chi\left(
        \frac{WI}{r\lambda^{3/2}\eta^{-1/2}}
        \right)\I_{W^d}.
\end{equation}
Here and below, \(I\in\tilde{\mathbb Z}_n^d\) is represented in the centered fundamental domain \([-n/2,n/2)^d\). This source is constant inside each block.  The next lemma is the only spatial resolvent estimate needed for the dynamical argument.

\begin{lemma}
\label{lemma:BA-radial-resolvent-profile}
Let \(x=(0,u)\), let \(0<\lambda\eta\leq1\), and suppose that \(E/\lambda\) belongs to a compact regular energy set in the applicable regime.  Assume \cref{eq:BA-profile-finite-volume-condition} and also that
\begin{equation}
\label{eq:BA-radial-support-condition}
        \frac{2r\lambda^{3/2}\eta^{-1/2}}{W}\leq\frac n{10}.
\end{equation}
Uniformly in \(E\), there are nonnegative random errors \(\mathcal E_{r,\eta}^{\rm in}(E)\) and \(\mathcal E_{r,\eta}^{\rm out}(E)\) such that
\begin{equation}
\label{eq:BA-small-ball-resolvent}
\begin{aligned}
        \lambda\eta
        (R\mathcal A_{r,\eta}R^\star)_{xx}
        \leq C\bigg[
        r^2\left(1+\bm 1_{\{d=2\}}|\log r|\right)
        +\frac{\sqrt{\lambda\eta}}{r}
        +\lambda\eta
        \bigg]
        +\mathcal E_{r,\eta}^{\rm in}(E)
\end{aligned}
\end{equation}
for \(0<r\leq1\), whereas
\begin{equation}
\label{eq:BA-large-ball-resolvent}
\begin{aligned}
        \lambda\eta
        \left(
        R(\I-\mathcal A_{r,\eta})R^\star
        \right)_{xx}
        \leq C\bigg[
        e^{-cr}
        +\frac{\sqrt{\lambda\eta}}{r}
        +\lambda\eta
        \bigg]
        +\mathcal E_{r,\eta}^{\rm out}(E)
\end{aligned}
\end{equation}
for \(r\geq1\).  In the fixed block size regime the constants may depend on \(W\) and on the finite collection of folded bulk windows.
\begin{equation}
\label{eq:BA-radial-random-error}
        \sup_{E/\lambda\in\mathcal G}
        \left(
        \mathcal E_{r,\eta}^{\rm in}(E)
        +\mathcal E_{r,\eta}^{\rm out}(E)
        \right)
        \prec
        \sup_{E/\lambda\in\mathcal G}
        \Phi_d^{\rm BA}(E+\mathrm i\eta)+\lambda^{-D}
\end{equation}
on every compact regular set \(\mathcal G\), for every fixed \(D>0\).
\end{lemma}

\begin{proof}
Define a function on the continuous block coordinate by
\[
        f_{r,\eta}(y)
        \coloneqq
        \chi\left(
        \frac{Wy}{r\lambda^{3/2}\eta^{-1/2}}
        \right).
\]
The chain rule and the scaling rule for the Fourier transform give
\begin{equation}
\label{eq:BA-radial-test-derivatives}
\begin{aligned}
        \|\nabla f_{r,\eta}\|_{L^\infty}
        \leq
        \frac{CW}{r\lambda^{3/2}\eta^{-1/2}},\quad
        \|\widehat f_{r,\eta}\|_{L^1}=
        \|\widehat\chi\|_{L^1}.
\end{aligned}
\end{equation}
Indeed, if \(a=r\lambda^{3/2}\eta^{-1/2}/W\), then \(\widehat f_{r,\eta}(k)=a^d\widehat\chi(ak)\), and the second identity follows after rescaling the integration variable by \(a\).

The function \(f_{r,\eta}\) is radial and therefore even, so Proposition~\ref{prop:BA-diagonal-elliptic-approximation} applies directly.  In the growing block size regime, that proposition and \cref{eq:BA-radial-test-derivatives} imply
\[
\begin{aligned}
        \left|
        p_u(0)-j_u(0)\mathcal L_{\rm BA}^{-1}f_{r,\eta}(0)
        \right|
        &\leq
        C\left[
        \frac{\lambda^2}{Ws(0)}
        \frac{W}{r\lambda^{3/2}\eta^{-1/2}}+1
        \right]\\
        &=
        C\left[
        \frac{\sqrt{\lambda\eta}}{r s(0)}+1
        \right].
\end{aligned}
\]
Multiplication by \(\lambda\eta\), followed by \(s(0)\asymp\lambda\eta\), yields
\begin{equation}
\label{eq:BA-radial-elliptic-error}
        \lambda\eta
        \left|
        p_u(0)-j_u(0)\mathcal L_{\rm BA}^{-1}f_{r,\eta}(0)
        \right|
        \leq
        C\left(
        \frac{\sqrt{\lambda\eta}}r+\lambda\eta
        \right).
\end{equation}
For fixed \(W\), the factor \(W\) resulting from \cref{eq:BA-radial-test-derivatives} is absorbed into the constant, and the same bound follows from \cref{eq:BA-fixed-W-elliptic-approximation}.

We next derive the two estimates for the continuum response.  The heat kernel representation of the massive Green function is
\begin{equation}
\label{eq:BA-massive-heat-kernel}
\begin{aligned}
        s(0)\mathcal L_{\rm BA}^{-1}f(0)
        ={}&
        \int_0^\infty e^{-t}
        \int_{\mathbb R^d}
        \frac{e^{-|z|^2/(4t)}}{(4\pi t)^{d/2}}
        f\left(\sqrt{\frac{d_{\rm BA}}{s(0)}}z\right)
        \,\rmd z\,\rmd t.
\end{aligned}
\end{equation}
To verify it, start from
\[
        \mathcal L_{\rm BA}^{-1}
        =
        \int_0^\infty e^{-s(0)v}e^{v d_{\rm BA}\Delta}\,\rmd v,
\]
and use \(t=s(0)v\) and \(y=\sqrt{d_{\rm BA}/s(0)}z\).  Moreover,
\begin{equation}
\label{eq:BA-radial-scale-ratio}
        \frac{
        r\lambda^{3/2}\eta^{-1/2}/W
        }{
        \sqrt{d_{\rm BA}/s(0)}}
        \asymp r
\end{equation}
by \(s(0)\asymp\lambda\eta\) and \(d_{\rm BA}\asymp\lambda^4/W^2\).

For \(0<r\leq1\), the integrand in \cref{eq:BA-massive-heat-kernel} vanishes unless \(|z|\leq Cr\). The Gaussian integral over this ball is bounded by \(\min\{1,Cr^dt^{-d/2}\}\).  Splitting the \(t\)-integral at \(t=r^2\) therefore gives
\begin{align}
        s(0)\mathcal L_{\rm BA}^{-1}f_{r,\eta}(0)
        &\leq
        Cr^2+Cr^d\int_{r^2}^1t^{-d/2}\,\rmd t+Cr^d
        \nonumber\\
        &\leq
        Cr^2\left(1+\bm 1_{\{d=2\}}|\log r|\right).
\label{eq:BA-massive-small-ball}
\end{align}
For \(r\geq1\), the fact that \(\chi=1\) on \(B_1\) and \cref{eq:BA-radial-scale-ratio} show that the integrand for \(1-s(0)\mathcal L_{\rm BA}^{-1}f_{r,\eta}(0)\) is supported in \(|z|\geq cr\); here we also used \(s(0)\mathcal L_{\rm BA}^{-1}1=1\).  A Gaussian tail estimate gives
\begin{align}
        1-s(0)\mathcal L_{\rm BA}^{-1}f_{r,\eta}(0)
        &\leq
        C\int_0^\infty
        \exp\left(-t-\frac{cr^2}{t}\right)\,\rmd t
        \leq Ce^{-c'r}.
\label{eq:BA-massive-large-ball}
\end{align}
In the last step we used \(t+cr^2/t\geq2\sqrt c\,r\), and hence
\[
        e^{-t-cr^2/t}
        \leq e^{-\sqrt c\,r}e^{-t/2},
\]
after discarding the nonpositive term \(-cr^2/(2t)\).

Lemmas~\ref{lemma:BA-fixed-row-moments} and \ref{lemma:BA-Ward-mass} imply
\[
        c\leq j_u(0)\leq C,
        \qquad
        \frac{\lambda\eta}{s(0)}\asymp1.
\]
Combining these bounds with \cref{eq:BA-radial-elliptic-error,eq:BA-massive-small-ball} proves the deterministic part of \cref{eq:BA-small-ball-resolvent}.  For the complement we use the exact linearity identity
\begin{equation}
\label{eq:BA-complement-response}
        p_u[1-f_{r,\eta}](0)
        =
        \frac{j_u(0)}{s(0)}-p_u[f_{r,\eta}](0).
\end{equation}
Together with \cref{eq:BA-radial-elliptic-error,eq:BA-massive-large-ball}, this proves the deterministic part of \cref{eq:BA-large-ball-resolvent}.

Finally, \cref{eq:BA-radial-support-condition} allows us to apply Corollary~\ref{coro:BA-fixed-row-random-profile} to either radial source. For $\mathcal A\in\{\mathcal A_{r,\eta},\I-\mathcal A_{r,\eta}\}$, let $f$ be its scalar block profile and write $p_u[f]$ for the corresponding infinite-volume response.  Then
\[
\begin{aligned}
        \lambda\eta
        \left|
        (R\mathcal A R^\star)_{xx}-p_u[f](0)
        \right|\prec
        \lambda\eta
        \left(1+(\lambda\eta)^{-1}\right)
        \Phi_d^{\rm BA}(E+\mathrm i\eta)+\lambda^{-D}
        \leq
        2\Phi_d^{\rm BA}(E+\mathrm i\eta)+\lambda^{-D},
\end{aligned}
\]
where we used \(\lambda\eta\leq1\).  After enlarging the deterministic constant in \cref{eq:BA-small-ball-resolvent,eq:BA-large-ball-resolvent}, take \(\mathcal E_{r,\eta}^{\rm in}\) and \(\mathcal E_{r,\eta}^{\rm out}\) to be the positive parts of the two remaining random discrepancies.  The preceding display and the energy net following this proof give \cref{eq:BA-radial-random-error}.  This completes the proof.
\end{proof}

Uniformity in the preceding lemma follows from the following standard net argument, included because the energy integration requires a single event. For \(0\preceq\mathcal A\preceq\I\), set
\[
        F_{\mathcal A}(E)
        =
        (R(E+\mathrm i\eta)\mathcal A
        R(E+\mathrm i\eta)^\star)_{xx}.
\]
Since \(\partial_E R=R^2\),
\begin{equation}
\label{eq:BA-energy-Lipschitz}
\begin{aligned}
        |\partial_EF_{\mathcal A}(E)|
        &\leq
        \|(R^\star)^2e_x\|_2\|R^\star e_x\|_2
        +\|R^\star e_x\|_2\|(R^\star)^2e_x\|_2
        \leq2\eta^{-3}.
\end{aligned}
\end{equation}
Throughout \cref{eq:BA-dynamical-time-volume-range}, \(\eta^{-1}\) is bounded by a fixed power of \(\lambda\).  A grid of polynomial size and mesh \(\lambda^{-C}\) changes \(\lambda\eta F_{\mathcal A}\) by at most \(2\lambda\eta^{-2}\lambda^{-C}\).  Given any fixed \(D>0\), choose \(C\) so large that this is at most \(\lambda^{-D}\).  The grid still has polynomial size.  The union bound in the definition of stochastic domination therefore upgrades the pointwise estimates to \cref{eq:BA-radial-random-error} on one event uniform over every compact regular energy set.

The profile has only been proved away from the critical energies.  The following lemma removes this restriction without introducing a microscopic test function.  For every sufficiently small \(\varepsilon>0\), let \(\mathcal G_\varepsilon\) be a finite union of open regular subintervals of \((0,4d)\setminus\Sigma_d^\Psi\).  We choose these intervals so that their complement in \((0,4d)\) is contained in shrinking neighborhoods of the critical level sets.  Choose $h_\varepsilon\in C_c^\infty(\mathcal G_\varepsilon)$, with $0\leq h_\varepsilon\leq1$, such that $h_\varepsilon(e)=1$ whenever $e\in(0,4d)$ and $\operatorname{dist}(e,\mathbb R\setminus\mathcal G_\varepsilon) \geq2\varepsilon$, and such that
\[
        \operatorname{dist}(\operatorname{supp}h_\varepsilon,
        \mathbb R\setminus\mathcal G_\varepsilon)\geq\varepsilon.
\]
In the growing block size regime these intervals are obtained by deleting neighborhoods of the finite set \(\Sigma_d^\Psi\).  In the fixed-block regime the same conclusion follows from Proposition~\ref{prop:BA-regular-energies-folded-bulk}.  In either case the intervals and cutoff may be chosen so that
\begin{equation}
\label{eq:BA-free-exceptional-weight}
\begin{aligned}
        \omega(\varepsilon)
        &\coloneqq
        \frac1{(2\pi)^d}
        \int_{\mathbb T^d}
        |1-h_\varepsilon(\Omega(p))|^2
        \,\rmd p
        \longrightarrow0
        \qquad(\varepsilon\downarrow0).
\end{aligned}
\end{equation}
Here \(\Omega\) is the Fourier symbol of \(\Psi\).  The last assertion follows because the critical level sets have Lebesgue measure zero. Moreover, the elementary error estimate for the Riemann sum of the smooth periodic function \(g_\varepsilon(p)=|1-h_\varepsilon(\Omega(p))|^2\) is
\begin{equation}
\label{eq:BA-free-exceptional-Riemann-error}
        \left|
        \frac1{L^d}\sum_{p\in\tilde{\mathbb T}_L^d}g_\varepsilon(p)
        -\frac1{(2\pi)^d}\int_{\mathbb T^d}g_\varepsilon(p)\,\rmd p
        \right|
        \leq\frac{C_\varepsilon}{L}.
\end{equation}
Indeed, partitioning \(\mathbb T^d\) into the \(L^d\) Fourier cells and applying the mean value theorem on each cell bounds the error by \(CL^{-1}\|\nabla g_\varepsilon\|_{L^\infty}\).

\begin{lemma}
\label{lemma:BA-exceptional-energy-removal}
For every deterministic \(0\preceq\mathcal A\preceq\I\),
\begin{equation}
\label{eq:BA-energy-removal}
\begin{aligned}
        \eta\int_{\mathbb R}
        \bigl(R\mathcal A R^\star\bigr)_{xx}\,\rmd E
        \leq
        C\omega(\varepsilon)
        +\frac{C_\varepsilon}{L}
        +C\frac{\eta}{\lambda\varepsilon}
        +C\sup_{E/\lambda\in\mathcal G_\varepsilon}
        \lambda\eta
        \bigl(R\mathcal A R^\star\bigr)_{xx}
        +\mathcal E_\varepsilon
\end{aligned}
\end{equation}
where \(\mathcal E_\varepsilon\geq0\) satisfies
\begin{equation}
\label{eq:BA-exceptional-energy-random-error}
        \mathcal E_\varepsilon
        \prec
        \frac{C_\varepsilon(\log\lambda)^C}{\lambda^2}.
\end{equation}
The estimate is uniform for the two radial observables in \cref{eq:BA-small-ball-resolvent,eq:BA-large-ball-resolvent}.
\end{lemma}

\begin{proof}
We first verify that the random Hamiltonian has asymptotically the same local good-energy weight as the free operator.  Set
\[
        P_{\rm g}=h_\varepsilon(H_L/\lambda),
        \qquad
        P_{\rm b}=\I-P_{\rm g}.
\]
Fourier inversion and Duhamel's formula give
\begin{align}
        h_\varepsilon(H_L/\lambda)-h_\varepsilon(\Psi_L)
        &=
        \frac1{2\pi}\int_{\mathbb R}
        \widehat h_\varepsilon(t)
        \left(e^{\mathrm itH_L/\lambda}
        -e^{\mathrm it\Psi_L}\right)\,\rmd t,\nonumber\\
        e^{\mathrm itH_L/\lambda}-e^{\mathrm it\Psi_L}
        &=
        \frac{\mathrm i}{\lambda}
        \int_0^t
        e^{\mathrm i(t-s)H_L/\lambda}
        \mathcal V e^{\mathrm is\Psi_L}\,\rmd s .
\end{align}
Both propagators are unitary.  Hence
\begin{equation}
\label{eq:BA-functional-calculus-comparison}
\begin{aligned}
        \|h_\varepsilon(H_L/\lambda)
        -h_\varepsilon(\Psi_L)\|_{2\to2}
        &\leq
        \frac{\|\mathcal V\|_{2\to2}}{2\pi\lambda}
        \int_{\mathbb R}|t|
        |\widehat h_\varepsilon(t)|\,\rmd t
        \prec
        \frac{C_\varepsilon(\log\lambda)^C}{\lambda}.
\end{aligned}
\end{equation}
For the last bound we used the GOE spectral-radius tail and a union bound over the polynomially many blocks, which give \(\|\mathcal V\|_{2\to2}\prec(\log\lambda)^C\).

The free Fourier transform gives the exact identity
\begin{equation}
\label{eq:BA-free-local-spectral-weight}
        \|(
        \I-h_\varepsilon(\Psi_L))e_x\|_2^2
        =
        \frac1{L^d}
        \sum_{p\in\tilde{\mathbb T}_L^d}
        |1-h_\varepsilon(\Omega(p))|^2.
\end{equation}
Combining \cref{eq:BA-free-exceptional-weight,eq:BA-free-exceptional-Riemann-error,eq:BA-functional-calculus-comparison,eq:BA-free-local-spectral-weight} with \((a+b)^2\leq2a^2+2b^2\) yields
\begin{equation}
\label{eq:BA-random-exceptional-weight}
        \|P_{\rm b}e_x\|_2^2
        \leq
        2\omega(\varepsilon)
        +\frac{C_\varepsilon}{L}
        +\mathcal R_\varepsilon,
        \qquad
        \mathcal R_\varepsilon
        \prec
        \frac{C_\varepsilon(\log\lambda)^C}{\lambda^2}.
\end{equation}

We now estimate the energy integral.  Since \(P_{\rm g}\) and \(P_{\rm b}\) are functions of \(H_L\), they commute with the resolvent.  For \(E/\lambda\notin\mathcal G_\varepsilon\),
\begin{align}
        (R\mathcal AR^\star)_{xx}
        &=
        \|\mathcal A^{1/2}R^\star e_x\|_2^2
        \nonumber\\
        &\leq
        2\|R^\star P_{\rm g}e_x\|_2^2
        +2\|R^\star P_{\rm b}e_x\|_2^2.
\label{eq:BA-good-bad-resolvent-splitting}
\end{align}
The spectral theorem and \(\int_{\mathbb R}\eta((E-\xi)^2+\eta^2)^{-1}\,\rmd E=\pi\) give the exact normalization
\begin{equation}
\label{eq:BA-bad-energy-spectral-normalization}
        \eta\int_{\mathbb R}
        \|R(E+\mathrm i\eta)^\star
        P_{\rm b}e_x\|_2^2\,\rmd E
        =
        \pi\|P_{\rm b}e_x\|_2^2.
\end{equation}
On the other hand, the spectral support of \(P_{\rm g}\) has distance at least \(\lambda\varepsilon\) from \(\lambda(\mathbb R\setminus\mathcal G_\varepsilon)\).  Therefore
\begin{align}
        &\eta\int_{\lambda(\mathbb R\setminus\mathcal G_\varepsilon)}
        \|R(E+\mathrm i\eta)^\star
        P_{\rm g}e_x\|_2^2\,\rmd E
        \nonumber\\
        &\quad\leq
        \sup_{\xi\in\lambda\operatorname{supp}h_\varepsilon}
        \eta\int_{|E-\xi|\geq\lambda\varepsilon}
        \frac{\rmd E}{(E-\xi)^2+\eta^2}
        \leq C\frac{\eta}{\lambda\varepsilon}.
\label{eq:BA-good-vector-outside-good-energy}
\end{align}
Finally, \(\lambda\mathcal G_\varepsilon\) has length \(O(\lambda)\), and hence
\begin{equation}
\label{eq:BA-integral-over-good-energy}
\begin{aligned}
        \eta\int_{\lambda\mathcal G_\varepsilon}
        (R\mathcal AR^\star)_{xx}\,\rmd E
        \leq
        C\sup_{E/\lambda\in\mathcal G_\varepsilon}
        \lambda\eta(R\mathcal AR^\star)_{xx}.
\end{aligned}
\end{equation}
Substituting \cref{eq:BA-random-exceptional-weight}, \cref{eq:BA-bad-energy-spectral-normalization}, \cref{eq:BA-good-vector-outside-good-energy} and \cref{eq:BA-integral-over-good-energy} into \cref{eq:BA-good-bad-resolvent-splitting} proves \cref{eq:BA-energy-removal}, with \(\mathcal E_\varepsilon=2\pi\mathcal R_\varepsilon\).
\end{proof}

\begin{proof}[Proof of Theorem~\ref{thm:BA-dynamical-delocalization}]
We first consider fixed \(W\) and growing blocks satisfying \(W\leq\lambda^{2-\gamma_0}\), with fixed \(\gamma_0>0\). The remaining growing-block range follows from the uniform argument in the proof of Theorem~\ref{thm:BA-genuine-large-block-dynamics} below.  The time-conversion identities proved here are used in both arguments. Set \(\eta=T^{-1}\).  We first derive the resolvent representation of the time average.  For \(t\geq0\), the Laplace transform gives
\begin{equation}
\label{eq:BA-resolvent-Laplace-transform}
        R(E+\mathrm i\eta)
        =
        \mathrm i\int_0^\infty
        e^{-\eta t}e^{\mathrm iEt}e^{-\mathrm itH_L}\,\rmd t.
\end{equation}
Indeed, on a spectral value \(\xi\) the right-hand side is \(\mathrm i(\eta+\mathrm i(\xi-E))^{-1} =(\xi-E-\mathrm i\eta)^{-1}\).  Applying Fourier Plancherel in \(E\), first to each matrix entry and then summing against a positive block-diagonal \(\mathcal A\), gives
\begin{equation}
\label{eq:BA-resolvent-propagator-Plancherel}
\begin{aligned}
        &\int_0^\infty e^{-2\eta t}
        (e^{-\mathrm itH_L}\mathcal A
        e^{\mathrm itH_L})_{xx}\,\rmd t =
        \frac1{2\pi}\int_{\mathbb R}
        (R(E+\mathrm i\eta)\mathcal A
        R(E+\mathrm i\eta)^\star)_{xx}\,\rmd E.
\end{aligned}
\end{equation}
For \(\mathcal A=\operatorname{diag}(a_y)_y\), the integrand on the left is exactly
\[
        (e^{-\mathrm itH_L}\mathcal A
        e^{\mathrm itH_L})_{xx}
        =
        \sum_ya_y|(e^{-\mathrm itH_L})_{xy}|^2.
\]
Since \(H_L^\top=H_L\), its propagator is symmetric, so the same row is the transition profile of the initially localized vector \(e_x\). Since \(e^{-2t/T}\geq e^{-2}\) on \([0,T]\), positivity implies the one-sided estimate
\begin{equation}
\label{eq:BA-finite-time-from-resolvent}
\begin{aligned}
        &\frac1T\int_0^T
        (e^{-\mathrm itH_L}\mathcal A
        e^{\mathrm itH_L})_{xx}\,\rmd t \leq
        \frac{e^2\eta}{2\pi}
        \int_{\mathbb R}
        (R(E+\mathrm i\eta)\mathcal A
        R(E+\mathrm i\eta)^\star)_{xx}\,\rmd E.
\end{aligned}
\end{equation}

We next check all error terms in the claimed time range.  From \cref{eq:BA-dynamical-time-volume-range},
\begin{equation}
\label{eq:BA-dynamical-parameter-conversion}
        \lambda\eta=\frac\lambda T\leq C_\delta^{-1},
        \qquad
        1+(\lambda\eta)^{-1}
        =1+\frac T\lambda
        \leq C\lambda^{\kappa_d^{\rm BA}-\delta}.
\end{equation}
Consequently,
\begin{equation}
\label{eq:BA-dynamical-Phi-smallness}
        \sup_{E/\lambda\in\mathcal G_\varepsilon}
        \Phi_d^{\rm BA}(E+\mathrm i\eta)
        \leq
        \begin{cases}
        C\lambda^{-13\delta/4},
        &d\in\{2,4,5,6\},\\[1mm]
        C\lambda^{-27\delta/8},
        &d=3,\\[1mm]
        C\lambda^{-1+2\kappa_d^{\rm BA}-2\delta}
        \leq C\lambda^{-2\delta},
        &d\geq7.
        \end{cases}
\end{equation}
For example, the first exponent is \(-1/2+(13/4)(2/13-\delta)=-13\delta/4\); the other two follow in the same way.  The definition of \(\kappa_d^{\rm BA}\) also places \(\eta\) inside the range of the local law.  Moreover, the arbitrarily small exponent allowed in stochastic domination is chosen smaller than the positive powers in \cref{eq:BA-dynamical-Phi-smallness}.  Thus the smallness hypothesis of Lemma~\ref{lemma:BA-T-error-and-kernel} also holds.  Finally, \cref{eq:BA-dynamical-time-volume-range} is precisely \cref{eq:BA-profile-finite-volume-condition}, while
\[
        \frac{\eta}{\lambda\varepsilon}
        =\frac1{\lambda T\varepsilon}=o(1)
\]
for fixed \(\varepsilon\). For every fixed \(r>0\), the volume condition also gives
\begin{equation}
\label{eq:BA-dynamical-radial-support}
        \frac{2r\lambda^{3/2}\eta^{-1/2}}{Wn}
        \leq
        \frac{2r\sqrt{\lambda\eta}}
        {W(\log\lambda)^2}
        =o(1).
\end{equation}
Thus \cref{eq:BA-radial-support-condition} holds for each of the two fixed radii chosen below.

We first prove the lower bound on the second moment.  Apply Lemma~\ref{lemma:BA-exceptional-energy-removal} with \(\mathcal A=\mathcal A_{r_-,\eta}\), followed by \cref{eq:BA-small-ball-resolvent}.  This gives
\begin{align}
        &\eta\int_{\mathbb R}
        (R\mathcal A_{r_-,\eta}R^\star)_{xx}\,\rmd E\leq C\bigg[
        \omega(\varepsilon)
        +r_-^2\left(1+\bm 1_{\{d=2\}}|\log r_-|\right)
        +\frac{\sqrt{\lambda\eta}}{r_-}
        +\lambda\eta
        +\frac{C_\varepsilon}{L}
        +\frac{\eta}{\lambda\varepsilon}
        \bigg]
        +\mathcal R_\lambda.
\label{eq:BA-integrated-small-ball}
\end{align}
Here we may take
\[
        \mathcal R_\lambda
        =
        C\mathcal E_\varepsilon
        +C\sup_{E/\lambda\in\mathcal G_\varepsilon}
        \mathcal E_{r_-,\eta}^{\rm in}(E).
\]
Thus \(\mathcal R_\lambda\geq0\) and
\begin{equation}
\label{eq:BA-integrated-small-ball-random-error}
        \mathcal R_\lambda
        \prec
        \sup_{E/\lambda\in\mathcal G_\varepsilon}
        \Phi_d^{\rm BA}(E+\mathrm i\eta)
        +\lambda^{-D}
        +\frac{C_\varepsilon(\log\lambda)^C}{\lambda^2}.
\end{equation}
Choose first \(\varepsilon\) so that \(\omega(\varepsilon)\) is small, then choose \(0<r_-<1\) so that the second term is small.  Next enlarge the constant in \(T\geq C_\delta\lambda\).  By \cref{eq:BA-dynamical-parameter-conversion}, this makes both \(r_-^{-1}\sqrt{\lambda\eta}\) and \(\lambda\eta\) small.  Finally take \(\lambda\) large.  Equation \cref{eq:BA-dynamical-Phi-smallness} shows that the right-hand side of \cref{eq:BA-integrated-small-ball-random-error} tends to zero even after the arbitrarily small power allowed by stochastic domination. The last two deterministic terms in \cref{eq:BA-integrated-small-ball} also tend to zero.  Hence
\begin{equation}
\label{eq:BA-integrated-small-ball-small}
        \eta\int_{\mathbb R}
        (R\mathcal A_{r_-,\eta}R^\star)_{xx}\,\rmd E
        \leq\frac{\pi}{2e^2}
\end{equation}
with probability at least \(1-C\lambda^{-1000}\).  By \cref{eq:BA-finite-time-from-resolvent},
\begin{equation}
\label{eq:BA-time-averaged-block-small-ball}
        \frac1T\int_0^T
        (e^{-\mathrm itH_L}\mathcal A_{r_-,\eta}
        e^{\mathrm itH_L})_{xx}\,\rmd t
        \leq\frac14.
\end{equation}

We now compare the block cutoff with a microscopic ball.  Write \(y=(I,v)\) and identify the internal labels \(u,v\) with their natural coordinates in \(\{0,\ldots,W-1\}^d\).  In the next display both distances are periodic Euclidean distances, represented in the centered microscopic and block fundamental domains.  The triangle inequality gives
\begin{equation}
\label{eq:BA-cell-geometry}
        \left|
        |y-x|-W|I|
        \right|
        \leq|v-u|\leq\sqrt d\,W.
\end{equation}
At the diffusive length in the theorem,
\begin{equation}
\label{eq:BA-block-boundary-negligible}
        \frac{W}{\lambda^{3/2}\sqrt T}
        \leq
        \begin{cases}
        C_W\lambda^{-2},&W\text{ is fixed},\\[1mm]
        C\lambda^{-\gamma_0},&W\leq\lambda^{2-\gamma_0},
        \end{cases}
\end{equation}
where we used \(T\geq C_\delta\lambda\).  Thus, for large \(\lambda\), the indicator of \(\{|y-x|\leq r_-\lambda^{3/2}\sqrt T/2\}\) is bounded by \(\mathcal A_{r_-,\eta}\).  By \cref{eq:BA-dynamical-radial-support}, the cutoff is supported inside \(\{|I|\leq n/10\}\), and the ball does not wrap around the torus.

It follows from \cref{eq:BA-time-averaged-block-small-ball} that at least three quarters of the time-averaged mass lies outside the microscopic ball of radius \(r_-\lambda^{3/2}\sqrt T/2\).  Since unitarity gives \(\sum_y|(e^{-\mathrm itH_L})_{xy}|^2=1\),
\begin{align}
        &\frac1T\int_0^T
        \sum_y|y-x|^2
        |(e^{-\mathrm itH_L})_{xy}|^2\,\rmd t
        \nonumber\\
        &\quad\geq
        \frac{r_-^2\lambda^3T}{4}
        \frac1T\int_0^T
        \sum_{|y-x|\geq r_-\lambda^{3/2}\sqrt T/2}
        |(e^{-\mathrm itH_L})_{xy}|^2\,\rmd t
        \geq\frac{3r_-^2}{16}\lambda^3T.
\label{eq:BA-second-moment-final-derivation}
\end{align}
This proves \cref{eq:BA-dynamical-second-moment}.

For the upper tail, apply Lemma~\ref{lemma:BA-exceptional-energy-removal} to \(\I-\mathcal A_{r_+,\eta}\) and use \cref{eq:BA-large-ball-resolvent}.  Exactly as above,
\begin{align}
        &\eta\int_{\mathbb R}
        \left(R(\I-\mathcal A_{r_+,\eta})R^\star\right)_{xx}\,\rmd E
        \nonumber\\
        &\quad\leq
        C\left[
        \omega(\varepsilon)+e^{-cr_+}
        +\frac{\sqrt{\lambda\eta}}{r_+}+\lambda\eta
        +\frac{C_\varepsilon}{L}
        +\frac{\eta}{\lambda\varepsilon}
        \right]
        +\mathcal R_\lambda^{\rm out},
\label{eq:BA-integrated-outer-tail}
\end{align}
where
\[
        \mathcal R_\lambda^{\rm out}
        =
        C\mathcal E_\varepsilon
        +C\sup_{E/\lambda\in\mathcal G_\varepsilon}
        \mathcal E_{r_+,\eta}^{\rm out}(E)
        \prec
        \sup_{E/\lambda\in\mathcal G_\varepsilon}
        \Phi_d^{\rm BA}(E+\mathrm i\eta)
        +\lambda^{-D}
        +\frac{C_\varepsilon(\log\lambda)^C}{\lambda^2}.
\]
Choose first \(\varepsilon\) and then \(r_+>1\) so that the first two terms, including the constants in \cref{eq:BA-energy-removal,eq:BA-finite-time-from-resolvent}, contribute at most \(\delta/2\) to the time average.  Enlarging \(C_\delta\) in the lower time bound and then taking \(\lambda\) large makes the sum of all remaining terms at most \(\delta/2\).  The same two equations then give
\begin{equation}
\label{eq:BA-time-averaged-block-outer-tail}
        \frac1T\int_0^T
        (e^{-\mathrm itH_L}(\I-\mathcal A_{r_+,\eta})
        e^{\mathrm itH_L})_{xx}\,\rmd t
        \leq\delta.
\end{equation}
By \cref{eq:BA-cell-geometry,eq:BA-block-boundary-negligible}, if \(|y-x|\geq3r_+\lambda^{3/2}\sqrt T\), then \(W|I|\geq2r_+\lambda^{3/2}\sqrt T\), and hence \(\chi(WI/(r_+\lambda^{3/2}\sqrt T))=0\).  Thus the microscopic outer indicator is bounded by \(\I-\mathcal A_{r_+,\eta}\), and \cref{eq:BA-time-averaged-block-outer-tail} proves \cref{eq:BA-dynamical-outer-tail} after enlarging \(C_\delta\) so that \(C_\delta\geq3r_+\).

Only two radial sources and polynomially many energy-net points have been used.  Taking the exponent in stochastic domination larger than \(1000\) and applying the union bound proves the stated probability. This completes the proof.
\end{proof}

\subsection{Growing blocks across the transition scale}
\label{subsec:BA-genuine-large-blocks}

The growing-block local law, the \(T\)-equation error bound, and the Ward mass do not require \(W\ll\lambda^2\).  We now keep both terms in \cref{eq:BA-cell-coercivity-crossover}.  In this subsection we will shorthand \(h=1+\lambda^2/W\); thus \(d_{\rm BA}\asymp(\lambda^2/W)h\).  The same estimates will apply on both sides of \(W\asymp\lambda^2\).

\begin{lemma}
\label{lemma:BA-genuine-large-block-fourth-moment}
Under the growing-block hypotheses of Proposition~\ref{prop:BA-cell-coercivity},
\begin{equation}
\label{eq:BA-genuine-large-block-fourth-moment}
        \frac1{W^d}\sum_{X\in\mathbb Z^d}
        |X|^4\|\tilde M_\infty(X)\|_{\rm HS}^2
        \leq C d_{\rm BA}h^2.
\end{equation}
For every multi-index \(\gamma\), with \(1\leq|\gamma|\leq4\),
\begin{equation}
\label{eq:BA-genuine-large-block-symbol-derivatives}
        \sup_{k\in\mathbb T^d}
        \|\partial^\gamma\widehat{\tilde K}_\infty(k)\|
        \leq C_\gamma d_{\rm BA}h^{|\gamma|-2}.
\end{equation}
The constants are independent of \(W\).
\end{lemma}

\begin{proof}
Use the notation in the growing-block proof of Proposition~\ref{prop:BA-cell-coercivity}. For \(|t|=qW+s\), with \(q\in\mathbb N_0\) and \(0\leq s<W\), the endpoint counting used in \cref{eq:BA-exact-cell-crossing-second-moment} gives
\[
        Q_W^{(4)}(t)
        =\left(1-\frac{s}{W}\right)q^4
          +\frac{s}{W}(q+1)^4
        \leq C\left(\frac{|t|}{W}+\frac{|t|^4}{W^4}\right).
\]
For \(q=0\) the first term on the right is exact, and for \(q\geq1\) we use \((q+1)^4\leq16q^4\).  Therefore \cref{eq:BA-scalar-microscopic-flight-moments} implies
\begin{equation}
\label{eq:BA-scalar-cell-fourth-moment}
        \frac1{W^d}\sum_X X_i^4\|M_g(X)\|_{\rm HS}^2
        \leq C\left(\frac{\lambda^2}{W}
                    +\frac{\lambda^8}{W^4}\right)
        \leq C d_{\rm BA}h^2.
\end{equation}

To transfer this bound to the matrix self-energy, diagonalize the scalar fiber using
\[
        U_{um}(\xi)=W^{-d/2}
        \exp\!\left(\mathrm i\frac{\xi+2\pi m}{W}\cdot u\right),
        \qquad
        G_g(\xi)=U(\xi)Q_g(\xi)U(\xi)^\star,
\]
where \(Q_g(\xi)=\operatorname{diag}_m (\Omega((\xi+2\pi m)/W)-\zeta)^{-1}\). Since \(\partial_iU=\mathrm iD_iU\), with \(D_i=\operatorname{diag}(u_i/W)\) of norm at most one,
\[
        \partial_iG_g=\mathrm i[D_i,G_g]+U(\partial_iQ_g)U^\star.
\]
The two terms have operator norms at most \(C\lambda^2\) and \(C\lambda^4/W\), respectively.  Hence
\begin{equation}
\label{eq:BA-scalar-cell-first-derivative-operator-bound}
        \sup_\xi\|\partial_iG_g(\xi)\|
        \leq C\lambda^2h.
\end{equation}
Differentiating \cref{eq:BA-cell-derivative-relative-factorization} once more gives
\[
        \partial_i^2G_{\tilde\theta}
        =A(\partial_i^2G_g)B
          +2A(\partial_iG_g)VA(\partial_iG_g)B.
\]
The uniform bounds for \(A,B\), the estimate \(\|V\|\leq C\lambda^{-2}\), and \cref{eq:BA-scalar-cell-first-derivative-operator-bound} imply
\[
        \|\partial_i^2G_{\tilde\theta}\|_{\rm HS}
        \leq C\left(\|\partial_i^2G_g\|_{\rm HS}
                      +h\|\partial_iG_g\|_{\rm HS}\right).
\]
Square, integrate, and divide by \(W^d\lambda^2\). With normalized Haar measure, Parseval gives exactly
\[
        \frac1{W^d}\sum_X X_i^4
        \|\tilde M_\infty(X)\|_{\rm HS}^2
        =\frac1{\lambda^2W^d}\int_{\mathbb T^d}
        \|\partial_i^2G_{\tilde\theta}(\xi)\|_{\rm HS}^2
        \,\rmd\mu(\xi).
\]
The scalar second moment is \(O(d_{\rm BA})\), and its fourth moment satisfies \cref{eq:BA-scalar-cell-fourth-moment}. This proves the desired bound for \(X_i^4\). Now use \(|X|^4\leq d\sum_iX_i^4\).

For the symbol derivatives, use
\[
        \|\tilde K_\infty(X)\|
        \leq\frac2{W^d}\|\tilde M_\infty(X)\|_{\rm HS}^2.
\]
Its zeroth moment is bounded, its second moment is \(O(d_{\rm BA})\), and its fourth moment is \(O(d_{\rm BA}h^2)\).  Cauchy--Schwarz bounds the third moment by \(Cd_{\rm BA}h\).  For the first moment, Cauchy--Schwarz and \(|X|\leq|X|^2\) for nonzero \(X\in\mathbb Z^d\) give
\[
        \sum_X|X|\|\tilde K_\infty(X)\|
        \leq C\min\{\sqrt{d_{\rm BA}},d_{\rm BA}\}
        \leq C d_{\rm BA}/h.
\]
The last inequality follows separately for \(\lambda^2/W\leq1\) and \(\lambda^2/W\geq1\). Termwise differentiation proves \cref{eq:BA-genuine-large-block-symbol-derivatives}.
\end{proof}

\begin{coro}
\label{coro:BA-genuine-large-block-multiplier}
Under the hypotheses of Lemma~\ref{lemma:BA-genuine-large-block-fourth-moment}, there is a constant \(c>0\), independent of \(W\), such that for \(|k|\leq c/h\),
\begin{equation}
\label{eq:BA-genuine-large-block-scalar-expansion}
        s(k)=s(0)+d_{\rm BA}|k|^2+O(d_{\rm BA}h^2|k|^4),
        \qquad
        |s(k)|\geq c\bigl(s(0)+d_{\rm BA}|k|^2\bigr).
\end{equation}
Moreover, uniformly in \(u\),
\begin{equation}
\label{eq:BA-genuine-large-block-row-multiplier}
        \left|\frac{j_u^{\rm ev}(k)}{s(k)}
        -\frac{j_u(0)}{s(0)+d_{\rm BA}|k|^2}\right|
        \leq C\frac{h^2}{d_{\rm BA}}.
\end{equation}
For real even \(f\in C_c^\infty(\mathbb R^d)\) and \(a_I=f(I)\),
\begin{equation}
\label{eq:BA-genuine-large-block-elliptic-approximation}
        \left|p_u[a](0)-j_u(0)\mathcal L_{\rm BA}^{-1}f(0)\right|
        \leq C\left(
        \frac{h\|\nabla f\|_\infty}{s(0)}
        +\frac{h^2\|\widehat f\|_1}{d_{\rm BA}}\right).
\end{equation}
\end{coro}

\begin{proof}
The hard inverse at zero is uniformly bounded.  By \cref{eq:BA-genuine-large-block-symbol-derivatives}, its perturbation on \(|k|\leq c/h\) is \(O(cd_{\rm BA}/h^2)=O(c)\). Thus the inverse remains bounded after decreasing \(c\). Differentiating the inverse and the normalized cell equation gives
\[
        \|\partial^\gamma H(k)\|
        +\|\partial^\gamma C(k)\|
        +|\partial^\gamma s(k)|
        \leq C_\gamma d_{\rm BA}h^{|\gamma|-2},
        \qquad 1\leq|\gamma|\leq4.
\]
Indeed, a term containing \(m\) differentiated kernel factors is bounded by \(Cd_{\rm BA}^m h^{|\gamma|-2m}\), which is at most \(Cd_{\rm BA}h^{|\gamma|-2}\), since \(d_{\rm BA}/h^2\leq C\).  Cubic symmetry and Taylor's formula give \cref{eq:BA-genuine-large-block-scalar-expansion}. The row moments in \cref{eq:BA-fixed-row-moments} similarly give \(j_u^{\rm ev}(k)=j_u(0)+O(h^2|k|^2)\).  Subtracting the two fractions bounds their difference by
\[
        C\frac{h^2|k|^2}{s(0)+d_{\rm BA}|k|^2}
        +C\frac{d_{\rm BA}h^2|k|^4}
                   {(s(0)+d_{\rm BA}|k|^2)^2}
        \leq C\frac{h^2}{d_{\rm BA}}.
\]
Finally, repeat the smoothing and Fourier argument of Proposition~\ref{prop:BA-diagonal-elliptic-approximation}, using smoothing length \(Ch\) so that the Fourier support lies inside \(|k|\leq c/h\).  Positivity bounds the smoothing error by \(Ch\|\nabla f\|_\infty/s(0)\), and the multiplier estimate gives the remaining term in \cref{eq:BA-genuine-large-block-elliptic-approximation}.
\end{proof}

\begin{proof}[Proof of Theorem~\ref{thm:BA-genuine-large-block-dynamics}]
All estimates below are uniform on a fixed compact regular energy set.  The assumed lower bound for \(W\) permits application of the growing-block local law and the \(T\)-equation error bound. Corollary~\ref{coro:BA-genuine-large-block-multiplier} gives
\begin{equation}
\label{eq:BA-large-block-elliptic-comparison}
        \left|p_u[(f(I))_I](0)
        -j_u(0)\mathcal L_{\rm BA}^{-1}f(0)\right|
        \leq C\left[
        \frac{h\|\nabla f\|_\infty}{s(0)}
        +\frac{h^2\|\widehat f\|_1}{d_{\rm BA}}\right]
\end{equation}
for every real even \(f\in C_c^\infty(\mathbb R^d)\).

The finite-volume comparison in Corollary~\ref{coro:BA-fixed-row-random-profile} also extends to every growing block size.  Its exponential-moment scale is \(\mu_0=cW/(\lambda^2+W)=c/h\).  For \(\mu=c'\eta/\lambda\),
\[
        C\frac\mu{\mu_0}
        \leq C'c'\frac\eta\lambda\left(1+\frac{\lambda^2}{W}\right)
        \leq \frac c2\lambda\eta
\]
after choosing \(c'\) sufficiently small.  Thus the weighted-kernel perturbation is smaller than half its Ward contraction gap.  The weighted inverse, source-periodization, and cross-image estimates in that proof follow unchanged.  Under \(n\eta/\lambda\geq(\log\lambda)^2\), the normalized random error remains \(\Phi_d^{\rm BA}(E+\mathrm i\eta)+\lambda^{-D}\) in stochastic domination, for every fixed \(D>0\).

Set
\begin{equation}
\label{eq:BA-large-block-radial-observable}
        R_\eta\coloneqq\sqrt{\frac{\lambda^3+\lambda W}{\eta}},
        \qquad
        f_{r,\eta}(y)\coloneqq\chi\left(\frac{Wy}{rR_\eta}\right),
        \qquad
        \mathcal A^{\rm lg}_{r,\eta}
        \coloneqq\bigoplus_I f_{r,\eta}(I)\I_{W^d}.
\end{equation}
Then \(rR_\eta/W\asymp r\sqrt{d_{\rm BA}/s(0)}\), and multiplication of \cref{eq:BA-large-block-elliptic-comparison} by \(\lambda\eta\) gives the normalized error
\begin{equation}
\label{eq:BA-large-block-radial-error}
        C\left[\frac{hW}{rR_\eta}
                    +\frac{\lambda\eta h^2}{d_{\rm BA}}\right]
        \leq C\left[
        \frac1r\sqrt{\left(\lambda+\frac W\lambda\right)\eta}
        +\left(\lambda+\frac W\lambda\right)\eta\right].
\end{equation}
The massive heat-kernel formula \cref{eq:BA-massive-heat-kernel} therefore gives
\begin{equation}
\label{eq:BA-large-block-small-ball}
\begin{aligned}
 &\lambda\eta(R\mathcal A^{\rm lg}_{r,\eta}R^\star)_{xx}\\
 &\quad\leq C\left[
 r^2\bigl(1+\bm1_{\{d=2\}}|\log r|\bigr)
 +\frac1r\sqrt{\left(\lambda+\frac W\lambda\right)\eta}
 +\left(\lambda+\frac W\lambda\right)\eta\right]
 +\mathcal E^{\rm lg,in}_{r,\eta}(E),
 \qquad 0<r\leq1,
\end{aligned}
\end{equation}
and
\begin{equation}
\label{eq:BA-large-block-large-ball}
\begin{aligned}
 &\lambda\eta(R(\I-\mathcal A^{\rm lg}_{r,\eta})R^\star)_{xx}\\
 &\quad\leq C\left[
 e^{-cr}+\frac1r\sqrt{\left(\lambda+\frac W\lambda\right)\eta}
 +\left(\lambda+\frac W\lambda\right)\eta\right]
 +\mathcal E^{\rm lg,out}_{r,\eta}(E),
 \qquad r\geq1.
\end{aligned}
\end{equation}
These bounds apply provided \(2rR_\eta/W\leq n/10\). For the complement use \cref{eq:BA-complement-response}. The same polynomial energy net as before gives
\[
        \sup_{E/\lambda\in\mathcal G}
        \bigl(|\mathcal E^{\rm lg,in}_{r,\eta}(E)|
              +|\mathcal E^{\rm lg,out}_{r,\eta}(E)|\bigr)
        \prec\sup_{E/\lambda\in\mathcal G}
        \Phi_d^{\rm BA}(E+\mathrm i\eta)+\lambda^{-D}.
\]

Take \(\eta=T^{-1}\).  The upper time bound is the one used in \cref{eq:BA-dynamical-Phi-smallness}, so the random errors tend to zero.  The lower bound in \cref{eq:BA-large-block-time-range} gives
\[
        \left(\lambda+\frac W\lambda\right)\eta
        \leq C_\delta^{-1},
        \qquad \lambda\eta\leq C_\delta^{-1}.
\]
For fixed \(r>0\), the volume condition and \(T\geq C_\delta\lambda\) imply
\[
        \frac{R_\eta}{Wn}
        \leq\frac1{(\log\lambda)^2}
        \sqrt{\frac{\lambda}{W^2T}+\frac1{\lambda WT}}
        =o(1).
\]
Thus the cutoffs have the required support. The exceptional-energy removal lemma and \cref{eq:BA-finite-time-from-resolvent} apply without change. Choose the exceptional-energy cutoff first, then choose \(0<r_-<1<r_+\) so that the heat-kernel small-ball and large-ball terms are sufficiently small.  Enlarge \(C_\delta\) to control \cref{eq:BA-large-block-radial-error} for these fixed radii. Taking \(\lambda\) large gives
\[
\begin{aligned}
        \frac1T\int_0^T
        (e^{-\mathrm itH_L}\mathcal A^{\rm lg}_{r_-,T^{-1}}
        e^{\mathrm itH_L})_{xx}\,\rmd t&\leq\frac14,\\
        \frac1T\int_0^T
        (e^{-\mathrm itH_L}(\I-\mathcal A^{\rm lg}_{r_+,T^{-1}})
        e^{\mathrm itH_L})_{xx}\,\rmd t&\leq\delta.
\end{aligned}
\]
Finally,
\[
        \frac W{\sqrt{(\lambda^3+\lambda W)T}}
        \leq\sqrt{\frac W{\lambda T}}
        \leq C_\delta^{-1/2}.
\]
Enlarging the same constant makes the block-coordinate error \(\sqrt d\,W\) in \cref{eq:BA-cell-geometry} smaller than \(r_-\sqrt{(\lambda^3+\lambda W)T}/2\). At least three quarters of the averaged mass then lies outside the corresponding microscopic ball, proving \cref{eq:BA-large-block-second-moment}.  The outer indicator at radius \(3r_+\sqrt{(\lambda^3+\lambda W)T}\) is bounded by \(\I-\mathcal A^{\rm lg}_{r_+,T^{-1}}\), proving \cref{eq:BA-large-block-outer-tail}. The union bound over the two radial sources and the polynomial energy net gives the stated probability. For \(W\leq\lambda^2\), we have \(\lambda^3+\lambda W\asymp\lambda^3\) and \(\lambda+W/\lambda\asymp\lambda\); hence this also proves part \textnormal{(b)} of Theorem~\ref{thm:BA-dynamical-delocalization} throughout its stated range, including \(W\asymp\lambda^2\).
\end{proof}

\subsection{Scaling limit of the resolvent}
\label{subsec:BA-resolvent-scaling-proof}

\begin{proof}[Proof of Theorem~\ref{thm:BA-resolvent-scaling}]
We use the block scale \(\rho=r_\lambda/W\).  All deterministic coefficients below are evaluated at \(z=\lambda e+\mathrm i\lambda^{-1-\kappa}\). The Ward mass, the row normalization, and the crossover coercivity give
\begin{equation}
\label{eq:BA-scaling-coefficients}
        s(0)\asymp\lambda\eta,
        \qquad j_u(0)\asymp1,
        \qquad
        d_{\rm BA}\asymp
        \frac{\lambda^4}{W^2}+\frac{\lambda^2}{W}.
\end{equation}
Consequently
\begin{equation}
\label{eq:BA-scaling-mass}
        \frac{s(0)\rho^2}{d_{\rm BA}}
        \asymp\lambda^{-2a}.
\end{equation}
Moreover \(\rho\to\infty\), and the volume assumption implies \(\rho\ll n\).  Thus the block source \(a_I=f(I/\rho)\) is supported in \(|I|\leq n/10\) for all sufficiently large \(\lambda\).

We first extend the elliptic comparison to functions that need not be even.  In this proof only, denote
\[
        h=
        \begin{cases}
        \lambda^2,&W\text{ fixed},\\
        1+\lambda^2/W,&W\text{ growing}.
        \end{cases}
\]
Lemma~\ref{lemma:BA-fixed-row-moments}, differentiation of the normalized cell equation, and the bounded hard inverse imply
\[
        |j_u(k)-j_u^{\rm ev}(k)|\leq Ch|k|
\]
on the corresponding low-frequency window.  For large blocks this uses the uniformly bounded first row moment and the derivative estimate in \cref{eq:BA-genuine-large-block-symbol-derivatives}. Choose the even smoothing kernel used in Proposition~\ref{prop:BA-diagonal-elliptic-approximation}, with smoothing length \(Ch\).  The odd part of the numerator now contributes at most
\[
        \frac{Ch}{s(0)}
        \int |k|\rho^d|\widehat f(\rho k)|\,\rmd k
        \leq\frac{C_fh}{s(0)\rho},
\]
since \(|s(k)|\geq cs(0)\).  The smoothing error has the same bound.  The even part is controlled by Lemma~\ref{lemma:BA-row-elliptic-multiplier} and Corollary~\ref{coro:BA-genuine-large-block-multiplier}. We therefore obtain, for arbitrary real \(f\in C_c^\infty\),
\begin{equation}
\label{eq:BA-scaling-general-test}
\begin{aligned}
        \left|p_u[a](0)
        -j_u(0)\mathcal L_{\rm BA}^{-1}
        \bigl[f(\cdot/\rho)\bigr](0)\right|
        \leq C_f\left[
        \frac{h}{s(0)\rho}
        +\frac{h^2}{d_{\rm BA}}\right].
\end{aligned}
\end{equation}
Complex-valued functions follow by linearity.

Apply Corollary~\ref{coro:BA-fixed-row-random-profile} to the positive and negative parts of the block source, after division by \(\|f\|_\infty\).  Its large-block extension was verified in the proof of Theorem~\ref{thm:BA-genuine-large-block-dynamics}. For \(\mathcal A=\bigoplus_I f(I/\rho)\I_{W^d}\), this gives
\begin{equation}
\label{eq:BA-scaling-random-error}
        \left|(R\mathcal AR^\star)_{xx}-p_u[a](0)\right|
        \prec s(0)^{-1}\Phi_d^{\rm BA}(z)+\lambda^{-D}
\end{equation}
for every fixed \(D>0\).  The same result for the constant source, together with \(p_u[1]=j_u(0)/s(0)\), gives \((RR^\star)_{xx}\prec s(0)^{-1}\).  Since \(|(y-x)-WI|\leq C_dW\) for \(y\) in block \(I\), replacing the block source by the microscopic source changes the response by at most
\[
        \frac{C_f}{\rho}(RR^\star)_{xx}
        \prec\frac{C_f}{s(0)\rho}.
\]
The supports stay away from the boundary of the centered fundamental domain, so this comparison also applies on the torus.

Multiplying \cref{eq:BA-scaling-general-test,eq:BA-scaling-random-error} and the last estimate by \(d_{\rm BA}/(j_u(0)\rho^2)\), we get
\begin{equation}
\label{eq:BA-scaling-normalized-error}
\begin{aligned}
        &\left|
        \frac{W^2d_{\rm BA}}{j_u(0)r_\lambda^2}
        \sum_y f\left(\frac{y-x}{r_\lambda}\right)|R(z)_{xy}|^2
        -\left[\left(-\Delta_{\mathbb R^d}
        +\frac{s(0)\rho^2}{d_{\rm BA}}\right)^{-1}f\right](0)
        \right|\\
        &\qquad\prec
        \lambda^{-a}
        +\lambda^{-6a}
        +\lambda^{2a}\Phi_d^{\rm BA}(z)+\lambda^{-D}.
\end{aligned}
\end{equation}
Indeed, for fixed \(W\) and growing \(W\leq\lambda^2\), the deterministic errors are bounded by \(C_f(\lambda^{3a-\kappa/2}+\lambda^{2a-\kappa})\). In the large-block regime they are bounded by
\[
        C_f\left(
        \lambda^{3a}\sqrt{\frac{W\eta}{\lambda}}
        +\lambda^{2a}\frac{W\eta}{\lambda}\right),
\]
and \(W\eta/\lambda\leq\lambda^{-8a}\). The microscopic sampling error is bounded by the same expression because \(h\geq1\).

By \cref{eq:BA-T-Q-Phi},
\[
        \Phi_d^{\rm BA}(z)\leq C
        \begin{cases}
        \lambda^{-1/2+13\kappa/4},&d\in\{2,4,5,6\},\\
        \lambda^{-3/4+27\kappa/8},&d=3,\\
        \lambda^{-1+2\kappa},&d\geq7.
        \end{cases}
\]
Each exponent is at most \(-2(\kappa_d^{\rm BA}-\kappa)\).  The choice of \(a\) therefore bounds the right-hand side of \cref{eq:BA-scaling-normalized-error} by \(C_f\lambda^{-a}\) in stochastic domination.

It remains to remove the mass from the continuum resolvent.  For \(d\geq3\), the Fourier representation and \cref{eq:BA-scaling-mass} give
\[
\begin{aligned}
        \left[\left(-\Delta+\frac{s(0)\rho^2}{d_{\rm BA}}\right)^{-1}f\right](0)
        &\longrightarrow
        \frac1{(2\pi)^d}\int_{\mathbb R^d}
        \frac{\widehat f(k)}{|k|^2}\,\rmd k\\
        &=\frac1{(d-2)|\mathbb S^{d-1}|}
        \int_{\mathbb R^d}f(y)|y|^{2-d}\,\rmd y.
\end{aligned}
\]
For \(d=2\), the assumption \(\int f=0\) implies \(\widehat f(k)=O_f(|k|)\) near zero.  The same dominated convergence argument applies and the limit is \(-\frac1{2\pi}\int f(y)\log|y|\,\rmd y\). In both cases, splitting the Fourier integral at \(|k|=\sqrt{s(0)\rho^2/d_{\rm BA}}\) bounds this continuum error by \(C_f\lambda^{-a}\).  This proves the stated rate in stochastic domination.  Choosing an arbitrarily large probability exponent and applying Borel--Cantelli along \(\lambda=m\) proves the almost-sure assertion, without any independence assumption between different values of \(m\). The additional limit in the mesoscopic regime follows from Lemma~\ref{lemma:BA-mesoscopic-coefficient-limit} below.
\end{proof}

\begin{lemma}
\label{lemma:BA-mesoscopic-coefficient-limit}
Fix a regular bulk energy \(e\), \(0<\kappa<\kappa_d^{\rm BA}\), and \(\gamma_0>0\).  Suppose
\[
 W^{d/p}\geq C_0(1+\log\lambda)^4,
 \qquad W\leq\lambda^{2-\gamma_0},
 \qquad z=\lambda e+\mathrm i\eta,
 \quad\eta=\lambda^{-1-\kappa},
 \qquad n\eta/\lambda\geq(\log\lambda)^2,
\]
with the remaining growing-block and polynomial-volume assumptions of the local law.  With \(\beta_e\) as in \cref{eq:BA-mesoscopic-beta},
\[
 \sup_{u\in[W]^d}|j_u(0)-1|\longrightarrow0,
 \qquad
 \frac{W^2d_{\rm BA}}{\lambda^4}\longrightarrow\beta_e^{-1}.
\]
\end{lemma}

\begin{proof}
Write \(g=\langle\tilde\theta\rangle_W\) and \(b=\pi\rho(e)>0\).  The growing-block scalarization bounds allow us to apply Proposition~\ref{prop:stability_renorm_resolv} to \(\tilde M\) with a fixed small off-scalar constant. Using \(\tilde\theta=\mathscr S_W[\operatorname{bl}_0\mathbb ER]\) and subtracting normalized traces gives
\[
 \|\tilde\theta-g\I\|
 \leq 2W^{-d}\|\operatorname{bl}_0\tilde M\|
       +4\|\mathbb ER-\tilde M\|_{1\to\infty}.
\]
Consequently \cref{eq:offscalar_small,eq:BA-T-averaged-max-bound} and the bulk bounds imply
\[
 \|\tilde\theta-g\I\|\leq\varepsilon_*/\lambda,
 \qquad |g|\leq C/\lambda,
 \qquad \Im g\geq c/\lambda.
\]
Here we may take
\[
 \varepsilon_*=C(1+\log\lambda)^2W^{-d/p}
       +4\lambda\|\mathbb ER-\tilde M\|_{1\to\infty}=o(1).
\]
All finite-volume errors below are smaller than every inverse power of \(\lambda\), by the assumed volume condition and the same-self-energy comparison.  Let \(M_g\) denote the infinite-volume scalar resolvent with block entries \(m(WX+u-v)\), as in the proof of Proposition~\ref{prop:BA-cell-coercivity}. Its row and column norms are bounded uniformly, by the scalar Ward identity and regularity of the free density of states.  The same bounds for \(\tilde M_\infty\) follow from \cref{eq:BA-row-moment-inputs}.  The resolvent identity therefore improves the diagonal comparison to
\[
 \sup_u\left|
 (\operatorname{bl}_{0,0}\tilde M_\infty)_{uu}
 -(M_g)_{(0,u),(0,u)}\right|
 \leq C\varepsilon_*/\lambda.
\]
Since \((z+g)/\lambda\to e\) from the upper half-plane, \(\lambda\Im(M_g)_{(0,u),(0,u)}\to b\). Moreover, \cref{eq:BA-T-remainder-weighted} and the finite-to-infinite comparison give \(\|\mathcal R_{K,\infty}\|=o(\lambda^{-1})\). Taking normalized traces in \cref{eq:BA-infinite-Ward-remainder}, using \(\langle\mathscr S_W[B]\rangle_W =(1+W^{-d})\langle B\rangle_W\), yields
\[
 \lambda\Im g\longrightarrow b,
 \qquad
 \Gamma=(b+o(1))\lambda^{-1}\I,
 \qquad
 \lambda\Im(\operatorname{bl}_{0,0}\tilde M_\infty)_{uu}
 =b+o(1)
\]
uniformly in \(u\).

We next maintain the vanishing size of the internal part of the renewal.  Estimate~\cref{eq:offscalar_small}, which also holds for \(\tilde\theta\), gives
\[
 W^{-d}\|\operatorname{bl}_{0,0}\tilde M_\infty\|
 \leq C\varepsilon_*/\lambda.
\]
The two generalized Ward identities and \(\Gamma\succeq c\lambda^{-1}\I\) imply
\[
 \|\mathcal E(0)[\I]\|
 +\|\mathcal E(0)^*[\I]\|
 \leq C\varepsilon_*.
\]
Here the adjoint is taken with respect to the ordinary trace pairing.  Positivity and duality show that both the operator-norm and trace-norm operator bounds of \(\mathcal E(0)\) are at most \(t=C\varepsilon_*=o(1)\).  Thus \(U_0=\I+O(t)\) in operator norm.

For completeness, this also sharpens the Hessian comparison. Use \(b_Z,p_Z\) from Lemma~\ref{lemma:BA-positive-one-flight} and shorthand \(\mathcal E_{2,i}=\sum_X X_i^2\mathcal E_X\). Here the ordinary trace, instead of the normalized trace, satisfies
\[
 \varphi_X(B)=\Tr\mathcal E_X[B],
 \qquad \Tr\mathcal E_{2,i}[\I]=\mathcal M_i.
\]
In a contribution to \(p_Z\) with \(m\) internal steps, \((X_{0,i}+\cdots+X_{m,i})^2 \leq(m+1)\sum_{a=0}^mX_{a,i}^2\). For each position of the squared step, the unweighted maps to its right cost \(t\) in operator norm, while those to its left cost \(t\) in trace norm.  There are \(m\) such maps altogether. Consequently the contribution with \(m\geq1\) is bounded by \((m+1)^2t^m\mathcal M_i\).  The contribution with \(m=0\) is exactly \(\mathcal M_i\), and every contribution is nonnegative. The same argument for \(b_Z\), with its normalized trace, gives
\[
 \sum_ZZ_i^2p_Z=(1+O(t))\mathcal M_i,
 \qquad
 \sum_ZZ_i^2b_Z
 \leq W^{-d}\mathcal M_i\sum_{m\geq1}m^2t^{m-1}
 \leq CW^{-d}\mathcal M_i.
\]
Since \(\langle U_0\rangle_W=1+O(t)\) and \(s(0)\asymp\lambda\eta\to0\), \cref{eq:BA-short-positive-Hessian} now gives
\[
 d_{\rm BA}=(1+o(1))\mathcal M_i/2.
\]
Also \(\I+C(0)=U_0/\langle U_0\rangle_W=\I+o(1)\). The first part of the proof and the generalized Ward identity give \(J_u(0)[\I]=1+o(1)\) uniformly in \(u\), and therefore \(j_u(0)=1+o(1)\).

Finally, the factors \(A,B\) in \cref{eq:BA-cell-derivative-relative-factorization} equal \(\I+O(\varepsilon_*)\).  Hence \(\mathcal M_i\) differs relatively by \(o(1)\) from its scalar counterpart. By \cref{eq:BA-exact-cell-crossing-second-moment,eq:BA-exact-cell-crossing-weight}, that counterpart is
\[
 \frac{\mu_{2,i}}{W^2}+O(\mu_{1,i}/W)
 =\frac{\mu_{2,i}}{W^2}(1+o(1)),
\]
because \(\mu_{1,i}=O(\lambda^2)\), \(\mu_{2,i}\asymp\lambda^4\), and \(W/\lambda^2\to0\). With \(\zeta=(z+g)/\lambda\), the first part of the proof gives \(\lambda^2\Im\zeta\to b\) and \(\Re\zeta\to e\). The coarea formula, smoothness at the regular level, and \(\int_{\mathbb R}(s^2+a^2)^{-2}\,\rmd s=\pi/(2a^3)\) give
\[
 \lambda^{-4}\mu_{2,i}
 \longrightarrow
 \frac{\pi}{2b^3}(2\pi)^{-d}
 \int_{\{\Omega=e\}}
 \frac{|\partial_i\Omega(p)|^2}{|\nabla\Omega(p)|}
 \,\rmd\mathcal{H}^{d-1}(p).
\]
Cubic symmetry makes this integral independent of \(i\). Combining the last three displays proves the asserted limit.
\end{proof}

\subsection{Delocalization of most eigenfunctions}
\label{subsec:BA-most-eigenfunctions}

We first record a deterministic consequence of positivity of the imaginary part of the resolvent.

\begin{lemma}
\label{lemma:BA-eigenfunction-compression}
Let \(A\) be a self-adjoint matrix with eigenpairs \((E_j,\psi_j)_j\) and orthonormal eigenvectors.  Let \(P\) be an arbitrary orthogonal projection, and denote \(G=(A-E-\mathrm i\eta)^{-1}\), where \(\eta>0\). Then
\[
        \card\left\{j:|E_j-E|\leq\eta,\quad
                   \|P\psi_j\|_2^2\geq\frac12\right\}
        \leq16\eta^2\|PGP\|_{\mathrm{HS}}^2.
\]
\end{lemma}

\begin{proof}
The spectral decomposition gives
\[
\begin{aligned}
        \Tr\bigl[(P\Im G P)^2\bigr]
        &=\sum_{j,k}
        \frac{\eta^2\,|\langle\psi_j,P\psi_k\rangle|^2}
        {((E_j-E)^2+\eta^2)((E_k-E)^2+\eta^2)}.
\end{aligned}
\]
Every summand is nonnegative. For each index in the set to be counted, the summand with \(k=j\) is at least \(1/(16\eta^2)\). Finally,
\[
        \Tr\bigl[(P\Im G P)^2\bigr]
        =\|\Im(PGP)\|_{\mathrm{HS}}^2
        \leq\|PGP\|_{\mathrm{HS}}^2,
\]
which proves the claim.
\end{proof}

\begin{proof}[Proof of Theorem~\ref{thm:BA-most-eigenfunctions}]
Take
\[
        \eta=\lambda^{-1-\kappa_d^{\rm BA}+\delta}.
\]
We claim that, on an event of probability at least \(1-C_{\varepsilon,\delta}\lambda^{-1000}\),
\begin{equation}
\label{eq:BA-eigenfunction-small-ball}
        \sup_{\substack{x\in\mathbb Z_L^d\\
                E/\lambda\in[e-2\varepsilon,e+2\varepsilon]}}
        \lambda\eta\sum_{y\in B_{4r}(x)}
                |R(E+\mathrm i\eta)_{xy}|^2
        \leq C_{\varepsilon,\delta}\lambda^{-\delta/2}.
\end{equation}
To check this, write \(x=(I_x,u)\) and use the block-diagonal source
\[
        \mathcal A_x
        =\bigoplus_I
        \chi\left(\frac{W(I-I_x)}{8r}\right)\I_{W^d},
\]
with \(\chi\) as in \cref{eq:BA-radial-cutoff}. In cases \textnormal{(a)} and \textnormal{(b)}, \(\lambda^2/r\leq\lambda^{-\delta}\), whereas in case \textnormal{(c)}, \(W/r\leq\lambda^{-\delta}\). In particular, \(W=o(r)\) in every case. Thus \(\mathcal A_x\) dominates the indicator of \(B_{4r}(x)\) for all sufficiently large \(\lambda\). The finite-volume assumption implies \(r=o(L)\) and ensures the support condition for this source.

The dimensionless radius of \(\mathcal A_x\) is comparable with \(8\lambda^{-\delta/2}\).  Apply Lemma~\ref{lemma:BA-radial-resolvent-profile} in case \textnormal{(a)} and \cref{eq:BA-large-block-small-ball} in cases \textnormal{(b)} and \textnormal{(c)}. The first deterministic term in either bound is therefore at most
\[
        C\lambda^{-\delta}
        \left(1+\bm1_{\{d=2\}}\log\lambda\right).
\]
The gradient error is bounded by \(C\lambda^2/r\) in case \textnormal{(a)} and by \(C(\lambda^2+W)/r\) in the growing-block cases.  These are at most \(C\lambda^{-\delta}\) throughout the stated ranges. The remaining deterministic error satisfies
\[
        \lambda\eta\leq\lambda^{-3\delta}
        \quad\text{in cases (a) and (b)},
        \qquad
        \frac{W\eta}{\lambda}\leq\lambda^{-3\delta}
        \quad\text{in case (c)}.
\]
By \cref{eq:BA-dynamical-Phi-smallness}, the normalized random error is stochastically dominated by \(C\lambda^{-2\delta}\). The energy net following Lemma~\ref{lemma:BA-radial-resolvent-profile}, followed by a union bound over the polynomially many sites, makes this estimate uniform in \(E\) and \(x\). This proves \cref{eq:BA-eigenfunction-small-ball}, including its stated probability. Translation here is by whole blocks; the existing row estimates already hold uniformly over the internal coordinate \(u\).

Choose an \(r\)-net of \(\mathbb Z_L^d\) with at most \(C(L/r)^d\) points. If \(E_j\in\mathfrak L(r)\), its localization ball is contained in \(B_{2r}(x_0)\) for some point \(x_0\) of this net. Consequently, for a specific choice \(P=\mathbf1_{B_{2r}(x_0)}\), \(\|P\psi_j\|_2^2\geq1/2\) for all sufficiently large \(\lambda\). For any \(E_0\) in the energy interval under consideration, \cref{eq:BA-eigenfunction-small-ball} gives
\[
\begin{aligned}
        \|PR(E_0+\mathrm i\eta)P\|_{\mathrm{HS}}^2
        &=\sum_{x,y\in B_{2r}(x_0)}
                    |R(E_0+\mathrm i\eta)_{xy}|^2\\
        &\leq
        \frac{C_{\varepsilon,\delta}r^d}{\lambda\eta}
                    \lambda^{-\delta/2}.
\end{aligned}
\]
Lemma~\ref{lemma:BA-eigenfunction-compression} thus bounds the number of these eigenfunctions with \(|E_j-E_0|\leq\eta\) by $C_{\varepsilon,\delta}\lambda^{-\delta/2} \eta r^d/\lambda.$ Covering \([\lambda(e-\varepsilon),\lambda(e+\varepsilon)]\) by at most \(C\lambda\varepsilon/\eta\) such intervals and summing over the spatial net yields
\begin{equation}
\label{eq:BA-localized-eigenvalue-count}
        \card\bigl(\mathfrak L(r)\cap
              [\lambda(e-\varepsilon),\lambda(e+\varepsilon)]\bigr)
        \leq C_{\varepsilon,\delta}
                    \lambda^{-\delta/2}\varepsilon L^d.
\end{equation}

It remains to bound the denominator. The GOE norm estimate used in \cref{eq:BA-functional-calculus-comparison} implies \(\|\mathcal V\|_{2\to2}=o(\lambda\varepsilon)\) with probability at least \(1-C_{\varepsilon}\lambda^{-1000}\). By Weyl's inequality, on this event every eigenvalue of \(\Psi_L\) in \([e-\varepsilon/2,e+\varepsilon/2]\) gives, with multiplicity, an eigenvalue of \(H_L/\lambda\) in \([e-\varepsilon,e+\varepsilon]\). The free Fourier transform and the Riemann sum give
\[
\begin{aligned}
        \frac1{L^d}\card\bigl(\spec(\Psi_L)\cap
                    [e-\varepsilon/2,e+\varepsilon/2]\bigr)
        &\longrightarrow
        \frac1{(2\pi)^d}
        \int_{\mathbb T^d}
        \mathbf1_{\{|\Omega(p)-e|\leq\varepsilon/2\}}\,\rmd p
        \geq c_{\varepsilon}\varepsilon.
\end{aligned}
\]
The last lower bound follows from the positive free density of states in the regular bulk window. Therefore
\[
        \card\bigl(\spec(H_L)\cap
              [\lambda(e-\varepsilon),\lambda(e+\varepsilon)]\bigr)
        \geq c_{\varepsilon}\varepsilon L^d.
\]
Combining this with \cref{eq:BA-localized-eigenvalue-count} proves the theorem.
\end{proof}

\section{Common deterministic tools}
\label{sec:common-deterministic-tools}

This section records the deterministic inputs shared by the growing- and fixed-block analyses.  The partial Bloch transform reduces a block-translation-invariant operator to finite-dimensional cell fibers, while the numerical-range and Combes--Thomas estimates control the resulting non-self-adjoint resolvents without invoking a spectral theorem.  We first recall the Bloch representation and then prove the dissipative resolvent bounds used in the finite-volume comparisons.

\subsection{Partial Bloch transform}

For orientation, let
\[
        M_s=(\lambda\Psi-z-s\Theta)^{-1},
        \qquad 0\leq s\leq1.
\]
The resolvent identity, or equivalently the fundamental theorem of calculus, gives the exact relation
\[
        M_1-M_0=\int_0^1M_s\Theta M_s\,\rmd s.
\]
Estimating this formula directly requires column norms of \(M_s\), which at the level of the bare spectral parameter carry an unfavorable dependence on \(\eta\).  The deterministic stability arguments below instead use the dissipation generated self-consistently by the scalar part of \(\Theta\).

Write the $d$-dimensional discrete torus $\widetilde{\mathbb T}_n^d
        =\left(\frac{2\pi}{n}\mathbb Z\right)^d/(2\pi\mathbb Z)^d.$
The partial Bloch transform in the block coordinate is defined as
\[
        (\mathcal U\psi)(\xi,u)
        =n^{-d/2}\sum_{b\in\widetilde{\mathbb Z}_n^d}
        e^{-\mathrm i\xi\cdot b}\psi(Wb+u),
        \qquad
        \xi\in\widetilde{\mathbb T}_n^d,\qquad u\in[W]^d.
\]
If \(A\) is block-translation invariant with block kernel \(A(b)\), then
\[
        \mathcal U A\mathcal U^\star
        =\bigoplus_{\xi\in\widetilde{\mathbb T}_n^d}\widehat A(\xi),
        \qquad
        \widehat A(\xi)=
        \sum_{b\in\widetilde{\mathbb Z}_n^d}e^{-\mathrm i\xi\cdot b}A(b).
\]
Moreover, \((\mathcal U\delta_{(0,u)})(\xi)=n^{-d/2}e_u\).  Consequently, if \(\Theta\) has constant cell symbol \(\theta_W\), then
\begin{equation}
\label{equ:Bloch_decomp}
        \bigl[(\lambda\Psi_L-z-\Theta)^{-1}\bigr]_{(0,u),(0,u')}
        =\frac1{n^d}\sum_{\xi\in\widetilde{\mathbb T}_n^d}
        \left[
        (\lambda\widehat\Psi_W(\xi)-z-\theta_W)^{-1}
        \right]_{u,u'}.
\end{equation}
Clearly, for the free operator the microscopic Fourier symbol is
\begin{equation}
\label{equ:micro_fourier_symbol}
    \begin{aligned}
         \Omega(p)=2d-2\sum_{j=1}^d\cos p_j,
    \end{aligned}
\end{equation}
and hence
\[
        (\lambda\Psi_L-z)^{-1}_{(0,0)}
        =\frac1{L^d}\sum_{p\in\widetilde{\mathbb T}_L^d}
        \frac1{\lambda\Omega(p)-z}.
\]

\subsection{Numerical-range and Combes--Thomas estimates}

The following elementary estimate replaces the spectral-theorem input for the non-self-adjoint operator \(\lambda\Psi-z-\Theta\).
\begin{prop}[Numerical-range bound and generalized Ward identity]
\label{prop:num_range}
Let $A$ be a finite-dimensional operator with $\Im A \coloneqq \frac{A - A^\star}{2i} \preceq -\delta\,\I$ for some $\delta > 0$.  Then $A$ is invertible,
\begin{equation}
\label{eq:num_range_bd}
    \|A^{-1}\|_{2\to2}\le \delta^{-1},
\end{equation}
and the generalized Ward identity holds in both orderings:
\begin{equation}
\label{eq:gen_ward}
    \Im (A^{-1})
    =
    A^{-1}(-\Im A)(A^{-1})^\star
    =
    (A^{-1})^\star(-\Im A)A^{-1}
    \succeq
    \delta\, A^{-1}(A^{-1})^\star
    \ \text{ and }\
    \succeq \delta\,(A^{-1})^\star A^{-1}.
\end{equation}
In particular, for $\Theta$ block-diagonal with $\Im\Theta \succeq \delta_0\I$, $\delta_0 \ge 0$, and $\Im z = \eta > 0$, the operator $M = (\lambda\Psi_L - z - \Theta)^{-1}$ exists, satisfies $\|M\|_{2\to2} \le (\eta + \delta_0)^{-1}$, and for every $x$,
\begin{equation}
\label{eq:ward_column}
    \big\|M\rvert x\rangle\big\|_2^2
    \;\vee\;
    \big\|M^\star\rvert x\rangle\big\|_2^2
    \;\le\;
    \frac{\Im M_{(x,x)}}{\eta + \delta_0}.
\end{equation}
\end{prop}

\begin{proof}
For any vector $v$, $\Im\langle v, Av\rangle \le -\delta\|v\|_2^2$, hence $\|Av\|_2\|v\|_2 \ge |\langle v, Av\rangle| \ge \delta\|v\|_2^2$.  Thus $A$ is injective, hence invertible in finite dimension, and substituting $v = A^{-1}u$ gives~\cref{eq:num_range_bd}.  The identity $A^{-1}-(A^\star)^{-1} = A^{-1}(A^\star - A)(A^\star)^{-1} = (A^\star)^{-1}(A^\star - A)A^{-1}$, divided by $2i$, gives~\cref{eq:gen_ward}; note $(A^\star)^{-1} = (A^{-1})^\star$.  For the application take $A = \lambda\Psi_L - z - \Theta$, so that $-\Im A = \eta\I + \Im\Theta \succeq (\eta+\delta_0)\I$; evaluating the two orderings of~\cref{eq:gen_ward} at $\langle x|\cdot|x\rangle$ yields~\cref{eq:ward_column}.
\end{proof}

The same argument applies to bounded operators on an infinite-dimensional Hilbert space; we record this form for the Combes--Thomas argument.
\begin{prop}[Infinite-dimensional numerical-range bound]
\label{prop:num_range_inf_d}
Let \(A_\infty\in\mathcal B(\mathcal H)\) satisfy
\[
        \Im A_\infty
        \coloneqq\frac{A_\infty-A_\infty^\star}{2\mathrm i}
        \preceq-\delta\I
\]
for some \(\delta>0\).  Then \(A_\infty\) is invertible,
\begin{equation}
    \|A_\infty^{-1}\|_{2\to2}\le \delta^{-1},
\end{equation}
and the generalized Ward identity holds in both orderings:
\begin{equation}
    \Im (A_\infty^{-1})
    =
    A_\infty^{-1}(-\Im A_\infty)(A_\infty^{-1})^\star
    =
    (A_\infty^{-1})^\star(-\Im A_\infty)A_\infty^{-1}
    \succeq
    \delta\, A_\infty^{-1}(A_\infty^{-1})^\star
    \ \text{ and }\
    \succeq \delta\,(A_\infty^{-1})^\star A_\infty^{-1}.
\end{equation}
\end{prop}

\begin{proof}
For every \(u\in\mathcal H\),
\[
        \|A_\infty u\|\,\|u\|
        \geq|\langle u,A_\infty u\rangle|
        \geq\delta\|u\|^2.
\]
Thus \(A_\infty\) is bounded below and has closed range.  The same estimate for \(A_\infty^\star\) gives \(\ker A_\infty^\star=\{0\}\), hence \(\overline{\operatorname{ran}A_\infty} =(\ker A_\infty^\star)^\perp=\mathcal H\).  Its range is therefore both closed and dense, so \(A_\infty\) is invertible with \(\|A_\infty^{-1}\|\leq\delta^{-1}\).  The two Ward identities follow from the same resolvent identities as in Proposition~\ref{prop:num_range}.
\end{proof}

\begin{lemma}[Non-self-adjoint Combes--Thomas comparison]
\label{lemma:CT_est_block_renorm_resolv}
Fix \(K<\infty\), let \(z=E+\mathrm i\eta\) with \(0<\eta\leq K\lambda\), and write \(L=nW\).  Suppose that \(\Theta\) is block-diagonal and block-translation invariant, with constant cell symbol \(\theta_W\) and \(\Im\Theta\succeq0\).  Denote by
\[
        R_\Theta(z)=(\lambda\Psi-z-\Theta)^{-1},
        \qquad
        R_{\Theta,n}(z)=(\lambda\Psi_L-z-\Theta)^{-1}
\]
the infinite-volume resolvent and its periodization to the block torus \(\widetilde{\mathbb Z}_n^d\), respectively.  Choose representatives of \(a,b\) such that \(|a-b|_\infty\leq n/2\).  If \(n\eta/\lambda\geq1\), then
\[
        \left\|\Pi_{[a]}R_\Theta(z)\Pi_{[b]}
        -\Pi_{[a]}R_{\Theta,n}(z)\Pi_{[b]}\right\|_{2\to2}
        \leq
        \frac{C_{d,K}}{\eta}
        \exp\left(-c_{d,K}\frac{n\eta}{\lambda}\right).
\]
\end{lemma}

\begin{proof}
Define the corresponding block kernels by
\[
        G_\Theta(r,z)
        =
        \frac1{(2\pi)^d}\int_{\mathbb T^d}
        e^{\mathrm i r\cdot\xi}
        (\lambda\widehat\Psi_W(\xi)-z-\theta_W)^{-1}\,\rmd\xi
\]
and
\[
        G_{\Theta,n}(r,z)
        =
        \frac1{n^d}\sum_{\xi\in\widetilde{\mathbb T}_n^d}
        e^{\mathrm i r\cdot\xi}
        (\lambda\widehat\Psi_W(\xi)-z-\theta_W)^{-1}.
\]
Since
\[
        \frac1{n^d}\sum_{\xi\in\widetilde{\mathbb T}_n^d}
        e^{-\mathrm i(l-r)\cdot\xi}
        =
        \mathbf 1_{\{l-r\in n\mathbb Z^d\}},
\]
Fourier inversion gives the exact periodization identity
\begin{equation}
\label{eq:CT-periodization}
        G_{\Theta,n}(r,z)
        =
        \sum_{k\in\mathbb Z^d}G_\Theta(r+nk,z).
\end{equation}

It remains to estimate the infinite-volume kernel.  Fix a block \(b\), set
\[
        \rho_N(x)=\min\{d_\infty(x,b),N\},
        \qquad
        Q_{\mu,N}=\sum_{x\in\mathbb Z^d}e^{\mu\rho_N(x)}\Pi_x,
\]
and note that \(Q_{\mu,N}\) commutes with \(\Theta\).  At block level the hopping can only have range one: \(\Pi_x\Psi\Pi_y=0\) whenever \(d_\infty(x,y)>1\).  Schur's test therefore gives, uniformly in \(N\),
\[
        \left\|
        \lambda(Q_{\mu,N}\Psi Q_{\mu,N}^{-1}-\Psi)
        \right\|
        \leq 4d\lambda(e^\mu-1)
        \eqqcolon s_\Psi(\mu).
\]
The numerical-range estimate applied to the conjugated operator yields
\[
        \left\|
        Q_{\mu,N}R_\Theta(z)Q_{\mu,N}^{-1}
        \right\|
        \leq\frac1{\eta-s_\Psi(\mu)}
\]
whenever \(s_\Psi(\mu)<\eta\).  With
\[
        \mu=\log\left(1+\frac{\eta}{9d\lambda}\right)
\]
we have \(s_\Psi(\mu)=4\eta/9<\eta/2\).  Taking \(N\) beyond \(d_\infty(a,b)\) gives
\begin{equation}
\label{eq:CT-infinite-kernel}
        \|\Pi_{[a]}R_\Theta(z)\Pi_{[b]}\|_{2\to2}
        \leq
        \frac2\eta e^{-\mu d_\infty(a,b)}.
\end{equation}
Moreover, since \(0<\eta\leq K\lambda\),
\[
        \mu
        \geq
        \frac{\eta/(9d\lambda)}{1+\eta/(9d\lambda)}
        \geq c_{d,K}\frac{\eta}{\lambda}.
\]

Finally, let \(r=a-b\).  For \(k\neq0\), the choice \(|r|_\infty\leq n/2\) gives
\[
        |r+nk|_\infty
        \geq n\left(|k|_\infty-\frac12\right).
\]
Combining \cref{eq:CT-periodization,eq:CT-infinite-kernel} and using \(\mu n\geq c_{d,K}n\eta/\lambda\geq c_{d,K}\), we obtain
\[
        \|G_{\Theta,n}(r,z)-G_\Theta(r,z)\|_{2\to2}
        \leq
        \frac2\eta\sum_{k\neq0}e^{-\mu|r+nk|_\infty}
        \leq
        \frac{C_{d,K}}{\eta}
        e^{-c_{d,K}n\eta/\lambda}.
\]
This proves the claimed comparison.
\end{proof}

\section{Growing blocks: scalar reduction and deterministic stability}
\label{sec:large-block-stability}
\label{sec:scalar_reduction}

This section proves the deterministic estimates required by the local-law bootstrap when the block size grows.  For the MDE solution
\[
        M=(\lambda\Psi_L-z-\Theta)^{-1},
        \qquad \Theta=S[M],
\]
we obtain, uniformly in \(\eta>0\), diagonal bounds of order \(\lambda^{-1}\), a block-trace lower bound of the same order, and \(\eta\)-uniform column bounds.  These serve as the main inputs that allow the fluctuation estimate to continue from the kinetic scale to the diffusive scale.

The main driving mechanism is a scalar reduction.  We decompose the cell self-energy as \(\theta_W=g\I+\vartheta_W\), where \(g\) is its normalized trace. The normalized block map and a Schur--H\" older estimate produce a shrinking factor \(W^{-d/p}\), which makes the traceless component \(\vartheta_W\) perturbative under the standing block-size condition.  The scalar term is absorbed into the spectral parameter, \(z\mapsto z+g\); its positive imaginary part supplies dissipation of order \(\lambda^{-1}\).  The remaining matrix-valued part is then controlled by a resolvent expansion around this shifted free operator.

In the proof, we first collect the free-resolvent bounds and define an invariant convex class for the self-energy.  The self-consistent map preserves this class, so Brouwer's theorem gives a fixed point, which is identified with the MDE solution by uniqueness. We then transfer the resulting resolvent estimates to the local-law comparison object.  The final subsection proves stability under the residual in the approximate self-consistent equation by interpolating between the two self-energies.  The \(W^{-d/p}\) gain controls the non-scalar component, while free density-of-states regularity controls the scalar component.

\subsection{Setting, free-resolvent inputs, and the a priori class}

Throughout this section $d \ge 2$, $z = E + i\eta$ with $0 < \eta \le 1$, $E/\lambda\in(0,4d)$, and $d(E/\lambda, \Sigma_d^\Psi) \ge \varepsilon$.  We fix an exponent $p > p_d$.  We abbreviate
\[
    \ell \coloneqq 1 + \log\lambda ,
    \qquad
    R_0(w) \coloneqq (\lambda\Psi_L - w)^{-1}
    \quad (\Im w > 0),
\]
and we write $\langle B \rangle \coloneqq W^{-d}\Tr B$ for the normalized trace of a $W^d\times W^d$ matrix $B$.  For a $W^d\times W^d$ matrix $\theta$ we denote by $\Theta(\theta) \coloneqq \sum_{b\in\tilde{\mathbb Z}_n^d} \iota_b\theta\iota_b^*$ its block-diagonal, block-translation-invariant lift to $\ell^2(\mathbb Z_L^d)$; note $\|\Theta(\theta)\|_{2\to2} = \|\theta\|_{2\to2}$ and $\Theta(\I_{W^d}) = \I_{L^d}$.  In particular, decomposing
\begin{equation}
\label{eq:scalar_split}
    \theta = g\,\I + \vartheta,
    \qquad
    g \coloneqq \langle\theta\rangle,
    \quad
    \vartheta \coloneqq \theta - \langle\theta\rangle\I ,
\end{equation}
we have the exact identity
\begin{equation}
\label{eq:shifted_param}
    \lambda\Psi_L - z - \Theta(\theta)
    =
    \lambda\Psi_L - \tilde z - \Theta(\vartheta),
    \qquad
    \tilde z \coloneqq z + g .
\end{equation}

We will use two facts about the \emph{free} operator only.  The first is the density-of-states estimation of \citep[App.~A.2--A.3]{black2025selfconsistentequationsquantumdiffusion}, transported to $\Psi = 2d\I - \Delta$ (a recentering of the dispersion relation): writing $\Omega(\xi) = 2d - 2\sum_{i}\cos\xi_i$ and
\[
    \varphi(\zeta)
    \coloneqq
    \frac{1}{(2\pi)^d}\int_{\mathbb T^d}\frac{\rmd \xi}{\Omega(\xi)-\zeta}
    =
    \int_{\mathbb R}\frac{\rho(\sigma)}{\sigma - \zeta}\,\rmd\sigma,
    \qquad
    \zeta \in \mathbb H,
\]
where $\rho$ is the free density of states, we have:

\begin{prop}[Free DOS regularity; Lemmas A.4--A.5 of \citep{black2025selfconsistentequationsquantumdiffusion}]
\label{prop:free_dos_package}
The density $\rho$ is smooth on $\mathbb R\setminus\Sigma_d^\Psi$ and there exist constants $c_\rho, C_\varphi$ depending only on $(d,\varepsilon)$ such that on the domain $D_{\varepsilon/2} \coloneqq \{\zeta\in\mathbb H : d(\Re\zeta, \Sigma_d^\Psi)\ge \varepsilon/2, \ \Im\zeta \le 1\}$,
\begin{equation}
\label{eq:phi_reg}
    \|\varphi\|_{C^1(D_{\varepsilon/2})} \le C_\varphi,
    \qquad
    \Im\varphi(\zeta) \ge c_\rho
    \quad\text{whenever }\Re\zeta \in (0, 4d).
\end{equation}
The same bounds hold for the finite-volume function $\varphi_L(\zeta) \coloneqq \lambda\,(\lambda\Psi_L - \lambda\zeta)^{-1}_{(0,0)}$ on the subdomain $\{\Im\zeta \ge \lambda^{-3}\}$, with error $|\varphi_L - \varphi| \le e^{-cL\Im\zeta}$ (\citep[Lemma A.2]{black2025selfconsistentequationsquantumdiffusion}), under the standing polynomial-volume assumption.
\end{prop}

The second input is the uniform-Sobolev-type bound for the free resolvent, which follows by composing the two halves of Proposition~\ref{prop:tomas_stein} (equivalently \citep[Prop. 3.7]{black2025selfconsistentequationsquantumdiffusion}; the underlying surface-measure decay is due to \citep{erdos2006quantumdiffusionrandomschrodinger,cuenin2021fouriertransformsurfacecarriedmeasures}). Since $\lambda\Psi_L-w$ is normal, for the dimension-dependent pair $p,p'$ fixed above and every $w$ with $w/\lambda\in D_{\varepsilon/2}$,
\begin{equation}
\label{eq:large-block-unif-sobolev}
    \|R_0(w)\|_{p'\to p}
    \le
    \|R_0(w)^{1/2}\|_{p'\to 2}\,\|R_0(w)^{1/2}\|_{2\to p}
    \lesssim
    \frac{1}{\lambda}\Big(1 + \big|\log(\Im w/\lambda)\big|\Big)^{2}.
\end{equation}

Finally, we fix, with constants $K_1, c_1, \varepsilon_*$ to be chosen in Proposition~\ref{prop:self_map}, the a priori class of self-energy symbols
\begin{equation}
\label{eq:Ksharp}
    \mathcal K_\sharp
    \coloneqq
    \Big\{\theta\in\mathbb C^{W^d\times W^d}
    :\ \Im\theta \succeq 0,\ \
    |\langle\theta\rangle| \le \tfrac{K_1}{\lambda},\ \
    \Im\langle\theta\rangle \ge \tfrac{c_1}{\lambda},\ \
    \|\theta - \langle\theta\rangle\I\|_{2\to2}\le \tfrac{\varepsilon_*}{\lambda}
    \Big\}.
\end{equation}
$\mathcal K_\sharp$ is convex, compact, and nonempty (it contains $\theta = 2ic_1\lambda^{-1}\I$ provided $K_1 \ge 2c_1$).

\subsection{Dissipative resolvent bounds without Hermiticity}
\begin{lemma}
\label{lemma:free_seed}
Let $w$ satisfy $w/\lambda \in D_{\varepsilon/2}$ and $\Im w \ge c_1/\lambda$.  Then, with constants depending on $(d,\varepsilon,p,c_1)$ and \emph{uniformly in $\eta$},
\begin{align}
    \|R_0(w)\|_{1\to\infty} &\le C\,\ell^{2}\,\lambda^{-1},
    \label{eq:free_1inf}\\
    \|R_0(w)\|_{1\to2} = \|R_0(w)\|_{2\to\infty} &\le C\,\ell,
    \label{eq:free_12}\\
    \|R_0(w)\|_{2\to p} &\le C\,\ell,
    \label{eq:free_2p}\\
    \|R_0(w)\|_{2\to2} &\le \lambda/c_1,
    \label{eq:free_22}\\
    W^{-d}\sup_{u\in[W]^d}\sum_{v\in[W]^d}\big|R_0(w)_{(u,v)}\big|
    &\le C\,\ell^{2}\,\lambda^{-1}\,W^{-d/p}.
    \label{eq:free_seed}
\end{align}
\end{lemma}

\begin{proof}
Since $\Im(w/\lambda) \ge c_1/\lambda^2$, the logarithmic factor in~\cref{eq:large-block-unif-sobolev} is $\lesssim \ell$, so $\|R_0(w)\|_{p'\to p}\lesssim \ell^2/\lambda$.  Then~\cref{eq:free_1inf} follows from H\"older, $|\langle\delta_x, R_0(w)\delta_y\rangle| \le \|\delta_x\|_{p'}\|R_0(w)\delta_y\|_{p} \le \|R_0(w)\|_{p'\to p}$.  For~\cref{eq:free_12}, the Ward identity~\cref{eq:ward_column} for the (normal, dissipative) operator $\lambda\Psi_L - w$ gives $\|R_0(w)\delta_x\|_2^2 \le \Im R_0(w)_{(x,x)}/\Im w \le (C\ell^2/\lambda)\cdot(\lambda/c_1)$, and the same for $R_0(w)^\star\delta_x$. For~\cref{eq:free_2p}, factor through the square root, which is again defined by the functional calculus of the normal operator $\lambda\Psi_L - w$: $\|R_0(w)\|_{2\to p}\le\|R_0(w)^{1/2}\|_{2\to p}\,\|R_0(w)^{1/2}\|_{2\to2} \lesssim \lambda^{-1/2}\ell\cdot(\Im w)^{-1/2} \lesssim \ell$, by Proposition~\ref{prop:tomas_stein}. The bound~\cref{eq:free_22} is~\cref{eq:num_range_bd}. Finally, for~\cref{eq:free_seed}, apply H\"older's inequality over the $W^d$ sites of a block:
\[
    \sum_{v\in[W]^d}\big|R_0(w)_{(u,v)}\big|
    \le
    \big(W^d\big)^{1/p'}\,
    \big\|R_0(w)\delta_u\big\|_{\ell^p}
    \le
    W^{d/p'}\,\|R_0(w)\|_{p'\to p}
    \lesssim
    W^{d/p'}\,\frac{\ell^{2}}{\lambda},
\]
and multiply by $W^{-d}$.
\end{proof}

\subsection{The stability estimate in the large block regime}

The next proposition is the core of the section.  We emphasize that it does not rely on self-consistency: it holds for any block self-energy whose symbol has a dissipative scalar part of size $\asymp\lambda^{-1}$ and an off-scalar part of comparably shrinking size.  This robustness is what allows it to be applied both to the MDE self-energy $\Theta$ and, once scalarization has been established, to the solution of approximate self-consistent equation with solution given by self-energy $\tilde\Theta$. The fixed-neighborhood estimate used to start the bootstrap is then proved in Corollary~\ref{coro:robust_app_bootstrap}.

\begin{prop}
\label{prop:stability_renorm_resolv}
Let $\theta = g\I + \vartheta$ as in~\cref{eq:scalar_split} with
\[
    \Im g \ge \frac{c_1}{\lambda},
    \qquad
    |g| \le \frac{K_1}{\lambda},
    \qquad
    \|\vartheta\|_{2\to2}\le\frac{\varepsilon_*}{\lambda},
    \qquad
    \varepsilon_* \le \frac{c_1}{2},
\]
and let $z = E + i\eta$, $\eta\in(0,1]$, $E/\lambda\in(0,4d)$, $d(E/\lambda,\Sigma_d^\Psi)\ge\varepsilon$, and $\lambda \ge \lambda_0(d,\varepsilon,p,K_1)$.  Set $M \coloneqq (\lambda\Psi_L - z - \Theta(\theta))^{-1}$ and $\tilde z \coloneqq z + g$.  Then, uniformly in $\eta$:
\begin{enumerate}[(a)]
    \item $M = \sum_{k\ge0}\big(R_0(\tilde z)\,\Theta(\vartheta)\big)^k R_0(\tilde z)$, and
    \begin{equation}
    \label{eq:neumann_1inf}
        \|M - R_0(\tilde z)\|_{1\to\infty}
        \le
        C_\flat\,\ell^{2}\,\frac{\varepsilon_*}{\lambda}.
    \end{equation}
    \item $\|M\|_{1\to\infty}\le C\ell^{2}/\lambda$,\ \ $\|M\|_{2\to2}\le 2\lambda/c_1$,\ \ and $\|M\|_{1\to2}\vee\|M\|_{2\to\infty} \le C\ell$.
    \item For every $x \in \mathbb Z_L^d$,
    \begin{equation}
    \label{eq:diag_close}
        \Big|M_{(x,x)} - \frac{1}{\lambda}\varphi\,\!\Big(\frac z\lambda\Big)\Big|
        \le
        \frac{C_\flat}{\lambda}\Big(\ell^{2}\varepsilon_*
        + \frac{C_\varphi K_1}{\lambda}\Big)
        + \lambda^{-1}e^{-cL/\lambda^2} .
    \end{equation}
    \item  With $m_W \coloneqq \Pi_{[0]}M\,\Pi_{[0]}$ viewed as a $W^d\times W^d$ matrix,
    \begin{equation}
    \label{eq:offscalar_small}
        W^{-d}\big\|m_W\big\|_{2\to2}
        =
        W^{-d}\big\|m_W^\top\big\|_{2\to2}
        \le
        C_\sharp\,\ell^{2}\,\frac{W^{-d/p}}{\lambda}\,(1+\varepsilon_*).
    \end{equation}
\end{enumerate}
The constants $C_\flat, C_\sharp, C$ depend only on $(d,\varepsilon,p,c_1,K_1)$.
\end{prop}

\begin{proof}
Write $R_0 = R_0(\tilde z)$ and $\widehat\Theta = \Theta(\vartheta)$ throughout; by~\cref{eq:shifted_param}, $M^{-1} = R_0^{-1} - \widehat\Theta$.  Note that $\tilde z/\lambda \in D_{\varepsilon/2}$: indeed $|\Re\tilde z/\lambda - E/\lambda| \le K_1/\lambda^2 \le \varepsilon/2$ for $\lambda \ge \lambda_0$, and $0 < \Im\tilde z/\lambda \le (1 + K_1/\lambda)/\lambda \le 1$.  Also $\Im\tilde z \ge \Im g \ge c_1/\lambda$, so all of Lemma~\ref{lemma:free_seed} is available at $\tilde z$.

\emph{(a).}  By~\cref{eq:free_22}, $\|R_0\widehat\Theta\|_{2\to2} \le (\lambda/c_1)(\varepsilon_*/\lambda) = \varepsilon_*/c_1 \le 1/2$, so the Neumann series converges in operator norm and equals $M$. For $k\ge1$, factor the $k$-th term as $1 \to 2 \to 2 \to \cdots \to 2 \to \infty$:
\[
    \big\|(R_0\widehat\Theta)^kR_0\big\|_{1\to\infty}
    \le
    \|R_0\|_{2\to\infty}\,\|\vartheta\|
    \,\big(\|R_0\|_{2\to2}\|\vartheta\|\big)^{k-1}\,
    \|R_0\|_{1\to2}
    \le
    C\ell^{2}\,\frac{\varepsilon_*}{\lambda}\,2^{-(k-1)},
\]
using~\cref{eq:free_12} and $\|\Theta(\vartheta)\|_{2\to2}=\|\vartheta\|_{2\to2}$; summing over $k\ge1$ gives~\cref{eq:neumann_1inf}.

\emph{(b).}  The first bound is~\cref{eq:free_1inf} plus~\cref{eq:neumann_1inf}.  For the second, $\Im\theta \succeq (\Im g - \|\vartheta\|)\I \succeq \frac{c_1}{2\lambda}\I$, so Proposition~\ref{prop:num_range} gives $\|M\|_{2\to2}\le(\eta+c_1/(2\lambda))^{-1}\le 2\lambda/c_1$. For the third, apply the Ward bound~\cref{eq:ward_column} with $\delta_0=c_1/(2\lambda)$ and the first bound: $\|M\delta_x\|_2^2 \le (C\ell^2/\lambda)(2\lambda/c_1) = C'\ell^2$, and likewise for $M^\star\delta_x$.

\emph{(c).}  Decompose $M_{(x,x)} = R_0(\tilde z)_{(x,x)} + (M - R_0(\tilde z))_{(x,x)}$.  The second term is bounded by~\cref{eq:neumann_1inf}.  Since $\lambda\Psi_L$ is translation invariant on the torus, all diagonal entries of $R_0(\tilde z)$ are equal to $L^{-d}\Tr R_0(\tilde z) = \lambda^{-1}\varphi_L(\tilde z/\lambda)$. By Proposition~\ref{prop:free_dos_package}, this differs from $\lambda^{-1}\varphi(\tilde z/\lambda)$ by at most $\lambda^{-1}e^{-cL/\lambda^{2}}$. Finally, the segment from $z/\lambda$ to $\tilde z/\lambda$ lies in $D_{\varepsilon/2}$ (the real part moves by at most $K_1/\lambda^2$, the imaginary part increases from $\eta/\lambda$), so the $C^1$ bound~\cref{eq:phi_reg} gives $|\varphi(\tilde z/\lambda) - \varphi(z/\lambda)| \le C_\varphi\,|g|/\lambda \le C_\varphi K_1/\lambda^{2}$. Combining the three errors yields~\cref{eq:diag_close}.  We emphasize that this step --- the matrix analogue of \citep[Prop.~A.7]{black2025selfconsistentequationsquantumdiffusion} --- uses only the regularity of the \emph{free} density of states, uniformly up to the real axis.

\emph{(d).}  Transposition preserves the operator norm.  Split $m_W = \Pi_{[0]}R_0(\tilde z)\Pi_{[0]} + \Pi_{[0]}D\,\Pi_{[0]}$ with $D \coloneqq M - R_0(\tilde z)$, and bound each piece by the Schur test
\[
    \|B\|_{2\to2}
    \le
    \Big(\sup_{u}\sum_{v}|B_{(u,v)}|\Big)^{1/2}
    \Big(\sup_{v}\sum_{u}|B_{(u,v)}|\Big)^{1/2},
    \qquad u,v\in[W^d].
\]
For the free piece, both row and column sums are controlled by~\cref{eq:free_seed} (the kernel $R_0(\tilde z)_{(u,v)}$ depends only on $u - v$ and is symmetric, since $\Psi_L^\top = \Psi_L$), giving $W^{-d}\|\Pi_{[0]}R_0(\tilde z)\Pi_{[0]}\| \le C\ell^{2}\lambda^{-1}W^{-d/p}$. For the correction, by H\"older over the block coordinate,
\[
    \sum_{v\in[W]^d}|D_{(u,v)}|
    \le
    W^{d/p'}\Big(\sum_{v\in[W]^d}|D_{(u,v)}|^p\Big)^{1/p}
    \le
    W^{d/p'}\,\big\|D^\top\delta_u\big\|_{\ell^p}
    \le
    W^{d/p'}\,\|D^\top\|_{1\to p},
\]
and, symmetrically,
\[
        \sum_u|D_{(u,v)}|
        \leq W^{d/p'}\|D\|_{1\to p}.
\]
Since $R_0^\top = R_0$, we have $D^\top = \sum_{k\ge1}(R_0\,\Theta(\vartheta)^\top)^kR_0$ with $\|\Theta(\vartheta)^\top\| = \|\vartheta\|$, so both $D$ and $D^\top$ obey the same estimate, obtained by factoring $1\to 2 \to 2 \to\cdots\to 2 \to p$:
\[
    \|D\|_{1\to p}
    \le
    \sum_{k\ge1}
    \|R_0\|_{2\to p}\,\|\vartheta\|\,
    \big(\|R_0\|_{2\to2}\|\vartheta\|\big)^{k-1}\,\|R_0\|_{1\to2}
    \le
    2\,C\ell\cdot\frac{\varepsilon_*}{\lambda}\cdot C\ell
    =
    C'\ell^{2}\,\frac{\varepsilon_*}{\lambda},
\]
using~\cref{eq:free_2p}, \cref{eq:free_12}, and $\|R_0\|_{2\to2}\|\vartheta\|\le1/2$.  Hence $W^{-d}\|\Pi_{[0]}D\Pi_{[0]}\| \le W^{-d}\cdot W^{d/p'}\cdot C'\ell^{2}\varepsilon_*/\lambda = C'\ell^{2}\varepsilon_*\,W^{-d/p}/\lambda$, and summing the two pieces gives~\cref{eq:offscalar_small}.
\end{proof}

\begin{rmk}
The decisive point in part (d) is the prefactor $W^{-d}$ from the superoperator competing against a Schur sum over only $W^{d/p'}$ effective sites: the block average gains $W^{-d/p}$.  Heuristically, for large blocks, the random-matrix fluctuations self-average and the block self-energy behaves like a scalar. Here the mechanism is purely deterministic, driven by the off-diagonal decay of the free resolvent encoded in the Stein--Tomas bound~\cref{eq:large-block-unif-sobolev}.
\end{rmk}

\subsection{The self-consistent fixed-point map}

We write $\Phi$ for the self-consistent map at the Bloch symbol level,
\[
    \Phi(\theta)
    \coloneqq
    \mathscr S\big[\Pi_{[0]}\,(\lambda\Psi_L - z - \Theta(\theta))^{-1}\,\Pi_{[0]}\big],
    \qquad
    \mathscr S[m] \coloneqq W^{-d}\big(m^\top + \Tr(m)\,\I\big),
\]
so that $\Theta$ solves the self-consistent equation~\cref{equ:self_consist_eq} if and only if $\theta_W = \Phi(\theta_W)$.  (The map is well defined on $\{\Im\theta\succeq0\}$ by Proposition~\ref{prop:num_range}, and $\Pi_{[b]}M\Pi_{[b]}$ is independent of $b$ by block-translation invariance, so the equation is consistently posed at the symbol level.)

\begin{prop}
\label{prop:self_map}
Set
\[
    c_1 \coloneqq \frac{c_\rho}{2},
    \qquad
    K_1 \coloneqq 2\sup_{D_{\varepsilon/2}}|\varphi| + 1,
    \qquad
    \varepsilon_* \coloneqq 4C_\sharp\,\ell^{2}\,W^{-d/p}.
\]
There exist $\lambda_0$ and $C_0$, depending only on $(d,\varepsilon,p)$, such that if
\begin{equation}
\label{eq:W_condition}
    \lambda \ge \lambda_0
    \qquad\text{and}\qquad
    W^{d/p} \ \ge\ C_0\,\ell^{4},
\end{equation}
then $\Phi(\mathcal K_\sharp)\subseteq\mathcal K_\sharp$, and $\Phi$ is continuous on $\mathcal K_\sharp$.
\end{prop}

\begin{proof}
Note first that~\cref{eq:W_condition} guarantees $\varepsilon_* \le c_1/2$ and, more stringently, $C_\flat\,\ell^{2}\varepsilon_* \le c_\rho/8$ once $C_0$ is large enough; all conditions of Proposition~\ref{prop:stability_renorm_resolv} are then met for every $\theta\in\mathcal K_\sharp$ (with $g = \langle\theta\rangle$, $\vartheta = \theta - \langle\theta\rangle\I$).  Write $m = m_W(\theta)$ and $\Phi = \Phi(\theta)$.

\emph{Positivity.}  $\Im\Theta(\theta)\succeq0$ implies $\Im M = M(\eta + \Im\Theta)M^\star\succeq0$ by~\cref{eq:gen_ward}, hence $\Im m \succeq 0$; both operations in $\mathscr S$ preserve positive semidefiniteness (the transpose of a Hermitian positive-semidefinite matrix is its complex conjugate and remains positive semidefinite), so $\Im\Phi\succeq0$.

\emph{Normalized trace.}  A direct computation gives $\langle\Phi\rangle = (1 + W^{-d})\,\langle m\rangle$ with $\langle m\rangle = W^{-d}\Tr_{[0]}M$, which is an average of diagonal entries of $M$.  By~\cref{eq:diag_close},
\begin{equation}
\label{eq:trace_close}
    \Big|\langle m\rangle - \frac1\lambda\varphi\,\!\Big(\frac z\lambda\Big)\Big|
    \le
    \frac{1}{\lambda}\Big(C_\flat\ell^2\varepsilon_*
    + \frac{C_\varphi K_1}{\lambda}\Big)
    +\lambda^{-1}e^{-cL/\lambda^2}
    \le
    \frac{c_\rho}{4\lambda},
\end{equation}
for $\lambda \ge \lambda_0$ and $C_0$ large.  The standing bulk assumptions $E/\lambda\in(0,4d)$ and $d(E/\lambda,\Sigma_d^\Psi)\ge\varepsilon$, together with \cref{eq:phi_reg}, give $\Im\varphi(z/\lambda) \ge c_\rho$ and $|\varphi(z/\lambda)|\le (K_1-1)/2$; hence $\Im\langle\Phi\rangle \ge \Im\langle m\rangle \ge (c_\rho - c_\rho/4)/\lambda \ge c_1/\lambda$ and $|\langle\Phi\rangle| \le 2|\langle m\rangle| \le K_1/\lambda$.

\emph{Off-scalar part.}  We have $\Phi - \langle\Phi\rangle\I = W^{-d}m^\top - W^{-d}\langle m\rangle\I$, so by~\cref{eq:offscalar_small} and~\cref{eq:trace_close},
\[
    \big\|\Phi - \langle\Phi\rangle\I\big\|
    \le
    C_\sharp\ell^{2}\frac{W^{-d/p}}{\lambda}(1+\varepsilon_*)
    + W^{-d}\,\frac{K_1}{\lambda}
    \le
    \frac{2C_\sharp\ell^2 W^{-d/p} + K_1 W^{-d}}{\lambda}
    \le
    \frac{\varepsilon_*}{\lambda},
\]
using $W^{-d}\le W^{-d/p}$ and the definition of $\varepsilon_*$ (enlarging $C_0$ if necessary so that $K_1 W^{-d(1-1/p)} \le C_\sharp\ell^2$, which holds for all $W\ge2$ after enlarging $C_\sharp$).  Thus $\Phi(\theta)\in\mathcal K_\sharp$.

\emph{Continuity.}  On $\mathcal K_\sharp$ the map $\theta\mapsto M(\theta)$ is continuous because $\lambda\Psi_L - z - \Theta(\theta)$ is invertible with $\|M\|_{2\to2}\le\eta^{-1}$ uniformly, and $\mathscr S$, $\Pi_{[0]}$ are linear.
\end{proof}

\begin{thm}[A priori description of the self-energy]
\label{thm:theta_apriori}
Assume $d\ge2$, $z = E+i\eta$ with $\eta\in(0,1]$, $E/\lambda\in(0,4d)$, $d(E/\lambda,\Sigma_d^\Psi)\ge\varepsilon$, and condition~\cref{eq:W_condition}.  Then the self-consistent equation~\cref{equ:self_consist_eq} has a unique solution $\Theta$ with $\Im\Theta\succeq0$ in the class of block-diagonal, block-translation-invariant matrices, its symbol $\theta_W$ lies in $\mathcal K_\sharp$, and
\begin{equation}
\label{eq:theta_scalar}
    \Big\|\theta_W - \frac{1}{\lambda}\varphi\,\!\Big(\frac z\lambda\Big)\I\,
    \Big\|_{2\to2}
    \;\le\;
    \frac{C}{\lambda}\Big(\ell^{4}\,W^{-d/p} + \frac1\lambda\Big),
    \qquad
    \Im\theta_W \succeq \frac{c_\rho}{4\lambda}\,\I,
    \qquad
    \|\theta_W\|_{2\to2}\le\frac{C}{\lambda},
\end{equation}
uniformly in $\eta\in(0,1]$.
\end{thm}

\begin{proof}
\emph{Existence in $\mathcal K_\sharp$.}  $\mathcal K_\sharp$ is a nonempty convex compact subset of $\mathbb C^{W^d\times W^d}$ and $\Phi:\mathcal K_\sharp\to\mathcal K_\sharp$ is continuous (Proposition~\ref{prop:self_map}), so Brouwer's fixed-point theorem yields $\theta_*\in\mathcal K_\sharp$ with $\Phi(\theta_*)=\theta_*$.

\emph{Uniqueness.}  Suppose $\Theta$ is any solution with $\Im\Theta\succeq0$.  Then $M(\Theta)$ satisfies the matrix Dyson equation $-M^{-1} = z - \lambda\Psi_L + S[M]$ with $\Im M(\Theta)\succ0$ (by Proposition~\ref{prop:num_range}, since $\eta>0$), where $S$ is the positivity-preserving, self-adjoint covariance map already analyzed in the proof of Lemma~\ref{lemma:theta_22_est_kinet}.  By \citep{helton2007operatorvaluedsemicircularelementssolving} the matrix Dyson equation has a unique solution with positive imaginary part for each $z\in\mathbb H$ (see also the stability theory of \citep{Ajanki_2018}); hence $M(\Theta) = M(\Theta(\theta_*))$ and $\Theta = S[M(\Theta)] = \Theta(\theta_*)$.  Note that, unlike the contraction argument of Lemma~\ref{lemma:theta_22_est_kinet}, this step does not require $\lambda\eta\gtrsim1$.

\emph{Bounds.}  Living in $\mathcal K_\sharp$ gives $\Im\theta_W \succeq (\Im\langle\theta_W\rangle - \varepsilon_*/\lambda)\I \succeq (c_1 - \varepsilon_*)\lambda^{-1}\I \succeq c_\rho\lambda^{-1}\I/4$ and $\|\theta_W\| \le (K_1+\varepsilon_*)/\lambda$. For~\cref{eq:theta_scalar}, write $\theta_W - \lambda^{-1}\varphi(z/\lambda)\I = (\langle\theta_W\rangle - \lambda^{-1}\varphi(z/\lambda))\I + (\theta_W - \langle\theta_W\rangle\I)$; the second term has norm at most $\varepsilon_*/\lambda \lesssim \ell^2 W^{-d/p}/\lambda$, and the first is controlled by~\cref{eq:trace_close} applied at the fixed point, where the explicit errors are $\lesssim \lambda^{-1}(\ell^2\varepsilon_* + \lambda^{-1}) \lesssim \lambda^{-1}(\ell^4 W^{-d/p} + \lambda^{-1})$.
\end{proof}

\subsection{Resulting two-sided bounds for the fluctuation lemma}
The fluctuation control in Lemma~\ref{lemma:flut_exp_bdd} requires three ingredients: the renormalized resolvent and the free resolvent are close on the diagonal, with two-sided bounds, uniformly in $\eta$. Following theorem establishes the estimation needed in the bootstrap argument.

\begin{thm}
\label{thm:main_closeness}
Under the assumptions of Theorem~\ref{thm:theta_apriori}, let $M = (\lambda\Psi_L - z - \Theta)^{-1}$ and set $\delta_W \coloneqq \ell^{4}W^{-d/p} + \lambda^{-1}$.  Then, uniformly in $\eta\in(0,1]$ and in $x\in\mathbb Z_L^d$:
\begin{enumerate}[(a)]
    \item \emph{(Diagonal closeness.)} $\displaystyle \Big|M_{(x,x)} - \frac1\lambda\varphi\,\!\Big(\frac z\lambda\Big)\Big| \le \frac{C\,\delta_W}{\lambda}$. If moreover $\eta \ge \lambda^{-10}$, the same bound (with an additional error $e^{-c\eta L}$) holds with $\lambda^{-1}\varphi(z/\lambda)$ replaced by the finite-volume free resolvent entry $(\lambda\Psi_L - z)^{-1}_{(x,x)}$.
    \item \emph{(Two-sided diagonal bounds.)} $\displaystyle \frac{c_\rho}{2\lambda} \le \Im M_{(x,x)} \le \frac{C}{\lambda}$, and consequently $\displaystyle W^{-d}\Tr_{[0]}\big(\Im M\big) \ge \frac{c_\rho}{2\lambda}$ and $\displaystyle W^{-d}\Tr_{[0]}\big(\Im S[M]\big) \ge \frac{c_\rho}{2\lambda}$.
    \item \emph{(Entry and column bounds.)} $\|M\|_{1\to\infty}\le C\ell^{2}/\lambda$ and $\|M\rvert x\rangle\|_2^2 \,\vee\, \|M^\star\rvert x\rangle\|_2^2 \le C$.
\end{enumerate}
All constants depend only on $(d,\varepsilon,p)$.
\end{thm}

\begin{proof}
By Theorem~\ref{thm:theta_apriori}, $\theta_W\in\mathcal K_\sharp$, so Proposition~\ref{prop:stability_renorm_resolv} applies with $\varepsilon_* = 4C_\sharp\ell^2W^{-d/p}$.  Part (a) is then~\cref{eq:diag_close}; the finite-volume variant follows from Lemma~A.2 of \citep{black2025selfconsistentequationsquantumdiffusion} (or Lemma~\ref{lemma:CT_est_block_renorm_resolv} with $\Theta=0$), which requires $\eta L \gg 1$ and is the only place where a lower bound on $\eta$ enters.  Part (b): by (a) and~\cref{eq:phi_reg}, $\Im M_{(x,x)} \ge \lambda^{-1}(c_\rho - C\delta_W) \ge c_\rho/(2\lambda)$ under~\cref{eq:W_condition} (enlarging $C_0$, $\lambda_0$), while the upper bound follows from $|M_{(x,x)}|\le|\varphi|/\lambda + C\delta_W/\lambda$; averaging the lower bound over the block gives the trace bound, and $\Im S[M]\rvert_{[0]} = W^{-d}\big((\Pi_{[0]}\Im M\,\Pi_{[0]})^\top + \Tr_{[0]}(\Im M)\,\Pi_{[0]}\big) \succeq W^{-d}\Tr_{[0]}(\Im M)\,\Pi_{[0]}$ since the transposed term is positive semidefinite.  Part (c) is Proposition~\ref{prop:stability_renorm_resolv}(b) together with~\cref{eq:ward_column} and $\lambda_{\min}(\Im\theta_W)\ge c_\rho/(4\lambda)$.
\end{proof}

\begin{corollary}
\label{coro:robust_app_bootstrap}
Under the assumptions of Theorem~\ref{thm:theta_apriori}, there is a constant \(c_*>0\), independent of \(\lambda,W,L\), with the following property.  Let \(\widetilde\Theta=\Theta(\widetilde\theta_W)\) be block-diagonal and block-translation-invariant, and suppose
\[
        \|\widetilde\theta_W-\theta_W\|_{2\to2}
        \leq \frac{c_*}{\lambda}.
\]
Then \(\widetilde M=(\lambda\Psi_L-z-\widetilde\Theta)^{-1}\) exists and, uniformly in \(\eta\in(0,1]\),
\begin{equation}
\label{eq:growing-fixed-radius-bounds}
\begin{gathered}
        \frac c\lambda\I\preceq\Im\widetilde\theta_W
        \preceq\frac C\lambda\I,
        \qquad
        \frac c\lambda\leq\Im\widetilde M_{xx}\leq\frac C\lambda,\\
        \|\widetilde M\|_{2\to2}\leq C\lambda,
        \qquad
        \sup_x\bigl(\|\widetilde M e_x\|_2
                   +\|\widetilde M^*e_x\|_2\bigr)\leq C,\\
        \|\widetilde M-M\|_{1\to\infty}
        \leq C\|\widetilde\theta_W-\theta_W\|_{2\to2},
        \qquad
        \|\widetilde M\|_{1\to\infty}\leq\frac{C\ell^2}{\lambda}.
\end{gathered}
\end{equation}
If, in addition, \(\|\widetilde\theta_W-\theta_W\|\leq\varepsilon_*/(2\lambda)\), the conclusions of Proposition~\ref{prop:stability_renorm_resolv} and Theorem~\ref{thm:main_closeness} also hold for \(\widetilde M\), with adjusted constants.
\end{corollary}

\begin{proof}
Theorem~\ref{thm:theta_apriori} gives \(\Im\theta_W\succeq c\lambda^{-1}\I\) and \(\|\theta_W\|\leq C/\lambda\), so the first two-sided operator bound persists when \(c_*\) is sufficiently small.  Set \(D=\Theta(\widetilde\theta_W-\theta_W)\). The numerical-range bound gives \(\|M\|\leq C\lambda\), and hence \(\|MD\|\vee\|DM\|\leq1/2\) after decreasing \(c_*\).  The identities
\[
        \widetilde M=(\I-MD)^{-1}M
                    =M(\I-DM)^{-1}
\]
and Theorem~\ref{thm:main_closeness}(c) give the uniform column and row bounds.  Consequently the resolvent identity yields
\[
        |(\widetilde M-M)_{xy}|
        =|e_x^*\widetilde MDM e_y|
        \leq C\|D\|.
\]
The diagonal bounds follow from Theorem~\ref{thm:main_closeness}(b), on decreasing \(c_*\) once more; the entry bound follows from part (c) of that theorem. For the last assertion, the smaller perturbation preserves the hypotheses of Proposition~\ref{prop:stability_renorm_resolv}, with \((c_1,K_1,\varepsilon_*)\) replaced by \((c_1/2,K_1+1,2\varepsilon_*)\).
\end{proof}

\begin{rmk}
\label{rmk:feed_lemma}
In the proof of Lemma~\ref{lemma:flut_exp_bdd}, the chain of estimates following the generalized Ward identity now closes as follows, uniformly in $\eta$:
\[
\begin{aligned}
    \big|\mathbb ER_{(x,y)} - \tilde M_{(x,y)}\big|
    &\le
    \|\tilde M\rvert x\rangle\|_2
    \big(\mathbb E\|S[R-\mathbb ER]\|_{2\to2}^2\big)^{1/2}
    \\
    &\qquad{}\times
    \big(\mathbb E\|R\rvert y\rangle\|_2^2\big)^{1/2}
    \\
    &\le
    C\,\eta^{-1/2}
    \big(\mathbb E\|R-\mathbb ER\|_{1\to\infty}^2\big)^{1/2}
    \\
    &\qquad{}\times
    \big(\mathbb E\|R\|_{1\to\infty}\big)^{1/2},
\end{aligned}
\]
where we used Corollary~\ref{coro:robust_app_bootstrap} for $\|\tilde M\rvert x\rangle\|_2\le C$ (in place of the bound $\Im\tilde M_{(x,x)}/\eta$, which carried the kinetic restriction) and the standard Ward identity for the self-adjoint $H$ in the last factor.  The prefactor improves from $\lambda^{-1/2}\eta^{-1}$ to $\eta^{-1/2}$, as required to enter the diffusive regime $\lambda^{-1}\ll\lambda\eta\ll1$. The lower bound $\big(\tilde M^\star\Im\tilde\Theta\tilde M\big)_{(x,x)} \ge \tfrac{c_\rho}{8\lambda}\,\|\tilde M^\star\rvert x\rangle\|_2^2$ required there is immediate from $\Im\tilde\theta_W \succeq \tfrac{c_\rho}{8\lambda}\I$.
\end{rmk}

\subsection{Stability under the approximate self-consistent equation}
\label{subsec:connect_resol_self_energy_large}
It remains to compare the MDE resolvent \(M\) with the approximate self-consistent resolvent \(\widetilde M\) used in Lemma~\ref{lemma:flut_exp_bdd}. To do this, we interpolate between their cell self-energies in a fixed neighborhood of the MDE solution.  The comparison of normalized trace and the traceless component are handled separately, as the factor \(W^{-d/p}\) shrinks only the non-scalar output of the covariance map.  In the scalar direction, stability follows instead from the regularity of the free density of states.

\begin{lemma}
\label{lemma:large_block_stability}
Under the assumptions of Theorem~\ref{thm:theta_apriori}, let \(\theta\) denote the MDE cell self-energy.  After choosing \(c_*\) in Corollary~\ref{coro:robust_app_bootstrap} sufficiently small and \(C_0,\lambda_0\) in \cref{eq:W_condition} sufficiently large, every \(\widetilde\theta\) satisfying
\[
        \|\widetilde\theta-\theta\|_{2\to2}
        \leq\frac{c_*}{\lambda}
\]
obeys
\begin{equation}
\label{eq:growing-fixed-radius-stability}
        \|\widetilde\theta-\theta\|_{2\to2}
        \leq C\|\widetilde\theta-\Phi(\widetilde\theta)\|_{2\to2}.
\end{equation}
The constants are independent of \(\lambda,W,L\).
\end{lemma}

\begin{proof}
Write
\[
        \delta=\widetilde\theta-\theta=a\I+B,
        \qquad a=\langle\delta\rangle,
        \qquad \langle B\rangle=0,
\]
and set \(\theta_t=\theta+t\delta\) and \(M_t=(\lambda\Psi_L-z-\Theta(\theta_t))^{-1}\) for \(0\leq t\leq1\). Corollary~\ref{coro:robust_app_bootstrap} applies along this segment, so
\begin{equation}
\label{eq:growing-interpolation-mapping}
        \|M_t\|_{2\to2}\leq C\lambda,
        \qquad
        \|M_t\|_{1\to2}+\|M_t\|_{2\to\infty}\leq C.
\end{equation}
Decompose \(\theta_t=g_t\I+\vartheta_t\), where \(g_t=\langle\theta_t\rangle\), and write \(R_{0,t}=(\lambda\Psi_L-z-g_t)^{-1}\). By Theorem~\ref{thm:theta_apriori},
\[
        \Im g_t\geq\frac c\lambda,
        \qquad |g_t|\leq\frac C\lambda,
        \qquad
        \lambda\|\vartheta_t\|\leq\varepsilon_*+2c_*.
\]
We choose \(c_*\) small and \(C_0\) large enough that \(\|\Theta(\vartheta_t)R_{0,t}\|\leq1/2\). Lemma~\ref{lemma:free_seed} and the identity
\[
        M_t=R_{0,t}
        (\I-\Theta(\vartheta_t)R_{0,t})^{-1}
\]
then give \(\|M_t\|_{2\to p}+\|M_t^\top\|_{2\to p}\leq C\ell\). Here the transposed resolvent satisfies the same estimates because \(R_{0,t}^\top=R_{0,t}\).  The scalar Ward identity and Proposition~\ref{prop:free_dos_package} also give
\[
        \|R_{0,t}\|_{1\to2}+\|R_{0,t}\|_{2\to\infty}\leq C.
\]

Differentiating the resolvent gives
\[
        \frac{\rmd}{\rmd t}M_t=M_t\Theta(\delta)M_t.
\]
For any operator \(A\), H\"older's inequality on a cell followed by the Schur test yields
\begin{equation}
\label{equ:schur_holder_bdd}
        W^{-d}\|\Pi_{[0]}A\Pi_{[0]}\|_{2\to2}
        \leq W^{-d/p}
        \bigl(\|A\|_{1\to p}\|A^\top\|_{1\to p}\bigr)^{1/2}.
\end{equation}
The preceding mapping bounds imply
\[
        \|M_t\Theta(\delta)M_t\|_{1\to p}
        +\|(M_t\Theta(\delta)M_t)^\top\|_{1\to p}
        \leq C\ell\|\delta\|.
\]
Since
\[
        \mathscr S[A]-\langle\mathscr S[A]\rangle\I
        =W^{-d}(A^\top-\langle A\rangle\I),
\]
we obtain
\begin{equation}
\label{eq:growing-traceless-derivative}
\begin{aligned}
        &\|\Phi(\widetilde\theta)-\Phi(\theta)
          -\langle\Phi(\widetilde\theta)-\Phi(\theta)\rangle\I\|\\
        &\hspace{25mm}\leq
        C\ell W^{-d/p}\bigl(|a|+\|B\|\bigr).
\end{aligned}
\end{equation}

To estimate the trace, first use \cref{eq:growing-interpolation-mapping} to obtain
\[
        |\langle\Pi_{[0]}M_t\Theta(B)M_t\Pi_{[0]}\rangle|
        \leq C\|B\|.
\]
For the scalar direction, translation invariance and Proposition~\ref{prop:free_dos_package} give
\[
        \langle\Pi_{[0]}R_{0,t}^2\Pi_{[0]}\rangle
        =\frac1{\lambda^2}
          \varphi_L'\!\left(\frac{z+g_t}{\lambda}\right)
        =O(\lambda^{-2}).
\]
The finite-volume regularity applies since \(\Im((z+g_t)/\lambda)\geq c\lambda^{-2}\). Set \(V_t=\Theta(\vartheta_t)\).  The identity
\[
        M_t^2-R_{0,t}^2
        =R_{0,t}V_tM_t^2+R_{0,t}^2V_tM_t
\]
shows, by the uniform row and column bounds and \(\|M_t\|+\|R_{0,t}\|\leq C\lambda\), that
\[
        |(M_t^2-R_{0,t}^2)_{xx}|
        \leq C\lambda\|\vartheta_t\|
        \leq C(\varepsilon_*+2c_*).
\]
Finally, \(\langle\mathscr S[A]\rangle=(1+W^{-d})\langle A\rangle\). Integrating in \(t\), we conclude that
\begin{equation}
\label{eq:growing-scalar-derivative}
        |\langle\Phi(\widetilde\theta)-\Phi(\theta)\rangle|
        \leq C(\lambda^{-2}+\varepsilon_*+2c_*)|a|+C\|B\|.
\end{equation}

Let \(e=\widetilde\theta-\Phi(\widetilde\theta)\). As \(\theta=\Phi(\theta)\), \cref{eq:growing-traceless-derivative,eq:growing-scalar-derivative} imply, after fixing \(c_*\), \(C_0\), and \(\lambda_0\),
\[
        |a|\leq\frac14|a|+C\|B\|+\|e\|,
        \qquad
        \|B\|\leq C\ell W^{-d/p}(|a|+\|B\|)+2\|e\|.
\]
The first inequality gives \(|a|\leq C\|B\|+C\|e\|\).  Substituting this in the second and using \cref{eq:W_condition} absorbs its \(\|B\|\) term. Thus \(|a|+\|B\|\leq C\|e\|\), proving the assertion.
\end{proof}

\begin{corollary}
\label{cor:tildeM_M_compar_large_block}
Under the hypotheses of Lemma~\ref{lemma:large_block_stability}, suppose
\[
        \widetilde\theta=\Phi(\widetilde\theta)+e,
        \qquad
        e=\mathscr S\!\left[
        \Pi_{[0]}(\mathbb ER-\widetilde M)\Pi_{[0]}\right].
\]
Then
\[
        \|\widetilde M-M\|_{1\to\infty}
        \leq C\|\mathbb ER-\widetilde M\|_{1\to\infty}.
\]
\end{corollary}

\begin{proof}
The explicit covariance map and the Schur test give
\[
        \|e\|_{2\to2}
        \leq2\|\mathbb ER-\widetilde M\|_{1\to\infty}.
\]
Lemma~\ref{lemma:large_block_stability} therefore bounds \(\|\widetilde\theta-\theta\|\) by the right-hand side, and \cref{eq:growing-fixed-radius-bounds} transfers this estimate to \(\widetilde M-M\).
\end{proof}

\section{Fixed-block deterministic stability}
\label{sec:fixed-block-stability}
\label{sec:fiber_reduction}
For fixed \(W\), unfortunately, the matrix structure inside a cell cannot be replaced by a scalar error. Indeed, previous research on the self-consistent equation tolerate the block-preserving superoperator \citep{Ajanki_2018}, yet they require the superoperator to be ``flat", which is missing in our setting. Instead, for the periodic nature of the self-energies, we proceed by passing them to Bloch fibers and retain the full \(W^d\times W^d\) self-energy.  The argument has three steps.  First, a averaged free eigenvalue counting near the given energy level, or we call free shell bound for the free fibers give uniform upper bounds for the fibered self-energy.  Second, we show the existence of a persistent eigenvalue, or we call the regular isolated band, which supplies a patch of positive Ward mass and hence a \emph{matching} lower bound on the dissipation.  Third, the same patch separates the stability map into a strictly contracting complement and a one-dimensional pole whose square is small after integration by parts. This proves continuum stability; a Combes--Thomas comparison then transfers the result to the finite torus.

Throughout this section \(W<\infty\) is fixed.  We first treat \(0<\lambda\eta\leq1\), which contains the diffusive scale, and then extend the estimates to every fixed bounded interval of \(\lambda\eta\) in Lemma~\ref{lemma:fh-bounded-kinetic-extension}.  This overlaps the kinetic range from Section~\ref{sec:kinetic-initialization}.
\subsection{Regular channels and dissipative bounds}
\label{subsec:fixed_block_MDE_pert}
For notational simplicity, normalize Haar measure on $\mathbb T^d$ by $\int_{\mathbb T^d}\rmd\xi=1$, and write $x=(b,u)$, with $b\in\tilde{\mathbb Z}_n^d$ and $u\in[W^d]$.  The partial Bloch transform in the block variable is again denoted by $\mathcal U$. For a $W^d\times W^d$ matrix $B$, recall the single-block Bloch superoperator
\begin{equation}
        \mathscr S_W[B]
        \coloneqq
        \frac{1}{W^d}\Bigl(B^{\top}+\Tr B\I_{W^d}\Bigr).
\end{equation}
If $B\succeq0$, then
\begin{equation}
\label{equ:block_S_trace_flat}
        \mathscr S_W[B]
        \succeq
        \frac{1}{W^d}\Tr B\I_{W^d} .
\end{equation}
This elementary trace term is the key mechanism that replaces the unavailable scalar flatness condition on the full lattice, where we maintain the flatness on the fiber level. The cost is that, we need to later analyze the spectrum of the Bloch symbol parametrized by the partial momentum space.

Let
\[
        \Omega(p)
        \coloneqq
        2d-2\sum_{j=1}^d\cos p_j,
        \qquad p\in\mathbb T^d,
\]
and recall the definition of the set of critical values of $\Omega$ given by $\Sigma_d^\Psi$. We fix once and for all a compact bulk window
\begin{equation}
\label{equ:fixed_block_bulk_window_I}
        I\Subset (0,4d)\setminus\Sigma_d^\Psi .
\end{equation}
\begin{lemma}
\label{lemma:fh-two-sided-flatness}
For every $B\succeq0$,
\begin{equation}
\label{eq:fh-two-sided-flatness}
        \langle B\rangle_W\I_{W^d}
        \preceq
        \mathscr S_W[B]
        \preceq
        2\langle B\rangle_W\I_{W^d}.
\end{equation}
Moreover,
\begin{equation}
\label{eq:fh-trace-S}
        \left\langle\mathscr S_W[B]\right\rangle_W
        =\left(1+\frac{1}{W^d}\right)\langle B\rangle_W.
\end{equation}
\end{lemma}
\begin{proof}
If \(B\succeq0\), then \(0\preceq B^\top\preceq(\Tr B)\I_{W^d}\). Substitution in the definition of \(\mathscr S_W\) proves \cref{eq:fh-two-sided-flatness}; taking normalized traces proves \cref{eq:fh-trace-S}.
\end{proof}
\begin{lemma}[Unfolding of the fixed-block Bloch symbol]
\label{lemma:fixed-block-unfolding}
There is a measurable family of unitary matrices $U_W(\xi)$ such that diagonalized the Bloch symbol
\begin{equation}
\label{equ:fold_dispersion_fixed_W}
        U_W(\xi)^*\hat{\Psi}_W(\xi)U_W(\xi)
        =
\operatorname{diag}\left\{
        \Omega\left(\frac{\xi+2\pi m}{W}\right)
        \right\}_{m\in\{0,\ldots,W-1\}^d} .
\end{equation}
Consequently, for every bounded Borel function $F$,
\begin{equation}
\label{equ:unfold_trace_fixed_W}
        \int_{\mathbb T^d}\Tr F(\hat\Psi_W(\xi))\,\rmd\xi
        =
        W^d\int_{\mathbb T^d}F(\Omega(p))\,\rmd p .
\end{equation}
Moreover, for $u,u'\in[W^d]$,
\begin{equation}
\label{equ:unfold_matrix_entry_fixed_W}
        \left[\int_{\mathbb T^d}(\hat{\Psi}_W(\xi)-\zeta)^{-1}\,\rmd\xi\right]_{u,u'}
        =
        \int_{\mathbb T^d}
        \frac{e^{ip\cdot(u-u')}}{\Omega(p)-\zeta}\,\rmd p,
        \qquad \Im \zeta\neq0 .
\end{equation}
\end{lemma}

\begin{proof}
The full Fourier transform diagonalizes $\Psi$ with dispersion $\Omega(p)$. The partial Bloch transform only fixes the block momentum $\xi=Wp\mod 2\pi$, so each fiber consists of the folded momenta $(\xi+2\pi m)/W$, $m\in\{0,\ldots,W-1\}^d$.  This gives \cref{equ:fold_dispersion_fixed_W}.  Summing over the folded momenta and using the change of variables $p=(\xi+2\pi m)/W$ gives \cref{equ:unfold_trace_fixed_W}.  Formula \cref{equ:unfold_matrix_entry_fixed_W} is the same unfolding applied to the kernel of the resolvent.
\end{proof}
The next definition describes the regularity of the folded energy surface needed in the fixed-block argument.
 \begin{defn}[Finite-height folded bulk window]
\label{def:fh-bulk-window}
A compact interval $I\Subset\mathbb R$ is called a finite-height folded bulk window for the self-adjoint operator $h$ parametrized by $\xi$ if the following three conditions hold.

\smallskip
\noindent \textnormal{(C1) Free shell bound.} There are $s_0>0$ and $C_{\rm sh}<\infty$ such that, uniformly for $e\in I$ and $0<s\le s_0$,
\begin{equation}
\label{eq:fh-free-shell}
        \mathcal N_e(s)
        \coloneqq
        \int_{\mathbb T^d}
        \Tr\,\mathbf 1_{[-s,s]}\bigl(h(\xi)-e\bigr)\,\rmd\xi
        \le C_{\rm sh}s.
\end{equation}

\smallskip
\noindent \textnormal{(C2) Free matrix Cauchy-transform bound.} There is $M_0<\infty$ such that
\begin{equation}
\label{eq:fh-free-Cauchy}
        \sup_{e\in I}\sup_{0<\nu\le1}
        \left\|
        \int_{\mathbb T^d}
        \bigl(h(\xi)-e-i\nu\bigr)^{-1}\,\rmd\xi
        \right\|_{2\to2}
        \le M_0.
\end{equation}

\smallskip
\noindent \textnormal{(C3) One uniformly isolated regular channel.} There exist constants
\[
        \delta_{\rm iso}>0,
        \qquad
        g_{\rm iso}>0,
        \qquad
        v_{\rm iso}>0,
        \qquad
        j_{\rm iso}>0,
        \qquad
        C_{\rm iso}<\infty,
\]
such that, for every $e\in I$, there exist a band index $j_e$, a bounded open set $V_e\subset\mathbb R^{d-1}$, and an open set $U_e\subset\mathbb T^d$ satisfying $|V_e|_{\mathbb R^{d-1}}\geq v_{\rm iso},$
together with a $C^2$ diffeomorphism onto its image
\[
        \Xi_e:
        (-2\delta_{\rm iso},2\delta_{\rm iso})\times V_e
        \longrightarrow
         U_e \subseteq \mathbb T^d.
\]
The coordinate maps satisfy the uniform nondegeneracy and regularity bounds
\begin{equation}
\label{eq:fh-isolated-chart-bounds}
        \inf_{\substack{|s|<2\delta_{\rm iso}\\ \omega\in V_e}}
        \left|\det D\Xi_e(s,\omega)\right|
        \geq j_{\rm iso},
        \qquad
        \sup_{\substack{|s|<2\delta_{\rm iso}\\ \omega\in V_e}}
        \left(
        \|D\Xi_e(s,\omega)\|
        +
        \|D^2\Xi_e(s,\omega)\|
        \right)
        \leq C_{\rm iso}.
\end{equation}
Moreover, $h(\xi)$ has a simple $C^2$ eigenvalue branch $\varepsilon_{j_e}(\xi)$ on $U_e$ such that
\begin{equation}
\label{eq:fh-isolated-coordinate}
        \varepsilon_{j_e}\bigl(\Xi_e(s,\omega)\bigr)-e=s
\end{equation}
for every
\[
        (s,\omega)
        \in
        (-2\delta_{\rm iso},2\delta_{\rm iso})\times V_e,
\]
and this eigenvalue is uniformly separated from the remaining spectrum:
\begin{equation}
\label{eq:fh-isolated-gap}
        \operatorname{dist}\!\left(
        \varepsilon_{j_e}\bigl(\Xi_e(s,\omega)\bigr),
        \spec h\bigl(\Xi_e(s,\omega)\bigr)
        \setminus
        \left\{
        \varepsilon_{j_e}\bigl(\Xi_e(s,\omega)\bigr)
        \right\}
        \right)
        \geq g_{\rm iso}.
\end{equation}
All constants are uniform in $e\in I$.
\end{defn}
\begin{remark}
Conditions~\textnormal{(C1)}--\textnormal{(C2)} are the usual shell and Cauchy-transform bounds away from the exceptional energies. Condition~\textnormal{(C3)} asks only for one nonresonant regular patch on each energy shell, where we allow such choice to depend on the energy $e$. The folded-band crossings elsewhere are absorbed into the resonant set and instead controlled by Ward's identity.  For this special choice of the chart \(\Xi_e(s,\omega)\), we expect to be $s$ transverse to the energy surface and $\omega$ parametrizes its $(d-1)$ tangential directions.  The Jacobian lower bound then prevents degeneration of this change of variables.
\end{remark}
The following example verifies condition~\textnormal{(C3)} for $\widehat\Psi$ when $W=2$.
\begin{example}
\label{ex:W2-isolated-channel}
Let $W=2$, and consider the free lattice dispersion
\[
        \Omega(p)
        =
        2d-2\sum_{\ell=1}^d\cos p_\ell,
        \qquad p\in\mathbb T^d .
\]
After folding to the block Brillouin zone, the Bloch symbol $\widehat\Psi(\xi)$ has the eigenvalue branches
\[
        \varepsilon_m(p)
        =
        2d-2\sum_{\ell=1}^d (-1)^{m_\ell}\cos p_\ell,
        \qquad
        m\in\{0,1\}^d,
\]
where $p=\xi/2$.  Equivalently, shifting one folded momentum by $\pi$ flips the sign of the corresponding cosine.

Fix a compact interval \(I\Subset(0,4d)\setminus\{2d\}\).  We verify condition~\textnormal{(C3)} on \(I\).  For \(e\in I\), choose \(r_e\in(0,\pi)\) from \(2d-2d\cos r_e=e\), equivalently \(\cos r_e=1-e/(2d)\), and set $p_e=(r_e,\ldots,r_e)$. For the branch $m=0$, we have $\varepsilon_0(p_e)=e.$ Moreover, for every $m\neq0$,
\[
\begin{aligned}
        \varepsilon_m(p_e)-\varepsilon_0(p_e)
        &=
        4|m|\cos r_e,
\end{aligned}
\]
where $|m|=m_1+\cdots+m_d$.  Hence
\[
        \operatorname{dist}
        \left(
        \varepsilon_0(p_e),
        \{\varepsilon_m(p_e):m\neq0\}
        \right)
        \ge
        4\inf_{e\in I}
        \left|1-\frac{e}{2d}\right|
        >0.
\]
Thus the selected branch is uniformly isolated at $p_e$.  The same branch is also uniformly regular.  Indeed,
\[
        \nabla_p\varepsilon_0(p)
        =
        2(\sin p_1,\ldots,\sin p_d),
\]
and therefore
\[
        |\nabla_p\varepsilon_0(p_e)|
        =
        2\sqrt d\,\sin r_e .
\]
Since $I\Subset(0,4d)$, the numbers $r_e$ stay uniformly away from $0$ and $\pi$, and hence $\inf_{e\in I}|\nabla_p\varepsilon_0(p_e)|>0.$
By the quantitative implicit function theorem \citep{KrantzParks_IFT}, there exist constants
\[
        \delta_{\rm iso}>0,
        \qquad
        g_{\rm iso}>0,
        \qquad
        v_{\rm iso}>0,
        \qquad
        j_{\rm iso}>0,
        \qquad
        C_{\rm iso}<\infty,
\]
depending only on $I,d$, such that for every $e\in I$ there is a $C^2$ coordinate chart
\[
        \tilde\Xi_e:
        (-2\delta_{\rm iso},2\delta_{\rm iso})\times V_e
        \longrightarrow \mathbb T^d_p
\]
which is a diffeomorphism onto its image, with
\[
        |V_e|_{\mathbb R^{d-1}}\ge v_{\rm iso},
        \qquad
        \inf |\det D\tilde\Xi_e|\ge j_{\rm iso},
        \qquad
        \|\tilde\Xi_e\|_{C^2}\le C_{\rm iso},
\]
and such that
\[
        \varepsilon_0(\tilde\Xi_e(s,\omega))-e=s .
\]
Furthermore, by shrinking $\delta_{\rm iso}$ if necessary, the isolation persists throughout the coordinate tube:
\[
        \operatorname{dist}
        \left(
        \varepsilon_0(\tilde\Xi_e(s,\omega)),
        \{\varepsilon_m(\tilde\Xi_e(s,\omega)):m\neq0\}
        \right)
        \ge g_{\rm iso}.
\]
Returning to the block momentum variable $\xi=2p$ gives the desired coordinate map $\Xi_e=2\tilde\Xi_e$ for the folded symbol $\widehat\Psi(\xi)$, with modified constants.  Thus condition \textnormal{(C3)} holds uniformly for every $I\Subset(0,4d)\setminus\{2d\}$.
\end{example}
Proposition~\ref{prop:BA-regular-energies-folded-bulk} proves the corresponding statement for every $I\Subset(0,4d)\setminus\Sigma_d^\Psi$.  We first verify conditions~\textnormal{(C1)}--\textnormal{(C2)} using the coarea formula.
\begin{lemma}
\label{lemma:BA-folded-shell-Cauchy}
Let \(d\geq2\), let \(W<\infty\) be fixed, and suppose that $I\Subset(0,4d)\setminus\Sigma_d^\Psi.$ Then \(h(\xi)=\widehat\Psi_W(\xi)\) satisfies conditions \textnormal{\rm(C1)} and \textnormal{\rm(C2)} of Definition~\ref{def:fh-bulk-window}.  The corresponding constants may depend on \(W,d\), and \(I\).
\end{lemma}

\begin{proof}
Recall that
\[
        \Omega(p)
        =
        2d-2\sum_{j=1}^d\cos p_j,
        \qquad
        p\in\mathbb T^d.
\]
The critical points of \(\Omega\) satisfy \(p_j\in\{0,\pi\}\) for every \(j\), and hence its critical values are $\{0,4,\ldots,4d\}.$ Choose an open interval \(J\) such that $I\Subset J\Subset(0,4d)\setminus\{0,4,\ldots,4d\}.$ Then
\begin{equation}
\label{eq:BA-regular-shell-gradient}
        \inf_{p\in\Omega^{-1}(\overline J)}
        |\nabla\Omega(p)|>0.
\end{equation}
Identify the internal coordinates with $u\in\{0,\ldots,W-1\}^d.$ For \(\alpha\in\mathbb Z^d\) and \(\sigma\in J\), define the weighted free density of states by
\begin{equation}
\label{eq:BA-weighted-free-DOS}
        \rho_\alpha(\sigma)
        \coloneqq
        \frac1{(2\pi)^d}
        \int_{\Omega^{-1}(\sigma)}
        \frac{e^{\mathrm ip\cdot\alpha}}
        {|\nabla\Omega(p)|}
        \,\rmd\mathcal H^{d-1}(p).
\end{equation}
Here \(\mathcal H^{d-1}\) denotes the \((d-1)\)-dimensional Hausdorff measure on \(\mathbb T^d\). Since $\Omega:\Omega^{-1}(J)\longrightarrow J$ is a proper smooth submersion, integration along its level sets gives $\rho_\alpha\in C^\infty(J).$ For completeness, one may use on \(\Omega^{-1}(J)\) the smooth vector field
\[
        V(p)
        \coloneqq
        \frac{\nabla\Omega(p)}{|\nabla\Omega(p)|^2},
        \qquad
        V\Omega=1.
\]
Its flow identifies nearby level sets, and after transporting the surface measure by this flow, the integrand and its Jacobian depend smoothly on the energy variable.  In particular, for every \(J'\Subset J\),
\begin{equation}
\label{eq:BA-weighted-DOS-C1}
        \sup_{\sigma\in J'}
        \left(
        |\rho_\alpha(\sigma)|+|\rho_\alpha'(\sigma)|
        \right)<\infty.
\end{equation}

We first prove \textnormal{\rm(C1)}.  Choose \(s_0>0\) and \(J'\Subset J\) such that
\[
        [e-s_0,e+s_0]\subset J'
        \qquad
        \text{for every }e\in I.
\]
By \cref{equ:unfold_trace_fixed_W} and the coarea formula,
\begin{align*}
        \mathcal N_e(s)
        &=
        \int_{\mathbb T^d}
        \Tr\,
        \mathbf 1_{[-s,s]}
        \bigl(\widehat\Psi_W(\xi)-e\bigr)
        \,\rmd\xi
        \\
        &=
        W^d
        \int_{\mathbb T^d}
        \mathbf 1_{\{|\Omega(p)-e|\leq s\}}
        \,\rmd p
        \\
        &=
        W^d
        \int_{e-s}^{e+s}\rho_0(\sigma)\,\rmd \sigma
        \\
        &\leq
        2W^d\|\rho_0\|_{L^\infty(J')}s
\end{align*}
for every \(e\in I\) and \(0<s\leq s_0\).  This proves \textnormal{\rm(C1)}.

We next prove \textnormal{\rm(C2)}.  Choose \(\chi\in C_c^\infty(J)\) such that \(\chi=1\) on an open neighborhood of \(I\).  By \cref{equ:unfold_matrix_entry_fixed_W}, the \((u,u')\)-entry of the matrix in \textnormal{\rm(C2)} is
\begin{equation}
\label{eq:BA-folded-Cauchy-entry}
\begin{aligned}
        A_{u,u'}(e,\nu)
        &\coloneqq
        \left[
        \int_{\mathbb T^d}
        \bigl(
        \widehat\Psi_W(\xi)-e-\mathrm i\nu
        \bigr)^{-1}
        \,\rmd\xi
        \right]_{u,u'}
        \\
        &=
        \int_{\mathbb T^d}
        \frac{
        e^{\mathrm i p\cdot(u-u')}
        }{
        \Omega(p)-e-\mathrm i\nu
        }
        \,\rmd p.
\end{aligned}
\end{equation}

We first consider the part away from \(I\).  Since \(\chi=1\) on a neighborhood of \(I\), there is \(c_{I,J}>0\) such that
\[
        1-\chi(\Omega(p))\neq0
        \quad\Longrightarrow\quad
        |\Omega(p)-e|\geq c_{I,J},
        \qquad e\in I.
\]
Consequently,
\begin{equation}
\label{eq:BA-Cauchy-far-energy}
\left|
        \int_{\mathbb T^d}
        \frac{
        (1-\chi(\Omega(p)))
        e^{\mathrm i p\cdot(u-u')}
        }{
        \Omega(p)-e-\mathrm i\nu
        }
        \,\rmd p
\right|
\leq c_{I,J}^{-1}.
\end{equation}
For the remaining part, set
\[
        \alpha=u-u',
        \qquad
        F_\alpha(\sigma)\coloneqq\chi(\sigma)\rho_\alpha(\sigma),
\]
and extend \(F_\alpha\) by zero outside \(J\).  Then \(F_\alpha\in C_c^1(\mathbb R)\), and the coarea formula gives
\begin{equation}
\label{eq:BA-Cauchy-near-energy}
        \int_{\mathbb T^d}
        \frac{
        \chi(\Omega(p))e^{\mathrm i p\cdot\alpha}
        }{
        \Omega(p)-e-\mathrm i\nu
        }
        \,\rmd p
        =
        \int_{\mathbb R}
        \frac{F_\alpha(\sigma)}
        {\sigma - e - \mathrm i\nu}
        \,\rmd \sigma.
\end{equation}
Choose \(a<b\) such that $\operatorname{supp} F_\alpha\cup I \Subset(a,b)\Subset J.$ The same \(a,b\) may be used for all the finitely many differences \(\alpha=u-u'\).  Subtracting \(F_\alpha(e)\), we obtain
\begin{align*}
 \left|
 \int_a^b
 \frac{F_\alpha(\sigma)}
 {\sigma-e-\mathrm i\nu}
 \,\rmd \sigma
 \right|
 &\leq
 \int_a^b
 \frac{|F_\alpha(\sigma)-F_\alpha(e)|}
 {|\sigma-e-\mathrm i\nu|}
 \,\rmd \sigma +
 |F_\alpha(e)|
 \left|
 \int_a^b
 \frac{\rmd \sigma}
 {\sigma-e-\mathrm i\nu}
 \right|.
\end{align*}
The first term is bounded by $(b-a)\|F_\alpha'\|_\infty$, because
\[
        |F_\alpha(\sigma)-F_\alpha(e)|
        \leq
        \|F_\alpha'\|_\infty|\sigma-e|
        \leq
        \|F_\alpha'\|_\infty
        |\sigma-e-\mathrm i\nu|.
\]
For the second term,
\[
        \sup_{e\in I}\sup_{0<\nu\leq1}
        \left|
        \int_a^b
        \frac{\rmd \sigma}
        {\sigma-e-\mathrm i\nu}
        \right|
        <\infty.
\]
Indeed, the real part of this integral is $\frac12\log\frac{(b-e)^2+\nu^2}{(a-e)^2+\nu^2},$ while its imaginary part is uniformly bounded by \(\pi\).  The real part is uniformly bounded because \(I\Subset(a,b)\).

We conclude that
\[
        \sup_{e\in I}\sup_{0<\nu\leq1}
        |A_{u,u'}(e,\nu)|
        \leq C_{W,I}.
\]
There are only finitely many pairs \(u,u'\).  Since the matrix size is \(W^d\),
\[
\begin{aligned}
        \left\|
        \int_{\mathbb T^d}
        \bigl(
        \widehat\Psi_W(\xi)-e-\mathrm i\nu
        \bigr)^{-1}
        \,\rmd\xi
        \right\|_{2\to2}
        &\leq
        W^d
        \max_{u,u'}|A_{u,u'}(e,\nu)|
        \\
        &\leq C_{W,I}.
\end{aligned}
\]
This proves \textnormal{\rm(C2)}.
\end{proof}

We first record an elementary algebraic lemma.
\begin{lemma}
\label{lemma:BA-separable-hyperplane-identity}
Let \(d\geq2\), let \(c\neq0\), and let
\[
        a_j,b_j\in\mathbb R,
        \qquad
        \sigma_j\in\{\pm1\},
        \qquad
        j=1,\ldots,d.
\]
Suppose that
\begin{equation}
\label{eq:BA-separable-hyperplane-identity}
        \sum_{j=1}^d
        \left(
        a_jx_j+
        b_j\sigma_j\sqrt{1-x_j^2}
        \right)
        =0
\end{equation}
on a nonempty relatively open subset of
\[
        \left\{
        x\in(-1,1)^d:
        \sum_{j=1}^dx_j=c
        \right\}.
\]
Then $a_j=b_j=0$ for every $j$.
\end{lemma}

\begin{proof}
For each \(j\), set $f_j(x)\coloneqq a_jx+b_j\sigma_j\sqrt{1-x^2}.$ We first suppose that \(d\geq3\).  Locally use \(x_1,\ldots,x_{d-1}\) as independent variables and write $x_d=c-\sum_{j=1}^{d-1}x_j.$ Then \cref{eq:BA-separable-hyperplane-identity} becomes $\sum_{j=1}^{d-1}f_j(x_j)+f_d(x_d)=0$. Differentiating with respect to \(x_i\), \(i<d\), gives
\begin{equation}
\label{eq:BA-separable-first-derivative}
        f_i'(x_i)-f_d'(x_d)=0.
\end{equation}
Choose \(k\neq i,d\) and differentiate this identity with respect to \(x_k\).  Since \(\partial_{x_k}x_d=-1\), we obtain $f_d''(x_d)=0$ on a nonempty interval.  But $f_d''(x) = -(b_d\sigma_d){(1-x^2)^{-3/2}},$ so \(b_d=0\).  Consequently \(f_d'=a_d\), and \cref{eq:BA-separable-first-derivative} gives $f_i'(x_i)=a_d$ for every $i<d$. Differentiating once more shows \(b_i=0\), and then \(a_i=a_d\). Thus all \(b_j\) vanish and all \(a_j\) have a common value \(a\). Equation \cref{eq:BA-separable-hyperplane-identity} now reduces to $a\sum_{j=1}^dx_j=ac=0$. Since \(c\neq0\), we have \(a=0\).

It remains to consider \(d=2\).  Write \(x_2=c-x_1\).  Then $f_1(x)+f_2(c-x)=0$ on a nonempty interval.  Differentiating twice gives
\[
        f_1''(x)+f_2''(c-x)=0,
\]
or
\begin{equation}
\label{eq:BA-separable-d2-second}
        \frac{b_1\sigma_1}{(1-x^2)^{3/2}}
        +
        \frac{b_2\sigma_2}
        {(1-(c-x)^2)^{3/2}}
        =0.
\end{equation}
If one of \(b_1,b_2\) vanishes, then \cref{eq:BA-separable-d2-second} shows that both vanish.  Otherwise, the quotient $(1-x^2)/(1-(c-x)^2)$ would be constant on an interval.  Comparing the coefficients of the two quadratic polynomials would then force \(c=0\), contrary to the hypothesis.  Hence \(b_1=b_2=0\).

The original identity is now $a_1x+a_2(c-x)=0$ on an interval.  Therefore \(a_1=a_2\), and then \(a_2c=0\).  Since \(c\neq0\), both coefficients vanish.
\end{proof}

\begin{lemma}[Existence of a nonresonant folded channel]
\label{lemma:BA-fold-nonresonance}
Let \(d\geq2\), let \(W<\infty\) be fixed, and let $e\in(0,4d)\setminus\Sigma_d^\Psi$. Then there exists \(p_e\in\mathbb T^d\) such that
\begin{equation}
\label{eq:BA-fold-nonresonant-point}
\begin{aligned}
        \Omega(p_e)=e,\quad
        \nabla\Omega(p_e)\neq0,\quad
        \Omega(p_e+q)\neq\Omega(p_e)
        \quad
        \text{for every }q\in\mathcal Q_W\setminus\{0\},
\end{aligned}
\end{equation}
where
\[
        \mathcal Q_W
        \coloneqq
        \left\{
        \frac{2\pi m}{W}:
        m\in\{0,\ldots,W-1\}^d
        \right\}.
\]
\end{lemma}

\begin{proof}
Let $F_e\coloneqq\{p\in\mathbb T^d:\Omega(p)=e\}.$ The set \(F_e\) is nonempty because the diagonal path
\[
        p=(r,\ldots,r),
        \qquad
        \Omega(p)=2d(1-\cos r),
        \qquad
        0<r<\pi,
\]
covers all energies in \((0,4d)\).

Moreover,
\[
        \nabla\Omega(p)
        =
        2(\sin p_1,\ldots,\sin p_d).
\]
The critical values of \(\Omega\) are \(\{0,4,\ldots,4d\}\).  Since \(e\notin\Sigma_d^\Psi\), it follows that $\nabla\Omega(p)\neq0$ for every $\in F_e$. Thus \(F_e\) is a compact smooth hypersurface. For \(q\in\mathcal Q_W\setminus\{0\}\), define the collision set $Z_q\coloneqq\left\{p\in F_e:\Omega(p+q)=\Omega(p)\right\}.$ It is closed in \(F_e\), and its relative interior is empty.  To see this, argue by contradiction.

Suppose otherwise.  Then there is a nonempty relatively open set \(\mathcal O\subset F_e\) on which
\begin{equation}
\label{eq:BA-fold-collision-open-patch}
        \Omega(p+q)-\Omega(p)=0.
\end{equation}
After shrinking \(\mathcal O\), we may assume that
\begin{equation}
\label{eq:BA-nonzero-coordinate-sines}
        \sin p_j\neq0
        \qquad
        \text{for every }j
        \text{ and every }p\in\mathcal O.
\end{equation}
Indeed, if \(\sin p_j=0\) on a relatively open subset of \(F_e\), then \(p_j\) is locally constant there, so every tangent vector has vanishing \(j\)-th component.  But $T_pF_e = \{v:\nabla\Omega(p)\cdot v=0\},$ and $\partial_j\Omega(p)=2\sin p_j=0$, so $e_j\in T_pF_e$, a contradiction.  On the patch \(\mathcal O\), set
\[
        x_j=\cos p_j,
        \qquad
        \sin p_j=\sigma_j\sqrt{1-x_j^2},
        \qquad
        \sigma_j\in\{\pm1\}.
\]
The signs are fixed after shrinking the connected patch.  By \cref{eq:BA-nonzero-coordinate-sines}, the map $p\longmapsto x=(\cos p_1,\ldots,\cos p_d)$ is a local diffeomorphism.  Since \(\Omega(p)=e\), its image in the energy shell satisfies
\begin{equation}
\label{eq:BA-cosine-affine-shell}
        \sum_{j=1}^dx_j = d-\frac e2
        \eqqcolon c_e.
\end{equation}
The definition of \(\Sigma_d^\Psi\) includes \(2d\), so \(c_e\neq0\). Using $\cos(p_j+q_j)=\cos p_j\cos q_j-\sin p_j\sin q_j,$ the collision identity \cref{eq:BA-fold-collision-open-patch} becomes
\[
        \sum_{j=1}^d
        \left[
        (1-\cos q_j)x_j
        +
        \sin q_j\,\sigma_j\sqrt{1-x_j^2}
        \right]
        =0
\]
on a nonempty relatively open subset of the hyperplane \cref{eq:BA-cosine-affine-shell}.  Applying Lemma~\ref{lemma:BA-separable-hyperplane-identity} with $a_j=1-\cos q_j,$ $b_j=\sin q_j$ gives $1-\cos q_j=0,$ $\sin q_j=0$ for every $j$. Therefore \(q_j=0\pmod{2\pi}\) for every \(j\), contradicting \(q\neq0\) in \(\mathbb T^d\).  Thus \(Z_q\) has empty relative interior.

It follows that \(F_e\setminus Z_q\) is open and dense in \(F_e\). Since \(\mathcal Q_W\setminus\{0\}\) is finite, $\bigcap_{q\in\mathcal Q_W\setminus\{0\}}(F_e\setminus Z_q)$ is open and dense in \(F_e\), and in particular is nonempty.  Choose \(p_e\) in this intersection.  Then $\Omega(p_e+q)\neq\Omega(p_e)$ for every $q\neq 0$ and \(\nabla\Omega(p_e)\neq0\) by regularity of the energy shell. This proves \cref{eq:BA-fold-nonresonant-point}.
\end{proof}

\begin{prop}
\label{prop:BA-regular-energies-folded-bulk}
Let \(d\geq2\), let \(W<\infty\) be fixed, and let \(I \Subset \mathbb R\) be a compact interval satisfying $I\Subset(0,4d)\setminus\Sigma_d^\Psi$. Then \(I\) is a finite-height folded bulk window for $h(\xi)=\widehat\Psi_W(\xi)$ in the sense of Definition~\ref{def:fh-bulk-window}.  All constants may depend on \(W,d\), and \(I\).
\end{prop}

\begin{proof}
Conditions \textnormal{\rm(C1)} and \textnormal{\rm(C2)} follow from Lemma~\ref{lemma:BA-folded-shell-Cauchy}.  It remains to prove \textnormal{\rm(C3)}.

Fix \(e_0\in I\).  By Lemma~\ref{lemma:BA-fold-nonresonance}, there is \(p_0\in\mathbb T^d\) such that $\Omega(p_0)=e_0, \nabla\Omega(p_0)\neq0,$ and $\Omega(p_0+q)\neq\Omega(p_0)$ for every $q\in\mathcal Q_W\setminus\{0\}.$ Set $\xi_0=Wp_0\pmod{2\pi}.$ Choose a local lift of \(\xi\) near \(\xi_0\).  For a suitable \(m_0\in\{0,\ldots,W-1\}^d\), the selected folded branch is locally $ \varepsilon_0(\xi) = \Omega\left( (\xi+2\pi m_0)/W \right),$ and $\varepsilon_0(\xi_0)=e_0.$ The remaining folded branches differ from this one at \(\xi_0\) by the shifts \(q\in\mathcal Q_W\setminus\{0\}\).  Therefore
\begin{equation}
\label{eq:BA-pointwise-folded-gap}
        \gamma_0
        \coloneqq
        \min_{q\in\mathcal Q_W\setminus\{0\}}
        |\Omega(p_0+q)-\Omega(p_0)|
        >0.
\end{equation}
Thus \(\varepsilon_0\) is a simple eigenvalue of \(\widehat\Psi_W(\xi_0)\).  Moreover,
\[
        \nabla_\xi\varepsilon_0(\xi_0)
        =
        \frac1W\nabla\Omega(p_0)
        \neq0.
\]

By continuity, after shrinking to a sufficiently small coordinate neighborhood of \(\xi_0\), the branch \(\varepsilon_0\) remains simple and \(C^\infty\), and its distance from the remaining spectrum is at least \(\gamma_0/2\).  The submersion theorem then gives \(r_0>0\), a bounded open set \(V_0\Subset\mathbb R^{d-1}\), and a \(C^2\) diffeomorphism onto its image
\[
        \Phi_0:
        (-4r_0,4r_0)\times V_0
        \longrightarrow U_0\subset\mathbb T^d
\]
such that
\begin{equation}
\label{eq:BA-master-folded-chart}
        \varepsilon_0(\Phi_0(t,\omega))
        =
        e_0+t.
\end{equation}
Shrinking \(r_0\) and \(V_0\) if necessary, there are constants $g_0>0, v_0>0, j_0>0, C_0<\infty$ such that, for \(|t|\leq3r_0\) and \(\omega\in V_0\),
\begin{align}
        |V_0|&\geq v_0,
        \label{eq:BA-master-chart-volume}
        \\
        |\det D\Phi_0(t,\omega)|
        &\geq j_0,
        \label{eq:BA-master-chart-Jacobian}
        \\
        \|D\Phi_0(t,\omega)\|
        +
        \|D^2\Phi_0(t,\omega)\|
        &\leq C_0,
        \label{eq:BA-master-chart-C2}
        \\
        \operatorname{dist}\left(
        \varepsilon_0(\Phi_0(t,\omega)),
        \spec\widehat\Psi_W(\Phi_0(t,\omega))
        \setminus
        \{\varepsilon_0(\Phi_0(t,\omega))\}
        \right)
        &\geq g_0.
        \label{eq:BA-master-chart-gap}
\end{align}

Let $J_0=(e_0-r_0,e_0+r_0).$ For every \(e\in J_0\), define
\begin{equation}
\label{eq:BA-shifted-folded-chart}
        \Xi_e(s,\omega)
        \coloneqq
        \Phi_0(e-e_0+s,\omega).
\end{equation}
If \(|s|<2r_0\), then $|e-e_0+s|<3r_0,$ so the preceding bounds apply.  Moreover,
\[
\begin{aligned}
        \varepsilon_0(\Xi_e(s,\omega))-e
        &=
        \varepsilon_0
        \bigl(
        \Phi_0(e-e_0+s,\omega)
        \bigr)-e
        \\
        &=
        e_0+(e-e_0+s)-e
        \\
        &=s.
\end{aligned}
\]

Repeating this construction for every \(e_0\in I\) produces an open cover of \(I\) by intervals \(J_0\).  By compactness, choose a finite subcover $I\subseteq\bigcup_{\alpha=1}^M J_\alpha,$ denote
\[
        \delta_{\rm iso}
        \coloneqq
        \min_{1\leq\alpha\leq M}r_\alpha,
        \qquad
        g_{\rm iso}
        \coloneqq
        \min_{1\leq\alpha\leq M}g_\alpha,
\]
\[
        v_{\rm iso}
        \coloneqq
        \min_{1\leq\alpha\leq M}v_\alpha,
        \qquad
        j_{\rm iso}
        \coloneqq
        \min_{1\leq\alpha\leq M}j_\alpha,
\]
and
\[
        C_{\rm iso}
        \coloneqq
        \max_{1\leq\alpha\leq M}C_\alpha.
\]
These constants are strictly positive or finite as appropriate.

For each \(e\in I\), choose an \(\alpha\) such that \(e\in J_\alpha\), take \(V_e=V_\alpha\), and restrict \cref{eq:BA-shifted-folded-chart} to $(-2\delta_{\rm iso},2\delta_{\rm iso})\times V_\alpha.$
Let \(U_e\) be its image.  Since the selected eigenvalue remains simple on this connected tube, it has a fixed local ordered index \(j_e\).  The resulting branch and chart satisfy \cref{eq:fh-isolated-chart-bounds,eq:fh-isolated-coordinate,eq:fh-isolated-gap} with the uniform constants just defined.  This proves \textnormal{\rm(C3)} and completes the proof.
\end{proof}

\begin{remark}
\label{rmk:BA-center-fold-obstruction}
The exclusion of \(2d\) in the definition of \(\Sigma_d^\Psi\) is natural for the folded problem.  If \(W\) is even, then $q_\star=(\pi,\ldots,\pi)\in\mathcal Q_W$ and $\Omega(p+q_\star)=4d-\Omega(p)$. Consequently, on the energy shell \(\Omega(p)=2d\),
\[
        \Omega(p+q_\star)=\Omega(p)=2d
\]
for every \(p\).  Hence every folded channel at this energy has a degenerate partner and condition \textnormal{\rm(C3)}, as presently formulated using a simple isolated band would fail. Instead, away from the critical energy, the natural next question to ask is that, if we can find at least a small open set of the torus with positive mass, such that we can always find a isolated eigenvalue for the Bloch symbol. If this holds true, this regular isolated band would protect the stability of the spectrum projection, and its positive volume permits the change of variable needed in proving the matching lower bound.
\end{remark}

\begin{lemma}[Non-commuting dissipation]
\label{lemma:fh-noncommuting-dissipation}
Let $B=B^*$ and let $\Gamma=\Gamma^*\succ0$.  Then
\begin{equation}
\label{eq:fh-imag-congruence}
        \Im(B-i\Gamma)^{-1}
        =\bigl(B\Gamma^{-1}B+\Gamma\bigr)^{-1}.
\end{equation}
If
\[
        \alpha\tau\I_{W^d}\preceq\Gamma\preceq\beta\tau\I_{W^d},
        \qquad 0<\alpha\le\beta,
\]
then
\begin{equation}
\label{eq:fh-Poisson-Loewner}
        \alpha\tau
        \bigl(B^2+\alpha\beta\tau^2\I_{W^d}\bigr)^{-1}
        \preceq
        \Im(B-i\Gamma)^{-1}
        \preceq
        \beta\tau
        \bigl(B^2+\alpha\beta\tau^2\I_{W^d}\bigr)^{-1}.
\end{equation}
\end{lemma}

\begin{proof}
Factor
\[
        B-i\Gamma
        =\Gamma^{1/2}
        \bigl(\Gamma^{-1/2}B\Gamma^{-1/2}-i\I_{W^d}\bigr)
        \Gamma^{1/2}.
\]
Taking the imaginary part of the inverse gives \cref{eq:fh-imag-congruence}.  Moreover,
\[
        \frac1{\beta\tau}B^2+\alpha\tau\I_{W^d}
        \preceq
        B\Gamma^{-1}B+\Gamma
        \preceq
        \frac1{\alpha\tau}B^2+\beta\tau\I_{W^d}.
\]
Inversion reverses the Loewner order and yields \cref{eq:fh-Poisson-Loewner}.
\end{proof}
\begin{lemma}
\label{lemma:fh-Poisson-shell}
Fix $e\in\mathbb R$ and assume \cref{eq:fh-free-shell}.  Let $X\in\mathbb C^{W^d\times W^d}$ be Hermitian, and let $\mu_1(\xi),\ldots,\mu_{W^d}(\xi)$ denote the eigenvalues of $h(\xi)-e-\tau X$.  Then there is $C<\infty$, depending only on $C_{\rm sh},s_0,W^d$, such that for every $w>0$ satisfying $w+\tau\|X\|_{2\to2}\leq s_0/2$,
\begin{equation}
\label{eq:fh-Poisson-shell-bound}
        \int_{\mathbb T^d}
        \sum_{j=1}^{W^d}\frac{w}{\mu_j(\xi)^2+w^2}\,\rmd\xi
        \le
        C\left(1+\frac{\tau}{w}\|X\|_{2\to 2}\right).
\end{equation}
\end{lemma}

\begin{proof}
Let $\varepsilon_1(\xi),\ldots,\varepsilon_{W^d}(\xi)$ denote the eigenvalues of $h(\xi)$. By Weyl's inequality, for every $s\geq0$,
\[
        \card\bigl\{
        j\in\{1,\ldots,W^d\}:|\mu_j(\xi)|\leq s
        \bigr\}
        \leq
        \card\bigl\{
        j\in\{1,\ldots,W^d\}:
        |\varepsilon_j(\xi)-e|\leq s+\tau \|X\|_{2\to 2}
        \bigr\}.
\]
After integration and using~\cref{eq:fh-free-shell}, we obtain
\begin{equation}
\label{eq:fh-pert-shell-count}
        \int_{\mathbb T^d}
        \card\bigl\{j:|\mu_j(\xi)|\leq s\bigr\}
        \,\mathrm d\xi
        \leq
        C_{\rm sh}(s+\tau\|X\|_{2\to 2}),
\end{equation}
provided that $s+\tau \|X\|_{2\to 2}\leq s_0$. Since $w>0$, the Poisson kernel admits the dyadic majorization
\begin{equation}
\label{eq:fh-dyadic-Poisson-majorization}
        \frac{w}{x^2+w^2}
        \leq
        \frac{1}{w}
        \mathbf{1}_{\{|x|\leq w\}}
        +
        \frac{1}{w}
        \sum_{\ell=0}^{\infty}
        2^{-2\ell}
        \mathbf{1}_{\{
        2^\ell w<|x|\leq 2^{\ell+1}w
        \}} .
\end{equation}
Indeed, on the annulus $2^\ell w<|x|\leq2^{\ell+1}w$, one has
\[
        \frac{w}{x^2+w^2}
        \leq
        \frac{w}{x^2}
        \leq
        \frac{2^{-2\ell}}{w}.
\]

Define $N(s)
        \coloneqq
        \int_{\mathbb T^d}
        \card\bigl\{j:|\mu_j(\xi)|\leq s\bigr\}
        \,\mathrm d\xi,$
apply~\cref{eq:fh-dyadic-Poisson-majorization} to each eigenvalue $\mu_j(\xi)$ gives
\begin{align}
\label{eq:fh-dyadic-Poisson-integral}
        \int_{\mathbb T^d}
        \sum_{j=1}^{W^d}
        \frac{w}{\mu_j(\xi)^2+w^2}
        \,\mathrm d\xi
        &\leq
        \frac{N(w)}{w}
        +
        \frac{1}{w}
        \sum_{\ell=0}^{\infty}
        2^{-2\ell}N(2^{\ell+1}w).
\end{align}
Let $L$ be the largest nonnegative integer such that $2^{L+1}w+\tau\|X\|_{2\to2}\leq s_0$.  For $0\leq\ell\leq L$, \cref{eq:fh-pert-shell-count} gives
\begin{align}
        \frac{N(w)}{w}
        +
        \frac{1}{w}
        \sum_{\ell=0}^{L}
        2^{-2\ell}N(2^{\ell+1}w)
        \nonumber&\leq
        C_{\rm sh}\left(1+\frac{\tau}{w}\|X\|_{2\to2}\right)
        +
        \frac{C_{\rm sh}}{w}
        \sum_{\ell=0}^{L}
        2^{-2\ell}
        \bigl(2^{\ell+1}w+\tau \|X\|_{2\to 2}\bigr)
        \nonumber\\
        &\leq
        C C_{\rm sh}
        \left(1+\frac{\tau}{w}\|X\|_{2\to 2}\right),
\label{eq:fh-small-dyadic-scales}
\end{align}
where we used the fact that
\[
        \sum_{\ell=0}^{\infty}2^{-\ell}<\infty,
        \qquad
        \sum_{\ell=0}^{\infty}2^{-2\ell}<\infty.
\]
For the remaining scales, use the elementary bound $N(s)\leq W^d$:
\begin{align}
        \frac{1}{w}
        \sum_{\ell=L+1}^{\infty}
        2^{-2\ell}N(2^{\ell+1}w)
        &\leq
        \frac{W^d}{w}
        \sum_{\ell=L+1}^{\infty}2^{-2\ell}
        \nonumber\\
        &\leq
        \frac{CW^d}{w}\,2^{-2L}
        \leq
        \frac{C'W^d\,w}{(s_0-\tau \|X\|_{2\to 2})^2}.
\label{eq:fh-large-dyadic-scales}
\end{align}
Since $\tau\|X\|_{2\to2}\leq s_0/2$,
\[
        \frac{C'W^d\,w}{(s_0-\tau \|X\|_{2\to 2})^2}
        \leq
        \frac{4C'W^d\,w}{s_0^2},
\]
and this term is absorbed by enlarging the constant.  Combining \cref{eq:fh-dyadic-Poisson-integral}, \cref{eq:fh-small-dyadic-scales}, and \cref{eq:fh-large-dyadic-scales}, we have proved~\cref{eq:fh-Poisson-shell-bound}.
\end{proof}
The following companion estimate controls the contribution of bands away from the rescaled energy level $e$.
\begin{lemma}
\label{lemma:fh-off-shell}
Assume \cref{eq:fh-free-shell} and define
\[
        r_e(\xi)
        \coloneqq
        \operatorname{dist}\bigl(e,\operatorname{Spec}h(\xi)\bigr).
\]
Then there is $C<\infty$, independent of $\delta$, such that for $0<\delta\le s_0/2$,
\begin{equation}
\label{eq:fh-off-shell-integral}
        \int_{\{r_e>\delta\}}
        \frac{\delta}{r_e(\xi)^2}\,\rmd\xi
        \le C.
\end{equation}
\end{lemma}

\begin{proof}
Decompose \(\{\delta<r_e\leq s_0\}\) into the shells \(2^\ell\delta<r_e\leq2^{\ell+1}\delta\).  On the \(\ell\)-th shell the integrand is at most \(2^{-2\ell}\delta^{-1}\), whereas \cref{eq:fh-free-shell} bounds its measure by \(C_{\rm sh}2^{\ell+1}\delta\).  Hence these shells contribute at most \(2C_{\rm sh}\sum_{\ell\geq0}2^{-\ell}\).  On \(\{r_e>s_0\}\), the integrand is at most \(\delta s_0^{-2}\leq(2s_0)^{-1}\).  Combining the two bounds proves \cref{eq:fh-off-shell-integral}.
\end{proof}
\subsection{The fibered matrix Dyson equation}
We now rescale the fixed-block MDE to quantities of order one.  Set
\[
        z=E+i\eta,
        \qquad
        e\coloneqq E/\lambda,
        \qquad
        \kappa\coloneqq \eta/\lambda,
        \qquad
        \tau\coloneqq \lambda^{-2}.
\]
For \(A\in\mathbb C^{W^d\times W^d}\) with \(\Im A\succeq0\), define
\begin{align}
\label{eq:fh-GA-def}
        G_A(\xi)\coloneqq
        \bigl(h(\xi)-e-i\kappa-\tau A\bigr)^{-1},\quad
        m_A\coloneqq
        \int_{\mathbb T^d}G_A(\xi)\,\rmd\xi.
\end{align}
The scaled fibered MDE is
\begin{equation}
\label{equ:scaled_MDE_fixed_block_pert}
        \sigma=\mathscr S_W[m_\sigma],
\end{equation}
for \(\sigma\in\mathbb C^{W^d\times W^d}\).  Equivalently, shorthand
\[
        \theta\coloneqq \lambda^{-1}\sigma,
        \qquad
        \Theta\coloneqq \I_{\tilde{\mathbb Z}_n^d}\otimes\theta,
\]
then \cref{equ:scaled_MDE_fixed_block_pert} is the unscaled equation
\begin{equation}
\label{equ:unscaled_MDE_fixed_block_pert}
        \theta
        =
        \mathscr S_W\left[
        \int_{\mathbb T^d}
        (\lambda \hat{\Psi}(\xi)-z-\theta)^{-1}\,\rmd\xi
        \right].
\end{equation}
\subsection{Uniform fixed-block MDE bounds}
The upper and lower bounds below have different sources.  Shell counting controls the total imaginary mass from above, whereas one isolated regular channel prevents that mass from vanishing.  Two-sided flatness then converts these trace estimates into matrix inequalities.

\begin{prop}[Upper bounds for the fixed-block MDE]
\label{prop:decoup-MDE-upper-bounds}
Assume that $I$ satisfies \textnormal{(C1)}--\textnormal{(C2)} of Definition~\ref{def:fh-bulk-window}.  There are constants $K_*,T_*<\infty$ and $\tau_*>0$, depending only on $(W^d,C_{\rm sh},s_0,M_0)$, such that the following holds.

Let $e\in I$, $0<\tau\le\tau_*$, $0<\kappa\le\tau$, and let $\sigma$ be any solution of~\cref{equ:scaled_MDE_fixed_block_pert} with \(\Im\sigma\succeq0\).  Then
\begin{equation}
\label{eq:fh-upper-final}
        \|\Re\sigma\|_{2\to 2}\le K_*,
        \qquad
        \Tr\Im m_\sigma\le T_*.
\end{equation}
\end{prop}

\begin{proof}
By the fibered MDE and Lemma~\ref{lemma:fh-two-sided-flatness},
\begin{equation}
\label{eq:fh-Y-vs-t}
        \frac {\Tr\Im m_\sigma}{W^d}\I_{W^d}\preceq \Im \sigma\preceq\frac {2\Tr\Im m_\sigma}{W^d}\I_{W^d}.
\end{equation}

Recall $\kappa/\tau\in(0,1]$.  The dissipative part $\Gamma\coloneqq\kappa\I_{W^d}+\tau \Im\sigma$ satisfies
\begin{equation}
\label{eq:fh-Gamma-alpha-beta}
        \left(\frac{\kappa}{\tau} + W^{-d}\Tr\Im m_\sigma\right)\tau\I_{W^d}
        \preceq\Gamma\preceq
        \left(\frac{\kappa}{\tau} + 2W^{-d}\Tr\Im m_\sigma\right)\tau\I_{W^d}.
\end{equation}

\emph{Step 1: Poisson-kernel type decoupled inequality.} The Ward identity gives
\[
        \Im G_\sigma(\xi)
        =G_\sigma(\xi)^*\Gamma G_\sigma(\xi).
\]
Using~\cref{eq:fh-Gamma-alpha-beta},
\[
\begin{split}
        \Tr\Im m_\sigma = \int_{\mathbb T^d}\Tr\Im G_\sigma(\xi)\,\rmd\xi\geq
        \frac{\tau \Tr\Im m_\sigma}{W^d}
        \int_{\mathbb T^d}
        \|G_\sigma(\xi)\|_{\rm HS}^2\,\rmd\xi.
\end{split}
\]
Since $\Tr\Im m_\sigma>0$, Jensen's inequality gives
\begin{equation}
\label{eq:fh-HS-cancellation}
        \Big\|\int_{\mathbb T^d}G_\sigma(\xi)\rmd \xi\Big\|_{\rm HS}^2\leq\int_{\mathbb T^d}
        \|G_\sigma(\xi)\|_{\rm HS}^2\,\rmd\xi
        \le\frac {W^d}{\tau}.
\end{equation}
The fibered MDE and the comparison of Hilbert--Schmidt and operator norms give
\begin{equation}
\label{eq:fh-preliminary-sigma}
        \|\sigma\|_{2\to2}
        \le C_{W^d}\|m_\sigma\|_{\rm HS}
        \le C_{W^d}\sqrt{W^d/\tau}.
\end{equation}
Moreover, from~\cref{eq:fh-imag-congruence}, $\Im G_\sigma(\xi)\preceq\Gamma^{-1}$, and therefore
\[
        \Tr\Im m_\sigma\le\Tr\Gamma^{-1}
        \le\frac{W^{2d}}{\tau \Tr\Im m_\sigma}.
\]
Thus
\begin{equation}
\label{eq:fh-preliminary-t}
        \Tr\Im m_\sigma\le W^d\tau^{-1/2}.
\end{equation}
These preliminary estimates imply $\tau \|\Re \sigma\|_{2\to 2}+\tau(\kappa/\tau + \Tr\Im m_\sigma)=O(\sqrt\tau)$, and notably, they do not rely on any bootstrap assumption.

\emph{Step 2: an intermediate trace inequality.} Write the real part
\[
        B(\xi)\coloneqq h(\xi)-e-\tau \Re \sigma
\]
and let $\mu_j(\xi)$ be its eigenvalues.  We denote
\[
        w\coloneqq\sqrt{(\kappa/\tau + W^{-d}\Tr\Im m_\sigma)
        (\kappa/\tau + 2W^{-d}\Tr\Im m_\sigma)}\,\tau,
\]
Lemma~\ref{lemma:fh-noncommuting-dissipation} gives
\[
        \Tr\Im m_\sigma\leq (\kappa/\tau + 2W^{-d}\Tr\Im m_\sigma)\tau
        \int_{\mathbb T^d}\sum_{j=1}^{W^d}\frac1{\mu_j(\xi)^2+w^2}\,\rmd\xi
        \leq
        \sqrt{2}\int_{\mathbb T^d}
        \sum_{j=1}^{W^d}\frac{w}{\mu_j(\xi)^2+w^2}\,\rmd\xi.
\]
For $\tau$ small enough, the preliminary estimates ensure $w+\tau \|\Re\sigma\|_{2\to2}\le s_0/2$.  Lemma~\ref{lemma:fh-Poisson-shell} yields
\begin{equation}
\label{eq:fh-t-vs-u}
        \Tr\Im m_\sigma
        \le C\left(1+\frac{\tau\|\Re\sigma\|_{2\to2}}{w}\right)
        \le C\left(1+\frac{\|\Re\sigma\|_{2\to2}}{\kappa/\tau + W^{-d}\Tr\Im m_\sigma}\right)
        \le C\left(1+\frac{W^d\|\Re\sigma\|_{2\to2}}{\Tr\Im m_\sigma}\right).
\end{equation}
Consequently,
\begin{equation}
\label{eq:fh-t-sqrt-u}
        \Tr\Im m_\sigma\le C\bigl(1+\sqrt {\|\Re\sigma\|_{2\to2}}\bigr).
\end{equation}

\emph{Step 3: bounding the real part of the fibered MDE solution.} Let
\[
        G_0(\xi)\coloneqq\bigl(h(\xi)-e-i\tau\bigr)^{-1},
        \qquad
        m_0\coloneqq\int_{\mathbb T^d}G_0(\xi)\,\rmd\xi.
\]
By~\cref{eq:fh-free-Cauchy}, $\|m_0\|\le M_0$.  Furthermore,
\[
        G_\sigma(\xi)^{-1}
        =G_0(\xi)^{-1}-Q,
        \qquad
        Q\coloneqq
        \tau \Re\sigma+i\bigl(\Gamma-\tau\I_{W^d}\bigr).
\]
By~\cref{eq:fh-Y-vs-t} and $0<\kappa/\tau\le1$,
\begin{equation}
\label{eq:fh-Q-bound}
        \|Q\|_{2\to2}
        \le C_{W^d}\tau (1+ \|\Re\sigma\|_{2\to 2} + \Tr\Im m_\sigma).
\end{equation}
Set
\[
        \delta\coloneqq4C_{W^d}\tau (1+ \|\Re\sigma\|_{2\to 2} + \Tr\Im m_\sigma),
        \qquad
        \mathcal R\coloneqq
        \{\xi:r_e(\xi)\le\delta\}.
\]
By~\cref{eq:fh-preliminary-sigma}--\cref{eq:fh-preliminary-t}, $\delta=O(\sqrt\tau)$, so $\delta\le s_0/2$ for $\tau$ small.  The free shell bound gives
\begin{equation}
\label{eq:fh-resonant-volume}
        |\mathcal R|\le C_{\rm sh}\delta.
\end{equation}
On $\mathcal R^c$, $\|G_0(\xi)Q\|\le1/4$.  Hence the Neumann resolvent identity and Lemma~\ref{lemma:fh-off-shell} imply
\begin{equation}
\label{eq:fh-off-resonant-difference}
\begin{split}
        \left\|
        \int_{\mathcal R^c}
        \bigl(G_\sigma(\xi)-G_0(\xi)\bigr)\,\rmd\xi
        \right\|
        &\le
        C\|Q\|
        \int_{\mathcal R^c}\frac{\rmd\xi}{r_e(\xi)^2}\\
        &\le C.
\end{split}
\end{equation}
On the resonant set, Cauchy--Schwarz, \cref{eq:fh-HS-cancellation}, and~\cref{eq:fh-resonant-volume} give
\begin{equation}
\label{eq:fh-resonant-G}
        \left\|\int_{\mathcal R}G_\sigma(\xi)\,\rmd\xi\right\|
        \le
        |\mathcal R|^{1/2}
        \left(\int_{\mathcal R}\|G_\sigma(\xi)\|_{\rm HS}^2\,\rmd\xi\right)^{1/2}
        \le C\sqrt{1+ \|\Re\sigma\|_{2\to 2} + \Tr\Im m_\sigma}.
\end{equation}
The free Ward identity and~\cref{eq:fh-free-shell} similarly yield
\begin{equation}
\label{eq:fh-resonant-G0}
        \left\|\int_{\mathcal R}G_0(\xi)\,\rmd\xi\right\|
        \le C\sqrt{1+ \|\Re\sigma\|_{2\to 2} + \Tr\Im m_\sigma} .
\end{equation}
Combining~\cref{eq:fh-free-Cauchy} and \cref{eq:fh-off-resonant-difference}--\cref{eq:fh-resonant-G0},
\begin{equation}
\label{eq:fh-matrix-average-ineq}
        \|m_\sigma\|_{2\to2}
        \le \|m_\sigma - m_0\|_{2\to 2} + \|m_0\|_{2\to 2} \leq C\bigl(1+\sqrt{1+\|\Re\sigma\|_{2\to 2}+\Tr\Im m_\sigma}\bigr).
\end{equation}
Since
\[
        \|\mathscr S_W[m_\sigma]\|_{2\to 2}
        \leq (1+W^{-d})\|m_\sigma\|_{2\to 2},
\]
the fibered MDE gives
\begin{equation}
\label{eq:fh-u-vs-t}
        \|\Re \sigma\|_{2\to 2}\le C\bigl(1+\sqrt{1+\|\Re\sigma\|_{2\to 2}+\Tr\Im m_\sigma}\bigr).
\end{equation}

\emph{Step 4: self-comparison inequality to reveal the upper bound.} By~\cref{eq:fh-t-sqrt-u}, $1+\|\Re\sigma\|_{2\to2}+ \Tr\Im m_\sigma\le C(1+\|\Re \sigma\|_{2\to2})$.  Hence~\cref{eq:fh-u-vs-t} implies
\[
        \|\Re\sigma\|_{2\to2}\le C\bigl(1+\sqrt{1+\|\Re\sigma\|_{2\to2}}\bigr),
\]
which forces $\|\Re\sigma\|_{2\to 2}\le K_*$ for a constant $K_*$.  Inserting this into \cref{eq:fh-t-sqrt-u} gives $\Tr\Im m_\sigma\le T_*$.  This proves \cref{eq:fh-upper-final}.
\end{proof}

\subsection{Lower bound from an isolated channel}

\begin{lemma}[Persistence of the isolated channel]
\label{lemma:fh-isolated-persistence}
Assume~\textnormal{(C3)} and fix $K<\infty$.  Then there exists $\tau_K>0$ such that, for every $e\in I$, every Hermitian $X$ with $\|X\|\le K$, and every $0<\tau\le\tau_K$, the matrix $h\bigl(\Xi_e(s,\omega)\bigr)-e-\tau X$ has a simple eigenvalue $\mu_X(s,\omega)$ satisfying
\begin{equation}
\label{eq:fh-mu-close}
        |\mu_X(s,\omega)-s|\le\tau K,
        \qquad
        |\partial_s\mu_X(s,\omega)-1|\le\frac12
\end{equation}
for $|s|\le\delta_{\rm iso}$ and $\omega\in V_e$.
\end{lemma}

\begin{proof}
Set
\[
        H_0(s,\omega)=h(\Xi_e(s,\omega))-e,
        \qquad
        H_X(s,\omega)=H_0(s,\omega)-\tau X.
\]
Weyl's inequality gives $|\mu_X-s|\leq\tau K$.  If $\tau K\leq g_{\rm iso}/4$, the contour
\[
        \Gamma_s=\{z\in\mathbb C:|z-s|=g_{\rm iso}/2\}
\]
encloses exactly one eigenvalue of both matrices.  Their rank-one Riesz projections are
\[
        P_Y(s,\omega)
        =\frac{1}{2\pi i}\int_{\Gamma_s}(z-H_Y(s,\omega))^{-1}\,\rmd z,
        \qquad Y\in\{0,X\}.
\]
On $\Gamma_s$ the two resolvents are bounded by $2/g_{\rm iso}$ and $4/g_{\rm iso}$, respectively.  The resolvent identity therefore gives
\[
        \|P_X(s,\omega)-P_0(s,\omega)\|
        \leq C\frac{\tau K}{g_{\rm iso}}.
\]
The perturbation $-\tau X$ is independent of $s$, while the unperturbed branch satisfies
\[
        \varepsilon_{j_e}(\Xi_e(s,\omega))=e+s.
\]
The Hellmann--Feynman formula therefore gives
\[
        \partial_s\mu_X
        =\Tr\bigl(P_X\,\partial_s h(\Xi_e(s,\omega))\bigr),
        \qquad
        1
        =\Tr\bigl(P_0\,\partial_s h(\Xi_e(s,\omega))\bigr).
\]
Condition~\textnormal{(C3)} now yields $|\partial_s\mu_X-1|\leq CC_{\rm iso}\tau K/g_{\rm iso}$, which is at most $1/2$ after decreasing $\tau_K$.
\end{proof}

\begin{prop}[Decoupled lower bound of the imaginary part]
\label{prop:decoup-MDE-lower-bound}
Assume that $I$ is a finite-height folded bulk window in the sense of Definition~\ref{def:fh-bulk-window}.  Let $K_*,T_*$ be the constants from Proposition~\ref{prop:decoup-MDE-upper-bounds}.  There are $t_*>0$ and $\tau_{**}>0$ such that every solution in Proposition~\ref{prop:decoup-MDE-upper-bounds} satisfies
\begin{equation}
\label{eq:fh-t-lower}
        \Tr\Im m_\sigma\ge t_*
\end{equation}
for $0<\tau\le\tau_{**}$.
\end{prop}

\begin{proof}
Proposition~\ref{prop:decoup-MDE-upper-bounds} gives
\[
        \|\Re \sigma\|\le K_*,
        \qquad
        \Tr\Im m_\sigma\le T_*,
\]
Note that $\kappa/\tau\in(0,1]$, then
\[
        \frac{\kappa}{\tau} + \frac{2\Tr\Im m_\sigma}{W^d}\leq 1+2T_*/{W^d}.
\]
For brevity set
\[
        \alpha=\frac{\kappa}{\tau}+W^{-d}\Tr\Im m_\sigma,
        \qquad
        \beta=\frac{\kappa}{\tau}+2W^{-d}\Tr\Im m_\sigma.
\]
Let $\mu_X(s,\omega)$ be the isolated perturbed eigenvalue from Lemma~\ref{lemma:fh-isolated-persistence}.  By the lower Loewner bound in \cref{eq:fh-Poisson-Loewner}, using again the choice $w = \sqrt{\alpha\beta}\,\tau$ from Proposition~\ref{prop:decoup-MDE-upper-bounds}, we have
\begin{equation}
\label{eq:fh-lower-patch-start}
\begin{split}
        \Tr\Im m_\sigma &\ge\int_{\mathbb T^d}\sum_{j=1}^{W^d}\frac{\alpha\tau}{\mu_j(\xi)^2+w^2}\,\rmd\xi\\
        &\ge
        \int_{V_e}\int_{-\delta_{\rm iso}}^{\delta_{\rm iso}}
        J_e(s,\omega)
        \frac{\alpha\tau}
        {\mu_X(s,\omega)^2+w^2}
        \,\rmd s\,\rmd\omega,
\end{split}
\end{equation}
where $J_e$ is the Jacobian of $\Xi_e$ and $J_e\ge j_{\rm iso}$.  By~\cref{eq:fh-mu-close}, for each $\omega$ the map $s\mapsto\mu_X(s,\omega)$ is strictly increasing, its derivative lies in $[1/2,3/2]$, and for $\tau$ small, its range contains $[-\delta_{\rm iso}/2,\delta_{\rm iso}/2]$.  Changing variables $r=\mu_X(s,\omega)$ in~\cref{eq:fh-lower-patch-start}, we have \(\rmd s/\rmd r=(\partial_s\mu_X)^{-1}\in[2/3,2]\), and hence
\begin{align*}
        \Tr\Im m_\sigma&\ge
        \frac{2v_{\rm iso}j_{\rm iso}}{3}\int_{-\delta_{\rm iso}/2}^{\delta_{\rm iso}/2}
        \frac{\alpha\tau}{r^2+w^2}\,\rmd r\\
        &=
        \frac{4v_{\rm iso}j_{\rm iso}}{3}\sqrt{\frac\alpha\beta}
        \arctan\!\left(
        \frac{\delta_{\rm iso}}{2\sqrt{\alpha\beta}\tau}
        \right).
\end{align*}
Since $\beta/\alpha\le2$ and $\beta\leq1+2T_*/W^d$, the last expression is bounded below by a positive constant $t_*$ for all sufficiently small $\tau$.  This proves~\cref{eq:fh-t-lower}.
\end{proof}
\subsection{The fixed-block MDE trap}
\begin{thm}
\label{thm:fh-fixed-block-MDE}
Assume that $W<\infty$ is fixed and that $I$ is a finite-height folded bulk window for $h=\widehat\Psi_W$.  There are constants $0<c<C<\infty$ and $\lambda_0<\infty$, depending only on the finite-height window parameters, such that whenever \(\lambda\geq\lambda_0\), \(e=E/\lambda\in I\), \(0<\eta\leq\lambda^{-1}\), \(\tau=\lambda^{-2}\), and \(\kappa=\eta/\lambda\), the scaled MDE \cref{equ:scaled_MDE_fixed_block_pert} admits a unique solution $\sigma$ with $\Im\sigma\succ0$, satisfying
\begin{equation}
\label{eq:fh-sigma-final-bounds}
        \|\sigma\|_{2\to2}\le C,
        \qquad
        c\I_{W^d}\preceq\Im\sigma\preceq C\I_{W^d}.
\end{equation}
Consequently, for the unscaled self-energy $\theta=\lambda^{-1}\sigma$,
\begin{equation}
\label{eq:fh-theta-final-bounds}
        \|\theta\|_{2\to2}\le\frac C\lambda,
        \qquad
        \frac c\lambda\I_{W^d}
        \preceq\Im\theta
        \preceq\frac C\lambda\I_{W^d}.
\end{equation}
\end{thm}

\begin{proof}
\emph{Existence and uniqueness.} The positivity-preserving covariance map satisfies the hypotheses of the standard MDE existence and uniqueness theorem \citep{Ajanki_2018,helton2007operatorvaluedsemicircularelementssolving}. For completeness, existence in the present finite-dimensional formulation also follows directly from the following invariant-ball argument.

\emph{Uniform bounds.} For $\kappa>0$, consider the compact convex set
\[
        \mathcal K_\kappa
        \coloneqq
        \left\{
        A\in\mathbb C^{W^d\times W^d}:
        \Im A\succeq0,
        \ \|A\|_{2\to2}\le (1 + \frac{1}{W^d})/\kappa
        \right\},
\]
Define $\Phi_{\tau,\kappa,e}(A)=\mathscr S_W[m_A]$.  If $A\in\mathcal K_\kappa$, then
\[
        \kappa\I_{W^d}+\tau\Im A\succeq\kappa\I_{W^d},
\]
so by Proposition~\ref{prop:num_range}, $\|G_A(\xi)\|\le\kappa^{-1}$ and $\|\Phi_{\tau,\kappa,e}(A)\|\le (1+1/W^d)/\kappa$.  Moreover, $\Im\Phi_{\tau,\kappa,e}(A)\succeq0$.  Hence $\Phi_{\tau,\kappa,e}$ maps $\mathcal K_\kappa$ continuously into itself, and Brouwer's theorem gives a fixed point.

Propositions~\ref{prop:decoup-MDE-upper-bounds} and \ref{prop:decoup-MDE-lower-bound} apply to every such fixed point.  Together with \cref{eq:fh-Y-vs-t}, they yield
\[
        \frac{t_*}{W^d}\I_{W^d}
        \preceq\Im\sigma
        \preceq\frac{2T_*}{W^d}\I_{W^d}.
\]
Recall that $\|\Re \sigma\| \leq K_*$, so~\cref{eq:fh-sigma-final-bounds} follows. Rescaling gives~\cref{eq:fh-theta-final-bounds}.
\end{proof}
\subsection{Stability and comparison of self-energies}
\label{subsec:connect_resol_self_energy_fixed}
We next compare an approximate self-energy with the exact MDE solution. Because the fibers are coupled through their average, the relevant stability operator acts on a single cell matrix rather than pointwise in \(\xi\). This prevents us from directly designing a pointwise stability operator at all $\xi \in \mathbb T^d$, but instead, we construct a torus-averaged stability operator. The fiber identity reads
\[
 \I+
 \bigl(z/\lambda-\widehat\Psi(\xi)+\tau\sigma\bigr)
 G_\sigma(\xi)=0.
\]
We first prove invertibility of the corresponding continuum one-point operator, extend it to the mixed two-point operator, and finally pass to the finite torus by the Combes--Thomas estimate.
\subsubsection{Uniform one-point stability from one isolated channel}
\label{subsec:fh-one-point-stability}
Let \(I\) be a finite-height folded bulk window in the sense of Definition~\ref{def:fh-bulk-window}.  Fix $e\in I, 0<\tau\leq\tau_0, 0<\kappa\leq\tau,$ and let \(\sigma\) be the exact solution of
\begin{equation}
\label{eq:fh-stab-MDE}
        \sigma
        =
        \mathscr S_W\!\left[\int_\xi G_\sigma(\xi)\right],
        \qquad
        G_\sigma(\xi)
        =
        \bigl(\hat{\Psi}(\xi)-e-i\kappa-\tau\sigma\bigr)^{-1}.
\end{equation}
By Theorem~\ref{thm:fh-fixed-block-MDE}, after decreasing \(\tau_0\) if necessary, there are constants \(0<c_\sigma<C_\sigma<\infty\), uniform in \((e,\tau,\kappa)\), such that
\begin{equation}
\label{eq:fh-stab-trap}
        \|\sigma\|_{2\to2}\leq C_\sigma,
        \qquad
        c_\sigma\I_{W^d}\preceq\Im\sigma\preceq C_\sigma\I_{W^d}.
\end{equation}

\begin{lemma}
\label{lem:fh-localized-Ward-Schwarz}
For a measurable set \(\mathscr E\subseteq\mathbb T^d\), define
\[
        H_{\mathscr E}\coloneqq\int_{\mathscr E}\Im G_\sigma(\xi)\,\rmd\xi, \quad \mathcal L_{\mathscr E}[X]
        \coloneqq
        \tau\mathscr S_W\!\left[
        \int_{\mathscr E}
        G_\sigma(\xi)XG_\sigma(\xi)\,\rmd\xi
        \right].
\]
Then, for arbitrary \(u,v\in\mathbb C^{W^d}\),
\begin{equation}
\label{eq:fh-localized-bilinear}
        \bigl|u^*\mathcal L_{\mathscr E}[X]v\bigr|
        \leq
        \|\mathcal C_{\Im\sigma}[X]\|_{2\to2}
        \bigl(u^*\mathscr S_W[H_{\mathscr E}]u\bigr)^{1/2}
        \bigl(v^*\mathscr S_W[H_\mathscr E]v\bigr)^{1/2}.
\end{equation}
Consequently, if $\mathscr S_W[H_\mathscr E]\preceq(1-\delta)\Im\sigma$ for some \(\delta\in[0,1]\), then
\begin{equation}
\label{eq:fh-localized-contraction}
        \|\mathcal C_{\Im \sigma}\mathcal L_{\mathscr E}[X]\|_{2\to2}
        \leq
        (1-\delta)\|\mathcal C_{\Im\sigma}[X]\|_{2\to2}.
\end{equation}
\end{lemma}

\begin{proof}
Recall that $\Gamma\coloneqq\kappa\I_{W^d}+\tau \Im\sigma.$ The generalized Ward identity gives
\begin{equation}
\label{eq:fh-stab-two-Ward}
        \Im G_\sigma(\xi)
        =
        G_\sigma(\xi)^*\Gamma G_\sigma(\xi)
        =
        G_\sigma(\xi)\Gamma G_\sigma(\xi)^*.
\end{equation}
Since \(\Gamma\succeq\tau\Im\sigma\), the definition of the sandwich operator and Cauchy--Schwarz in \(\xi\) imply that, for arbitrary \(a,b\in\mathbb C^{W^d}\),
\begin{equation}
\label{eq:fh-stab-bilinear-Ward}
    \begin{aligned}
        \tau\left|
        a^*\left(
        \int_{\mathscr E}G_\sigma XG_\sigma
        \right)b
        \right|
        \leq
        \|\mathcal{C}_{\Im\sigma}[X]\|_{2\to2}
        \bigl(a^*H_{\mathscr E}a\bigr)^{1/2}
        \bigl(b^*H_{\mathscr E}b\bigr)^{1/2}.
    \end{aligned}
\end{equation}
Indeed, for $\Im \sigma \in \mathbb C^{W^d\times W^d}$ we have
\begin{equation}
    \begin{aligned}
        |a^\star G_\sigma XG_\sigma b| &= \left|a^\star G_\sigma \Im\sigma^{1/2}\Im\sigma^{-1/2} X \Im\sigma^{-1/2} \Im\sigma^{1/2}G_\sigma b\right|\\
        &\leq
        \|\mathcal C_{\Im\sigma}[X]\|_{2\to 2}
        (a^\star G_\sigma \Im\sigma G_\sigma^\star a)^{1/2}
        (b^\star G_\sigma^\star \Im\sigma G_\sigma b)^{1/2}.
    \end{aligned}
\end{equation}
Cauchy--Schwarz in \(\xi\) gives
\begin{equation}
    \begin{aligned}
        \tau\left|a^\star\int_\mathscr E G_\sigma X G_\sigma\rmd \xi b\right|&\leq\tau\int_\mathscr E\left|a^\star G_\sigma X G_\sigma b\right|\rmd \xi\\
        &\leq \|\mathcal{C}_{\Im\sigma}[X]\|_{2\to 2}
        \left(\int_\mathscr E \tau a^\star G_\sigma \Im\sigma G_\sigma^\star a\rmd \xi\right)^{1/2}
        \left(\int_\mathscr E \tau b^\star G_\sigma^\star \Im\sigma G_\sigma b\rmd \xi\right)^{1/2}
    \end{aligned}
\end{equation}
Both \(\tau G_\sigma^*\Im\sigma G_\sigma\) and \(\tau G_\sigma\Im\sigma G_\sigma^*\) are bounded above by \(\Im G_\sigma\) in Loewner order, which proves \cref{eq:fh-stab-bilinear-Ward}.  Applying this estimate to an orthonormal basis gives
\[
        \tau\left|\Tr \int_\mathscr E G_\sigma X G_\sigma \rmd \xi\right|
        \leq
        \|\mathcal C_{\Im\sigma}[X]\|_{2\to2}\Tr H_{\mathscr E}.
\]
Using the explicit formula for \(\mathscr S_W\), followed by Cauchy--Schwarz, we obtain
\begin{align*}
        \bigl|u^*\mathcal L_{\mathscr E}[X]v\bigr|
        &\leq
        \frac{\|\mathcal C_{\Im\sigma}[X]\|_{2\to2}}{W^d}
        \left[
        \bigl(u^*H_{\mathscr E}^{\top}u\bigr)^{1/2}
        \bigl(v^*H_{\mathscr E}^{\top}v\bigr)^{1/2}
        +
        (\Tr H_{\mathscr E})\|u\|\|v\|
        \right]
        \\
        &\leq
        \|\mathcal C_{\Im\sigma}[X]\|_{2\to2}
        \bigl(u^*\mathscr S_W[H_\mathscr E]u\bigr)^{1/2}
        \bigl(v^*\mathscr S_W[H_\mathscr E] v\bigr)^{1/2}.
\end{align*}
This proves \cref{eq:fh-localized-bilinear}.  If \(\mathscr S_W[H_\mathscr E]\preceq(1-\delta)\Im\sigma\), divide by \((u^*\Im\sigma u)^{1/2}(v^*\Im\sigma v)^{1/2}\) and take the supremum over \(u,v\) to get \cref{eq:fh-localized-contraction}.
\end{proof}

We now isolate a fixed part of the regular channel from condition \textnormal{(C3)}.  Choose once and for all
\[
        0<\delta_{\rm p}
        <
        \min\{\delta_{\rm iso},g_{\rm iso}/8\}
\]
and set
\begin{equation}
\label{eq:fh-stab-patch}
        \mathcal P_e
        \coloneqq
        \Xi_e\bigl(
        (-\delta_{\rm p},\delta_{\rm p})\times V_e
        \bigr).
\end{equation}
Only this patch will be treated by cancellation.  The complementary resonances, including all folded-band crossings, will be treated solely by Lemma~\ref{lem:fh-localized-Ward-Schwarz}.

\begin{lemma}
\label{lem:fh-positive-patch-mass}
There are \(c_{\rm p}>0\) and \(\tau_{\rm p}>0\), depending only on the finite-height window parameters, \(W\), \(d\), and the trapped bounds in Theorem~\ref{thm:fh-fixed-block-MDE}, such that
\begin{equation}
\label{eq:fh-positive-patch-mass}
        \Tr H_{\mathcal P_e}
        =
        \Tr\int_{\mathcal P_e}\Im G_\sigma(\xi)\,\rmd\xi
        \geq c_{\rm p}
\end{equation}
for \(0<\tau\leq\tau_{\rm p}\), uniformly in \(e\in I\) and \(0<\kappa\leq\tau\).
\end{lemma}

\begin{proof}
This is the localized estimate already used in the proof of Proposition~\ref{prop:decoup-MDE-lower-bound}; we include it to make the stability argument self-contained.

Write
\[
        B(\xi)\coloneqq\hat{\Psi}(\xi)-e-\tau\Re\sigma.
\]
By Lemma~\ref{lemma:fh-isolated-persistence}, on \(\xi=\Xi_e(s,\omega)\) the matrix \(B(\xi)\) has a simple eigenvalue \(\mu(s,\omega)\) satisfying
\[
        |\mu(s,\omega)-s|\leq C\tau,
        \qquad
        \frac12\leq\partial_s\mu(s,\omega)\leq\frac32.
\]
The trapped bounds give
\[
        c\tau\I_{W^d}
        \preceq
        \Gamma=\kappa\I_{W^d}+\tau \Im\sigma
        \preceq
        C\tau\I_{W^d}.
\]
The lower Loewner bound from Lemma~\ref{lemma:fh-noncommuting-dissipation}, followed by the lower Jacobian bound in condition~\textnormal{(C3)}, yields
\[
        \Tr H_{\mathcal P_e}
        \geq
        c
        \int_{V_e}\int_{-\delta_{\rm p}}^{\delta_{\rm p}}
        \frac{\tau}
        {\mu(s,\omega)^2+C\tau^2}
        \,\rmd s\,\rmd\omega.
\]
Changing variables \(r=\mu(s,\omega)\), and using \(|V_e|\geq v_{\rm iso}\), gives a strictly positive lower bound independent of \(\tau\).
\end{proof}
The following Lemma quantifies the size of the double resolvent and the reveals the underlying cancellation on the isolated patch, which would be adopted to prove the invertibility of the one-point stability operator.
\begin{lemma}[Double-resolvent bound on the isolated patch]
\label{lemma:fh_patch_two_resolv}
There are \(C_{\rm p}<\infty\) and \(\tau_{\rm p}'>0\) such that
\begin{equation}
\label{eq:fh-patch-two-resolvent}
        \left\|
        \int_{\mathcal P_e}
        G_\sigma(\xi)XG_\sigma(\xi)\,\rmd\xi
        \right\|_{2\to2}
        \leq
        C_{\rm p}\bigl(1+|\log\tau|\bigr)\|X\|_{2\to2}
\end{equation}
for \(0<\tau\leq\tau_{\rm p}'\), uniformly in \(e,\kappa\), and \(X\).
\end{lemma}

\begin{proof}
Fix \(\omega\in V_e\) and suppress it from the notation.  Set
\begin{equation}
\label{eq:fh-patch-B-Gamma-clean}
        B(s)\coloneqq \hat{\Psi}(\Xi_e(s,\omega))-e-\tau\Re\sigma,\qquad
        \Gamma\coloneqq\kappa\I_{W^d}+\tau\Im\sigma .
\end{equation}
Thus
\[
        G_\sigma(\Xi_e(s,\omega))
        =
        (B(s)-i\Gamma)^{-1}.
\]
The matrix \(\Gamma\) is independent of \(s\), and the trapped bounds give
\begin{equation}
\label{eq:fh-patch-Gamma-clean}
        c\tau\I_{W^d}\preceq\Gamma\preceq C\tau\I_{W^d}.
\end{equation}

Let \(\mu(s)=\mu(s,\omega)\) be the persisted simple eigenvalue of \(B(s)\), let $u(s)$ be a corresponding normalized eigenvector, let \(P(s)=P(s,\omega)\) be its rank-one orthogonal Riesz projection, and set \(Q(s)=\I_{W^d}-P(s)\).  The isolated gap, the \(C^2\)-bounds on the chart, and the Riesz formula imply
\begin{equation}
\label{eq:fh-patch-channel-C2-clean}
        \|P\|+\|\partial_sP\|+\|\partial_s^2P\|\leq C,
        \qquad
        |\mu-s|\leq C\tau,
        \qquad
        |\partial_s\mu-1|\leq C\tau,
        \qquad
        |\partial_s^2\mu|\leq C.
\end{equation}
The projection bounds follow from the Riesz formula, while the estimates for $\mu$ and $\partial_s\mu$ are those of Lemma~\ref{lemma:fh-isolated-persistence}.  To bound the second derivative, choose a phase of $u$ so that $\langle u,\partial_su\rangle=0$. Differentiating $Bu=\mu u$ and restricting to $u^\perp$ gives
\[
        \partial_su
        =-\bigl(Q(B-\mu)Q\bigr)^{-1}_{u^\perp}
          Q(\partial_sB)u.
\]
The persistent gap bounds the reduced inverse by $2/g_{\rm iso}$, so $\|\partial_su\|\leq C$.  A second differentiation of the Hellmann--Feynman identity yields
\[
        \partial_s^2\mu
        =\langle u,(\partial_s^2B)u\rangle
         +2\Re\langle\partial_su,(\partial_sB)u\rangle,
\]
which proves the last estimate in \cref{eq:fh-patch-channel-C2-clean}.  For fixed \(\omega\), choose a \(C^2\) unitary frame $U(s)=\bigl(u(s),U_Q(s)\bigr)$ such that \(u(s)\) spans \(\operatorname{Ran}P(s)\) and the columns of \(U_Q(s)\) span \(\operatorname{Ran}Q(s)\).  Kato parallel transport constructs such a frame along the \(s\)-interval; the short argument is recalled in Appendix~\ref{appx:kato_parall_transport}.  The gap and Proposition~\ref{prop:kato_para_trans} give
\begin{equation}
\label{eq:fh-patch-frame-C2-clean}
        P=uu^*,
        \qquad
        Q=U_QU_Q^*,
        \qquad
        \sum_{j=0}^2\|\partial_s^jU\|\leq C.
\end{equation}

\medskip
\noindent \emph{Step 1: the reduced inverse on the nonresonant subspace.} Define the \((W^d-1)\times(W^d-1)\) coordinate matrices
\[
        B_Q\coloneqq U_Q^*(\hat{\Psi} - e - \tau\Re\sigma)U_Q,
        \qquad
        \Gamma_Q\coloneqq U_Q^*\Gamma U_Q,
        \qquad
        D_Q(s)
        \coloneqq
        B_Q-i\Gamma_Q.
\]
On the smaller patch, the complementary spectrum of \(B\) stays a fixed distance from zero.  Since \(\|\Gamma\|=O(\tau)\), a Neumann argument on \(\mathbb C^{W^d-1}\), followed by differentiation of \(D_QD_Q^{-1}=\I_{W^d-1}\), gives
\begin{equation}
\label{eq:fh-patch-DQ-C2-clean}
        \sum_{j=0}^2
        \|\partial_s^jD_Q^{-1}\|
        \leq C.
\end{equation}
Indeed, $|s|\leq\delta_{\rm p}\leq g_{\rm iso}/8$ and $\tau\|\Re\sigma\|=O(\tau)$ imply that every complementary eigenvalue $\widehat\mu$ of $B$ satisfies $|\widehat\mu|\geq g_{\rm iso}/2$ for small $\tau$.  Thus
\[
        \|B_Q^{-1}\|_{2\to2}\leq\frac{2}{g_{\rm iso}}.
\]
Since $D_Q=B_Q(\I-iB_Q^{-1}\Gamma_Q)$ and $\|\Gamma_Q\|=O(\tau)$, the Neumann series gives
\[
        D_Q^{-1}
        =\sum_{n=0}^\infty(iB_Q^{-1}\Gamma_Q)^nB_Q^{-1},
        \qquad
        \|D_Q^{-1}\|_{2\to2}\leq\frac{4}{g_{\rm iso}}.
\]
For the derivatives, the inverse rule gives \(\partial_sD_Q^{-1}=-D_Q^{-1}(\partial_sD_Q)D_Q^{-1}\), so it suffices to bound \(\partial_sD_Q\).  Explicitly,
\begin{equation}
    \begin{aligned}
        \partial_s D_Q &= \partial_s B_Q - i\partial_s \Gamma_Q\\
        &=(U_Q')^\star B U_Q + U_Q^\star B' U_Q + U_Q^\star B U_Q' - i \left((U_Q')^\star \Gamma U_Q + U_Q^\star \Gamma U_Q'\right),
    \end{aligned}
\end{equation}
For \(j=0,1,2\), \(\|\partial_s^jU_Q\|\leq\|\partial_s^jU\|\). Together with \cref{eq:fh-patch-frame-C2-clean} and boundedness of \(\Gamma,B,B'\), this proves the first-derivative bound.  Differentiating the inverse rule once more proves the second-derivative bound.

Here \(D_Q^{-1}\) is the \((W^d-1)\times(W^d-1)\) coordinate inverse. Its lift to the original \(W^d\)-dimensional space is, by definition,
\begin{equation}
\label{eq:fh-patch-lifted-RQ}
        \mathcal R_Q
        \coloneqq
        U_QD_Q^{-1}U_Q^*.
\end{equation}
Thus \(\mathcal R_Q\) is a \(W^d\times W^d\) matrix, \(\mathcal R_Q=Q\mathcal R_QQ\), and although \(D_Q\) depends on the chosen coordinate frame, its lift \(\mathcal R_Q\) does not.  Moreover, \cref{eq:fh-patch-frame-C2-clean}--\cref{eq:fh-patch-DQ-C2-clean} imply
\begin{equation}
\label{eq:fh-patch-RQ-C2}
        \|\mathcal R_Q\|
        +\|\partial_s\mathcal R_Q\|
        +\|\partial_s^2\mathcal R_Q\|
        \leq C.
\end{equation}

\medskip
\noindent \emph{Step 2: the scalar Schur complement.} We now isolate the only possible singularity in a scalar Schur complement.  Define
\begin{equation}
\label{eq:fh-patch-gamma-chi-clean}
        \gamma\coloneqq u^*\Gamma u=\Tr(P\Gamma),
        \qquad
        g\coloneqq U_Q^*\Gamma u,
        \qquad
        \chi\coloneqq g^*D_Q^{-1}g
        =\Tr(P\Gamma\mathcal R_Q\Gamma P).
\end{equation}
The coefficient $\gamma$ is the direct damping of the persistent channel, whereas $\chi$ records its coupling through the nonresonant complement. Since \(P\) has rank one,
\[
        P\Gamma P=\gamma P,
        \qquad
        P\Gamma\mathcal R_Q\Gamma P=\chi P.
\]
Because \(\Gamma\) is independent of \(s\),
\[
        \partial_s\gamma=\Tr((\partial_sP)\Gamma),
        \qquad
        \partial_s^2\gamma=\Tr((\partial_s^2P)\Gamma).
\]
Equations \cref{eq:fh-patch-Gamma-clean}, \cref{eq:fh-patch-channel-C2-clean}, and \cref{eq:fh-patch-RQ-C2} therefore imply
\begin{equation}
\label{eq:fh-patch-gamma-chi-bounds-clean}
        c\tau\leq\gamma\leq C\tau,
        \qquad
        |\partial_s\gamma|+|\partial_s^2\gamma|\leq C\tau,
        \qquad
        |\chi|+|\partial_s\chi|+|\partial_s^2\chi|
        \leq C\tau^2.
\end{equation}
Indeed, the frame bounds give
\[
        \|g\|+\|\partial_sg\|+\|\partial_s^2g\|
        \leq C\tau,
\]
and the last estimate in \cref{eq:fh-patch-gamma-chi-bounds-clean} follows by differentiating \(\chi=g^*D_Q^{-1}g\).

Since \(P\) is a spectral projection of \(B\) for the persistent eigenvalue $\mu$,
\[
        P(\hat\Psi - e - \tau \Re\sigma)P=\mu P,
        \qquad
        P(\hat\Psi - e - \tau \Re\sigma)Q=Q(\hat\Psi - e - \tau \Re\sigma)P=0.
\]
Thus, in the frame \(U=(u,U_Q)\),
\begin{equation}
\label{eq:fh-patch-coordinate-block-clean}
        U^*(B-i\Gamma)U
        =
        \begin{pmatrix}
        \mu-i\gamma&-ig^*\\
        -ig&D_Q
        \end{pmatrix}.
\end{equation}
The Schur complement of \(D_Q\) is
\[
        (\mu-i\gamma)
        -
        (-ig^*)D_Q^{-1}(-ig)
        =
        \mu-i\gamma+\chi.
\]
Set
\begin{equation}
        \zeta\coloneqq\mu-i\gamma+\chi.
\end{equation}
The exact coordinate block inverse is
\begin{equation}
\label{eq:fh-patch-coordinate-inverse-clean}
        U^*G_\sigma U
        =
        \begin{pmatrix}
        \zeta^{-1}
        &
        i\zeta^{-1} g^*D_Q^{-1}
        \\[1mm]
        iD_Q^{-1}g\zeta^{-1}
        &
        D_Q^{-1}-D_Q^{-1}g\zeta^{-1}g^*D_Q^{-1}
        \end{pmatrix}.
\end{equation}
Conjugating \cref{eq:fh-patch-coordinate-inverse-clean} by \(U\), and using
\[
        P=uu^*,
        \qquad
        g=U_Q^*\Gamma u,
        \qquad
        \mathcal R_Q=U_QD_Q^{-1}U_Q^*,
\]
gives the following exact \(W^d\times W^d\) identity:
\begin{align}
\label{eq:fh-patch-block-inverse-clean}
        G_\sigma
        &=
        \zeta^{-1}P
        +i\zeta^{-1}P\Gamma\mathcal R_Q
        +i\zeta^{-1}\mathcal R_Q\Gamma P
        +\mathcal R_Q
        -
        \zeta^{-1}\mathcal R_Q\Gamma P\Gamma\mathcal R_Q.
\end{align}
Define the regular remainder
\begin{equation}
\label{eq:fh-patch-Kreg-definition}
        \mathcal K_{\rm reg}
        \coloneqq
        i\zeta^{-1} P\Gamma\mathcal R_Q
        +i\zeta^{-1}\mathcal R_Q\Gamma P
        +\mathcal R_Q
        -\zeta^{-1}\mathcal R_Q\Gamma P\Gamma\mathcal R_Q.
\end{equation}
Then the exact decomposition is
\begin{equation}
\label{eq:fh-patch-clean-decomposition}
        G_\sigma(\Xi_e(s,\omega))
        =
        \zeta^{-1}(s,\omega)P(s,\omega)
        +
        \mathcal K_{\rm reg}(s,\omega).
\end{equation}

\medskip
\noindent \emph{Bounds on the scalar pole and on the regular remainder.} By \cref{eq:fh-patch-channel-C2-clean} and \cref{eq:fh-patch-gamma-chi-bounds-clean},
\begin{equation}
\label{eq:fh-patch-zeta-regular-clean}
        |\partial_s\zeta-1|\leq C\tau,
        \qquad
        |\partial_s^2\zeta|\leq C,
        \qquad
        c\tau\leq-\Im\zeta\leq C\tau,
        \qquad
        \Re\zeta=s+O(\tau).
\end{equation}
After decreasing \(\tau_0\), the function \(s\mapsto\Re\zeta(s,\omega)\) has a unique zero \(s_*(\omega)=O(\tau)\), and
\begin{equation}
\label{eq:fh-patch-a-size-clean}
        |\zeta^{-1}(s,\omega)|
        \leq
        \frac{C}{|s-s_*(\omega)|+\tau}.
\end{equation}
In particular, \(|\zeta^{-1}|\leq C/\tau\).  Set \(a\coloneqq\zeta^{-1}\).  Using \(\|\Gamma\|\leq C\tau\) and \cref{eq:fh-patch-RQ-C2} in \cref{eq:fh-patch-Kreg-definition}, we obtain
\begin{equation}
\label{eq:fh-patch-Kreg-bound}
        \|\mathcal K_{\rm reg}(s,\omega)\|
        \leq
        C\bigl(1+|a|\tau+|a|\tau^2\bigr)
        \leq C.
\end{equation}
This is the substantive content of the decomposition: all singular behavior is confined to the scalar factor $\zeta^{-1}P$.

\medskip
\noindent \emph{Step 3: cancellation of the squared pole.} Let $J_e(s,\omega)\coloneqq|\det D\Xi_e(s,\omega)|.$ The chart bounds imply
\begin{equation}
\label{eq:fh-patch-J-bounds-clean}
        \|J_e\|_\infty+\|\partial_sJ_e\|_\infty\leq C.
\end{equation}
Expanding \cref{eq:fh-patch-clean-decomposition},
\begin{align*}
        G_\sigma XG_\sigma
        &=
        \zeta^{-2}PXP
        +\zeta^{-1}PX\mathcal K_{\rm reg}
        +\zeta^{-1}\mathcal K_{\rm reg}XP
        +\mathcal K_{\rm reg}X\mathcal K_{\rm reg}.
\end{align*}
Equations \cref{eq:fh-patch-a-size-clean} and \cref{eq:fh-patch-Kreg-bound} give
\begin{equation}
\label{eq:fh-patch-a-L1-clean}
        \int_{-\delta_{\rm p}}^{\delta_{\rm p}}|\zeta^{-1}(s,\omega)|\,\rmd s
        \leq C(1+|\log\tau|).
\end{equation}
Hence the last three terms are absolutely integrable and contribute at most \(C(1+|\log\tau|)\|X\|\).

For the first term, set
\[
        \mathscr F(s,\omega)\coloneqq P(s,\omega)XP(s,\omega)J_e(s,\omega).
\]
The identity
\[
        \partial_s\zeta^{-1}=-(\partial_s\zeta)\zeta^{-2},
        \qquad
        \zeta^{-2}=-\frac{\partial_s\zeta^{-1}}{\partial_s\zeta},
\]
followed by integration by parts gives
\begin{equation}
\label{eq:fh-patch-integration-parts-clean}
    \begin{aligned}
        \int_{-\delta_{\rm p}}^{\delta_{\rm p}}\zeta^{-2}\mathscr F\,\rmd s
        &=
        -\left[
        \frac{\zeta^{-1}\mathscr F}{\partial_s\zeta}
        \right]_{-\delta_{\rm p}}^{\delta_{\rm p}}
        +
        \int_{-\delta_{\rm p}}^{\delta_{\rm p}}
        \zeta^{-1}\,
        \partial_s\!\left(
        \frac{\mathscr F}{\partial_s\zeta}
        \right)\rmd s.
    \end{aligned}
\end{equation}
By \cref{eq:fh-patch-channel-C2-clean}, \cref{eq:fh-patch-zeta-regular-clean}, and \cref{eq:fh-patch-J-bounds-clean},
\[
        |\partial_s\zeta|\geq\frac12,
        \qquad
        \left\|
        \partial_s\!\left(
        \frac{\mathscr F}{\partial_s\zeta}
        \right)
        \right\|
        \leq C\|X\|.
\]
The endpoints \(s=\pm\delta_{\rm p}\) remain a fixed distance from \(s_*(\omega)\), so the boundary term in \cref{eq:fh-patch-integration-parts-clean} is bounded by \(C\|X\|\). The remaining integral is bounded by \(C(1+|\log\tau|)\|X\|\) using \cref{eq:fh-patch-a-L1-clean}.

Finally, the \(\omega\)-volume is uniformly bounded above.  Indeed, injectivity of \(\Xi_e\) and the lower Jacobian bound imply
\[
        4\delta_{\rm iso}j_{\rm iso}|V_e|
        \leq|\mathbb T^d|.
\]
Integrating the preceding estimates over \(\omega\) proves \cref{eq:fh-patch-two-resolvent}.
\end{proof}
\begin{remark}
    The mechanism in the proof of Lemma~\ref{lemma:fh_patch_two_resolv} is the matrix analogue of
\[
        \frac1{(s-i\tau)^2}
        =
        -\partial_s\frac1{s-i\tau}.
\]
The product \(G_\sigma XG_\sigma\) contains the cancellative squared pole \(\zeta^{-2}\), whereas an adjoint product contains \(|\zeta|^{-2}\) and has no such cancellation.
\end{remark}
\begin{thm}[Uniform continuum one-point stability]
\label{thm:fh-one-point-reduced-stability}
Assume that \(W<\infty\) is fixed, that \(I\) is a finite-height folded bulk window for \(h=\widehat\Psi_W\), and that
\[
        e=E/\lambda\in I,
        \qquad
        \tau=\lambda^{-2},
        \qquad
        0<\kappa=\eta/\lambda\leq\tau.
\]
For all sufficiently large \(\lambda\), the one-point reduced stability operator
\begin{equation}
\label{eq:fh-one-point-B}
        \mathcal B_\sigma[X]
        \coloneqq
        X-
        \tau\mathscr S_W\!\left[
        \int_\xi
        G_\sigma(\xi)XG_\sigma(\xi)
        \right]
\end{equation}
is invertible and satisfies
\begin{equation}
\label{eq:fh-one-point-B-inverse}
        \bigl\|\mathcal B_\sigma^{-1}\bigr\|_{(2\to2)\to(2\to2)}
        \leq C_{W,I}.
\end{equation}
Particularly, for some \(c_{W,I}>0\),
\begin{equation}
\label{eq:fh-one-point-spectrum}
        \operatorname{Spec}(\mathcal L_\sigma)
        \subset
        \{z\in\mathbb C:|z|\leq1-c_{W,I}\},
\end{equation}
where
\[
        \mathcal L_\sigma[X]
        \coloneqq
        \tau\mathscr S_W\!\left[
        \int_\xi G_\sigma(\xi)XG_\sigma(\xi)
        \right].
\]
The constants are independent of \(\lambda,\eta\), and \(e\in I\).
\end{thm}

\begin{proof}
Split $\mathcal L_\sigma = \mathcal L_{\mathcal P_e} + \mathcal L_{\mathcal P_e^c}.$ By Lemma~\ref{lem:fh-positive-patch-mass} and two-sided flatness,
\begin{equation}
\label{eq:fh-patch-Y-lower}
        \mathscr S_W[H_{\mathcal P_e}]
        \succeq
        \frac{\Tr H_{\mathcal P_e}}{W^d}\I_{W^d}        \succeq
        \frac{c_{\rm p}}{W^d}\I_{W^d}.
\end{equation}
Since \(\Im\sigma\preceq C_\sigma\I_{W^d}\), this implies
\[
        \Im\sigma-\mathscr S_W[H_{\mathcal{P}_e}] = \mathscr S_W[H_{\mathcal{P}_e^c}]
        \preceq
        \left(1-\delta_0\right)\Im\sigma,
        \qquad
        \delta_0\coloneqq\frac{c_{\rm p}}{W^dC_\sigma}>0.
\]
Lemma~\ref{lem:fh-localized-Ward-Schwarz} therefore gives
\begin{equation}
\label{eq:fh-complement-strict}
        \|\mathcal{C}_{\Im\sigma}\mathcal L_{\mathcal P_e^c}[X]\|_{2\to 2}
        \leq
        (1-\delta_0)\|\mathcal{C}_{\Im\sigma}[X]\|_{2\to 2}.
\end{equation}
This estimate uses no regularity whatsoever on \(\mathcal P_e^c\).

On the isolated patch, Lemma~\ref{lemma:fh_patch_two_resolv}, the boundedness of \(\mathscr S_W\), and the equivalence of the two matrix norms give
\begin{equation}
\label{eq:fh-patch-small}
        \|\mathcal{C}_{\Im\sigma}\mathcal L_{\mathcal P_e}[X]\|_{2\to 2}
        \leq
        C\tau(1+|\log\tau|)\|\mathcal{C}_{\Im\sigma}[X]\|_{2\to 2}.
\end{equation}
After decreasing \(\tau_0\), the right-hand side is at most \(\delta_0\|\mathcal{C}_{\Im\sigma}[X]\|_{2\to2}/2\).  Combining \cref{eq:fh-complement-strict} and \cref{eq:fh-patch-small}, we obtain
\[
        \|\mathcal{C}_{\Im\sigma}\mathcal L_\sigma[X]\|_{2\to 2}
        \leq
        \left(1-\frac{\delta_0}{2}\right)\|\mathcal{C}_{\Im\sigma}[X]\|_{2\to 2}.
\]
Introduce the weighted norm $\|X\|_\sigma \coloneqq \|\mathcal C_{\Im\sigma}[X]\|_{2\to2}.$ The preceding estimate is a genuine contraction in this norm, so the Neumann series converges and
\[
        \|(\operatorname{Id}-\mathcal L_\sigma)^{-1}X\|_\sigma
        \leq
        \frac{2}{\delta_0}\|X\|_\sigma.
\]
By \cref{eq:fh-stab-trap}, \(\|X\|_{2\to2}\leq C_\sigma\|X\|_\sigma\) and \(\|X\|_\sigma\leq c_\sigma^{-1}\|X\|_{2\to2}\).  Hence
\[
        \|(\operatorname{Id}-\mathcal L_\sigma)^{-1}\|_{{(2\to2)}\to (2\to 2)}
        \leq \frac{2C_\sigma}{\delta_0c_\sigma}.
\]
This proves \cref{eq:fh-one-point-B-inverse}; the contraction in \(\|\cdot\|_\sigma\) also proves \cref{eq:fh-one-point-spectrum}.
\end{proof}

\subsubsection{Bootstrap passage to the continuum two-point operator}
\label{subsec:fh-two-point-stability}

We now pass from the exact one-point operator to the mixed operator needed for the comparison of the two self-energies.  For matrices $\sigma, \tilde{\sigma}\in \mathbb C^{W^d\times W^d},$ define
\begin{equation}
\label{eq:fh-continuum-two-point-B}
        \mathcal B_{\tilde\sigma,\sigma}[X]
        \coloneqq
        X-
        \tau\mathscr S_W\!\left[
        \int_\xi G_{\tilde\sigma}(\xi)XG_\sigma(\xi)
        \right].
\end{equation}
Thus \(\mathcal B_\sigma=\mathcal B_{\sigma,\sigma}\).

\begin{lemma}
\label{lemma:fh-two-point-Lipschitz}
There is a constant $C_{\rm Lip}>0$ depending only on \((W,I)\) such that, for every $\tilde\sigma\in\mathbb C^{W^d\times W^d}$ satisfying
\begin{equation}
\label{eq:fh-A-close-sigma}
        \|\tilde\sigma-\sigma\|_{2\to2}\leq c_\sigma/2,
\end{equation}
\(G_{\tilde\sigma}(\xi)\) is well defined for every \(\xi\) and
\begin{equation}
\label{eq:fh-two-point-Lipschitz-bound}
        \|\mathcal B_{\tilde\sigma,\sigma}-\mathcal B_\sigma\|_
        {(2\to2)\to(2\to2)}
        \leq
        C_{\rm Lip}\|\tilde\sigma-\sigma\|_{2\to2}.
\end{equation}
The constants are uniform in \(e\in I\), \(0<\kappa\leq\tau\), and sufficiently small \(\tau\).
\end{lemma}

\begin{proof}
From \cref{eq:fh-stab-trap} and \cref{eq:fh-A-close-sigma}, $\Im \tilde\sigma\succeq c_\sigma/2\I_{W^d}.$ The numerical-range estimate in Proposition~\ref{prop:num_range} therefore gives
\begin{equation}
\label{eq:fh-G-A-pointwise}
        \sup_{\xi\in\mathbb T^d}\|G_{\tilde\sigma}(\xi)\|_{2\to2}
        \leq
        \frac{1}{\kappa+\tau c_\sigma/2}
        \leq\frac{2}{\tau c_\sigma}.
\end{equation}
We also recall the exact Ward cancellation for the one-point solution from \cref{eq:fh-HS-cancellation}:
\begin{equation}
\label{eq:fh-G-sigma-L2}
        \int_\xi\|G_\sigma(\xi)\|_{2\to2}^2
        \leq
        \int_\xi\|G_\sigma(\xi)\|_{\rm HS}^2
        \leq\frac {W^d}\tau.
\end{equation}

The resolvent identity reads
\begin{equation}
\label{eq:fh-G-A-resolvent-identity}
        G_{\tilde\sigma}(\xi)-G_\sigma(\xi)
        =
        \tau G_{\tilde\sigma}(\xi)(\tilde\sigma - \sigma) G_\sigma(\xi).
\end{equation}
Consequently,
\begin{align}
\label{eq:fh-two-point-difference-exact}
        (\mathcal B_{\tilde\sigma,\sigma}-\mathcal B_\sigma)[X]
        &=-\tau\mathscr S_W\!\left[
        \int_\xi
        (G_{\tilde\sigma}-G_\sigma)XG_\sigma
        \right]
        \nonumber\\
        &=-\tau^2\mathscr S_W\!\left[
        \int_\xi
        G_{\tilde\sigma}(\tilde\sigma - \sigma) G_\sigma XG_\sigma
        \right].
\end{align}
The boundedness of \(\mathscr S_W\), together with \cref{eq:fh-G-A-pointwise,eq:fh-G-sigma-L2}, implies
\begin{align*}
        \| (\mathcal B_{\tilde\sigma,\sigma}-\mathcal B_\sigma)[X]\|_{2\to2}
        &\leq
        C_W\tau^2\|\tilde\sigma - \sigma\|_{2\to2}\|X\|_{2\to2}
        \sup_\xi\|G_{\tilde\sigma}(\xi)\|_{2\to2}
        \int_\xi\|G_\sigma(\xi)\|_{2\to2}^2
        \\
        &\leq
        C_{\rm Lip}\|\tilde\sigma - \sigma\|_{2\to2}\|X\|_{2\to2}.
\end{align*}
Both powers of \(\tau^{-1}\) have canceled the prefactor \(\tau^2\). This proves \cref{eq:fh-two-point-Lipschitz-bound}.
\end{proof}

\begin{lemma}
\label{lemma:fh-bounded-kinetic-extension}
Fix \(K_0<\infty\).  The bounds in Theorem~\ref{thm:fh-fixed-block-MDE}, the inverse estimate in Theorem~\ref{thm:fh-one-point-reduced-stability}, and Lemma~\ref{lemma:fh-two-point-Lipschitz} remain valid for \(0<\kappa\leq K_0\tau\).  Their constants may also depend on \(K_0\).
\end{lemma}
\begin{proof}
In the proof of Proposition~\ref{prop:decoup-MDE-upper-bounds}, replace the condition \(\kappa/\tau\leq1\) by \(\kappa/\tau\leq K_0\).  The Ward cancellation \cref{eq:fh-HS-cancellation} remains unchanged.  The preliminary bounds still give
\[
 \tau\|\Re\sigma\|+\tau(\kappa/\tau+\Tr\Im m_\sigma)
 \leq C_{K_0}\sqrt\tau.
\]
Moreover, \(\beta/\alpha\leq2\) remains valid, and the bound for \(Q\) in \cref{eq:fh-Q-bound} only changes its constant.  The same shell estimate and resonant-set decomposition therefore give
\[
 \Tr\Im m_\sigma\leq C_{K_0}(1+\sqrt{\|\Re\sigma\|}),
 \qquad
 \|\Re\sigma\|\leq
 C_{K_0}(1+\sqrt{1+\|\Re\sigma\|+\Tr\Im m_\sigma}).
\]
These inequalities again imply uniform upper bounds.  In the lower-bound argument, \(\beta\leq K_0+2T_*/W^d\), so the same arctangent expression is bounded below for sufficiently small \(\tau\).  This proves the MDE trap in the enlarged range.

Consequently, \(c\tau\I\preceq\kappa\I+\tau\Im\sigma\preceq C_{K_0}\tau\I\). These are precisely the bounds on \(\kappa\) needed in the positive-patch estimate and in the cancellation of the squared pole.  Their proofs, and hence the weighted contraction in Theorem~\ref{thm:fh-one-point-reduced-stability}, apply with the new constants.  Finally, the Lipschitz proof uses only the lower MDE trap, the numerical-range bound, and the exact Ward cancellation.  Its constant is therefore uniform throughout the enlarged range as well.
\end{proof}

\subsubsection{Quantitative passage from the continuum fixed point to finite volume}
\label{subsec:quantitative-finite-volume-passage}

Recall the scaled continuum and finite-volume maps
\[
        \Phi_\infty(A)
        \coloneqq 
        \mathscr S_W\!\left[\int_{\mathbb T^d}G_A(\xi)\,\rmd\xi\right],
        \qquad
        \Phi_L(A)
        \coloneqq 
        \mathscr S_W\!\left[
        \frac1{n^d}\sum_{\xi\in\tilde{\mathbb T}_n^d}
        G_A^L(\xi)
        \right],
\]
where
\[
        G_A(\xi)
        =\bigl(\hat{\Psi}(\xi)-e-i\kappa-\tau A\bigr)^{-1},
        \qquad
        G_A^L(\xi)
        =\bigl(\hat{\Psi}_L(\xi)-e-i\kappa-\tau A\bigr)^{-1}.
\]
We write
\[
        F_\infty(A)\coloneqq A-\Phi_\infty(A),
        \qquad
        F_L(A)\coloneqq A-\Phi_L(A).
\]
For clarity, the continuum two-point stability operator is denoted by
\begin{equation}
\label{eq:qfv-continuum-two-point-B}
        \mathcal B_{A,B}^{\infty}[X]
        \coloneqq 
        X-
        \tau\mathscr S_W\!\left[
        \int_{\mathbb T^d}G_A(\xi)XG_B(\xi)\,\rmd\xi
        \right].
\end{equation}

\begin{thm}[Quantitative finite-volume fixed point]
\label{thm:finite_volume_from_infinite_trap}
Let \(\sigma\) be the continuum fixed point $\sigma=\Phi_\infty(\sigma)$, and fix \(C_{\rm st}\geq1\) such that \(\|(\mathcal B_\sigma^\infty)^{-1}\|\leq C_{\rm st}\). There exists \(r_*>0\), depending only on \(W\), the folded bulk-window parameters, and the fixed upper bound for \(\lambda\eta\), such that the following holds.  We choose
\[
 r_*\leq c_\sigma/2,
 \qquad C_{\rm st}C_{\rm Lip}r_*\leq1/2.
\]
Define
\begin{equation}
\label{eq:qfv-map-error}
        \varepsilon_L
        \coloneqq 
        \sup_{\|\hat\sigma -\sigma\|_{2\to2}\le r_*}
        \|\Phi_L(\hat\sigma)-\Phi_\infty(\hat\sigma)\|_{2\to2}.
\end{equation}
If \(\varepsilon_L\) is sufficiently small, then the finite-volume MDE $\sigma_L=\Phi_L(\sigma_L)$ has a unique solution in the upper half-plane satisfying
\begin{equation}
\label{eq:qfv-scaled-closeness}
    \|\sigma_L-\sigma\|_{2\to2}
        \le 2C_{\rm st}\varepsilon_L.
\end{equation}
Equivalently, for \(\theta=\lambda^{-1}\sigma\) and \(\theta_L=\lambda^{-1}\sigma_L\),
\begin{equation}
\label{eq:qfv-unscaled-closeness}
        \|\theta_L-\theta\|_{2\to2}
        \le
        \frac{2C_{\rm st}}{\lambda}\varepsilon_L.
\end{equation}
\end{thm}

\begin{proof}
We first record a uniform quadratic expansion of \(F_\infty\) at \(\sigma\).  Choose \(r_*\leq c_\sigma/2\).  Then Lemma~\ref{lemma:fh-two-point-Lipschitz} implies
\begin{equation}
\label{eq:qfv-B-Lipschitz}
        \bigl\|
        \mathcal B_{\hat\sigma,\sigma}^{\infty}
        -\mathcal B_\sigma^\infty
        \bigr\|_{(2\to2)\to(2\to2)}
        \le C_{\rm Lip}\|\hat\sigma - \sigma\|_{2\to2},
\end{equation}
and hence
\begin{equation}
\label{eq:qfv-N-quadratic}
        \|(\mathcal B_{\hat\sigma,\sigma}^{\infty}
        -\mathcal B_\sigma^\infty)[\hat\sigma - \sigma]\|_{2\to2}
        \le C_{\rm Lip}\|\hat\sigma - \sigma\|_{2\to2}^2.
\end{equation}
We now use the stability inverse to construct a finite-volume zero locally.  Define
\begin{equation}
\label{eq:qfv-Newton-map}
        \mathcal T_L(\hat{\sigma})
        \coloneqq 
        \hat{\sigma}-(\mathcal B_\sigma^\infty)^{-1}F_L(\hat{\sigma}).
\end{equation}
Set $r_L\coloneqq 2C_{\rm st}\varepsilon_L.$ Decrease \(r_*\) if necessary, so that $C_{\rm st}C_{\rm Lip}r_*\leq 1/2, r_*\le{c_\sigma/2},$ and assume \(r_L\le r_*\).  If \(\|\hat{\sigma} - \sigma\|_{2\to2}\le r_L\), then
\[
        F_L(\hat{\sigma})
        =F_\infty(\hat{\sigma})-
        \bigl(\Phi_L(\hat{\sigma})-\Phi_\infty(\hat{\sigma})\bigr).
\]
The resolvent identity gives
\[
        \Phi_\infty(\hat\sigma)-\Phi_\infty(\sigma)
        =\mathscr S_W\!\left[
          \int_{\mathbb T^d}G_{\hat\sigma}(\xi)
          \tau(\hat\sigma-\sigma)G_\sigma(\xi)\,\rmd\xi
          \right],
\]
and hence
\begin{equation}
\label{eq:qfv-Newton-map-expansion}
\begin{aligned}
        \mathcal T_L(\hat\sigma)-\sigma
        &=-
        (\mathcal B_\sigma^\infty)^{-1} (\mathcal{B}_{\hat{\sigma}, \sigma}^\infty - \mathcal{B}_{\sigma}^\infty)[\hat{\sigma} - \sigma]
        +(\mathcal B_\sigma^\infty)^{-1}
        \bigl(\Phi_L(\hat{\sigma})-\Phi_\infty(\hat\sigma)\bigr).
\end{aligned}
\end{equation}
Therefore
\begin{align}
\label{eq:qfv-Newton-self-map}
        \|\mathcal T_L(\hat{\sigma})-\sigma\|_{2\to2}
        &\le
        C_{\rm st}
        \left(
        C_{\rm Lip}r_L^2+\varepsilon_L
        \right)\nonumber\\
        &\le \frac{r_L}{2}+\frac{r_L}{2}
        =r_L,
\end{align}
where we used \(C_{\rm st}\varepsilon_L=r_L/2\) and \(C_{\rm st}C_{\rm Lip}r_L\le1/2\).

Thus \(\mathcal T_L\) maps the closed ball \(\overline B_{r_L}(\sigma)\) into itself.  Regard the complex matrix space as the finite-dimensional real vector space \(\mathbb C^{W^d\times W^d}\cong\mathbb R^{2W^{2d}}\).  Brouwer's fixed-point theorem gives \(\sigma_L\in\overline B_{r_L}(\sigma)\) such that \(\mathcal T_L(\sigma_L)=\sigma_L\).  By \cref{eq:qfv-Newton-map} and the invertibility of \(\mathcal B_\sigma^\infty\), $F_L(\sigma_L)=0$, or equivalently, $\sigma_L=\Phi_L(\sigma_L)$. Moreover, \(r_L\le c_\sigma/2\) implies \(\Im \sigma_L\succeq c_\sigma\I_{W^d}/2\).  Uniqueness follows from the finite-volume MDE theorem for positivity-preserving covariance maps. The unscaled estimate \cref{eq:qfv-unscaled-closeness} follows by dividing \cref{eq:qfv-scaled-closeness} by $\lambda$; the trapped bounds follow from Theorem~\ref{thm:fh-fixed-block-MDE}.
\end{proof}
Now the Combes--Thomas estimate supplies the required map error.
\begin{corollary}
\label{cor:gap_inf_finit_map}
Assume \(n\eta/\lambda\geq1\).  Uniformly for \(\sigma\in\mathbb C^{W^d\times W^d}\) with \(\Im\sigma\succeq0\),
\[
 \|\Phi_\infty(\sigma)-\Phi_L(\sigma)\|_{2\to2}
 \leq C\frac{\lambda}{\eta}
 \exp\!\left(-\frac{cn\eta}{d\lambda}\right).
\]
The same estimate is therefore uniform on every subset of the upper half-plane.
\end{corollary}
\begin{proof}
Apply Lemma~\ref{lemma:CT_est_block_renorm_resolv} with the unscaled self-energy \(\theta=\sigma/\lambda\).  For the diagonal block it gives
\[
\begin{split}
 &\left\|
 \int_{\mathbb T^d}
 (\lambda\widehat\Psi(\xi)-z-\sigma/\lambda)^{-1}\,\rmd\xi
 -\frac1{n^d}\sum_{\xi\in\widetilde{\mathbb T}_n^d}
 (\lambda\widehat\Psi_L(\xi)-z-\sigma/\lambda)^{-1}
 \right\|_{2\to2}
 \\
 &\hspace{35mm}\leq
 \frac C\eta\exp\!\left(-\frac{cn\eta}{d\lambda}\right).
\end{split}
\]
Multiplication by \(\lambda\) converts the two resolvents into \(G_\sigma\) and \(G_\sigma^L\).  The boundedness of \(\mathscr S_W\) proves the claim.
\end{proof}
\begin{corollary}
\label{cor:qfv-finite-bootstrap-to-continuum}
Let \(r_*\) and \(\varepsilon_L\) be as in Theorem~\ref{thm:finite_volume_from_infinite_trap}, and assume \(\varepsilon_L\leq r_*/(4C_{\rm st})\).  Suppose that \(\tilde\sigma_L=\lambda\tilde\theta_L\) satisfies
\begin{equation}
\label{eq:qfv-finite-approximate-equation}
        \tilde\sigma_L
        =\Phi_L(\tilde\sigma_L)+\mathcal E_L
\end{equation}
and the fixed-neighborhood bootstrap hypothesis
\begin{equation}
\label{eq:qfv-finite-bootstrap}
        \|\tilde\sigma_L-\sigma_L\|_{2\to2}
        =\lambda\|\tilde\theta_L-\theta_L\|_{2\to2}
        \leq r_*/2.
\end{equation}
Then
\begin{equation}
\label{eq:qfv-bootstrap-close-to-continuum}
        \|\tilde\sigma_L-\sigma\|_{2\to2}\leq r_*,
\end{equation}
\begin{equation}
\label{eq:qfv-mixed-inverse}
        \bigl\|
        (\mathcal B_{\tilde\sigma_L,\sigma}^{\infty})^{-1}
        \bigr\|_{(2\to2)\to(2\to2)}
        \leq2C_{\rm st},
\end{equation}
and
\begin{equation}
\label{eq:qfv-finite-self-energy-comparison-scaled}
        \|\tilde\sigma_L-\sigma_L\|_{2\to2}
        \leq
        2C_{\rm st}\|\mathcal E_L\|_{2\to2}
        +4C_{\rm st}\varepsilon_L.
\end{equation}
The constants and the bootstrap radius are independent of \(\lambda\) and \(\eta\) in every fixed bounded range of \(\lambda\eta\).
\end{corollary}

\begin{proof}
Theorem~\ref{thm:finite_volume_from_infinite_trap} gives
\[
 \|\tilde\sigma_L-\sigma\|_{2\to2}
 \leq r_*/2+2C_{\rm st}\varepsilon_L\leq r_*.
\]
In particular, \(\Im\tilde\sigma_L\succeq c_\sigma\I/2\), and \cref{eq:qfv-B-Lipschitz} implies
\[
 \bigl\|(\mathcal B_\sigma^\infty)^{-1}
       (\mathcal B_{\tilde\sigma_L,\sigma}^\infty
                    -\mathcal B_\sigma^\infty)\bigr\|
 \leq C_{\rm st}C_{\rm Lip}r_*\leq1/2.
\]
The Neumann series proves \cref{eq:qfv-mixed-inverse}.  Subtracting \(\sigma=\Phi_\infty(\sigma)\) from \cref{eq:qfv-finite-approximate-equation} and using the resolvent identity gives
\begin{equation}
\label{eq:qfv-exact-continuum-mixed-equation}
 \mathcal B_{\tilde\sigma_L,\sigma}^{\infty}
       [\tilde\sigma_L-\sigma]
 =\mathcal E_L+
       \Phi_L(\tilde\sigma_L)-\Phi_\infty(\tilde\sigma_L).
\end{equation}
The last difference is bounded by \(\varepsilon_L\), since \(\tilde\sigma_L\) belongs to the ball in \cref{eq:qfv-map-error}. It follows that
\[
 \|\tilde\sigma_L-\sigma\|_{2\to2}
 \leq2C_{\rm st}(\|\mathcal E_L\|_{2\to2}+\varepsilon_L).
\]
Combining this with \(\|\sigma_L-\sigma\|_{2\to2}\leq2C_{\rm st}\varepsilon_L\) proves \cref{eq:qfv-finite-self-energy-comparison-scaled}. Uniformity in \(\eta\) follows from Lemma~\ref{lemma:fh-bounded-kinetic-extension}.
\end{proof}

We next verify the fixed-neighborhood hypothesis directly from the bootstrap assumption used in the local law.
\begin{corollary}
\label{cor:qfv-unscaled-random-application}
Let \(W<\infty\) be fixed, let \(e=E/\lambda\in I\), and suppose \(n\eta/\lambda\geq1\).  Assume that \(\lambda\chi_L(\eta)\) and the constant \(c_0>0\) are sufficiently small, where
\[
 \chi_L(\eta)=\eta^{-1}
       \exp\!\left(-\frac{cn\eta}{d\lambda}\right),
 \qquad
 \|\mathbb ER-M_L\|_{1\to\infty}\leq c_0/\lambda.
\]
Take
\begin{equation}
\label{eq:qfv-correct-scaled-residual}
 \mathcal E_L=
 \lambda\mathscr S_W\!\left[
        \Pi_{[0]}(\mathbb ER-\tilde M_L)\Pi_{[0]}
        \right].
\end{equation}
Then
\begin{equation}
\label{eq:qfv-final-unscaled-comparison}
 \|\tilde\theta_L-\theta_L\|_{2\to2}
 \leq C_{W,I}\|\mathbb ER-\tilde M_L\|_{1\to\infty}
        +C_{W,I}\chi_L(\eta).
\end{equation}
The column and entrywise comparison estimates in \cref{eq:local-law-deterministic-columns,eq:local-law-deterministic-comparison} also hold in the fixed-block regime.
\end{corollary}
\begin{proof}
First suppose \(\lambda\eta\leq C_K\), where \(C_K\) is the fixed kinetic threshold in Lemma~\ref{lemma:contra_self_map_kinetic}. Use Lemma~\ref{lemma:fh-bounded-kinetic-extension} with \(K_0=\max\{1,C_K\}\).  The covariance estimate \cref{equ:sup_op_ineq} gives
\[
 \|\tilde\sigma_L-\sigma_L\|_{2\to2}
 =\lambda\|\mathscr S_W[
       \Pi_{[0]}(\mathbb ER-M_L)\Pi_{[0]}]\|_{2\to2}
 \leq2\lambda\|\mathbb ER-M_L\|_{1\to\infty}
 \leq2c_0.
\]
Thus \(c_0\leq r_*/4\) implies \cref{eq:qfv-finite-bootstrap}. Corollary~\ref{cor:gap_inf_finit_map} gives \(\varepsilon_L\leq C\lambda\chi_L(\eta)\); the assumed smallness ensures \(\varepsilon_L\leq r_*/(4C_{\rm st})\). By the definition of \(\tilde\theta_L\), the residual in \cref{eq:qfv-correct-scaled-residual} satisfies the exact equation \cref{eq:qfv-finite-approximate-equation}.  Moreover,
\[
 \|\mathcal E_L\|_{2\to2}
 \leq2\lambda\|\mathbb ER-\tilde M_L\|_{1\to\infty}.
\]
Apply Corollary~\ref{cor:qfv-finite-bootstrap-to-continuum} and divide by \(\lambda\) to obtain \cref{eq:qfv-final-unscaled-comparison}. Note that this application requires no improved bootstrap estimate.

For the column estimates, the finite-volume MDE trap gives \(\Im\theta_L\succeq c\lambda^{-1}\I\), \(\|M_L\|_{2\to2}\leq C\lambda\), and, for the block containing \(x\),
\[
 \Im(M_L)_{xx}\leq\Tr_{[x]}\Im M_L
 \leq W^d\|\Im\theta_L\|_{2\to2}\leq C_{W,I}/\lambda.
\]
The generalized Ward identity therefore bounds both the columns and the rows of \(M_L\) in \(\ell^2\).  The original bootstrap gives \(\|\widetilde\Theta_L-\Theta_L\|_{2\to2}\leq2c_0/\lambda\). After decreasing \(c_0\), both resolvent factorizations
\[
 \tilde M_L
 =\bigl(\Id-M_L(\widetilde\Theta_L-\Theta_L)\bigr)^{-1}M_L
 =M_L\bigl(\Id-(\widetilde\Theta_L-\Theta_L)M_L\bigr)^{-1}
\]
have inverse factors bounded by two.  They give the same row and column bounds for \(\tilde M_L\).  The resolvent identity consequently implies
\[
 \|\tilde M_L-M_L\|_{1\to\infty}
 \leq C_{W,I}\|\tilde\theta_L-\theta_L\|_{2\to2},
\]
which proves the stated entrywise comparison.

For \(\lambda\eta\geq C_K\), the proof of Lemma~\ref{lemma:kinetic-M-tildeM-comparison} gives directly
\[
 \|\tilde\theta_L-\theta_L\|_{2\to2}
 \leq4\|\mathbb ER-\tilde M_L\|_{1\to\infty},
\]
as well as the required row, column, and entrywise bounds.  The two ranges cover the full interval of \(\eta\) in the local law.
\end{proof}
The remaining deterministic input for the local law is a uniform $1\to\infty$ bound on $M_L$.  We obtain it by combining the stability estimate above with the dyadic shell decomposition of Lemma~\ref{lemma:fh-Poisson-shell}, first for the continuum fibers and then for the finite torus.
\begin{lemma}
\label{lemma:fh-finite-M-max-bound}
Let $d\geq2$, let $W<\infty$ be fixed, and let $I$ be a finite-height folded bulk window.  Fix $K<\infty$, and suppose
\[
        e=E/\lambda\in I,
        \qquad
        0<\eta\leq K\lambda.
\]
Set
\[
        \tau=\lambda^{-2},
        \qquad
        \kappa=\eta/\lambda,
        \qquad
        \chi_L(\eta)
        =\frac1\eta
        \exp\!\left(-\frac{cn\eta}{d\lambda}\right).
\]
For a block symbol $\vartheta$, write \( \Theta(\vartheta)=\I_{\tilde{\mathbb Z}_n^d}\otimes\vartheta \), with the analogous constant block operator in infinite volume. If $\lambda\chi_L(\eta)$ is sufficiently small, then
\begin{equation}
\label{eq:fh-finite-M-max-bound}
        \|M_L\|_{1\to\infty}
        \leq
        \frac{C_{d,W,I,K}(1+\log\lambda)}{\lambda}
        +C_{d,W,I,K}\chi_L(\eta).
\end{equation}
\end{lemma}

\begin{proof}  
Introduce the notation
\[
        M_L(\vartheta)
        =\bigl(\lambda\Psi_L-z-\Theta(\vartheta)\bigr)^{-1},
        \qquad
        \Theta=\Theta(\theta),
        \qquad
        \Theta_L=\Theta(\theta_L),
\]
and the three resolvents
\[
\begin{aligned}
        M_\infty
        \coloneqq
        \bigl(\lambda\Psi-z-\Theta(\theta)\bigr)^{-1},\\
        \widehat M_L
        \coloneqq
        \bigl(\lambda\Psi_L-z-\Theta(\theta)\bigr)^{-1},\\
        M_L
        \coloneqq
        \bigl(\lambda\Psi_L-z-\Theta(\theta_L)\bigr)^{-1}.
\end{aligned}
\]
Thus Combes--Thomas will compare $\widehat M_L$ with $M_\infty$, since these two resolvents contain the same self-energy $\theta$.  The remaining comparison between $M_L$ and $\widehat M_L$ will follow from stability of the two fixed points.

\emph{Step 1: entrywise bound in infinite volume.} Write $h(\xi)=\widehat\Psi_W(\xi)$, denote
\begin{equation}
\label{eq:fh-Mmax-B-Gamma}
        B_\sigma(\xi)
        =h(\xi)-e-\tau\Re\sigma,
        \qquad
        \Gamma_\sigma
        =\kappa\I_{W^d}+\tau\Im\sigma,
        \qquad
        \varrho=\kappa+\tau.
\end{equation}
The fixed-block MDE estimates in the bounded-$\lambda\eta$ range, together with the kinetic contraction in the complementary range, give
\begin{equation}
\label{eq:fh-Mmax-trap}
        \|\sigma\|_{2\to2}\leq C,
        \qquad
        c\varrho\I_{W^d}
        \preceq\Gamma_\sigma
        \preceq C\varrho\I_{W^d}.
\end{equation}
Indeed, when $\lambda\eta$ is bounded, this follows from $c\I\preceq\Im\sigma\preceq C\I$; when $\lambda\eta$ is sufficiently large, $\kappa\gtrsim\tau$, $\Im\sigma\succeq0$, and $\|\sigma\|\leq C$ suffice.

Let
\[
        G_\sigma(\xi)
        =\bigl(B_\sigma(\xi)-i\Gamma_\sigma\bigr)^{-1},
        \qquad
        r_\sigma(\xi)
        =\operatorname{dist}
        \bigl(0,\operatorname{Spec}B_\sigma(\xi)\bigr).
\]
By Weyl's inequality and the free shell bound, for $s+\tau\|\Re\sigma\|\leq s_0$,
\begin{equation}
\label{eq:fh-Mmax-perturbed-shell}
\begin{split}
        N_\sigma(s)
        &\coloneqq
        \int_{\mathbb T^d}
        \Tr\mathbf 1_{[-s,s]}\bigl(B_\sigma(\xi)\bigr)\,\rmd\xi
        \\
        &\leq
        C\bigl(s+\tau\|\Re\sigma\|\bigr)
        \leq C(s+\varrho).
\end{split}
\end{equation}
The noncommuting-dissipation identity and \cref{eq:fh-Mmax-trap} imply
\begin{equation}
\label{eq:fh-Mmax-noncommuting-pointwise}
\begin{split}
        c\varrho G_\sigma(\xi)^*G_\sigma(\xi)
        &\preceq
        G_\sigma(\xi)^*\Gamma_\sigma G_\sigma(\xi)
        =\Im G_\sigma(\xi),
        \\
        \Im G_\sigma(\xi)
        &\preceq
        C\varrho
        \bigl(B_\sigma(\xi)^2+c'\varrho^2\I\bigr)^{-1}.
\end{split}
\end{equation}
Consequently,
\begin{equation}
\label{eq:fh-Mmax-pointwise}
        \|G_\sigma(\xi)\|_{2\to2}
        \leq
        \frac{C}{\sqrt{r_\sigma(\xi)^2+\varrho^2}}
        \leq\frac{C}{r_\sigma(\xi)+\varrho}.
\end{equation}

We integrate \cref{eq:fh-Mmax-pointwise} by dyadic shells.  On $\{r_\sigma\leq\varrho\}$, equations \cref{eq:fh-Mmax-perturbed-shell,eq:fh-Mmax-pointwise} give an $O(1)$ contribution.  Similarly, for every $j\geq0$ for which $2^{j+1}\varrho\leq s_0/2$,
\begin{equation}
\label{eq:fh-Mmax-one-dyadic-shell}
\begin{aligned}
     \int_{\{2^j\varrho<r_\sigma\leq2^{j+1}\varrho\}}
        \|G_\sigma(\xi)\|_{2\to2}\,\rmd\xi \leq
        \frac{C}{2^j\varrho}
        N_\sigma(2^{j+1}\varrho)
        \leq C.
\end{aligned}
\end{equation}
The region $r_\sigma\geq s_0/4$ also contributes $O(1)$.  Hence
\begin{equation}
\label{eq:fh-Mmax-fiber-L1}
        \int_{\mathbb T^d}
        \|G_\sigma(\xi)\|_{2\to2}\,\rmd\xi
        \leq
        C\bigl(1+\log_+(\varrho^{-1})\bigr)
        \leq C(1+\log\lambda),
\end{equation}
where the last inequality uses $\varrho=\kappa+\tau\geq\tau=\lambda^{-2}$. If $\varrho\geq s_0/4$, the same conclusion follows directly from the numerical-range estimate $\|G_\sigma\|\leq C/\varrho$; thus the preceding argument covers all values of $\varrho$.

Bloch inversion gives, for any two blocks $a,b$,
\begin{equation}
\label{eq:fh-Mmax-Bloch-inversion}
        \Pi_{[a]}M_\infty\Pi_{[b]}
        =\frac1\lambda
        \int_{\mathbb T^d}
        e^{i(a-b)\cdot\xi}G_\sigma(\xi)\,\rmd\xi.
\end{equation}
Since every matrix entry of a block is bounded by its operator norm, \cref{eq:fh-Mmax-fiber-L1,eq:fh-Mmax-Bloch-inversion} yield
\begin{equation}
\label{eq:fh-Mmax-infinite-bound}
        \|M_\infty\|_{1\to\infty}
        \leq\frac{C(1+\log\lambda)}{\lambda}.
\end{equation}

\emph{Step 2: comparison of the two fixed points.} Let
\[
        \varepsilon_L
        =\sup_{\|A-\sigma\|\leq r_*}
        \|\Phi_L(A)-\Phi_\infty(A)\|_{2\to2}.
\]
The same-self-energy Combes--Thomas estimate gives
\begin{equation}
\label{eq:fh-Mmax-map-error}
        \varepsilon_L
        \leq
        C\frac{\lambda}{\eta}
        \exp\!\left(-\frac{c_{\rm CT}n\eta}{d\lambda}\right)
        =C\lambda\chi_L(\eta).
\end{equation}
In the bounded-$\lambda\eta$ range, quantitative fixed-point stability gives $\|\sigma_L-\sigma\|\leq C\varepsilon_L$.  In the kinetic range, the same estimate follows directly by absorbing the contraction factor $C/(\lambda\eta)$.  Therefore, throughout the full range,
\begin{equation}
\label{eq:fh-Mmax-self-energy-gap}
        \|\theta_L-\theta\|_{2\to2}
        =\lambda^{-1}\|\sigma_L-\sigma\|_{2\to2}
        \leq C\chi_L(\eta).
\end{equation}

\emph{Step 3: uniform column estimations for the perturbation argument.} It remains to prove
\begin{equation}
\label{eq:fh-Mmax-column-bound}
\begin{split}
 \sup_x\bigl(&\|M_L|x\rangle\|_2+
                 \|M_L^*|x\rangle\|_2+
                 \|\widehat M_L|x\rangle\|_2+
                 \|\widehat M_L^*|x\rangle\|_2\bigr)
 \leq C.
\end{split}
\end{equation}
For bounded $\lambda\eta$, the infinite-volume MDE trap gives
\[
        \frac c\lambda\I\preceq\Im\theta_L
        \preceq\frac C\lambda\I.
\]
Moreover, since $\Theta_L= S[M_L]$ and $\Im M_L\succeq0$, the explicit formula for $S$ implies, for the block containing $x$,
\[
        0\leq\Im(M_L)_{xx}
        \leq\Tr_{[x]}\Im M_L
        \leq W^d\|\Im\theta_L\|_{2\to2}
        \leq\frac C\lambda.
\]
The generalized Ward identities
\[
        \Im M_L
        =M_L(\eta\I+\Im\Theta_L)M_L^*
        =M_L^*(\eta\I+\Im\Theta_L)M_L
\]
then give the first two column bounds in \cref{eq:fh-Mmax-column-bound}.  For $\widehat M_L$, use $\Im\theta\succeq c\lambda^{-1}\I$ and observe that
\[
\begin{split}
        \Im(M_\infty)_{xx}
        &\leq\frac1\lambda
        \int_{\mathbb T^d}\Tr\Im G_\sigma(\xi)\,\rmd\xi
        \\
        &\leq\frac{C\varrho}{\lambda}
        \int_{\mathbb T^d}
        \sum_{j=1}^{W^d}
        \frac{1}{\mu_j(\xi)^2+c'\varrho^2}\,\rmd\xi
        \leq\frac C\lambda,
\end{split}
\]
where $\mu_j(\xi)$ are the eigenvalues of $B_\sigma(\xi)$ and the last inequality follows from the perturbed Poisson-shell bound in \cref{eq:fh-Poisson-shell-bound}.  The same-self-energy Combes--Thomas comparison between $\widehat{M}_L$ and ${M}_\infty$, and the smallness of $\lambda\chi_L$ now give
\[
        \Im(\widehat M_L)_{xx}
        \leq\frac C\lambda+C\chi_L
        \leq\frac C\lambda,
\]
and the generalized Ward identities give the remaining two bounds.

When $\lambda\eta$ is sufficiently large, both pairs of column bounds follow instead from the kinetic Neumann estimate
\[
 \sup_x\bigl(
       \|M_L(\vartheta)|x\rangle\|_2
       +\|M_L(\vartheta)^*|x\rangle\|_2
       \bigr)
 \leq\frac{C}{\sqrt{\lambda\eta}}
 \leq C,
 \qquad \vartheta\in\{\theta,\theta_L\}.
\]
The two ranges overlap, proving \cref{eq:fh-Mmax-column-bound} uniformly.

\emph{Step 4: perturb the self-energy and conclude.} The resolvent identity gives
\begin{equation}
\label{eq:fh-Mmax-self-energy-resolvent-identity}
        M_L-\widehat M_L
        =M_L\bigl(\Theta(\theta_L)-\Theta(\theta)\bigr)\widehat M_L.
\end{equation}
Using \cref{eq:fh-Mmax-column-bound} entrywise and then \cref{eq:fh-Mmax-self-energy-gap}, we obtain
\begin{equation}
\label{eq:fh-Mmax-self-energy-resolvent-gap}
\begin{split}
        \|M_L-\widehat M_L\|_{1\to\infty}
        &\leq
        C\|\Theta(\theta_L)-\Theta(\theta)\|_{2\to2}
        \\
        &=C\|\theta_L-\theta\|_{2\to2}
        \leq C\chi_L(\eta).
\end{split}
\end{equation}
On the other hand, $\widehat M_L$ and $M_\infty$ contain the same self-energy, so the Combes--Thomas estimate applies directly:
\begin{equation}
\label{eq:fh-Mmax-same-self-energy-CT}
        \|\widehat M_L-M_\infty\|_{1\to\infty}
        \leq C\chi_L(\eta).
\end{equation}
Combining \cref{eq:fh-Mmax-infinite-bound,eq:fh-Mmax-self-energy-resolvent-gap,eq:fh-Mmax-same-self-energy-CT} proves \cref{eq:fh-finite-M-max-bound}.
\end{proof}

\appendix
\section{Kato parallel transport}
\label{appx:kato_parall_transport}
We record the standard moving-frame construction used in Lemma~\ref{lemma:fh_patch_two_resolv}.  It is included to fix the normalization and the uniform \(C^2\)-bounds.

\begin{proposition}[Kato transport for a rank-one projection]
\label{prop:kato_para_trans}
Let \(J\subset\mathbb R\) be an interval and let \(P\colon J\to\mathbb C^{W^d\times W^d}\) be a \(C^2\) family of rank-one orthogonal projections.  Fix \(s_0\in J\).  There is a \(C^2\) unitary \(V(s)\) such that
\[
        P(s)=V(s)P(s_0)V(s)^\star.
\]
Consequently, if \(u(s_0)\) spans \(\operatorname{Ran}P(s_0)\) and the columns of \(U_Q(s_0)\) span its orthogonal complement, then
\[
        U(s)\coloneqq
        V(s)\bigl(u(s_0),U_Q(s_0)\bigr)
        =
        \bigl(u(s),U_Q(s)\bigr)
\]
is a unitary frame adapted to \(P(s)\).  Moreover,
\[
        \|U'(s)\|+\|U''(s)\|
        \leq
        C\bigl(\|P'(s)\|+\|P'(s)\|^2+\|P''(s)\|\bigr).
\]
In particular, the right-hand side is uniform for the Riesz projections used in Lemma~\ref{lemma:fh_patch_two_resolv}.
\end{proposition}

\begin{proof}
Set \(K(s)=[P'(s),P(s)]\) and solve
\[
        V'(s)=K(s)V(s),
        \qquad
        V(s_0)=\I.
\]
Differentiating \(P^2=P\) gives \(PP'P=0\), hence \(K^\star=-K\) and \([K,P]=P'\).  Thus \(V(s)\) is unitary, while both \(P(s)\) and \(V(s)P(s_0)V(s)^\star\) solve the same linear ODE \(Q'=[K,Q]\) with the same initial condition.  They are therefore equal. The asserted frame follows by transporting an orthonormal frame at \(s_0\).

Finally,
\[
        U'=KU,
        \qquad
        U''=(K'+K^2)U,
        \qquad
        K'=[P'',P],
\]
because \([P',P']=0\).  The displayed derivative bound follows immediately.
\end{proof}
\section{Cubic symmetry of the scalar symbol}

Although the microscopic hopping operator is cubically symmetric, the decomposition into \(W\)-cells hides this symmetry: a coordinate reflection acts simultaneously on the block coordinate and on the ordering of the sites inside each block.  We now make this action explicit.  The only purpose of the construction is to show that the scalar response symbol \(s(k)\) inherits the reflection and coordinate-permutation symmetries of the cubic lattice.

Recall that the internal label \(u\in[W^d]\) is identified with its natural coordinate
\[
        u=(u_1,\ldots,u_d)\in\{0,\ldots,W-1\}^d.
\]
Let $g\in O(d,\mathbb Z)$ be a signed coordinate-permutation matrix.  Define the permutation matrix \(U_g\) on \(\mathbb C^{W^d}\) by
\begin{equation}
\label{eq:BA-internal-signed-permutation}
        U_ge_u
        \coloneqq
        e_{
        \frac{W-1}{2}\mathbf 1_d
        +
        g\left(
        u-\frac{W-1}{2}\mathbf 1_d
        \right)
        }.
\end{equation}
Every coordinate on the right-hand side is either \(u_j\) or \(W-1-u_j\) for some \(j\).  Hence it again belongs to \(\{0,\ldots,W-1\}^d\).  Thus \(U_g\) is a real orthogonal permutation matrix and
\[
        U_g^\top=U_g^{-1}=U_{g^{-1}}.
\]

The matrix \(U_g\) acts only on the internal coordinate \(u\).  The corresponding transformation of the full microscopic lattice is the permutation
\begin{equation}
\label{eq:BA-global-signed-permutation}
        \mathcal P_ge_{WI+u}
        \coloneqq
        e_{
        WgI
        +
        \frac{W-1}{2}\mathbf 1_d
        +
        g\left(
        u-\frac{W-1}{2}\mathbf 1_d
        \right)
        }.
\end{equation}
Indeed,
\[
\begin{aligned}
        &\frac{W-1}{2}\mathbf 1_d
        +
        g\left(
        WI+u-\frac{W-1}{2}\mathbf 1_d
        \right)
        \\
        &\qquad=
        WgI
        +
        \frac{W-1}{2}\mathbf 1_d
        +
        g\left(
        u-\frac{W-1}{2}\mathbf 1_d
        \right).
\end{aligned}
\]
Consequently, the full action is
\[
        (I,u)
        \longmapsto
        \left(
        gI,\,
        \frac{W-1}{2}\mathbf 1_d
        +
        g\left(
        u-\frac{W-1}{2}\mathbf 1_d
        \right)
        \right).
\]
The centering at the middle of the reference cell is important: it prevents the coordinate-dependent block carry that would arise from using the naive linear map \(x\mapsto gx\).

The centered affine transformation preserves nearest-neighbour adjacency.  Therefore $\mathcal P_g\Psi\mathcal P_g^\star=\Psi.$
Moreover, conjugation by \(\mathcal P_g\) only relabels the independent random blocks and conjugates each block by the real orthogonal matrix \(U_g\).  Since the GOE law is invariant under orthogonal conjugation, $\mathcal P_gH_L\mathcal P_g^\star
        \overset{\mathrm d}=H_L.$
It follows that $\mathcal P_g\,\mathbb E R(z)\,\mathcal P_g^\star = \mathbb E R(z).$
Since the block with index zero is preserved, it is easy to see that $U_g
        \blk_0[\mathbb E R(z)]
        U_g^\top
        =
        \blk_0[\mathbb E R(z)].$
Therefore, the conjugation map satisfies
\[
        \mathscr S_W[U_gBU_g^\top]
        =
        U_g\mathscr S_W[B]U_g^\top.
\]
Recalling that $\tilde\theta
        =
        \mathscr S_W\!\left[
        \blk_0[\mathbb E R(z)]
        \right],$
we obtain $U_g\tilde\theta U_g^\top=\tilde\theta.$
The hopping fibers and the renormalized resolvent fibers consequently satisfy
\begin{equation}
\label{eq:BA-fiber-signed-covariance}
        \hat\Psi(g\xi)
        =
        U_g\hat\Psi(\xi)U_g^\top,
        \qquad
        G_{\tilde\theta}(g\xi)
        =
        U_gG_{\tilde\theta}(\xi)U_g^\top.
\end{equation}

We next transfer this conjugation to the diffusion kernel.  Changing variables \(\xi=g\zeta\) in \cref{eq:BA-kernel-two-momentum}, and using \cref{eq:BA-fiber-signed-covariance}, gives
\[
\begin{aligned}
        \widehat{\tilde K}_\infty(gk)
        [U_gBU_g^\top]
        &=
        \frac1{\lambda^2}
        \mathscr S_W\!\left[
        \int_{\mathbb T^d}
        G_{\tilde\theta}(g\zeta)
        U_gBU_g^\top
        G_{\tilde\theta}(g(\zeta-k))^\star
        \,\rmd\zeta
        \right]
        \\
        &\quad=
        U_g
        \widehat{\tilde K}_\infty(k)[B]
        U_g^\top.
\end{aligned}
\]
Equivalently, for every cell matrix \(B\),
\[
\begin{aligned}
        \left(
        \Id-\widehat{\tilde K}_\infty(gk)
        \right)[U_gBU_g^\top] =
        U_g
        \left(
        \Id-\widehat{\tilde K}_\infty(k)
        \right)[B]
        U_g^\top.
\end{aligned}
\]

For \(k\) in the low-frequency region where \cref{eq:BA-scalar-response-definitions} defines \(s(k)\), the same covariance gives
\[
\begin{aligned}
        \left(
        \Id-\widehat{\tilde K}_\infty(gk)
        \right)^{-1}[\I_{W^d}]
        &=
        U_g
        \left(
        \Id-\widehat{\tilde K}_\infty(k)
        \right)^{-1}[\I_{W^d}]
        U_g^\top.
\end{aligned}
\]
Taking normalized traces, and using their invariance under orthogonal conjugation, yields
\begin{equation}
\label{eq:BA-scalar-symbol-cubic-invariance}
        s(gk)=s(k),
        \qquad
        g\in O(d,\mathbb Z).
\end{equation}

Coordinate reflections make \(s\) even in each component of \(k\), while ordinary coordinate permutations force the diagonal entries of \(D^2s(0)\) to be equal.  Together, these facts make the Hessian a scalar multiple of the identity.
\bibliography{BAM}
\bibliographystyle{plainnat}
\end{document}